\documentclass[12pt,a4paper]{article}
\usepackage{latexsym}
\usepackage{graphicx}
\usepackage{amsmath,amsthm}
\usepackage{amssymb}
\usepackage{epsfig}
\usepackage{float,lscape}
\usepackage{rotating}
\usepackage{graphicx}
\usepackage{sidecap}
\usepackage{subfigure}
\usepackage{epstopdf}
\usepackage{algorithm}
\usepackage{algpseudocode}
\newtheorem{thm}{Theorem}[section]

\newtheorem{prop}[thm]{Proposition}

\theoremstyle{definition}

\newtheorem{Exp}{Example}[section]

\theoremstyle{remark}
\newtheorem{rem}{Remark}[section]

\makeatletter \@addtoreset{equation}{section}
 \makeatother

\title{Families of non-linear subdivision schemes for scattered data fitting and their non-tensor product extensions}

\author {Ghulam Mustafa\thanks{Corresponding author: ghulam.mustafa@iub.edu.pk} \ and
Rabia Hameed\thanks{rabiahameedrazi@hotmail.com}  \\
Department of Mathematics\\
The Islamia University of Bahawalpur\\ }

\date{}
\begin{document}

\maketitle

\noindent\rule{15cm}{0.4pt}
\begin{abstract}
In this article, families of non-linear subdivision schemes are presented that are based on univariate polynomials up to degree three. Theses families of schemes are constructed by using dynamic iterative reweighed least squares method. These schemes are suitable for fitting scattered data with noise and outliers. Although these schemes are non-interpolatory, but have the ability to preserve the shape of the initial polygon in case of non-noisy initial data. The numerical examples illustrate that the schemes constructed by non-linear polynomials give better performance than the schemes that are constructed by linear polynomials (Computer-Aided Design, 58, 189-199). Moreover, the numerical examples show that these schemes have the ability to reproduce polynomials and do not cause over and under fitting of the data. Furthermore, families of non-linear bivariate subdivision schemes are also presented that are based on linear and non-linear bivariate polynomials.
\end{abstract}
\\
\renewcommand{\thefootnote}{\arabic{footnote}}

\section{Introduction}\label{Introduction}

Curve fitting is the method of constructing a curve, or a mathematical function, that has the best fit to a series of data points, probably subject to constraints. In statistics and machine learning, the common task is to fit a "model" to a set of training data, so as to be able to make consistent predictions on general untrained data. The cause of poor performance in statistics and machine learning is either overfitting or underfitting the data. In overfitting, a statistical model describes random error or noise instead of the underlying relationship. Overfitting occurs when a model is extremely complex, such as owning too many parameters relative to the number of observations. A model that has been overfit has poor predictive performance, as it overreacts to small variations in the training data. For example, the presence of outliers in the initial data points may increase the chance of overfitting in the fitted curves. Underfitting happens when a statistical model or machine learning algorithm cannot capture the underlying trend of the data. It occurs when the model or algorithm does not fit the data enough. It is often a result of an extremely simple model. For example, when fitting a linear model to non-linear data. Such a model would have poor predictive performance.

The method of least squares is one of the wonderful technique in statistics for curves fitting, whereas in the field of geometric modeling subdivision schemes are used for curves and  surfaces fitting. Subdivision schemes are recursive processes used for the quick generation of curves and surfaces. Subdivision schemes are also important ingredients in many multiscale algorithms used in data compression. In some applications, the given data need not to be retained at each step of the refinement process, which requires the use of approximating subdivision schemes. Nowadays, approximating subdivision schemes are further categorized as classical schemes and least squares subdivision schemes. Many authors in the field of geometric modeling have been constructed classical subdivision schemes for different type of initial data. Some of the recently constructed classical schemes along with their properties are available in \cite{Hameed, Hameed1, Jeong, Mustafa10, Novara}. The detail information about these schemes can be found in \cite{Dyn6}. Classical schemes are not suitable to fit scattered data and they are very sensitive to outliers. The presence of only one outlier in the initial data can badly damage the whole fitted model. This drawback of classical schemes has been inspired to construct least squares techniques based subdivision schemes.

The first attempt to construct subdivision schemes by using least squares technique is done by Dyn et al. \cite{Dyn5}. They introduced univariate linear binary subdivision schemes for refining noisy data by fitting local least squares polynomials. They also give several numerical experiments to show the performance of their schemes on initial noisy data. Mustafa et al. \cite{Mustafa9} proposed the $\ell_{1}$-regression based univariate and non-tensor product bivariate binary subdivision schemes to handle with noisy data with outliers. They also give many numerical experiments to determine the usability and practicality of their proposed schemes. Mustafa and Bari \cite{Mustafa11} proposed univariate and non-tensor product bivariate $N$-ary ($N \in \mathbb{N} \backslash \{1\}$) subdivision schemes to fit local least squares polynomials of degree 3. However, they have not given any numerical experiment to see the performance of their schemes on noisy data. The following motivation plays a vital role behind the construction of our families of schemes.
\subsection{Motivation and contribution}
\begin{figure}[htb] 
 \begin{center}
\begin{tabular}{ccccccccccc}
\epsfig{file=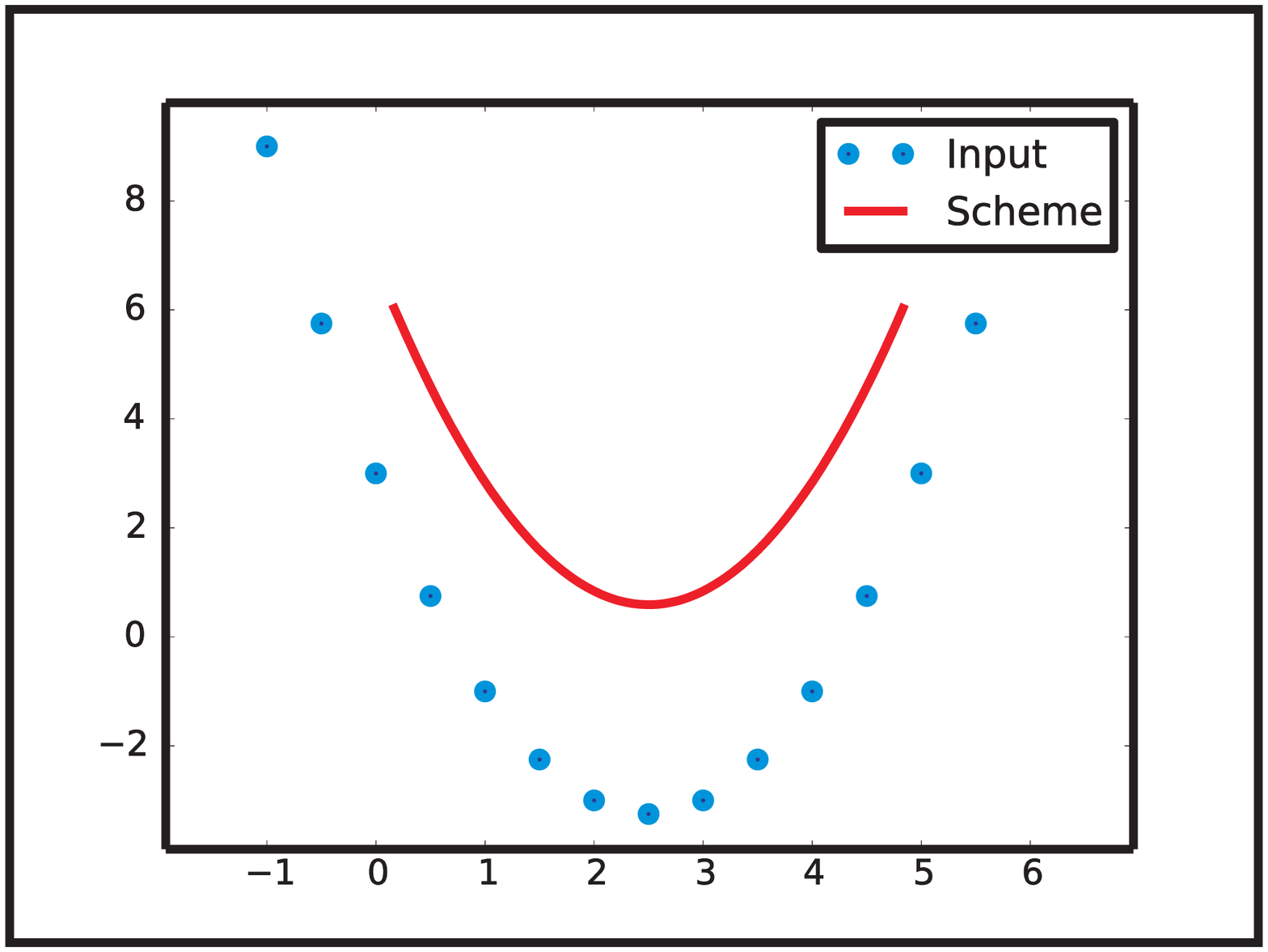, width=2.0 in}&
\epsfig{file=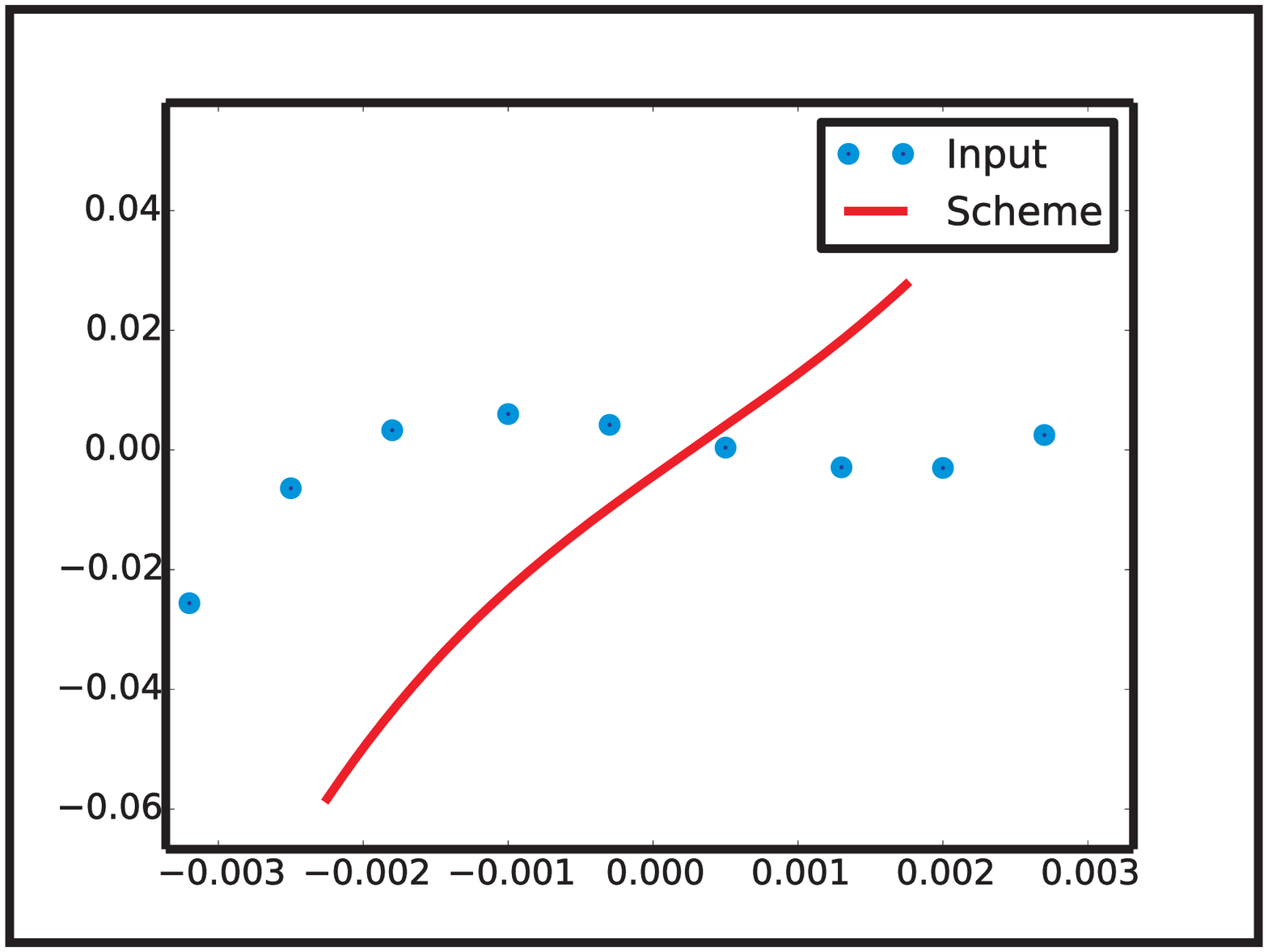, width=2.0 in}&\\
(a) & (b) \\
\epsfig{file=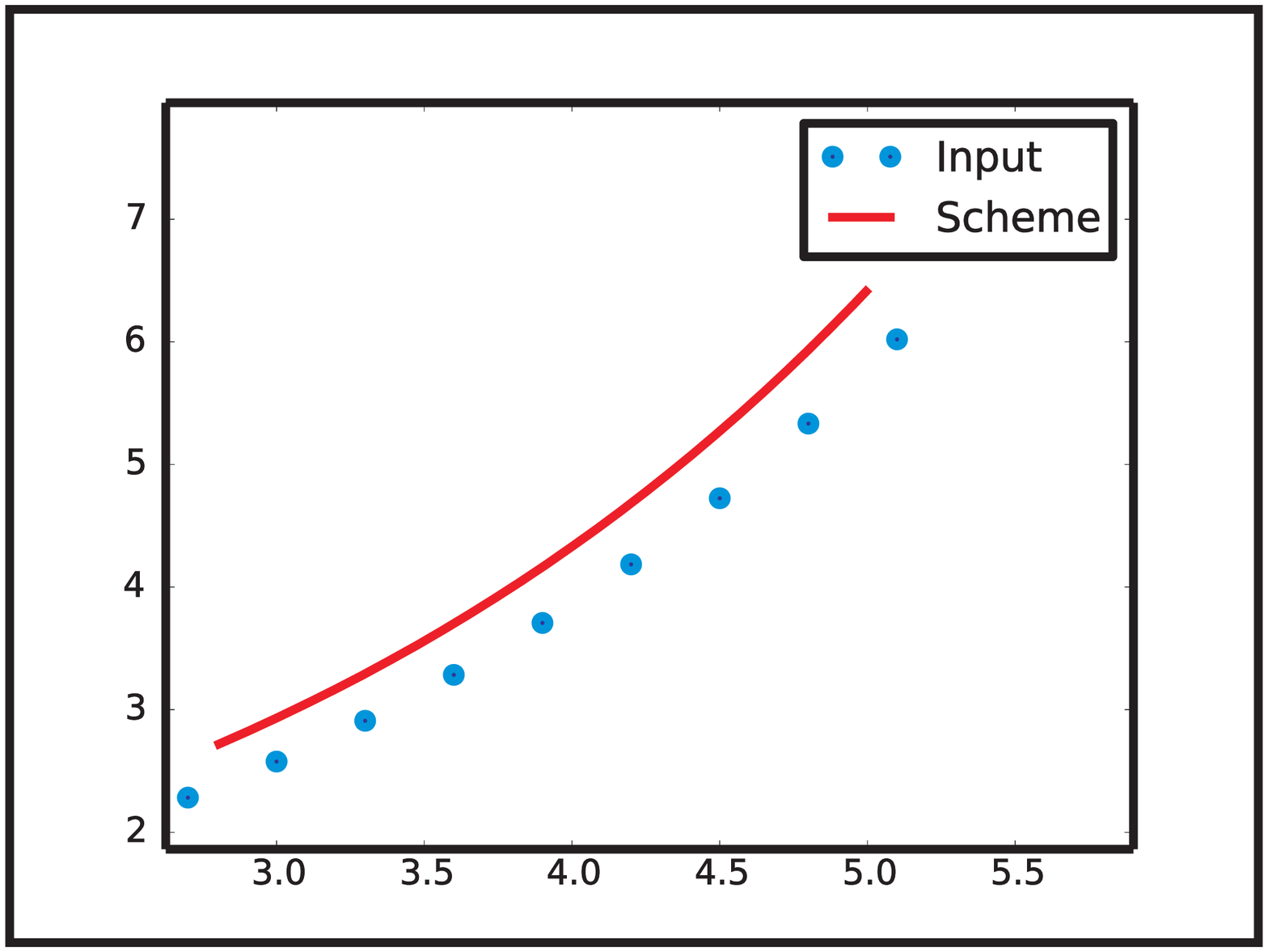, width=2.0 in}&
\epsfig{file=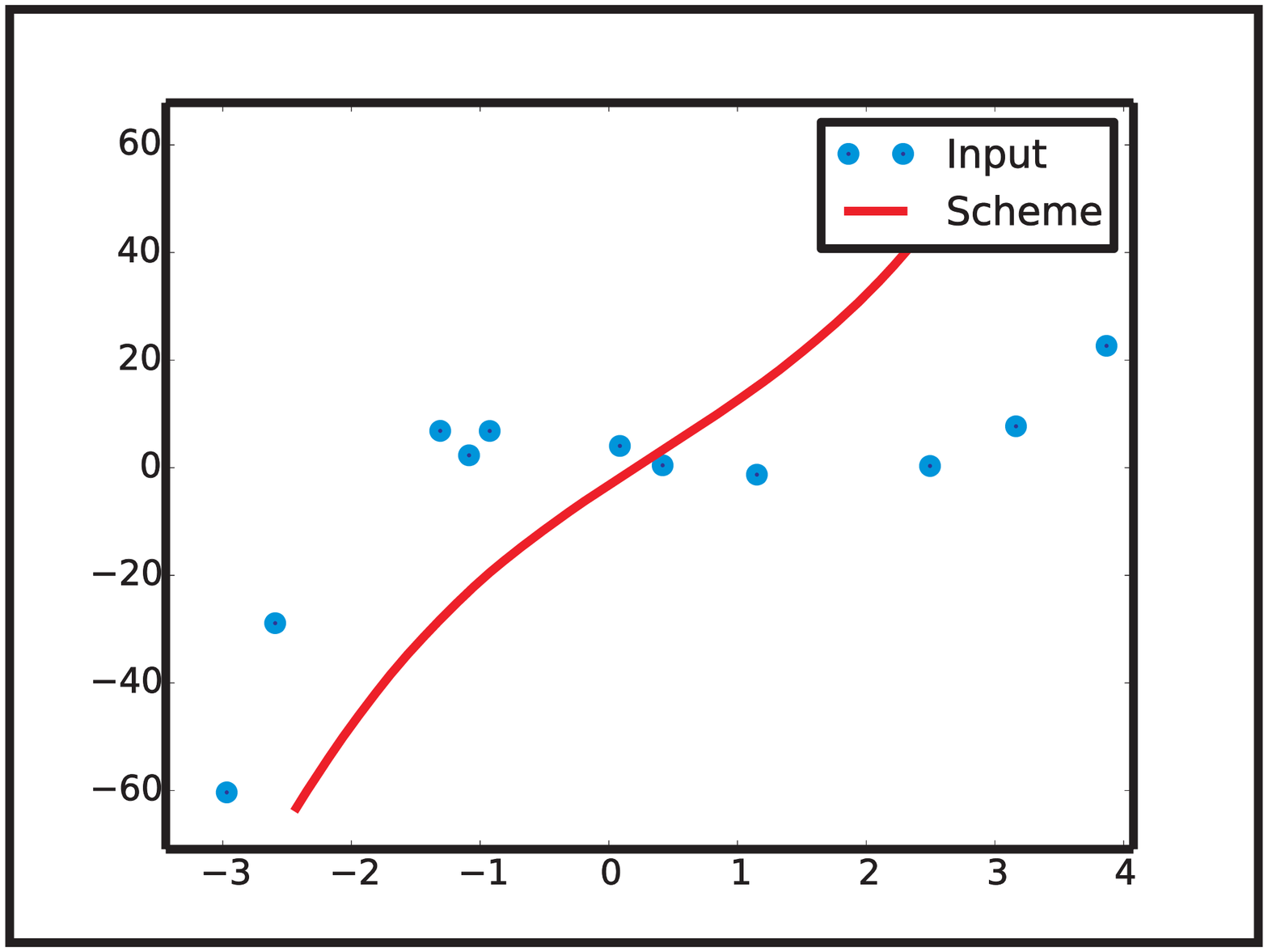, width=2.0 in}&\\
(c) & (d)
 \end{tabular}
\end{center}
\caption[Points are fitted by the schemes based on linear polynomial functions.]{\label{under-fit} \emph{Points are fitted by the schemes based on linear polynomial functions (proposed by \cite{Mustafa9}). Solid filled circles are initial data taken from some non-linear polynomial function (a)-(c) without noise and (d) with noise.}}
\end{figure}

Mustafa et al. \cite{Mustafa9} proposed subdivision schemes with dynamic iterative re-weighted weights by fitting 1D and 2D line functions to $2n$-observations and $(2n+1)$-observations in 2D and 3D spaces respectively. Their schemes have following merits:
\begin{itemize}
  \item The schemes remove noises and outliers without prior information about the input data.
  \item Overfitting of the data does not occur by using their schemes.
\end{itemize}
But the following de-merits have been noted in the schemes based on fitting 1D and 2D line functions:
\begin{itemize}
  \item The schemes exhibit underfitting when data comes from non-linear polynomial functions (see Figures \ref{under-fit}(b) and \ref{under-fit}(d)).
  \item The schemes do not reproduce non-linear polynomial functions (see Figures \ref{under-fit}(a)-\ref{under-fit}(c)).
\end{itemize}
The above issues motivate us to present the schemes based on fitting non-linear polynomial functions. We observe the following characteristics in the schemes that are based on fitting non-linear polynomial functions:
\begin{itemize}
  \item The proposed schemes have all the merits of the schemes \cite{Mustafa9}. The schemes introduced by \cite{Dyn5}, \cite{Mustafa11} and \cite{Mustafa9} are the special cases of our schemes.
  \item The proposed schemes do not cause overfit as well as underfit of the data (see Figures \ref{sharpe-corners} and \ref{outlier1}-\ref{Noise-outlier1}).
  \item These schemes also reproduce the polynomial data (see Figures \ref{quadpoly}-\ref{exppoly}).
  \item The schemes show interpolatory and approximating behaviors on non-noisy and noisy initial data respectively (see Figures \ref{sharpe-corners} and \ref{outlier1}-\ref{Noise-outlier1}).
  \item In surface case, our schemes produce smooth models comparative to the schemes of \cite{Mustafa9}.
\end{itemize}

For better understanding of this article, we have organized the material in the following fashion. In Section 2, we give the algorithms for the construction of families of binary univariate subdivision schemes. We also construct families of univariate schemes based on $\ell_{1}$-regression by using degree three, two and one polynomials respectively in this section. In Section 3, the numerical examples of the proposed families of schemes are presented. In Section 4, we construct families of bivariate subdivision schemes that are based on bivariate polynomials up to degree two. We give the conclusion of our work in Section 5.

\section{Framework for the construction of univariate schemes}
%
%
In this section, we present a framework for the construction of families of univariate approximating schemes. The framework has two steps. In first step, we fit the univariate $d$-degree polynomial to the $2n$-observations/data by using the IRLS. In second step, we derive the families of $2n$-point approximating subdivision schemes also called $\ell_{1}$-schemes.
\subsection{Step 1}
Consider the univariate polynomial function of degree $d$
\begin{eqnarray}\label{polyfunc}
f(x_{r})&=&\sum\limits_{a=0}^{d}\beta_{a}f_{a}(x_{r}),
\end{eqnarray}
where $f_{a}(x_{r})=x_{r}^{a}$, $a=0,1,\ldots,d$, are the monomial functions.\\
Now express the polynomial function (\ref{polyfunc}) with respect to the observations $(x_{r}=r,f_{r})$ for $r=-n+1,\ldots,n$ ($n \in \mathbb{N}\backslash \{1\}$) as
\begin{eqnarray}\label{poly}
f_{r}&=&f(r)=\sum\limits_{a=0}^{d}r^{a}\beta_{a}.
\end{eqnarray}
Let $\beta=\{\beta_{a}:a=0,1,\ldots.,d\}$, then use $\ell_{1}$ fitting function of \cite{Rey} for estimating $\beta$. The $\ell_{1}$-norm operator is defined as
\begin{eqnarray*}
&& \min_{\beta \in \mathbb{R}} \sum\limits_{r=-n+1}^{n}\left|f_{r}-\sum\limits_{a=0}^{d}r^{a}\beta_{a}\right|.
\end{eqnarray*}
The optimum of $\beta$ is expressed as
\begin{eqnarray}\label{optimization}
\nonumber \beta_{0},\beta_{1},\ldots,\beta_{d}&=& \mbox{arg} \,\ \min_{\beta \in \mathbb{R}} \sum\limits_{r=-n+1}^{n}\left|f_{r}-\sum\limits_{a=0}^{d}r^{a}\beta_{a}\right|\\
&&=\mbox{arg} \,\ \min_{\beta \in \mathbb{R}}F(\beta_{0},\beta_{1},\ldots,\beta_{d}).
\end{eqnarray}
Use the IRLS method \cite{Bissantz} which is based on the thought that, in a first step, the $\ell_{1}$ norm $F$, being a convex functional, can be approximated by a family of smooth convex functionals $F_{\delta}$, $\delta >0$. That is
\begin{eqnarray}\label{F-delta}
F_{\delta}(\beta_{0},\beta_{1},\ldots,\beta_{d})&=&\sum\limits_{r=-n+1}^{n}\left[\left(f_{r}-\sum\limits_{a=0}^{d}r^{a}\beta_{a}\right)^{2}+\delta\right]^{\frac{1}{2}}.
\end{eqnarray}
The regularization of the non-smooth functional (\ref{optimization}) by (\ref{F-delta}) has been done by \cite{Vogel}. So by (\cite{Bissantz}, Theorem 1), we get
\begin{eqnarray*}
\beta_{0,\delta},\beta_{1,\delta},\ldots,\beta_{d,\delta}&=& \mbox{arg} \,\ \min_{\beta \in \mathbb{R}}F_{\delta}(\beta_{0},\beta_{1},\ldots,\beta_{d}).
\end{eqnarray*}
This means that $\beta_{\delta}=\{\beta_{a,\delta}:a=0,1,\ldots,d\}$ is the approximation of $\beta$. \\
Thus, in order to compute $\beta_{\delta}$, use the iterative formula given below
\begin{eqnarray}\label{I-Formulae}
\beta^{(m+1)}_{0,\delta},\beta^{(m+1)}_{1,\delta},\ldots,\beta^{(m+1)}_{d,\delta}&=& \mbox{arg} \,\ \min_{\beta \in \mathbb{R}} \sum\limits_{r=-n+1}^{n}w^{(m)}_{r}\left(f_{r}-\sum\limits_{a=0}^{d}r^{a}\beta_{a}\right)^{2},
\end{eqnarray}
where
\begin{eqnarray}\label{weights}
w^{(m)}_{r}&=&\left[\left(f_{r}-\sum\limits_{a=0}^{d}r^{a}\beta^{(m)}_{a,\delta}\right)^{2}+\delta\right]^{-\frac{1}{2}}, \,\ m \geqslant 0 \,\ \mbox{and} \,\ n >0.
\end{eqnarray}
Here $w^{(m)}_{r}$ is called the dynamic iterative weights. When $w^{(m)}_{r}=1$, $\forall$ $r$ $\&$ $m$, then the $\ell_{1}$-regression problem converts into the $\ell_{2}$-regression problem.

The ordinary least squares method can be used to find out the starting values $\beta^{(0)}_{\delta}=\beta=\{\beta^{(0)}_{a,\delta}=\beta_{a}:a=0,1,\ldots,d\}$, i.e. determine the values of unknown parameters in (\ref{poly}) to make the sum of the squares of residuals minimum. The sum of the squares of residuals is defined as
\begin{eqnarray}\label{SSRM}
R&=&\sum\limits_{r=-n+1}^{n}\left(f_{r}-\sum\limits_{a=0}^{d}r^{a}\beta^{(0)}_{a,\delta}\right)^{2}.
\end{eqnarray}
Now differentiating (\ref{I-Formulae}) and (\ref{SSRM}) with respect to $\beta=\{\beta_{a}:a=0,1,\ldots,d\}$, then setting each system of the $(d+1)$-equations equal to zero and then by solving these systems for unknowns values, we get the values of $\beta^{(m+1)}_{\delta}=\{\beta^{(m+1)}_{a,\delta}:a=0,1,\ldots,d\}$ and $\beta^{(0)}_{\delta}=\{\beta^{(0)}_{a,\delta}:a=0,1,\ldots,d\}$ respectively.

Continue iterations until
\begin{eqnarray*}
&&\mbox{max}\left(\left|\beta^{(m+1)}_{a,\delta}-\beta^{(m)}_{a,\delta}\right|:a=0,1,\ldots,d\right)<\epsilon .
\end{eqnarray*}
Now by substituting optimum $\beta^{(m+1)}_{\delta}=\{\beta^{(m+1)}_{a,\delta}:a=0,1,\ldots,d\}$ of $\beta=\{\beta_{a}:a=0,1,\ldots,d\}$ into (\ref{poly}), we get the following best fitted $d$-degree polynomial to the $2n$-observations
\begin{eqnarray}\label{beta-a-delta}
f(r)&=&\sum\limits_{a=0}^{d}r^{a}\beta^{(m+1)}_{a,\delta}.
\end{eqnarray}
This step is briefly described in Algorithm \ref{parameters-Estimating}.
\subsection{Step 2}
In this step, first evaluate the best fitted $d$-degree polynomial (\ref{beta-a-delta}) at $r=\frac{1}{4}$ and $r=\frac{3}{4}$. Then set the following notations for iterative point of view. That is, set
\begin{eqnarray*}
&&f\left(\frac{1}{4}\right)=f_{2i}^{k+1}, \qquad
f\left(\frac{3}{4}\right)=f_{2i+1}^{k+1},\\
&&\beta^{(m+1)}_{a,\delta}=\beta^{(m+1)}_{a,\delta,i}, \qquad
w_{r}^{(m)}=w_{i+r}^{(m)},
\end{eqnarray*}
and also set the data values $f_{r}=f^{k}_{i+r}$, involved in $\beta^{(m+1)}_{a,\delta}$, where $f^{k+1}_{i}$ and $f^{k}_{i}$ are the control points at level $k+1$ and $k$ respectively. The refinement rules $f_{2i}^{k+1}$ and $f_{2i+1}^{k+1}$ make the families of $2n$-point $\ell_{1}$-schemes denoted by $D_{2n,d}$ to fit a curve to the set of data points. Step 2 is briefly described in Algorithm \ref{SS-A}.

\begin{algorithm}[htb!] 
   \caption[Estimation of the unknown parameters $\beta$ in (\ref{poly}).]{\label{parameters-Estimating} \emph{Estimation of the unknown parameters $\beta$ in (\ref{poly}).}}
    \begin{algorithmic}[1]
    \State \textbf{input:}  $d$, $n$, $\epsilon$, $\sigma$, $m^{**}$
    \For{$m^{*} = 0$ to ${m^{**}}$}
            \If {$m^{*}=0$}
                \State calculate: $\beta^{(m^{*})}_{\delta}=\{\beta^{(m^{*})}_{a,\delta}:a=0,1,\ldots,d\}$ by differentiating (\ref{SSRM}) and solving system of equations
            \Else
                \State calculate: $\beta^{(m^{*})}_{\delta}=\{\beta^{(m^{*})}_{a,\delta}:a=0,1,\ldots,d\}$ by differentiating (\ref{I-Formulae}) and solving system of equations
            \EndIf
            \State calculate: $w^{(m^{*})}_{r},r=-n+1,\ldots,n$ by (\ref{weights})
            \If {$m^{*} = 0$}
            \State \textbf{go to} 2
            \EndIf
            \State calculate: $M^{*}=\max \limits_{a}\left|\beta^{(m^{*})}_{a,\delta}-\beta^{(m^{*}-1)}_{a,\delta}\right|$, $a=0,1,\ldots,d$.
            \If {$M^{*}<\epsilon$}
            \State \textbf{go to} 17
            \EndIf
        \EndFor
        \State set: $\beta^{(m+1)}_{a,\delta}=\beta^{(m^{*})}_{a,\delta}$
        \State \textbf{output:} $\beta^{(m+1)}_{a,\delta},a=0,1,\ldots,d$
\end{algorithmic}
\end{algorithm}

\begin{algorithm}[htb!] 
\caption[Subdivision of 2D data points.]{\label{SS-A} \emph{Subdivision of 2D data points.}}
    \begin{algorithmic}[1]
    \State \textbf{input:}  $k^{**}$, $f^{0}=\{f^{0}_{i+r}:i \in \mathbb{Z};r=-n+1,\ldots,n\}$, $\beta^{(m+1)}_{\delta}=\{\beta^{(m+1)}_{a,\delta,i}:i \in \mathbb{Z};a=0,1,\ldots,d\}$
    \State calculate: $f\left(\frac{1}{4}\right)$, $f\left(\frac{3}{4}\right)$ from (\ref{beta-a-delta})
    \State set: $f^{1}_{2i}=f\left(\frac{1}{4}\right)$, $f^{1}_{2i+1}=f \left(\frac{3}{4}\right)$
    \State rewrite step 3: $f^{1}=S_{2n,d}f^{0}$ \Comment{$S_{2n,d}$ is the subdivision matrix/rule of the scheme $D_{2n,d}$}
    \For{$k^{*}=1$ to $k^{**}$} \Comment{$k^{**} \in \mathbb{N}$ number of subdivision steps}
        \State $f^{k^{*}}=S_{2n,d}f^{k^{*}-1}$
    \EndFor
    \State set $f^{k^{**}}=f^{k+1}$
    \State \textbf{output:} $f^{k+1}=\{f^{k+1}_{i}:i \in \mathbb{Z}\}$ \Comment{$(k+1)$-th level subdivided data}
\end{algorithmic}
\end{algorithm}

\begin{rem}
The $N$-ary subdivision scheme is a tool which maps a polygon/mesh $f^{k}=\{f^{k}_{i} \in \mathbb{R}^{\rho}, i \in \mathbb{Z}^{\rho}\}$ to a refined mesh $f^{k+1}=\{f^{k+1}_{i} \in \mathbb{R}^{\rho}, i \in \mathbb{Z}^{\rho}\}$ where $f^{k}_{i}$ denotes a sequence of points in $\mathbb{R}^{\rho}$ and $k$ is a non negative integer. For $\rho=1$ and $\rho=2$, the $N$-ary scheme is univariate and  bivariate respectively.\\
Furthermore, for $\rho=1,2$, the $D_{(2n)^\rho,d}$ and $D_{(2n+1)^\rho,d}$ schemes mean the binary $\rho$-variate $\ell_{1}$-schemes computed from $\rho$-variate $d$-degree polynomial. These schemes take $(2n)^\rho$ and $(2n+1)^\rho$ consecutive points from initial polygon/mesh respectively to compute a new point in order to get a refined polygon/mesh. In other words, $D_{h^\rho,d}$ schemes mean the binary $\rho$-variate schemes computed from $\rho$-variate $d$-degree polynomials where $h \in \left\{2n,2n+1:n \in \mathbb{N}\backslash \{1\}\right\}$. From numerical experiments, it is to be noted that small value of $\delta$ is a good choice, while good choice of $m$ is 5 or 6.
\end{rem}

\begin{prop}\label{2n-betas-OLS-cor}
If we substitute $d=3$ in (\ref{poly}), then from Algorithm \ref{parameters-Estimating} the starting values $\beta^{(0)}_{\delta}=\{\beta^{(0)}_{a,\delta}:a=0,1,\ldots,3\}$ are:
\begin{eqnarray}\label{ols-betas-3-2n}
&&\beta^{(0)}_{3,\delta}=\sum\limits_{r=-n+1}^{n}\frac{35(10r^{3}-15r^{2}+11r-6n^{2}r+3n^{2}-3)}{n(n^{2}-1)(4n^{2}-1)(4n^{2}-9)}f_{r},
\end{eqnarray}
\begin{eqnarray}\label{ols-betas-2-2n}
&&\beta^{(0)}_{2,\delta}=\sum\limits_{r=-n+1}^{n}\frac{15(3r^{2}-3r-n^{2}+1)}{2n(n^{2}-1)(4n^{2}-1)}f_{r}-\frac{3}{2}\beta^{(0)}_{3,\delta},
\end{eqnarray}
\begin{eqnarray}\label{ols-betas-1-2n}
&&\beta^{(0)}_{1,\delta}=\sum\limits_{r=-n+1}^{n}\frac{3(2r-1)}{n(4n^{2}-1)}f_{r}-\beta^{(0)}_{2,\delta}-\frac{1}{5}(3n^{2}+2)\beta^{(0)}_{3,\delta},
\end{eqnarray}
\begin{eqnarray}\label{ols-betas-0-2n}
&&\beta^{(0)}_{0,\delta}=\frac{1}{2n}\sum\limits_{r=-n+1}^{n}f_{r}-\frac{1}{2}\beta^{(0)}_{1,\delta}-\frac{1}{6}(2n^{2}+1)\beta^{(0)}_{2,\delta}-\frac{n^{2}}{2}\beta^{(0)}_{3,\delta}.
\end{eqnarray}
\end{prop}

\begin{prop}\label{betas-IRLS-cor}
If we substitute $d=3$ in (\ref{poly}), then from Algorithm \ref{parameters-Estimating} the optimum $\beta^{(m+1)}_{\delta}=\{\beta^{(m+1)}_{a,\delta}:a=0,1,\ldots,3\}$ of unknown parameters $\beta$ are
\begin{eqnarray}\label{IRLS-betas-3-2n}
\nonumber\beta^{(m+1)}_{3,\delta}&=&\frac{1}{\lambda^{(m)}_{0}\lambda^{(m)}_{1}-\left(\lambda^{(m)}_{2}\right)^{2}}\sum\limits_{r=-n+1}^{n}\left[\lambda^{(m)}_{0} \left\{\left(r^{3}\alpha^{(m)}_{0}-\alpha^{(m)}_{3}\right)\chi^{(m)}_{0}\right.\right.\\&&
\left.\left. \nonumber-\left(r\alpha^{(m)}_{0}-\alpha^{(m)}_{1}\right)\chi^{(m)}_{2}\right\}-\lambda^{(m)}_{2}\left\{\left(r^{2}\alpha^{(m)}_{0}-\alpha^{(m)}_{2}\right)\chi^{(m)}_{0}\right.\right.\\&&
\left.\left.-\left(r\alpha^{(m)}_{0}-\alpha^{(m)}_{1}\right)\chi^{(m)}_{1}\right\}\right]w^{(m)}_{r}f_{r},
\end{eqnarray}

\begin{eqnarray}\label{IRLS-betas-2-2n}
\nonumber\beta^{(m+1)}_{2,\delta}&=&\frac{1}{\lambda^{(m)}_{0}} \left[\sum\limits_{r=-n+1}^{n}\left\{\left(r^{2}\alpha^{(m)}_{0}-\alpha^{(m)}_{2}\right)\chi^{(m)}_{0}-\left(r\alpha^{(m)}_{0}-\alpha^{(m)}_{1}\right)\chi^{(m)}_{1}\right\}\right.\\ &&\left.
\times w^{(m)}_{r}f_{r} -\lambda^{(m)}_{2}\beta^{(m+1)}_{3,\delta}\right],
\end{eqnarray}

\begin{eqnarray}\label{IRLS-betas-1-2n}
\nonumber \beta^{(m+1)}_{1,\delta}&=&\frac{1}{\chi^{(m)}_{0}}\left[\sum\limits_{r=-n+1}^{n}\left(r\alpha^{(m)}_{0}-\alpha^{(m)}_{1}\right)w^{(m)}_{r}f_{r}-\sum\limits_{a=2}^{3}\chi_{a-1}\beta^{(m+1)}_{a,\delta}\right],\\
\end{eqnarray}

\begin{eqnarray}\label{IRLS-betas-0-2n}
\beta^{(m+1)}_{0,\delta}&=&\frac{1}{\alpha^{(m)}_{0}}\left[\sum\limits_{r=-n+1}^{n}w^{(m)}_{r}f_{r}-\sum\limits_{a=1}^{3}\alpha^{(m)}_{a}\beta^{(m+1)}_{a,\delta}\right],
\end{eqnarray}
where
\begin{eqnarray}\label{weights1}
\left\{\begin{array}{ccccccc}
&&\lambda^{(m)}_{0}=\chi^{(m)}_{0}\chi^{(m)}_{3}-\left(\chi^{(m)}_{1}\right)^{2}, \\ \\ &&\lambda^{(m)}_{1}=\chi^{(m)}_{0}\chi^{(m)}_{5}-\left(\chi^{(m)}_{2}\right)^{2}, \\ \\
&&\lambda^{(m)}_{2}=\chi^{(m)}_{0}\chi^{(m)}_{4}-\chi^{(m)}_{1}\chi^{(m)}_{2}, \,\
\end{array}\right.
\end{eqnarray}
\begin{eqnarray}\label{weights3}
\left\{\begin{array}{ccccccc}
&&\chi^{(m)}_{0}=\alpha^{(m)}_{0}\alpha^{(m)}_{2}-\left(\alpha^{(m)}_{1}\right)^{2}, \\ \\
&&\chi^{(m)}_{1}=\alpha^{(m)}_{0}\alpha^{(m)}_{3}-\alpha^{(m)}_{1}\alpha^{(m)}_{2}, \\ \\ &&\chi^{(m)}_{2}=\alpha^{(m)}_{0}\alpha^{(m)}_{4}-\alpha^{(m)}_{1}\alpha^{(m)}_{3},
\end{array}\right.
\end{eqnarray}

\begin{eqnarray}\label{weights5}
\left\{\begin{array}{ccccccc}
&&\chi^{(m)}_{3}=\alpha^{(m)}_{0}\alpha^{(m)}_{4}-\left(\alpha^{(m)}_{2}\right)^{2}, \\ \\ &&\chi^{(m)}_{4}=\alpha^{(m)}_{0}\alpha^{(m)}_{5}-\alpha^{(m)}_{2}\alpha^{(m)}_{3}, \\ \\
&&\chi^{(m)}_{5}=\alpha^{(m)}_{0}\alpha^{(m)}_{6}-\left(\alpha^{(m)}_{3}\right)^{2},
\end{array}\right.
\end{eqnarray}
\begin{eqnarray}\label{weights4}
\alpha^{(m)}_{b}=\sum\limits_{r=-n+1}^{n}r^{b}w^{(m)}_{r}\,\ \mbox{for} \,\ 0 \leqslant b \leqslant 6.
\end{eqnarray}
\end{prop}

\begin{prop}\label{2n+1-betas-OLS-cor}
If we put $d=3$ in (\ref{poly}) and replace $r=-n+1,\ldots,n$ by $r=-n,\ldots,n$, then from Algorithm \ref{parameters-Estimating} the starting values $\beta^{(0)}_{\delta}=\{\beta^{(0)}_{a,\delta}:a=0,1,\ldots,3\}$ are:
\begin{eqnarray}\label{ols-betas-3-2n+1}
\beta^{(0)}_{3,\delta}&=&\sum\limits_{r=-n}^{n}\frac{35(5r^{3}-3n^{2}r-3nr+r)f_{r}}{n(n^{2} -1)(n+2)(4n^{2} -1)(2n+3)},
\end{eqnarray}

\begin{eqnarray}\label{ols-betas-2-2n+1}
\beta^{(0)}_{2,\delta}&=&\sum\limits_{r=-n}^{n}\frac{15(3 r^{2}-n^{2}-n)f_{r}}{n(n+1)(4n^{2} -1)(2n+3)},
\end{eqnarray}

\begin{eqnarray}\label{ols-betas-1-2n+1}
\beta^{(0)}_{1,\delta}&=&\sum\limits_{r=-n}^{n}\frac{3 r f_{r}}{n(n+1)(2n+1)}-\frac{1}{5}(3 n^{2}+3 n-1) \beta^{(0)}_{3,\delta},
\end{eqnarray}

\begin{eqnarray}\label{ols-betas-0-2n+1}
\beta^{(0)}_{0,\delta}&=&\sum\limits_{r=-n}^{n}\frac{1}{2n+1} f_{r}-\frac{1}{3} n(n+1) \beta^{(0)}_{2,\delta},
\end{eqnarray}
\end{prop}
\begin{rem}
Since from iterative point of view, we have set $\beta^{(m+1)}_{a,\delta}=\beta^{(m+1)}_{a,\delta,i}$ and $w^{(m)}_{r}=w^{(m)}_{i+r}$. Furthermore, set the symbols used in Proposition \ref{betas-IRLS-cor} as $\lambda^{(m)}_{b}=\lambda^{(m)}_{b,i}$ for $b=0,1,2$; $\chi^{(m)}_{b}=\chi^{(m)}_{b,i}$ for $b=0,1,\ldots,5$ and $\alpha^{(m)}_{b}=\alpha^{(m)}_{b,i}$ for $b=0,1,\ldots,6$.
\end{rem}
\subsection{Family of $2n$-point schemes based on polynomial of degree 3}

From Propositions \ref{2n-betas-OLS-cor}-\ref{betas-IRLS-cor} and Algorithm \ref{SS-A} for $d=3$, we get the following family of $2n$-point subdivision schemes $D_{2n,3}$ based on fitting polynomial of degree 3
\begin{eqnarray}\label{2n-cubic-scheme}
\left\{\begin{array}{ccccccc}
f^{k+1}_{2i}&=&\frac{1}{64\gamma^{(m)}_{n,i}}\sum\limits_{r=-n+1}^{n}\left(p^{(m)}_{0,r,n,i}+p^{(m)}_{1,r,n,i}+p^{(m)}_{2,r,n,i}+p^{(m)}_{3,r,n,i}\right)w^{(m)}_{i+r}f^{k}_{i+r},\\ \\
f^{k+1}_{2i+1}&=&\frac{1}{64\gamma^{(m)}_{n,i}}\sum\limits_{r=-n+1}^{n}\left(p^{(m)}_{4,r,n,i}+p^{(m)}_{5,r,n,i}+p^{(m)}_{6,r,n,i}+p^{(m)}_{7,r,n,i}\right)w^{(m)}_{i+r}f^{k}_{i+r},
\end{array}\right.
\end{eqnarray}
where
\begin{eqnarray}\label{gamma-n}
\nonumber\gamma^{(m)}_{n,i}&=&-(\alpha^{(m)}_{5,i})^2 \alpha^{(m)}_{0,i} \alpha^{(m)}_{2,i}+\alpha^{(m)}_{0,i} \alpha^{(m)}_{2,i} \alpha^{(m)}_{6,i} \alpha^{(m)}_{4,i}+2 \alpha^{(m)}_{5,i} \alpha^{(m)}_{3,i} \alpha^{(m)}_{4,i} \alpha^{(m)}_{0,i}-\alpha^{(m)}_{0,i} \alpha^{(m)}_{6,i}\\&&
\nonumber\times (\alpha^{(m)}_{3,i})^2-\alpha^{(m)}_{0,i} (\alpha^{(m)}_{4,i})^3- \alpha^{(m)}_{6,i} (\alpha^{(m)}_{2,i})^3+2 \alpha^{(m)}_{3,i} (\alpha^{(m)}_{2,i})^2 \alpha^{(m)}_{5,i}+(\alpha^{(m)}_{2,i})^2 \times\\&&
\nonumber(\alpha^{(m)}_{4,i})^2-2 \alpha^{(m)}_{5,i} \alpha^{(m)}_{2,i} \alpha^{(m)}_{1,i} \alpha^{(m)}_{4,i}-3 (\alpha^{(m)}_{3,i})^2 \alpha^{(m)}_{2,i} \alpha^{(m)}_{4,i}+2 \alpha^{(m)}_{2,i} \alpha^{(m)}_{6,i} \alpha^{(m)}_{1,i} \alpha^{(m)}_{3,i}\\&&
\nonumber-2 \alpha^{(m)}_{5,i}(\alpha^{(m)}_{3,i})^2 \alpha^{(m)}_{1,i}+(\alpha^{(m)}_{5,i})^2 (\alpha^{(m)}_{1,i})^2+2 (\alpha^{(m)}_{4,i})^2 \alpha^{(m)}_{1,i} \alpha^{(m)}_{3,i}-\alpha^{(m)}_{6,i} (\alpha^{(m)}_{1,i})^2 \\&& \times
\alpha^{(m)}_{4,i}+(\alpha^{(m)}_{3,i})^4,
\end{eqnarray}

\begin{eqnarray}\label{p-0-n}
\nonumber p^{(m)}_{0,r,n,i}&=&4 \alpha^{(m)}_{2,i}r\alpha^{(m)}_{1,i}\alpha^{(m)}_{6,i}-64 \alpha^{(m)}_{2,i}\alpha^{(m)}_{4,i}r\alpha^{(m)}_{5,i}-64 \alpha^{(m)}_{1,i}r^3\alpha^{(m)}_{5,i}\alpha^{(m)}_{3,i}+64 r^3(\alpha^{(m)}_{4,i})^2\\&&
\nonumber\times \alpha^{(m)}_{1,i}- \alpha^{(m)}_{0,i}r^3(\alpha^{(m)}_{3,i})^2-64 r(\alpha^{(m)}_{3,i})^2\alpha^{(m)}_{5,i}-4 (\alpha^{(m)}_{3,i})^2r^3\alpha^{(m)}_{1,i}-4 \alpha^{(m)}_{4,i}\alpha^{(m)}_{1,i}\\&&
\nonumber\times \alpha^{(m)}_{5,i}+64 \alpha^{(m)}_{2,i}\alpha^{(m)}_{6,i}\alpha^{(m)}_{4,i}-16 \alpha^{(m)}_{4,i}\alpha^{(m)}_{1,i}\alpha^{(m)}_{6,i}+64 (\alpha^{(m)}_{2,i})^2r^3\alpha^{(m)}_{5,i}+(\alpha^{(m)}_{2,i})^2\\&&
\nonumber\times \alpha^{(m)}_{3,i}r^2+\alpha^{(m)}_{5,i}r^2(\alpha^{(m)}_{1,i})^2-16 \alpha^{(m)}_{4,i}r(\alpha^{(m)}_{3,i})^2-16 (\alpha^{(m)}_{2,i})^2r\alpha^{(m)}_{6,i}- \alpha^{(m)}_{5,i}\times \\&&
\nonumber\alpha^{(m)}_{3,i}\alpha^{(m)}_{1,i}-16 \alpha^{(m)}_{5,i}(\alpha^{(m)}_{3,i})^2- \alpha^{(m)}_{2,i}\alpha^{(m)}_{4,i}\alpha^{(m)}_{1,i}r^2+\alpha^{(m)}_{0,i}\alpha^{(m)}_{4,i}\alpha^{(m)}_{3,i}r^2+64 \alpha^{(m)}_{2,i}\\&&
\nonumber\times \alpha^{(m)}_{3,i}\alpha^{(m)}_{5,i}r^2+4 \alpha^{(m)}_{0,i}\alpha^{(m)}_{4,i}r\alpha^{(m)}_{5,i}-4 (\alpha^{(m)}_{2,i})^2\alpha^{(m)}_{6,i}+\alpha^{(m)}_{0,i}\alpha^{(m)}_{2,i}r^3\alpha^{(m)}_{4,i}+16 \\&&
\times \alpha^{(m)}_{0,i}r^3\alpha^{(m)}_{5,i}\alpha^{(m)}_{3,i}-4 \alpha^{(m)}_{2,i}\alpha^{(m)}_{4,i}r\alpha^{(m)}_{3,i}-4 \alpha^{(m)}_{2,i}\alpha^{(m)}_{4,i}r^3\alpha^{(m)}_{1,i},
\end{eqnarray}

\begin{eqnarray}\label{p-1-n}
\nonumber p^{(m)}_{1,r,n,i}&=&+16 \alpha^{(m)}_{2,i}r^2\alpha^{(m)}_{6,i}\alpha^{(m)}_{1,i}+16 \alpha^{(m)}_{0,i}\alpha^{(m)}_{5,i}r^2\alpha^{(m)}_{4,i}+64 \alpha^{(m)}_{2,i}r\alpha^{(m)}_{3,i}\alpha^{(m)}_{6,i}- \alpha^{(m)}_{2,i}\alpha^{(m)}_{5,i}\\&&
\nonumber \times \alpha^{(m)}_{1,i}r-16 \alpha^{(m)}_{2,i}r^3\alpha^{(m)}_{5,i}\alpha^{(m)}_{1,i}-16 \alpha^{(m)}_{0,i}r^2\alpha^{(m)}_{6,i}\alpha^{(m)}_{3,i}-64 \alpha^{(m)}_{4,i}\alpha^{(m)}_{1,i}\alpha^{(m)}_{5,i}r^2+\\&&
\nonumber 4 \alpha^{(m)}_{0,i}\alpha^{(m)}_{2,i}r^2 \alpha^{(m)}_{6,i}+16 r^3\alpha^{(m)}_{4,i}\alpha^{(m)}_{3,i}\alpha^{(m)}_{1,i}-16 \alpha^{(m)}_{2,i}\alpha^{(m)}_{4,i}\alpha^{(m)}_{3,i}r^2+32 \alpha^{(m)}_{2,i}\alpha^{(m)}_{5,i}\\&&
\nonumber \times \alpha^{(m)}_{3,i}r+4 \alpha^{(m)}_{0,i}\alpha^{(m)}_{4,i}r^3\alpha^{(m)}_{3,i}-64 (\alpha^{(m)}_{5,i})^2\alpha^{(m)}_{2,i}-64 \alpha^{(m)}_{6,i}(\alpha^{(m)}_{3,i})^2+(\alpha^{(m)}_{2,i})^2 \\&&
\nonumber \times \alpha^{(m)}_{5,i}+16 \alpha^{(m)}_{1,i}(\alpha^{(m)}_{5,i})^2+4 \alpha^{(m)}_{2,i}(\alpha^{(m)}_{4,i})^2+(\alpha^{(m)}_{4,i})^2\alpha^{(m)}_{1,i}+4 r(\alpha^{(m)}_{3,i})^3+16 \\&&
\nonumber \times (\alpha^{(m)}_{4,i})^2\alpha^{(m)}_{3,i}+16 (\alpha^{(m)}_{3,i})^3r^2+64 r^3(\alpha^{(m)}_{3,i})^3-16 \alpha^{(m)}_{0,i}r(\alpha^{(m)}_{5,i})^2+64 \alpha^{(m)}_{1,i}\\&&
\times r(\alpha^{(m)}_{5,i})^2- r^3 (\alpha^{(m)}_{1,i})^2\alpha^{(m)}_{4,i}- \alpha^{(m)}_{2,i}(\alpha^{(m)}_{3,i})^2r+4 (\alpha^{(m)}_{2,i})^2r^3\alpha^{(m)}_{3,i},
\end{eqnarray}

\begin{eqnarray}\label{p-2-n}
\nonumber p^{(m)}_{2,r,n,i}&=&-64 r^2(\alpha^{(m)}_{3,i})^2\alpha^{(m)}_{4,i}+128 \alpha^{(m)}_{5,i}\alpha^{(m)}_{3,i}\alpha^{(m)}_{4,i}+64 \alpha^{(m)}_{2,i}(\alpha^{(m)}_{4,i})^2r^2- \alpha^{(m)}_{0,i}(\alpha^{(m)}_{4,i})^2\\&&
\nonumber \times r+4 \alpha^{(m)}_{2,i}\alpha^{(m)}_{3,i}\alpha^{(m)}_{5,i}+16 \alpha^{(m)}_{2,i}\alpha^{(m)}_{6,i}\alpha^{(m)}_{3,i}- (\alpha^{(m)}_{3,i})^2\alpha^{(m)}_{1,i}r^2-4 r^2\alpha^{(m)}_{6,i}\times \\&&
\nonumber (\alpha^{(m)}_{1,i})^2+4 \alpha^{(m)}_{1,i}\alpha^{(m)}_{3,i}\alpha^{(m)}_{6,i}+16 (\alpha^{(m)}_{2,i})^2r^3\alpha^{(m)}_{4,i}+4 r^3(\alpha^{(m)}_{1,i})^2\alpha^{(m)}_{5,i}-4 \alpha^{(m)}_{2,i}r^2\\&&
\nonumber \times (\alpha^{(m)}_{3,i})^2-16 \alpha^{(m)}_{2,i}\alpha^{(m)}_{4,i}\alpha^{(m)}_{5,i}-16 \alpha^{(m)}_{0,i}r^3(\alpha^{(m)}_{4,i})^2-4 \alpha^{(m)}_{0,i}(\alpha^{(m)}_{4,i})^2r^2-2 \times \\&&
\nonumber \alpha^{(m)}_{2,i}\alpha^{(m)}_{4,i}\alpha^{(m)}_{3,i}+(\alpha^{(m)}_{2,i})^2\alpha^{(m)}_{4,i}r-64 (\alpha^{(m)}_{2,i})^2r^2\alpha^{(m)}_{6,i}-16 \alpha^{(m)}_{2,i}r^3(\alpha^{(m)}_{3,i})^2+64 \\&&
\times r\alpha^{(m)}_{3,i}(\alpha^{(m)}_{4,i})^2- r^3(\alpha^{(m)}_{2,i})^3-4 \alpha^{(m)}_{4,i}(\alpha^{(m)}_{3,i})^2+(\alpha^{(m)}_{3,i})^3-64 (\alpha^{(m)}_{4,i})^3,
\end{eqnarray}

\begin{eqnarray}\label{p-3-n}
\nonumber p^{(m)}_{3,r,n,i}&=&-4 r\alpha^{(m)}_{3,i}\alpha^{(m)}_{5,i}\alpha^{(m)}_{1,i}+16 \alpha^{(m)}_{0,i}\alpha^{(m)}_{4,i}r\alpha^{(m)}_{6,i}+2 \alpha^{(m)}_{2,i}r^3\alpha^{(m)}_{1,i}\alpha^{(m)}_{3,i}+64 \alpha^{(m)}_{1,i}r^2\times \\&&
\nonumber \alpha^{(m)}_{6,i}\alpha^{(m)}_{3,i}- \alpha^{(m)}_{0,i}\alpha^{(m)}_{2,i}\alpha^{(m)}_{5,i}r^2-16 \alpha^{(m)}_{5,i}r^2 \alpha^{(m)}_{3,i}\alpha^{(m)}_{1,i}+\alpha^{(m)}_{4,i}\alpha^{(m)}_{1,i}\alpha^{(m)}_{3,i}r-64 \\&&
\nonumber \times r\alpha^{(m)}_{1,i}\alpha^{(m)}_{6,i}\alpha^{(m)}_{4,i}-4 \alpha^{(m)}_{0,i}r\alpha^{(m)}_{3,i}\alpha^{(m)}_{6,i}-4 \alpha^{(m)}_{0,i} \alpha^{(m)}_{2,i}r^3\alpha^{(m)}_{5,i}+8 \alpha^{(m)}_{4,i}r^2\alpha^{(m)}_{3,i}\\&&
\times \alpha^{(m)}_{1,i}-128 \alpha^{(m)}_{2,i}\alpha^{(m)}_{4,i}r^3\alpha^{(m)}_{3,i}+\alpha^{(m)}_{0,i}\alpha^{(m)}_{5,i}\alpha^{(m)}_{3,i}r,
\end{eqnarray}

\begin{eqnarray}\label{p-4-n}
\nonumber p^{(m)}_{4,r,n,i}&=&48 (\alpha^{(m)}_{4,i})^{2} \alpha^{(m)}_{3,i}+64 r^{3} (\alpha^{(m)}_{3,i})^{3}-64 (\alpha^{(m)}_{5,i})^{2} \alpha^{(m)}_{2,i}+36 \alpha^{(m)}_{2,i} (\alpha^{(m)}_{4,i})^{2}+36 r \times \\&&
\nonumber(\alpha^{(m)}_{3,i})^{3}-48 \alpha^{(m)}_{5,i} r^{2} \alpha^{(m)}_{1,i} \alpha^{(m)}_{3,i}+48 \alpha^{(m)}_{0,i} r^{3} \alpha^{(m)}_{5,i} \alpha^{(m)}_{3,i}-36 \alpha^{(m)}_{0,i} r \alpha^{(m)}_{3,i} \alpha^{(m)}_{6,i}-\\&&
\nonumber64 \alpha^{(m)}_{4,i} \alpha^{(m)}_{1,i} \alpha^{(m)}_{5,i} r^{2}+64 \alpha^{(m)}_{1,i} r^{2} \alpha^{(m)}_{6,i} \alpha^{(m)}_{3,i}-36 \alpha^{(m)}_{0,i} \alpha^{(m)}_{2,i} r^{3} \alpha^{(m)}_{5,i}+36 \alpha^{(m)}_{0,i} \times \\&&
\nonumber\alpha^{(m)}_{4,i} r^{3} \alpha^{(m)}_{3,i}+48 \alpha^{(m)}_{0,i} \alpha^{(m)}_{4,i} r \alpha^{(m)}_{6,i}+48 \alpha^{(m)}_{0,i} \alpha^{(m)}_{5,i} r^{2} \alpha^{(m)}_{4,i}+36 \alpha^{(m)}_{0,i} \alpha^{(m)}_{4,i} r \alpha^{(m)}_{5,i}\\&&
\nonumber+96 \alpha^{(m)}_{2,i} \alpha^{(m)}_{5,i} \alpha^{(m)}_{3,i} r+27 \alpha^{(m)}_{0,i} \alpha^{(m)}_{2,i} r^{3} \alpha^{(m)}_{4,i}-27 \alpha^{(m)}_{2,i} \alpha^{(m)}_{5,i} \alpha^{(m)}_{1,i} r+27 \alpha^{(m)}_{0,i}\times\\&&
\nonumber\alpha^{(m)}_{4,i} \alpha^{(m)}_{3,i} r^{2}+27 \alpha^{(m)}_{0,i} \alpha^{(m)}_{5,i} \alpha^{(m)}_{3,i} r-27 \alpha^{(m)}_{0,i} \alpha^{(m)}_{2,i} \alpha^{(m)}_{5,i} r^{2}-64 r \alpha^{(m)}_{1,i} \alpha^{(m)}_{6,i} \alpha^{(m)}_{4,i},
\end{eqnarray}

\begin{eqnarray*}
\nonumber p^{(m)}_{5,r,n,i}&=&-48 \alpha^{(m)}_{2,i} \alpha^{(m)}_{4,i} \alpha^{(m)}_{3,i} r^{2}+72 \alpha^{(m)}_{4,i} r^{2} \alpha^{(m)}_{3,i} \alpha^{(m)}_{1,i}+54 \alpha^{(m)}_{2,i} r^{3} \alpha^{(m)}_{1,i} \alpha^{(m)}_{3,i}+27\\&& \times
\nonumber \alpha^{(m)}_{4,i}\alpha^{(m)}_{1,i} \alpha^{(m)}_{3,i} r-27 \alpha^{(m)}_{2,i} \alpha^{(m)}_{4,i} \alpha^{(m)}_{1,i} r^{2}+48 r^{3} \alpha^{(m)}_{4,i} \alpha^{(m)}_{3,i} \alpha^{(m)}_{1,i}-36 r \times\\&&
\nonumber\alpha^{(m)}_{3,i} \alpha^{(m)}_{5,i}\alpha^{(m)}_{1,i}+64 \alpha^{(m)}_{2,i} \alpha^{(m)}_{3,i} \alpha^{(m)}_{5,i} r^{2}-36 \alpha^{(m)}_{2,i} \alpha^{(m)}_{4,i} r \alpha^{(m)}_{3,i}-36 \alpha^{(m)}_{2,i}\alpha^{(m)}_{4,i} \\&& \times
\nonumber r^{3} \alpha^{(m)}_{1,i}+36\alpha^{(m)}_{2,i} r \alpha^{(m)}_{1,i} \alpha^{(m)}_{6,i}-48 \alpha^{(m)}_{2,i} r^{3} \alpha^{(m)}_{5,i} \alpha^{(m)}_{1,i}-128 \alpha^{(m)}_{2,i} \alpha^{(m)}_{4,i}\times\\&&
\nonumber r^{3} \alpha^{(m)}_{3,i}-64 \alpha^{(m)}_{2,i}\alpha^{(m)}_{4,i} r \alpha^{(m)}_{5,i}-64 \alpha^{(m)}_{1,i} r^{3} \alpha^{(m)}_{5,i} \alpha^{(m)}_{3,i}+48 \alpha^{(m)}_{2,i} r^{2} \alpha^{(m)}_{6,i} \alpha^{(m)}_{1,i}\\&&
\nonumber+64 \alpha^{(m)}_{2,i} r \alpha^{(m)}_{3,i} \alpha^{(m)}_{6,i}-48 \alpha^{(m)}_{0,i} r^{2} \alpha^{(m)}_{6,i} \alpha^{(m)}_{3,i}+36 \alpha^{(m)}_{0,i} \alpha^{(m)}_{2,i} r^{2} \alpha^{(m)}_{6,i}-27 r^{3}\\&&
\nonumber \times (\alpha^{(m)}_{2,i})^{3}+27 (\alpha^{(m)}_{4,i})^{2} \alpha^{(m)}_{1,i}+128 \alpha^{(m)}_{5,i} \alpha^{(m)}_{3,i} \alpha^{(m)}_{4,i}+27 \alpha^{(m)}_{5,i} r^{2} (\alpha^{(m)}_{1,i})^{2} +\\&&
\nonumber64 \alpha^{(m)}_{2,i} \alpha^{(m)}_{6,i} \alpha^{(m)}_{4,i}-48 \alpha^{(m)}_{4,i} r(\alpha^{(m)}_{3,i})^{2}+48 \alpha^{(m)}_{2,i} \alpha^{(m)}_{6,i} \alpha^{(m)}_{3,i}-27 (\alpha^{(m)}_{3,i})^{2} \times \\&&
\alpha^{(m)}_{1,i} r^{2}-36 r^{2} \alpha^{(m)}_{6,i} (\alpha^{(m)}_{1,i})^{2},
\end{eqnarray*}

\begin{eqnarray*}
\nonumber p^{(m)}_{6,r,n,i}&=&-36 \times \alpha^{(m)}_{0,i} (\alpha^{(m)}_{4,i})^{2} r^{2}-36 \alpha^{(m)}_{2,i} r^{2} (\alpha^{(m)}_{3,i})^{2}-36 \alpha^{(m)}_{4,i} (\alpha^{(m)}_{3,i})^{2}-48 \alpha^{(m)}_{5,i}\\&&
\times (\alpha^{(m)}_{3,i})^{2}-64 \alpha^{(m)}_{6,i} (\alpha^{(m)}_{3,i})^{2}+64 \alpha^{(m)}_{1,i} r (\alpha^{(m)}_{5,i})^{2}+36 (\alpha^{(m)}_{2,i})^{2} r^{3} \alpha^{(m)}_{3,i}-\\&&
64 r^{2} (\alpha^{(m)}_{3,i})^{2} \alpha^{(m)}_{4,i}+ 64 r \alpha^{(m)}_{3,i} (\alpha^{(m)}_{4,i})^{2}+36 r^{3} \alpha^{(m)}_{5,i} (\alpha^{(m)}_{1,i})^{2}-64 (\alpha^{(m)}_{2,i})^{2} r^{2}\\&&
\times \alpha^{(m)}_{6,i}+27 (\alpha^{(m)}_{2,i})^{2} \alpha^{(m)}_{3,i} r^{2} -48 \alpha^{(m)}_{4,i} \alpha^{(m)}_{1,i} \alpha^{(m)}_{6,i}-36 (\alpha^{(m)}_{3,i})^{2} r^{3} \alpha^{(m)}_{1,i}-\\&&
54 \alpha^{(m)}_{2,i} \alpha^{(m)}_{4,i} \alpha^{(m)}_{3,i}+36 \alpha^{(m)}_{2,i} \alpha^{(m)}_{5,i} \alpha^{(m)}_{3,i} +64 r^{3} (\alpha^{(m)}_{4,i})^{2} \alpha^{(m)}_{1,i}+27 (\alpha^{(m)}_{2,i})^{2}\\&&
\times \alpha^{(m)}_{5,i}-64 r (\alpha^{(m)}_{3,i})^{2} \alpha^{(m)}_{5,i}-48 \alpha^{(m)}_{2,i} r^{3} (\alpha^{(m)}_{3,i})^{2} +48 (\alpha^{(m)}_{2,i})^{2} r^{3} \alpha^{(m)}_{4,i},
\end{eqnarray*}

\begin{eqnarray*}
\nonumber p^{(m)}_{7,r,n,i}&=&64 \alpha^{(m)}_{2,i} (\alpha^{(m)}_{4,i})^{2} r^{2}-27 \alpha^{(m)}_{2,i} (\alpha^{(m)}_{3,i})^{2} r-27 r^{3} (\alpha^{(m)}_{1,i})^{2} \alpha^{(m)}_{4,i}-27 \alpha^{(m)}_{0,i} r^{3}\times \\&&
(\alpha^{(m)}_{3,i})^{2}-27 \alpha^{(m)}_{0,i} (\alpha^{(m)}_{4,i})^{2} r+27 (\alpha^{(m)}_{2,i})^{2} \alpha^{(m)}_{4,i} r+36 \alpha^{(m)}_{1,i} \alpha^{(m)}_{3,i} \alpha^{(m)}_{6,i}-48\\&&
\times (\alpha^{(m)}_{2,i})^{2} r \alpha^{(m)}_{6,i}+64 (\alpha^{(m)}_{2,i})^{2} r^{3} \alpha^{(m)}_{5,i}-48 \alpha^{(m)}_{0,i} r^{3} (\alpha^{(m)}_{4,i})^{2}-36 (\alpha^{(m)}_{2,i})^{2} \times \\&&
\alpha^{(m)}_{6,i}-64 (\alpha^{(m)}_{4,i})^{3}+27 (\alpha^{(m)}_{3,i})^{3}-48 \alpha^{(m)}_{0,i} r (\alpha^{(m)}_{5,i})^{2}-36 \alpha^{(m)}_{4,i} \alpha^{(m)}_{1,i} \alpha^{(m)}_{5,i}\\&&
-48\alpha^{(m)}_{2,i} \alpha^{(m)}_{4,i} \alpha^{(m)}_{5,i}-27 \alpha^{(m)}_{5,i} \alpha^{(m)}_{3,i} \alpha^{(m)}_{1,i}+48 (\alpha^{(m)}_{3,i})^{3} r^{2}+48 \alpha^{(m)}_{1,i} (\alpha^{(m)}_{5,i})^{2},
\end{eqnarray*}
for $k\geqslant 0$, $m\geqslant 0$, $\delta \geqslant 0$,
\begin{eqnarray}\label{weights2}
w^{(m)}_{i+r}&=&\left[\left(f^{k}_{i+r}-\sum\limits_{a=0}^{3}r^{a}\beta^{(m)}_{a,\delta,i}\right)^{2}+\delta\right]^{-\frac{1}{2}},
\end{eqnarray}
\begin{eqnarray}\label{weights6}
\alpha^{(m)}_{b,i}=\sum\limits_{r=-n+1}^{n}r^{b}w^{(m)}_{i+r}\,\ \mbox{for} \,\ 0 \leqslant b \leqslant 6,
\end{eqnarray}
and $\beta^{(m)}_{a,\delta,i}$ are defined in Propositions \ref{2n-betas-OLS-cor}-\ref{betas-IRLS-cor}.

\subsection{Family of $(2n+1)$-point schemes based on polynomial of degree 3}
We have introduced a family of schemes $D_{2n,3}$ in (\ref{2n-cubic-scheme}) that is based on the fitting of a cubic univariate polynomial to $2n$-observations in two-dimensional space. We can get an alternative of the schemes $D_{2n,3}$ by fitting polynomial function of degree 3 to $(2n+1)$-observations. This alternative can be obtained by replacing $r=-n+1,\ldots,n$ to $r=-n,\ldots,n$ in (\ref{IRLS-betas-3-2n})-(\ref{weights4}) and (\ref{2n-cubic-scheme})-(\ref{weights6}), $\frac{1}{4}$ by $-\frac{1}{4}$ and $\frac{3}{4}$ by $\frac{1}{4}$ in (\ref{2n-cubic-scheme}). Ultimately, we get the following family of schemes denoted by $D_{2n+1,3}$
\begin{eqnarray}\label{2n+1-cubic-scheme}
\left\{\begin{array}{ccccccc}
f^{k+1}_{2i}&=&\frac{1}{64\gamma^{(m)}_{n,i}}\sum\limits_{r=-n}^{n}\left(p^{(m)}_{8,r,n,i}+p^{(m)}_{9,r,n,i}+p^{(m)}_{10,r,n,i}+p^{(m)}_{11,r,n,i}\right)w^{(m)}_{i+r}f^{k}_{i+r},\\ \\
f^{k+1}_{2i+1}&=&\frac{1}{64\gamma^{(m)}_{n,i}}\sum\limits_{r=-n}^{n}\left(p^{(m)}_{0,r,n,i}+p^{(m)}_{1,r,n,i}+p^{(m)}_{2,r,n,i}+p^{(m)}_{3,r,n,i}\right)w^{(m)}_{i+r}f^{k}_{i+r},
\end{array}\right.
\end{eqnarray}
where $\gamma^{(m)}_{n,i}$, $p^{(m)}_{0,r,n,i}$, $p^{(m)}_{1,r,n,i}$, $p^{(m)}_{2,r,n,i}$ and $p^{(m)}_{3,r,n,i}$ are defined in (\ref{gamma-n})-(\ref{p-3-n}) and

\begin{eqnarray*}\label{p-8-n}
\nonumber p^{(m)}_{8,r,n,i}&=&64 \alpha^{(m)}_{2,i} (\alpha^{(m)}_{4,i})^2 r^2+16 \alpha^{(m)}_{2,i} \alpha^{(m)}_{4,i} \alpha^{(m)}_{5,i}+16 \alpha^{(m)}_{0,i} r^3 (\alpha^{(m)}_{4,i})^2+r^3 (\alpha^{(m)}_{1,i})^2 \\&&
\times \alpha^{(m)}_{4,i}+4 \alpha^{(m)}_{2,i} \alpha^{(m)}_{5,i} \alpha^{(m)}_{3,i}+16 \alpha^{(m)}_{4,i} \alpha^{(m)}_{1,i} \alpha^{(m)}_{6,i}+\alpha^{(m)}_{2,i} (\alpha^{(m)}_{3,i})^2 r+16 \times \\&&
\alpha^{(m)}_{2,i} r^3 (\alpha^{(m)}_{3,i})^2+64 r \alpha^{(m)}_{3,i} (\alpha^{(m)}_{4,i})^2+64 r^3 (\alpha^{(m)}_{4,i})^2 \alpha^{(m)}_{1,i}-16 \alpha^{(m)}_{2,i} \alpha^{(m)}_{6,i}\\&&
\times \alpha^{(m)}_{3,i}+64 \alpha^{(m)}_{2,i} \alpha^{(m)}_{6,i}\alpha^{(m)}_{4,i}+(\alpha^{(m)}_{3,i})^2 \alpha^{(m)}_{1,i} r^2+128 \alpha^{(m)}_{5,i} \alpha^{(m)}_{3,i} \alpha^{(m)}_{4,i}+\\&&
64(\alpha^{(m)}_{2,i})^2 r^3 \alpha^{(m)}_{5,i}-\alpha^{(m)}_{5,i} r^2 (\alpha^{(m)}_{1,i})^2+16 \alpha^{(m)}_{4,i} r (\alpha^{(m)}_{3,i})^2-64 r (\alpha^{(m)}_{3,i})^2\\&&
\times \alpha^{(m)}_{5,i}-4 \alpha^{(m)}_{0,i} (\alpha^{(m)}_{4,i})^2 r^2+4 r^3 (\alpha^{(m)}_{1,i})^2 \alpha^{(m)}_{5,i}-16 (\alpha^{(m)}_{2,i})^2 r^3 \alpha^{(m)}_{4,i},
\end{eqnarray*}

\begin{eqnarray*}
p^{(m)}_{9,r,n,i}&=&-(\alpha^{(m)}_{2,i})^2 \alpha^{(m)}_{4,i} r+\alpha^{(m)}_{0,i} r^3 (\alpha^{(m)}_{3,i})^2-4 r^2 \alpha^{(m)}_{6,i} (\alpha^{(m)}_{1,i})^2-(\alpha^{(m)}_{2,i})^2 \alpha^{(m)}_{3,i} r^2\\&&
+4 \alpha^{(m)}_{1,i} \alpha^{(m)}_{3,i} \alpha^{(m)}_{6,i}+16 (\alpha^{(m)}_{2,i})^2 r \alpha^{(m)}_{6,i}-4 \alpha^{(m)}_{2,i} r^2 (\alpha^{(m)}_{3,i})^2-64 r^2 (\alpha^{(m)}_{3,i})^2 \\&&
\times \alpha^{(m)}_{4,i}-16 \alpha^{(m)}_{0,i} r^3 \alpha^{(m)}_{5,i} \alpha^{(m)}_{3,i}+2 \alpha^{(m)}_{2,i} \alpha^{(m)}_{4,i} \alpha^{(m)}_{3,i}+16 \alpha^{(m)}_{0,i} r (\alpha^{(m)}_{5,i})^2-\\&&
4 (\alpha^{(m)}_{3,i})^2 r^3 \alpha^{(m)}_{1,i}-4 \alpha^{(m)}_{0,i} \alpha^{(m)}_{2,i} r^3 \alpha^{(m)}_{5,i}-64 \alpha^{(m)}_{2,i} \alpha^{(m)}_{4,i} r \alpha^{(m)}_{5,i}+4 (\alpha^{(m)}_{2,i})^2\\&&
\times r^3 \alpha^{(m)}_{3,i}-(\alpha^{(m)}_{4,i})^2 \alpha^{(m)}_{1,i}-4 \alpha^{(m)}_{4,i} (\alpha^{(m)}_{3,i})^2-16 \alpha^{(m)}_{1,i} (\alpha^{(m)}_{5,i})^2-64 \alpha^{(m)}_{6,i} \\&&
\times(\alpha^{(m)}_{3,i})^2-64 (\alpha^{(m)}_{5,i})^2 \alpha^{(m)}_{2,i}-16 (\alpha^{(m)}_{4,i})^2 \alpha^{(m)}_{3,i}+r^3 (\alpha^{(m)}_{2,i})^3-(\alpha^{(m)}_{2,i})^2 \\&&
\times \alpha^{(m)}_{5,i}-4 (\alpha^{(m)}_{2,i})^2 \alpha^{(m)}_{6,i}+4 r (\alpha^{(m)}_{3,i})^3+64 r^3 (\alpha^{(m)}_{3,i})^3+16 \alpha^{(m)}_{5,i} (\alpha^{(m)}_{3,i})^2\\&&
-\emph{}16 (\alpha^{(m)}_{3,i})^3 r^2+4 \alpha^{(m)}_{2,i} (\alpha^{(m)}_{4,i})^2-4 \alpha^{(m)}_{4,i} \alpha^{(m)}_{1,i} \alpha^{(m)}_{5,i}+\alpha^{(m)}_{0,i} (\alpha^{(m)}_{4,i})^2 r,
\end{eqnarray*}

\begin{eqnarray*}
p^{(m)}_{10,r,n,i}&=&-16 r^3 \alpha^{(m)}_{4,i} \alpha^{(m)}_{3,i} \alpha^{(m)}_{1,i}-64 r \alpha^{(m)}_{1,i} \alpha^{(m)}_{6,i} \alpha^{(m)}_{4,i}-\alpha^{(m)}_{0,i} \alpha^{(m)}_{2,i} r^3 \alpha^{(m)}_{4,i} -64 \times \\&&
(\alpha^{(m)}_{2,i})^2 r^2 \alpha^{(m)}_{6,i}+16 \alpha^{(m)}_{5,i} r^2 \alpha^{(m)}_{3,i} \alpha^{(m)}_{1,i}-64 \alpha^{(m)}_{4,i} \alpha^{(m)}_{1,i} \alpha^{(m)}_{5,i} r^2+8 \alpha^{(m)}_{4,i} r^2 \\&&
\nonumber\times \alpha^{(m)}_{3,i} \alpha^{(m)}_{1,i}+\alpha^{(m)}_{0,i} \alpha^{(m)}_{2,i} \alpha^{(m)}_{5,i} r^2-16 \alpha^{(m)}_{0,i} \alpha^{(m)}_{5,i} r^2 \alpha^{(m)}_{4,i}-4 r \alpha^{(m)}_{3,i} \alpha^{(m)}_{5,i} \times \\&&
\alpha^{(m)}_{1,i}-64 \alpha^{(m)}_{1,i} r^3 \alpha^{(m)}_{5,i} \alpha^{(m)}_{3,i}+64 \alpha^{(m)}_{1,i} r^2 \alpha^{(m)}_{6,i} \alpha^{(m)}_{3,i}+16 \alpha^{(m)}_{2,i} \alpha^{(m)}_{4,i} \alpha^{(m)}_{3,i} r^2 \\&&
+\alpha^{(m)}_{2,i} \alpha^{(m)}_{5,i} \alpha^{(m)}_{1,i} r+\alpha^{(m)}_{2,i} \alpha^{(m)}_{4,i} \alpha^{(m)}_{1,i} r^2+4 \alpha^{(m)}_{0,i} \alpha^{(m)}_{2,i} r^2 \alpha^{(m)}_{6,i}-\alpha^{(m)}_{4,i} \alpha^{(m)}_{1,i}\\&&
\times \alpha^{(m)}_{3,i} r+\alpha^{(m)}_{5,i} \alpha^{(m)}_{3,i} \alpha^{(m)}_{1,i}-2 \alpha^{(m)}_{2,i} r^3 \alpha^{(m)}_{1,i} \alpha^{(m)}_{3,i}+4 \alpha^{(m)}_{0,i} \alpha^{(m)}_{4,i} r^3 \alpha^{(m)}_{3,i},
\end{eqnarray*}

\begin{eqnarray*}
p^{(m)}_{11,r,n,i}&=&-4 \alpha^{(m)}_{0,i} r \alpha^{(m)}_{3,i} \alpha^{(m)}_{6,i}-(\alpha^{(m)}_{3,i})^3-64 (\alpha^{(m)}_{4,i})^3-\alpha^{(m)}_{0,i} \alpha^{(m)}_{5,i} \alpha^{(m)}_{3,i} r-16 \alpha^{(m)}_{2,i}\\&&
r^2 \alpha^{(m)}_{6,i} \alpha^{(m)}_{1,i}-\alpha^{(m)}_{0,i} \alpha^{(m)}_{4,i}\alpha^{(m)}_{3,i} r^2+16 \alpha^{(m)}_{2,i} r^3 \alpha^{(m)}_{5,i} \alpha^{(m)}_{1,i}-32 \alpha^{(m)}_{2,i} \alpha^{(m)}_{5,i} \alpha^{(m)}_{3,i}\\&&
\times r+64 \alpha^{(m)}_{2,i} r \alpha^{(m)}_{3,i} \alpha^{(m)}_{6,i} +4 \alpha^{(m)}_{2,i} r \alpha^{(m)}_{1,i} \alpha^{(m)}_{6,i}-128 \alpha^{(m)}_{2,i} \alpha^{(m)}_{4,i} r^3 \alpha^{(m)}_{3,i}+64 \\&&
\times \alpha^{(m)}_{2,i} \alpha^{(m)}_{3,i} \alpha^{(m)}_{5,i} r^2+16 \alpha^{(m)}_{0,i} r^2 \alpha^{(m)}_{6,i} \alpha^{(m)}_{3,i}-16 \alpha^{(m)}_{0,i} \alpha^{(m)}_{4,i} r \alpha^{(m)}_{6,i}+4 \alpha^{(m)}_{0,i} \times \\&&
\alpha^{(m)}_{4,i} r \alpha^{(m)}_{5,i}+64 \alpha^{(m)}_{1,i} r (\alpha^{(m)}_{5,i})^2-4 \alpha^{(m)}_{2,i} \alpha^{(m)}_{4,i} r^3 \alpha^{(m)}_{1,i}-4 \alpha^{(m)}_{2,i} \alpha^{(m)}_{4,i} r \alpha^{(m)}_{3,i}.
\end{eqnarray*}
The starting values $\beta^{(0)}_{a,\delta,i}$ are defined in Proposition \ref{2n+1-betas-OLS-cor}.

%
%
%

Now we introduce two more families of schemes that are based on fitting polynomial of degree 2.
\subsection{Family of $2n$-point schemes based on polynomial of degree 2}
From Propositions \ref{2n-betas-OLS-cor}-\ref{betas-IRLS-cor} with $\beta^{(m+1)}_{3,\delta,i}=0$ and Algorithm \ref{SS-A} with $d=2$, we get the following family of $2n$-points schemes $D_{2n,2}$ based on fitting polynomial of degree 2.
\begin{eqnarray}\label{2n-quadratic-scheme}
\left\{\begin{array}{ccccccc}
f^{k+1}_{2i}&=&\frac{\alpha^{(m)}_{0,i}}{16 \lambda^{(m)}_{0,i}}\sum\limits_{r=-n+1}^{n}q^{(m)}_{1,r,n,i}w^{(m)}_{i+r}f^{k}_{i+r},\\ \\
f^{k+1}_{2i+1}&=&\frac{\alpha^{(m)}_{0,i}}{16 \lambda^{(m)}_{0,i}}\sum\limits_{r=-n+1}^{n}q^{(m)}_{2,r,n,i}w^{(m)}_{i+r}f^{k}_{i+r},
\end{array}\right.
\end{eqnarray}
where $\lambda^{(m)}_{0,i}$ is defined in (\ref{weights1}) and
\begin{eqnarray}\label{q-1-n}
\nonumber q^{(m)}_{1,r,n,i}&=&r^{2} \alpha^{(m)}_{0,i} \alpha^{(m)}_{2,i}+4 \alpha^{(m)}_{0,i} r \alpha^{(m)}_{4,i}-4 \alpha^{(m)}_{0,i} \alpha^{(m)}_{3,i} r^{2}-\alpha^{(m)}_{0,i} \alpha^{(m)}_{3,i} r-(\alpha^{(m)}_{2,i})^{2}-4 \\&&
\nonumber \times (\alpha^{(m)}_{2,i})^{2} r-16 r^{2} (\alpha^{(m)}_{2,i})^{2}+\alpha^{(m)}_{1,i} \alpha^{(m)}_{2,i} r+4 \alpha^{(m)}_{2,i} \alpha^{(m)}_{3,i}+4 \alpha^{(m)}_{2,i} \alpha^{(m)}_{1,i} r^{2}+\\&&
\nonumber 16 \alpha^{(m)}_{2,i} \alpha^{(m)}_{3,i} r+16 \alpha^{(m)}_{4,i} \alpha^{(m)}_{2,i}-16 \alpha^{(m)}_{1,i} \alpha^{(m)}_{4,i} r+\alpha^{(m)}_{1,i} \alpha^{(m)}_{3,i}-16 (\alpha^{(m)}_{3,i})^{2}-\\&&
r^{2} (\alpha^{(m)}_{1,i})^{2}+16 r^{2} \alpha^{(m)}_{1,i} \alpha^{(m)}_{3,i}-4 \alpha^{(m)}_{1,i} \alpha^{(m)}_{4,i},
\end{eqnarray}

\begin{eqnarray*}\label{q-2-n}
\nonumber q^{(m)}_{2,r,n,i}&=&9 r^{2} \alpha^{(m)}_{0,i} \alpha^{(m)}_{2,i}+12 \alpha^{(m)}_{0,i} r \alpha^{(m)}_{4,i}-9 \alpha^{(m)}_{0,i} \alpha^{(m)}_{3,i} r-12 \alpha^{(m)}_{0,i} \alpha^{(m)}_{3,i} r^{2}-12 (\alpha^{(m)}_{2,i})^{2}\\&&
\nonumber \times r-9 (\alpha^{(m)}_{2,i})^{2}-16 r^{2} (\alpha^{(m)}_{2,i})^{2}+9 \alpha^{(m)}_{1,i} \alpha^{(m)}_{2,i} r+12 \alpha^{(m)}_{2,i} \alpha^{(m)}_{1,i} r^{2}+16 \alpha^{(m)}_{4,i} \\&&
\nonumber \times \alpha^{(m)}_{2,i}+12 \alpha^{(m)}_{2,i} \alpha^{(m)}_{3,i}+16 \alpha^{(m)}_{2,i} \alpha^{(m)}_{3,i} r+9 \alpha^{(m)}_{1,i} \alpha^{(m)}_{3,i}-16 \alpha^{(m)}_{1,i} \alpha^{(m)}_{4,i} r-12 \\&&
\times \alpha^{(m)}_{1,i} \alpha^{(m)}_{4,i}-9 r^{2} (\alpha^{(m)}_{1,i})^{2}+16 r^{2} \alpha^{(m)}_{1,i} \alpha^{(m)}_{3,i}-16 (\alpha^{(m)}_{3,i})^{2}.
\end{eqnarray*}
\subsection{Family of $(2n+1)$-point schemes based on polynomial of degree 2}
Similarly, if we fit a polynomial function of degree two to $(2n+1)$-observations, we get a further alternative of the schemes $D_{2n,2}$, i.e. replace $r=-n+1,\ldots,n$ by $r=-n,\ldots,n$ in Proposition \ref{betas-IRLS-cor} and (\ref{2n-quadratic-scheme})-(\ref{q-1-n}) and then use Propositions \ref{betas-IRLS-cor}-\ref{2n+1-betas-OLS-cor} for $\beta^{m+1}_{3,\delta,i}=0$ to get the following family of $(2n+1)$-point schemes $D_{2n+1,2}$
\begin{eqnarray}\label{2n+1-quadratic-scheme}
\left\{\begin{array}{ccccccc}
f^{k+1}_{2i}&=&\frac{\alpha^{(m)}_{0,i}}{16 \lambda^{(m)}_{0,i}}\sum\limits_{r=-n}^{n}q^{(m)}_{3,r,n,i}w^{(m)}_{i+r}f^{k}_{i+r},\\ \\
f^{k+1}_{2i+1}&=&\frac{\alpha^{(m)}_{0,i}}{16 \lambda^{(m)}_{0,i}}\sum\limits_{r=-n}^{n}q^{(m)}_{1,r,n,i}w^{(m)}_{i+r}f^{k}_{i+r},
\end{array}\right.
\end{eqnarray}
where $\lambda^{(m)}_{0,i}$ and $q^{(m)}_{1,r,n,i}$ are defined in (\ref{weights1}) and (\ref{q-1-n}) respectively and
\begin{eqnarray*}\label{q-3-n}
\nonumber q^{(m)}_{3,r,n,i}&=&-16 r^{2} (\alpha^{(m)}_{2,i})^{2}+16 \alpha^{(m)}_{2,i} \alpha^{(m)}_{4,i}+16 \alpha^{(m)}_{2,i} \alpha^{(m)}_{3,i} r-16 \alpha^{(m)}_{1,i} \alpha^{(m)}_{4,i} r+16 r^{2} \alpha^{(m)}_{1,i}\\&&
\nonumber \times \alpha^{(m)}_{3,i}-16 (\alpha^{(m)}_{3,i})^{2}-4 \alpha^{(m)}_{0,i} r \alpha^{(m)}_{4,i}+4 \alpha^{(m)}_{0,i} \alpha^{(m)}_{3,i} r^{2}+4 r (\alpha^{(m)}_{2,i})^{2}-4 \alpha^{(m)}_{2,i}\\&&
\nonumber \times \alpha^{(m)}_{3,i}-4 \alpha^{(m)}_{2,i} \alpha^{(m)}_{1,i} r^{2}+4 \alpha^{(m)}_{1,i} \alpha^{(m)}_{4,i}+r^{2} \alpha^{(m)}_{0,i} \alpha^{(m)}_{2,i}-r^{2} (\alpha^{(m)}_{1,i})^{2}-(\alpha^{(m)}_{2,i})^{2}\\&&
-\alpha^{(m)}_{0,i} \alpha^{(m)}_{3,i} r+\alpha^{(m)}_{1,i} \alpha^{(m)}_{3,i}+\alpha^{(m)}_{1,i} \alpha^{(m)}_{2,i} r.
\end{eqnarray*}

\subsection{Families of schemes based on polynomial of degree 1}

If we use Propositions \ref{2n-betas-OLS-cor}-\ref{betas-IRLS-cor} with $\beta^{(m+1)}_{3,\delta,i}=\beta^{(m+1)}_{2,\delta,i}=0$, $d=1$ and Algorithm \ref{SS-A}, we get the family of $2n$-point schemes $D_{2n,1}$ proposed by Mustafa et al. \cite{Mustafa9}. The refinement rules of this family of schemes are

\begin{eqnarray}\label{2n-linear-scheme}
\left\{\begin{array}{ccccccc}
f^{k+1}_{2i}&=&\frac{1}{4 \chi^{(m)}_{0,i}}\sum\limits_{r=-n+1}^{n}t^{(m)}_{1,r,n,i}w^{(m)}_{i+r}f^{k}_{i+r},\\ \\
f^{k+1}_{2i+1}&=&\frac{1}{4 \chi^{(m)}_{0,i}}\sum\limits_{r=-n+1}^{n}t^{(m)}_{2,r,n,i}w^{(m)}_{i+r}f^{k}_{i+r},
\end{array}\right.
\end{eqnarray}
where $\chi^{(m)}_{0,i}$ is defined in (\ref{weights3}) and
\begin{eqnarray}\label{t-1-n}
t^{(m)}_{1,r,n,i}&=&4 \alpha^{(m)}_{2,i}-4 r \alpha^{(m)}_{1,i}+r \alpha^{(m)}_{0,i}-\alpha^{(m)}_{1,i},
\end{eqnarray}

\begin{eqnarray*}\label{t-2-n}
t^{(m)}_{2,r,n,i}&=&4 \alpha^{(m)}_{2,i}-4 r \alpha^{(m)}_{1,i}+3 r \alpha^{(m)}_{0,i}-3 \alpha^{(m)}_{1,i}.
\end{eqnarray*}
%


If we replace $r=-n+1,\ldots,n$ by $r=-n,\ldots,n$ in Proposition \ref{betas-IRLS-cor}, then use Propositions \ref{betas-IRLS-cor}-\ref{2n+1-betas-OLS-cor} for $\beta^{(m+1)}_{3,\delta,i}=\beta^{(m+1)}_{2,\delta,i}=0$ and then apply Algorithm \ref{SS-A} for $d=1$ $\&$ $r=-\frac{1}{4},\frac{1}{4}$, we get the following family of schemes $D_{2n+1,1}$ proposed by Mustafa et al. \cite{Mustafa9}. The family of schemes is
\begin{eqnarray}\label{2n+1-linear-scheme}
\left\{\begin{array}{ccccccc}
f^{k+1}_{2i}&=&\frac{1}{4 \chi^{(m)}_{0,i}}\sum\limits_{r=-n}^{n}t^{(m)}_{3,r,n,i}w^{(m)}_{i+r}f^{k}_{i+r},\\ \\
f^{k+1}_{2i+1}&=&\frac{1}{4 \chi^{(m)}_{0,i}}\sum\limits_{r=-n}^{n}t^{(m)}_{1,r,n,i}w^{(m)}_{i+r}f^{k}_{i+r},
\end{array}\right.
\end{eqnarray}
where $\chi^{(m)}_{0,i}$ and $t^{(m)}_{1,r,n,i}$ are defined in (\ref{weights3}) and (\ref{t-1-n}) respectively and

\begin{eqnarray*}\label{t-3-n}
t^{(m)}_{3,r,n,i}&=&4 \alpha^{(m)}_{2,i}-4 r \alpha^{(m)}_{1,i}-r \alpha^{(m)}_{0,i}+ \alpha^{(m)}_{1,i}.
\end{eqnarray*}

\begin{rem}
If we take constant weights, i.e. $w^{(m)}_{i+r}=1$ for all $i$, $r$ and $m$ in (\ref{2n-cubic-scheme}), (\ref{2n+1-cubic-scheme}), (\ref{2n-quadratic-scheme}), (\ref{2n+1-quadratic-scheme}), (\ref{2n-linear-scheme}) and (\ref{2n+1-linear-scheme}), then for $n \in \mathbb{N} \backslash \{1\}$, these families of schemes reduce to the families of schemes introduced in \cite{Dyn5}.
\end{rem}

\subsection{Basic limit functions}
The basic limit functions of the families of schemes $D_{h,d}$, $h \in \{2n,2n+1:n \in \mathbb{N}\backslash\{1\}\}$, $1 \leqslant d \leqslant 3$ are defined as
\begin{eqnarray*}
\phi_{h,d}&=&D^{\infty}_{h,d}\Delta, \,\ h \in \{2n,2n+1:n \in \mathbb{N}\backslash\{1\}\}, \,\ 1 \leqslant d \leqslant 3,
\end{eqnarray*}
where $\Delta$ is the initial data such that $\Delta=\{(-p,0),(-p+1,0),\ldots,(0,1),\ldots,(p-1,0),(p,0)\}$ for sufficiently large positive integer $p$.

\begin{figure}[tbp] 
 \begin{center}
\begin{tabular}{ccccccccccc}
\epsfig{file=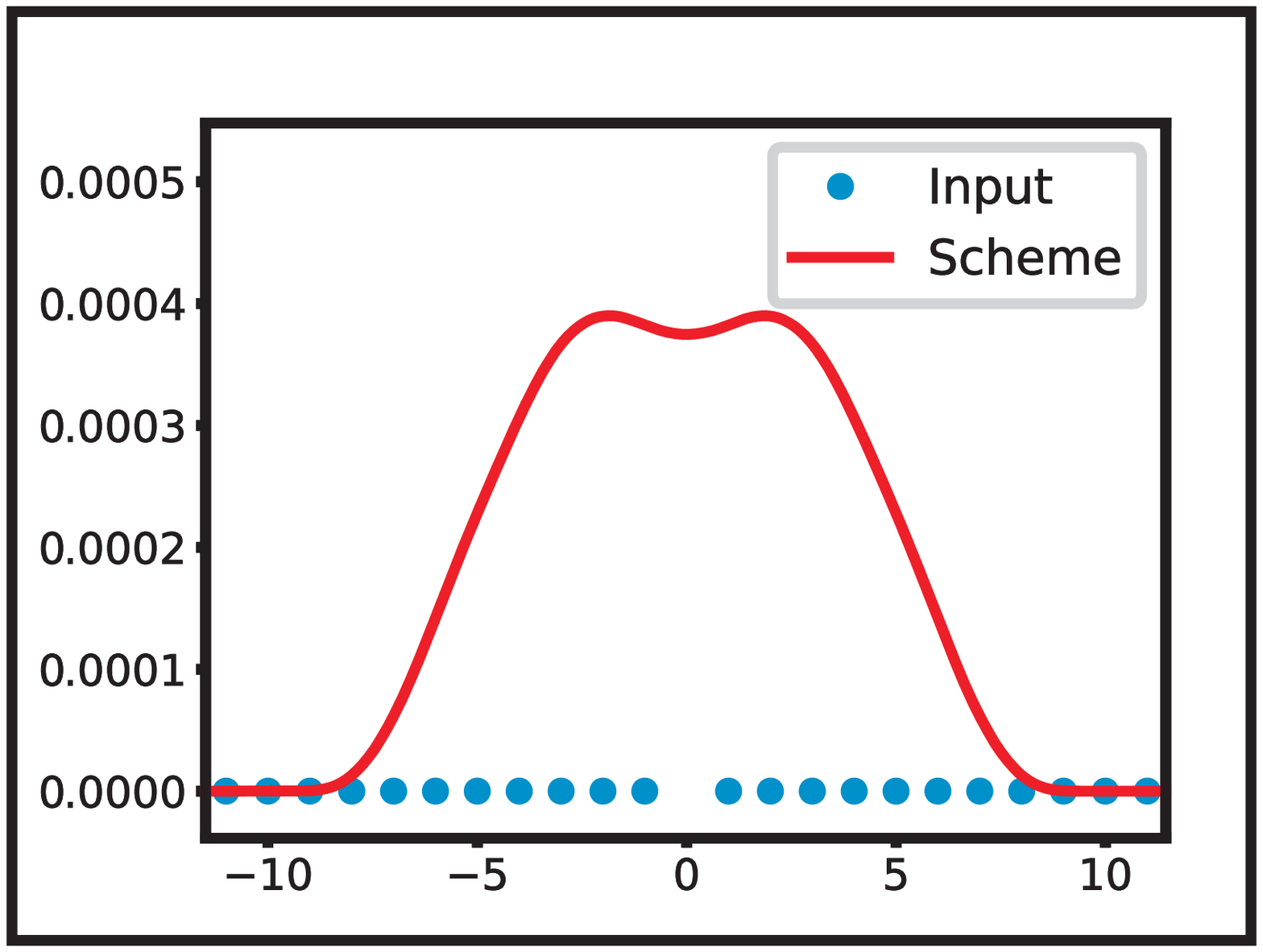, width=1.85 in}&
\epsfig{file=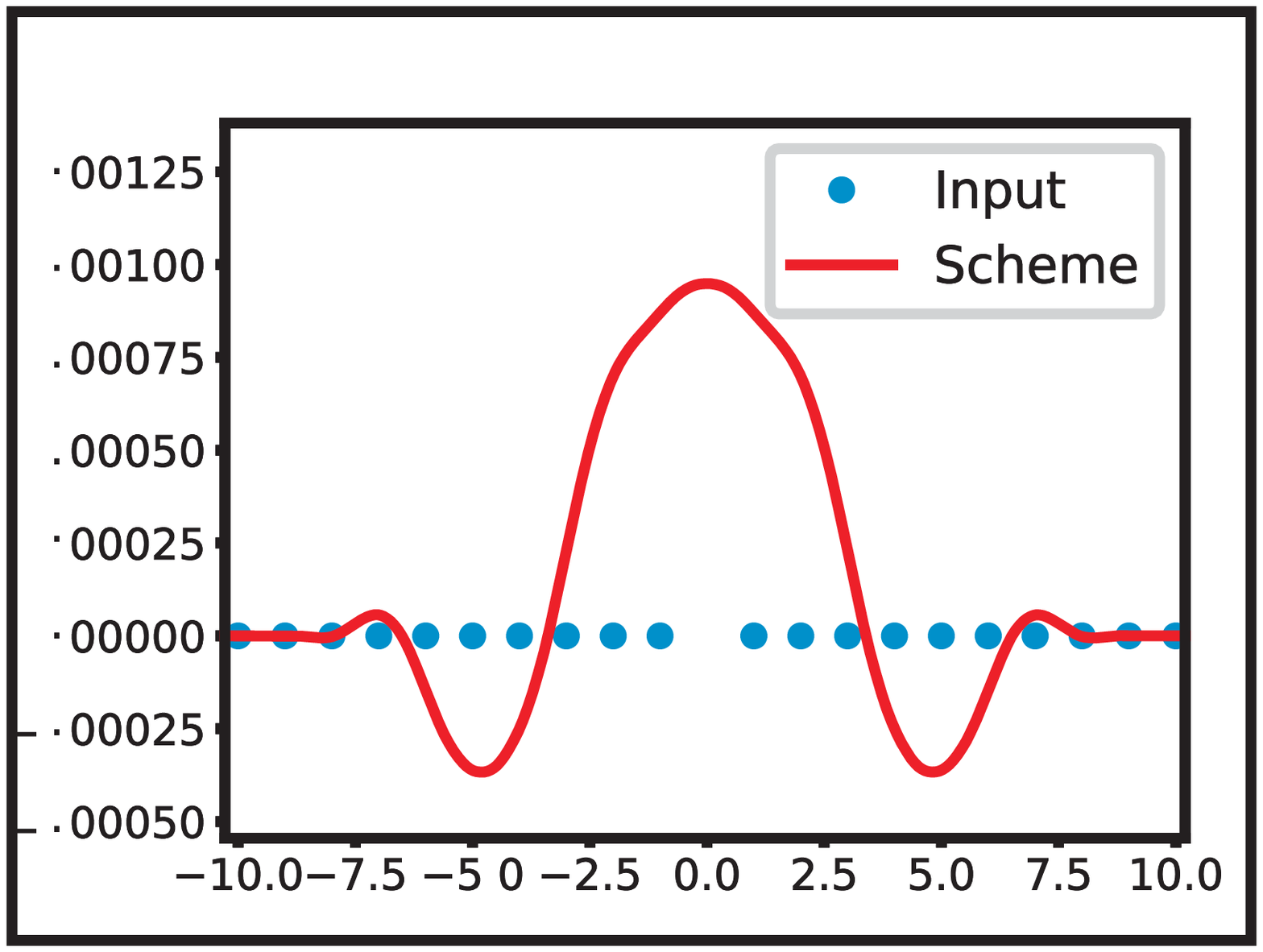, width=1.85 in}&
\epsfig{file=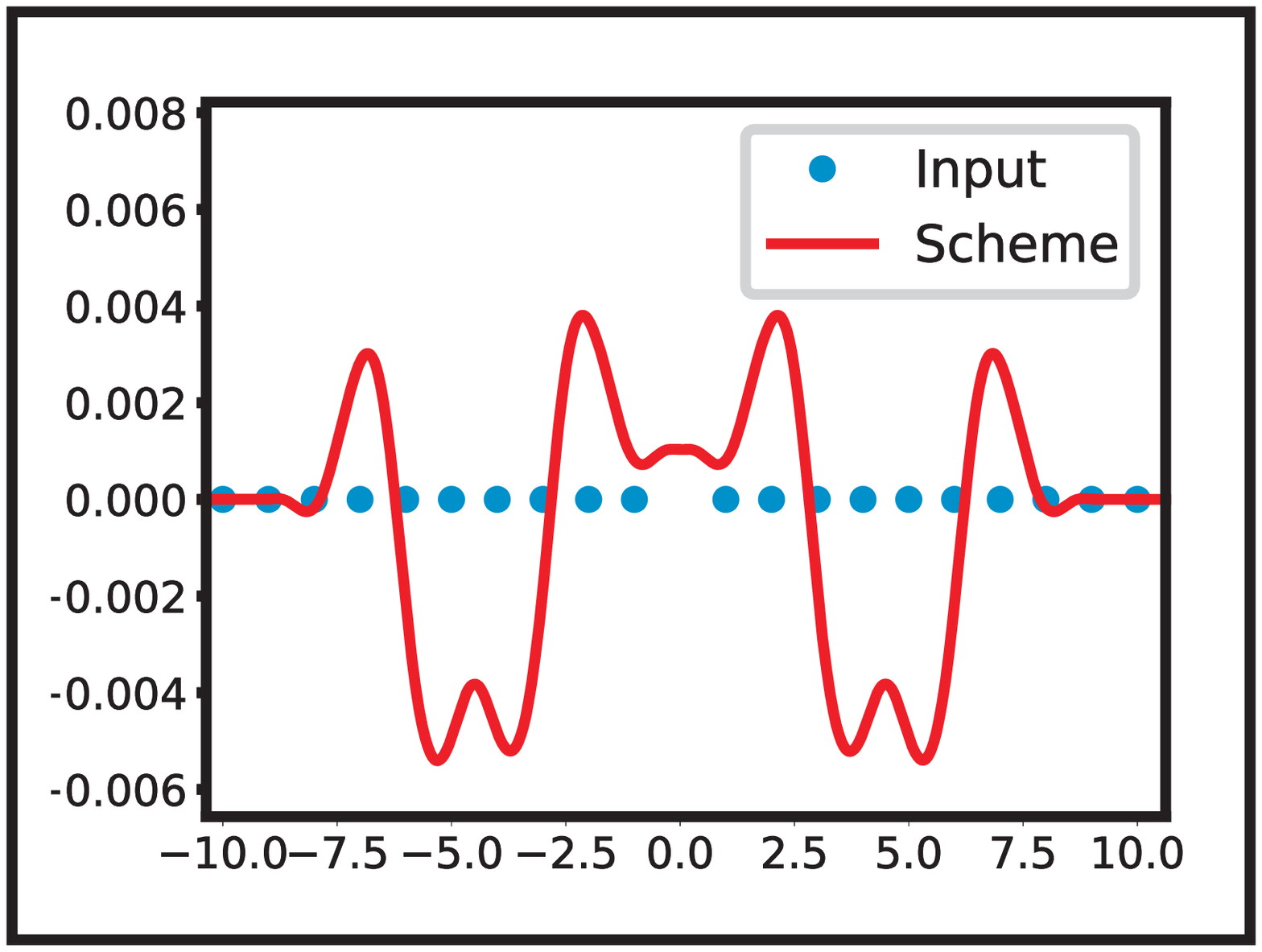, width=1.85 in}&\quad \\
(a) $D_{10,1}$ & (b) $D_{10,2}$ & (c) $D_{10,3}$
 \end{tabular}
\end{center}
 \begin{center}
\begin{tabular}{ccccccccccc}
\epsfig{file=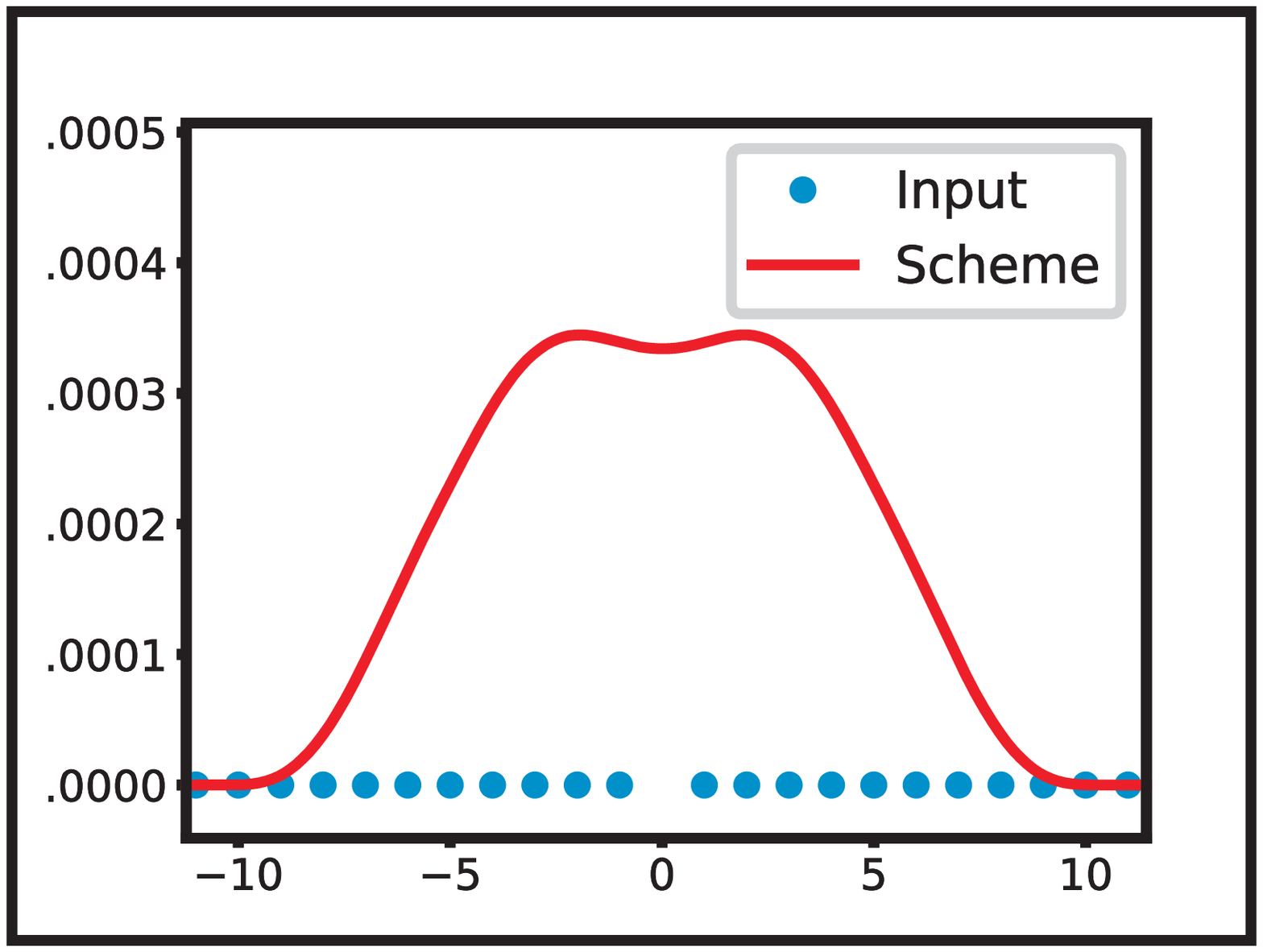, width=1.85 in}&
\epsfig{file=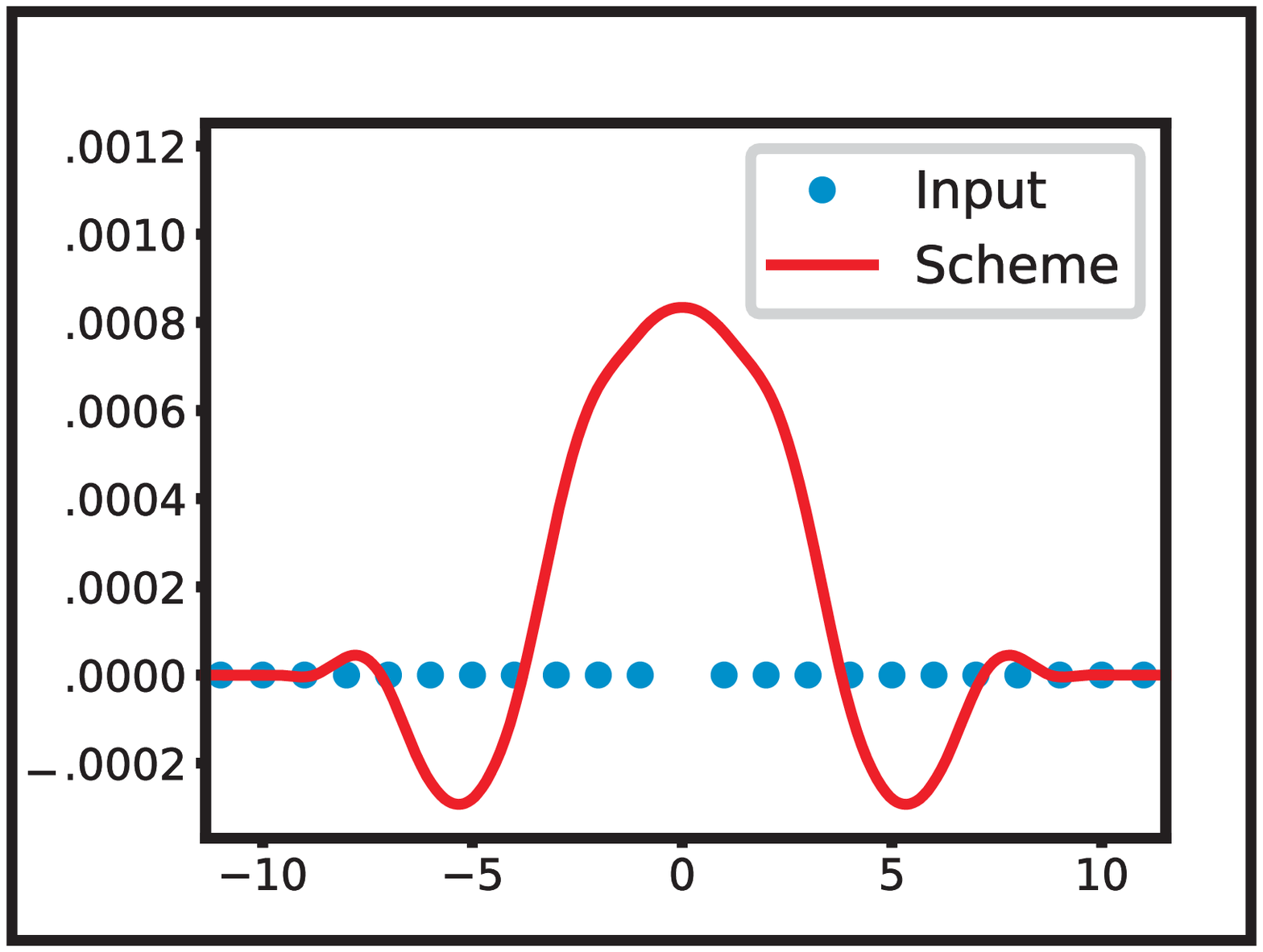, width=1.85 in}&
\epsfig{file=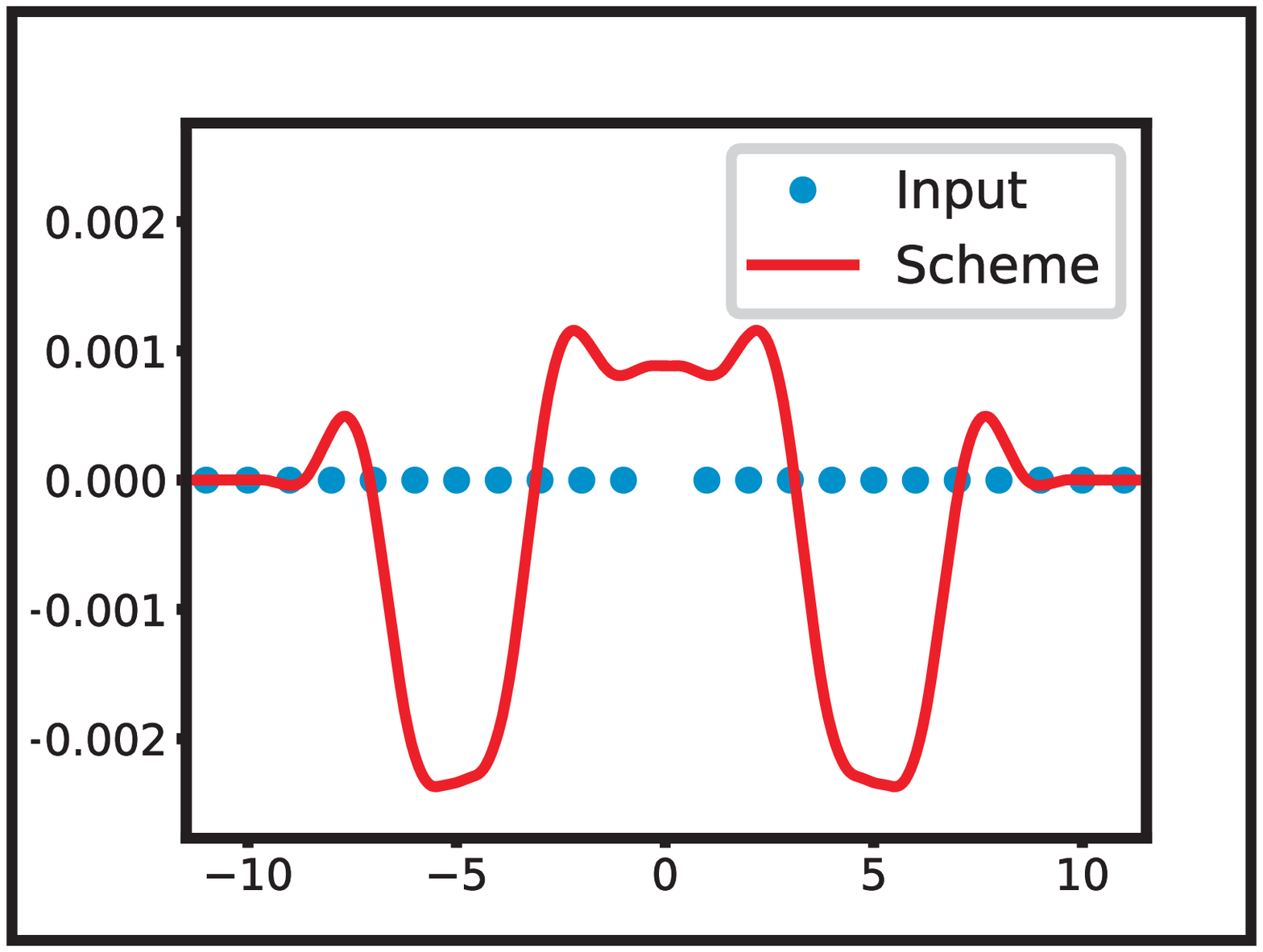, width=1.85 in}&\quad \\
(d) $D_{11,1}$ & (e) $D_{11,2}$ & (f) $D_{11,3}$
 \end{tabular}
\end{center}
 \begin{center}
\begin{tabular}{ccccccccccc}
\epsfig{file=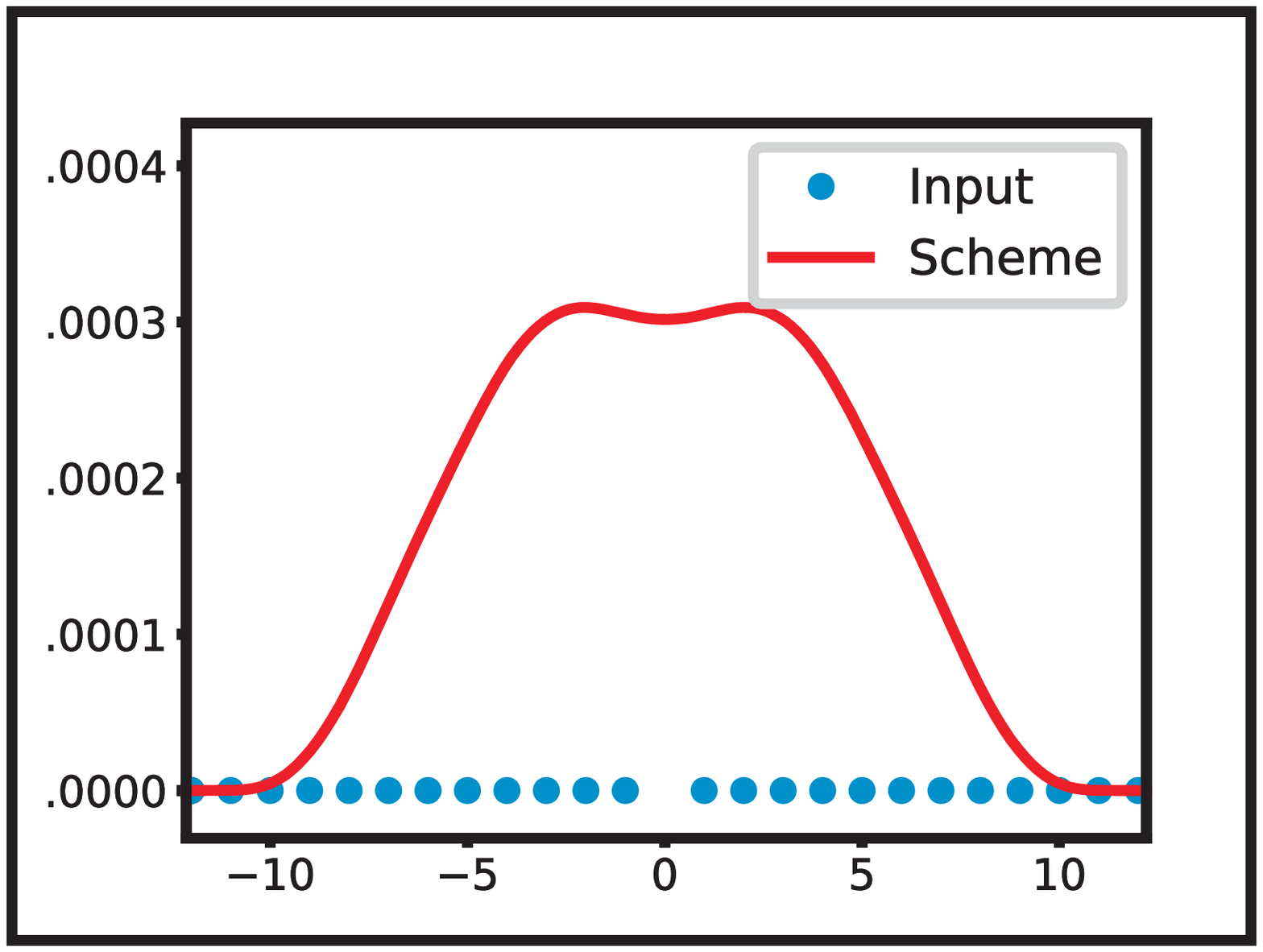, width=1.85 in}&
\epsfig{file=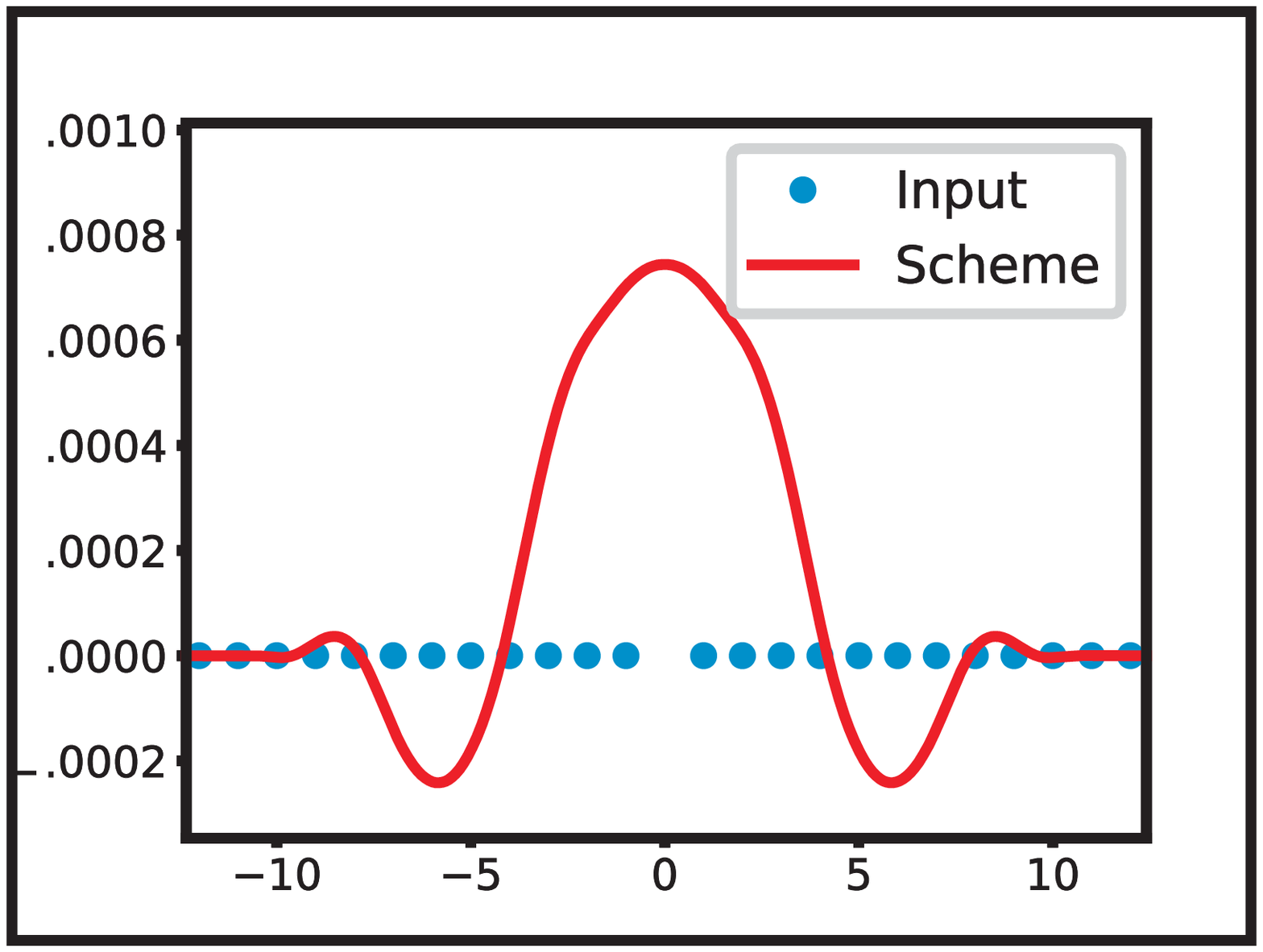, width=1.85 in}&
\epsfig{file=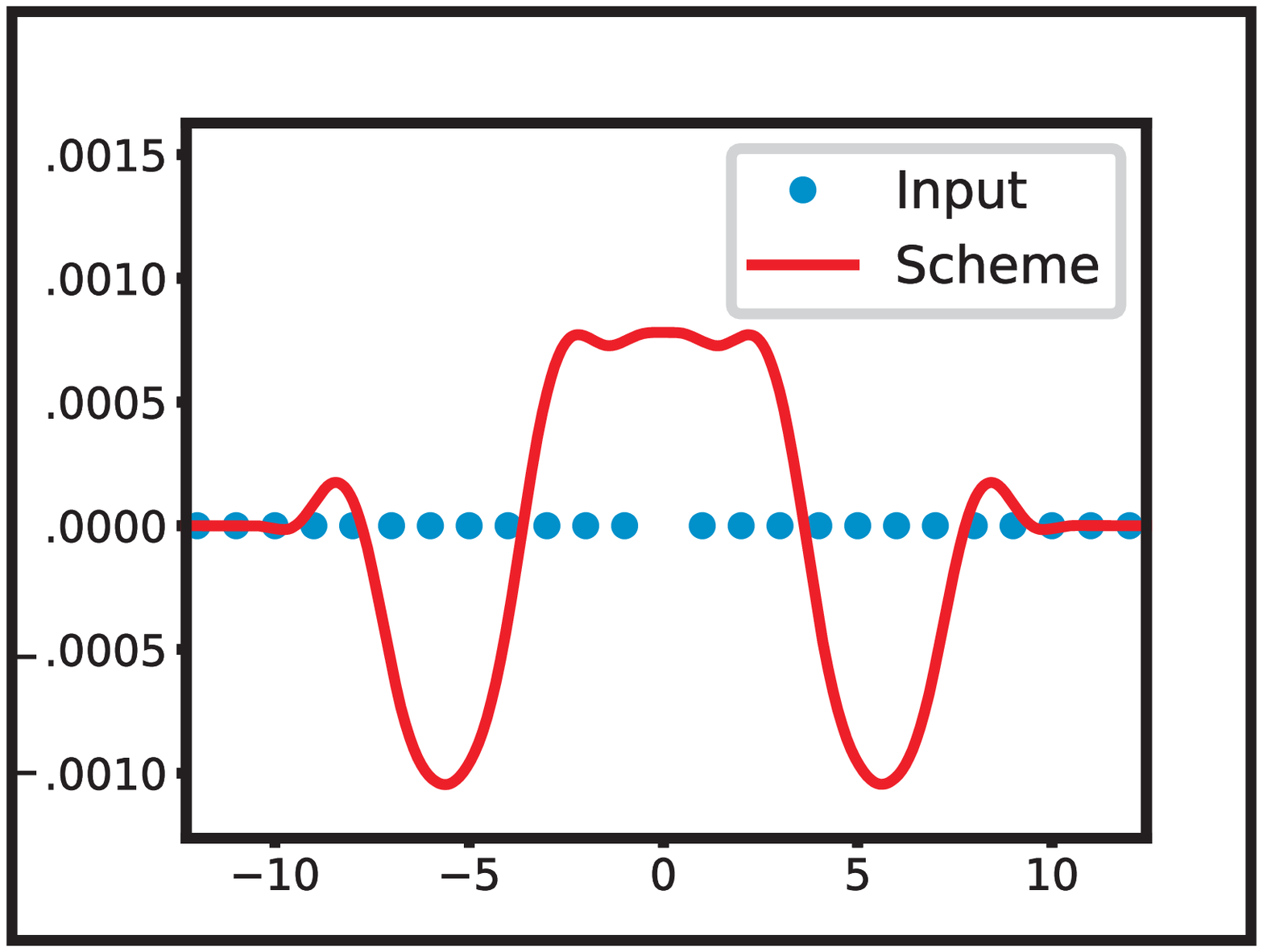, width=1.85 in}&\quad \\
(g) $D_{12,1}$ & (h) $D_{12,2}$ & (i) $D_{12,3}$
 \end{tabular}
\end{center}
 \begin{center}
\begin{tabular}{ccccccccccc}
\epsfig{file=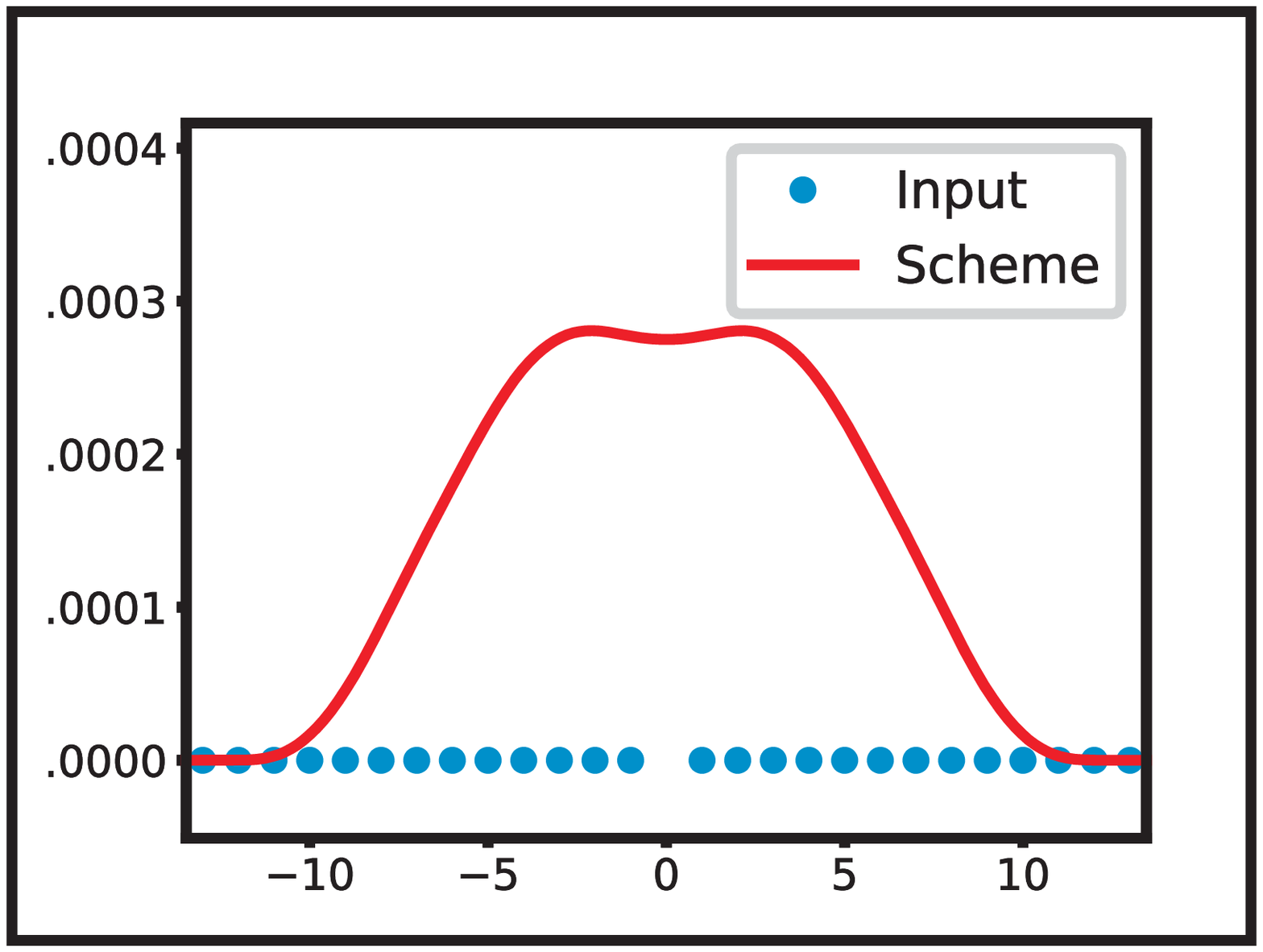, width=1.85 in}&
\epsfig{file=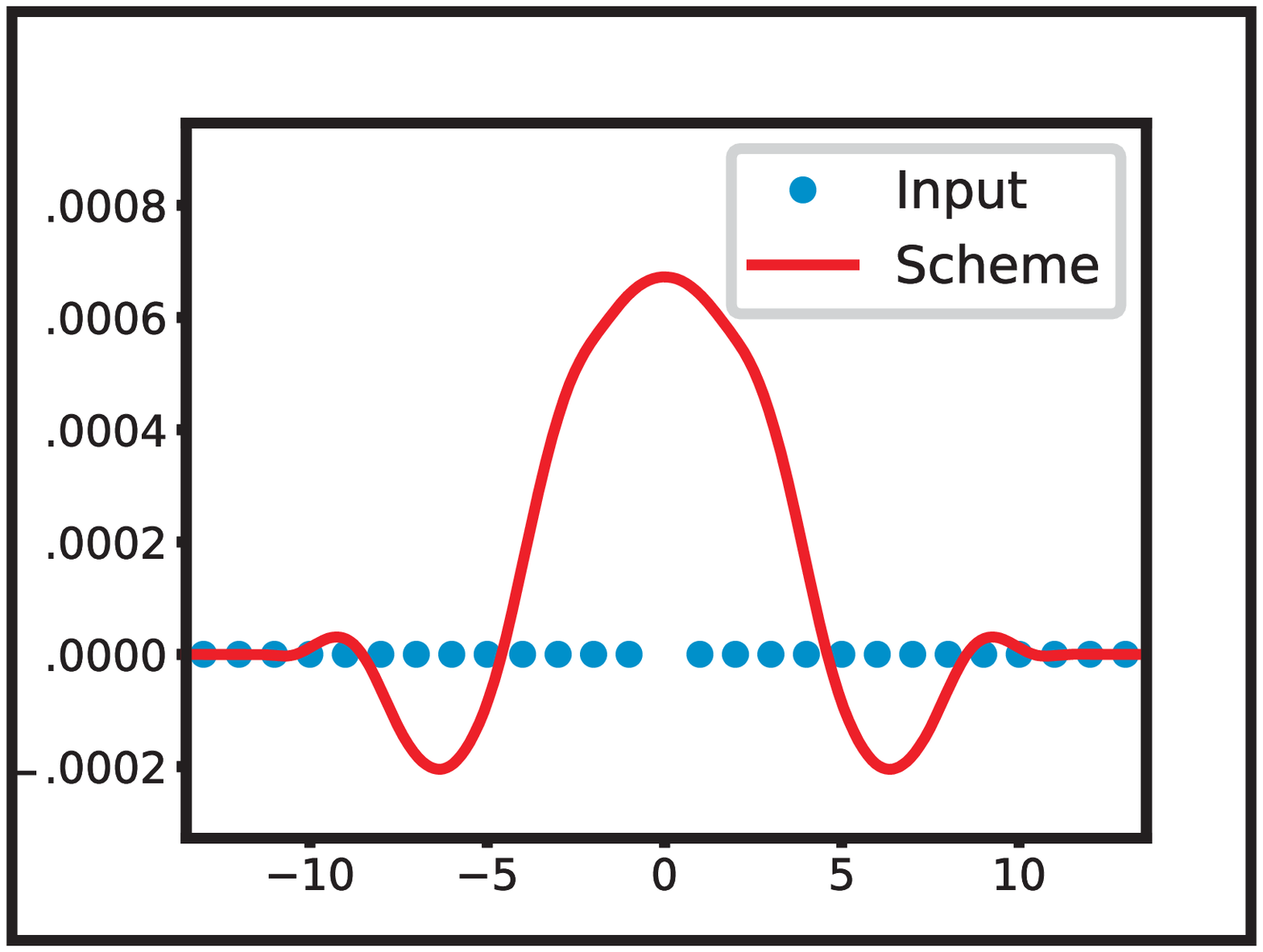, width=1.85 in}&
\epsfig{file=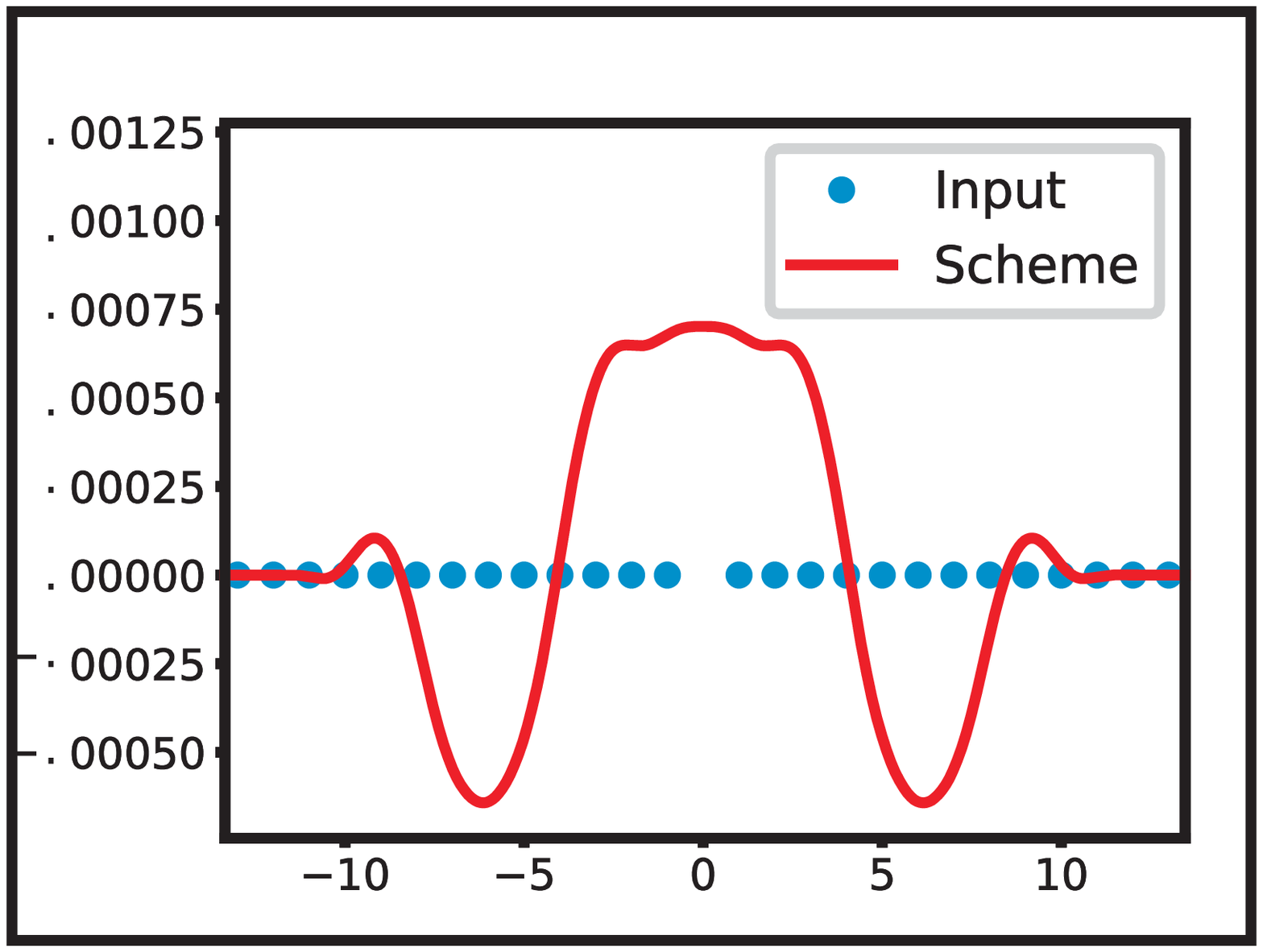, width=1.85 in}&\quad \\
(j) $D_{13,1}$ & (k) $D_{13,2}$ & (l) $D_{13,3}$
 \end{tabular}
\end{center}
\caption[Basic limit functions of the schemes $D_{h,d}$.]{\label{BLF} \emph{Basic limit functions of the schemes $D_{h,d}$, where $h \in \{2n,2n+1: 5 \leqslant n \leqslant 6\}$, $1 \leqslant d \leqslant 3$.}}
\end{figure}

The support of the basic limit functions and the subdivision schemes is the area of the limit curve that will be affected by the displacement of a single control point from its initial place. The part which is dependent on that given control point is called the support width of that control point. By using the approach of Beccari et al. \cite{Beccari}, it can be easily proved that the basic limit functions of families of scheme $D_{2n,d}$ have the support width $4n-1$, where $n \in \mathbb{N}\backslash \{1\}$ which implies that it vanishes out side the interval $[-\frac{4n-1}{2},\frac{4n-1}{2}]$. Similarly, the basic limit functions of the families of schemes $D_{2n+1,d}$ have the support width $4n+1$, where $n \in \mathbb{N}\backslash \{1\}$ which implies that it vanishes out side the interval $[-\frac{4n+1}{2},\frac{4n+1}{2}]$. Figure \ref{BLF}(a)-(l) shows the basic limit functions of the schemes $D_{10,1}$, $D_{10,2}$, $D_{10,3}$, $D_{11,1}$, $D_{11,2}$, $D_{11,3}$, $D_{12,1}$, $D_{12,2}$, $D_{12,3}$, $D_{13,1}$, $D_{13,2}$ and $D_{13,3}$ respectively. From this figure we see that the effects of the schemes constructed with different degree polynomials are different on initial data $\Delta$. The schemes $D_{h,1}$ and $D_{h,2}$ have two and one top peaks respectively $\forall$ $h \in \{2n,2n+1:n=5,6\}$, while the behavior of schemes $D_{h,3}$ is greatly effected with a small change in $n$.

\section{Numerical examples}
In this section, we give several numerical examples to see the performance of the proposed families of schemes. Different types of initial data have been used for this purpose. We give the comparisons of the schemes based on fitting non-linear and linear polynomials.
\begin{Exp}\textbf{Response to the polynomial data}\\
According to \cite{Hormann1}, three numerical quantities indicate the behavior of a subdivision scheme when all the initial control points lie on a polynomial. The limit curves on all polynomial data up to a certain degree passes through the same polynomial, a polynomial of the same degree or passes through the original control points.
These three behaviors show three maximal degrees, that are reproduction degree, generation degree, and interpolation degree respectively. Figures \ref{quadpoly}-\ref{exppoly} show the responses of the schemes $D_{h,1}$, $D_{h,2}$ and $D_{h,3}$, where $h=10,11$, on the data taken from quadratic polynomial function $g_{1}(x)=x^{2}-5 x+3$, cubic polynomial function $g_{2}(x)=x^{3}-x^{2}-5 x+3$ and exponential polynomial function $g_{3}(x)=(0.7670)e^{0.4040 x}$ respectively. We can see from these figures that fitted curves by schemes $D_{h,2}$, $D_{h,3}$, $h=10,11$ generate, reproduce and interpolate these polynomials while the limit curves of the schemes $D_{h,1}$ do not reproduce and interpolate these polynomials.
\begin{figure}[h] 
 \begin{center}
\begin{tabular}{ccccccccccc}
\epsfig{file=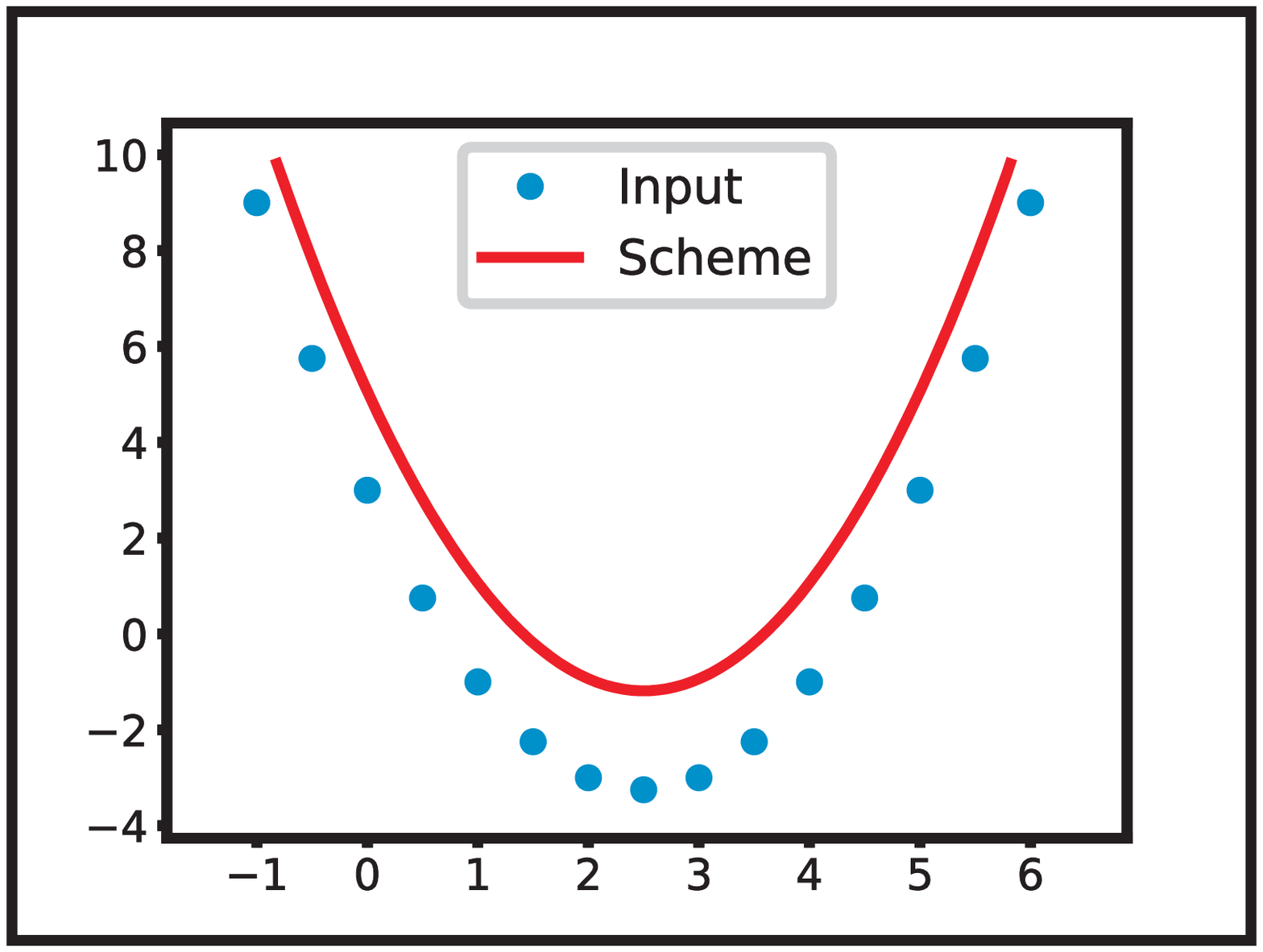, width=1.85 in}&
\epsfig{file=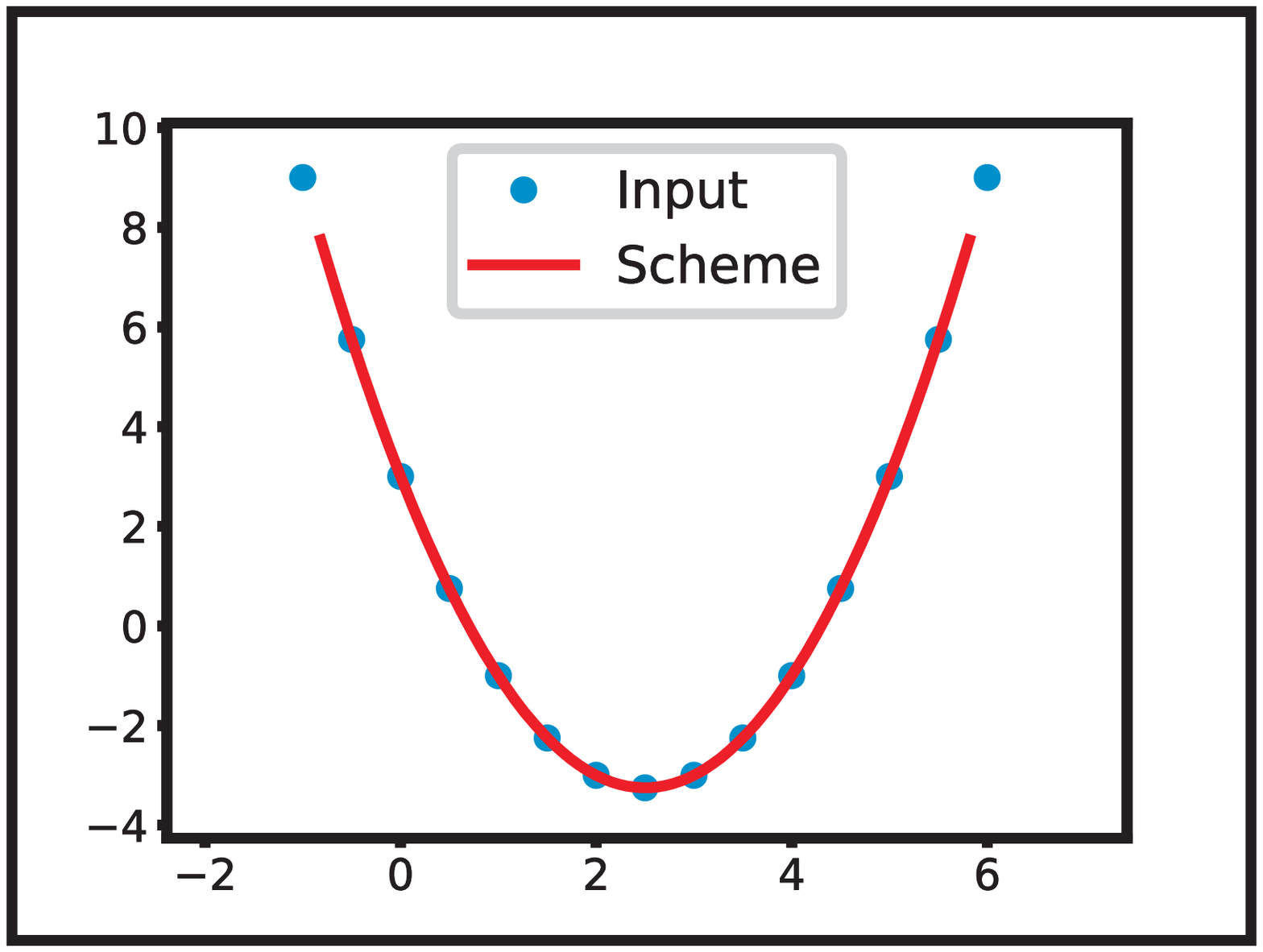, width=1.85 in}&
\epsfig{file=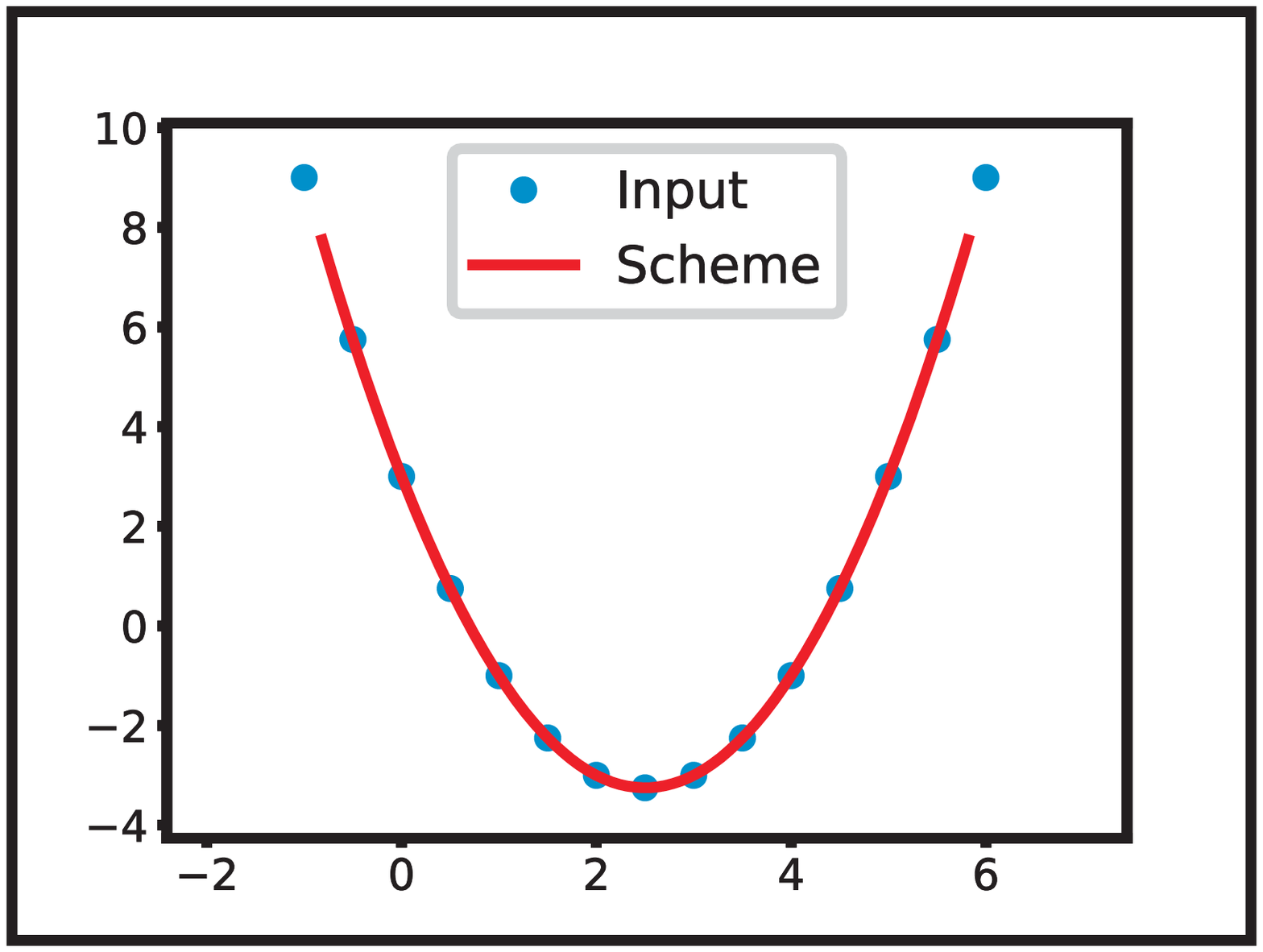, width=1.85 in}&\quad \\
(a) $D_{10,1}$ & (b) $D_{10,2}$ & (c) $D_{10,3}$
 \end{tabular}
\end{center}
 \begin{center}
\begin{tabular}{ccccccccccc}
\epsfig{file=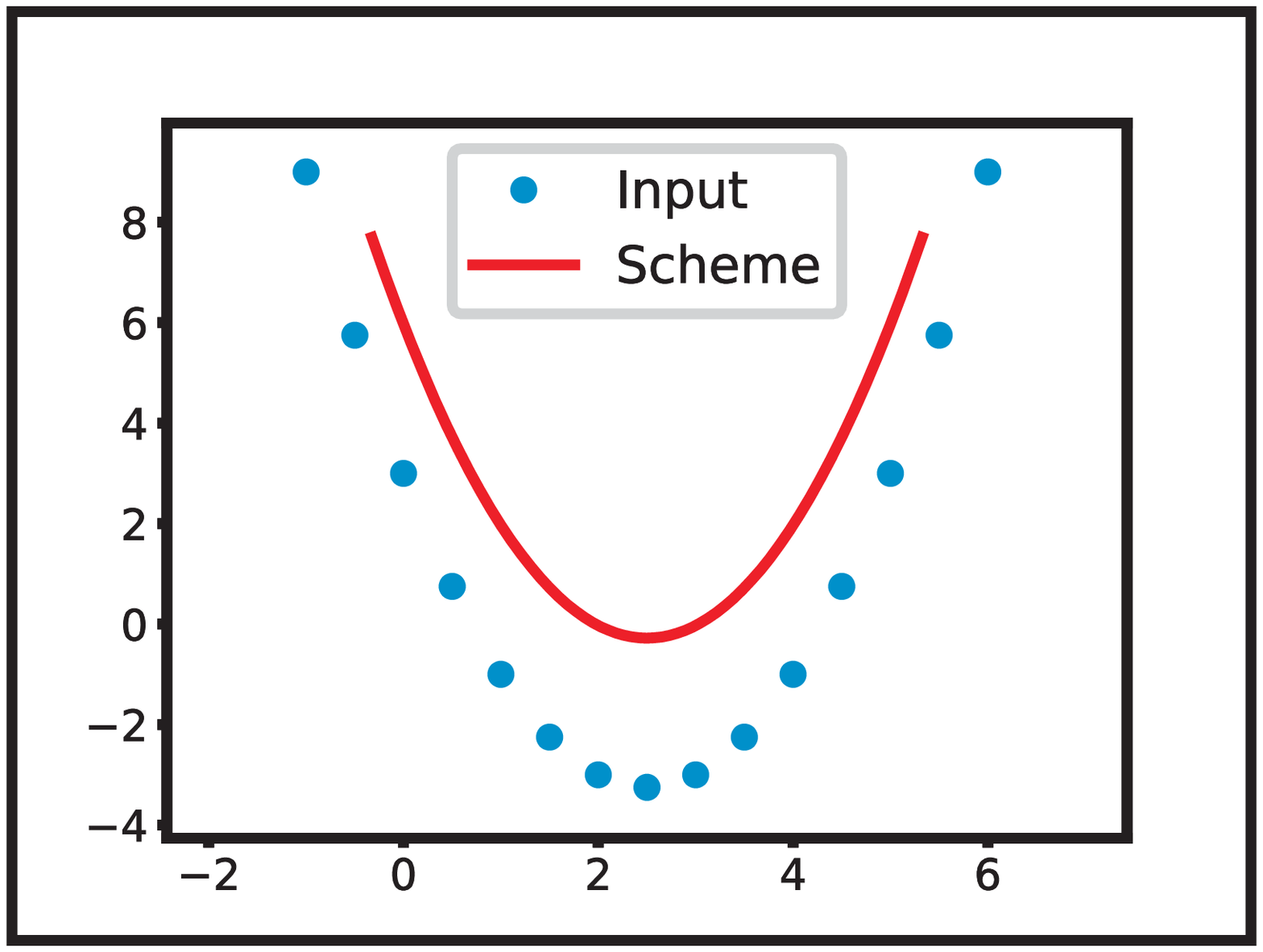, width=1.85 in}&
\epsfig{file=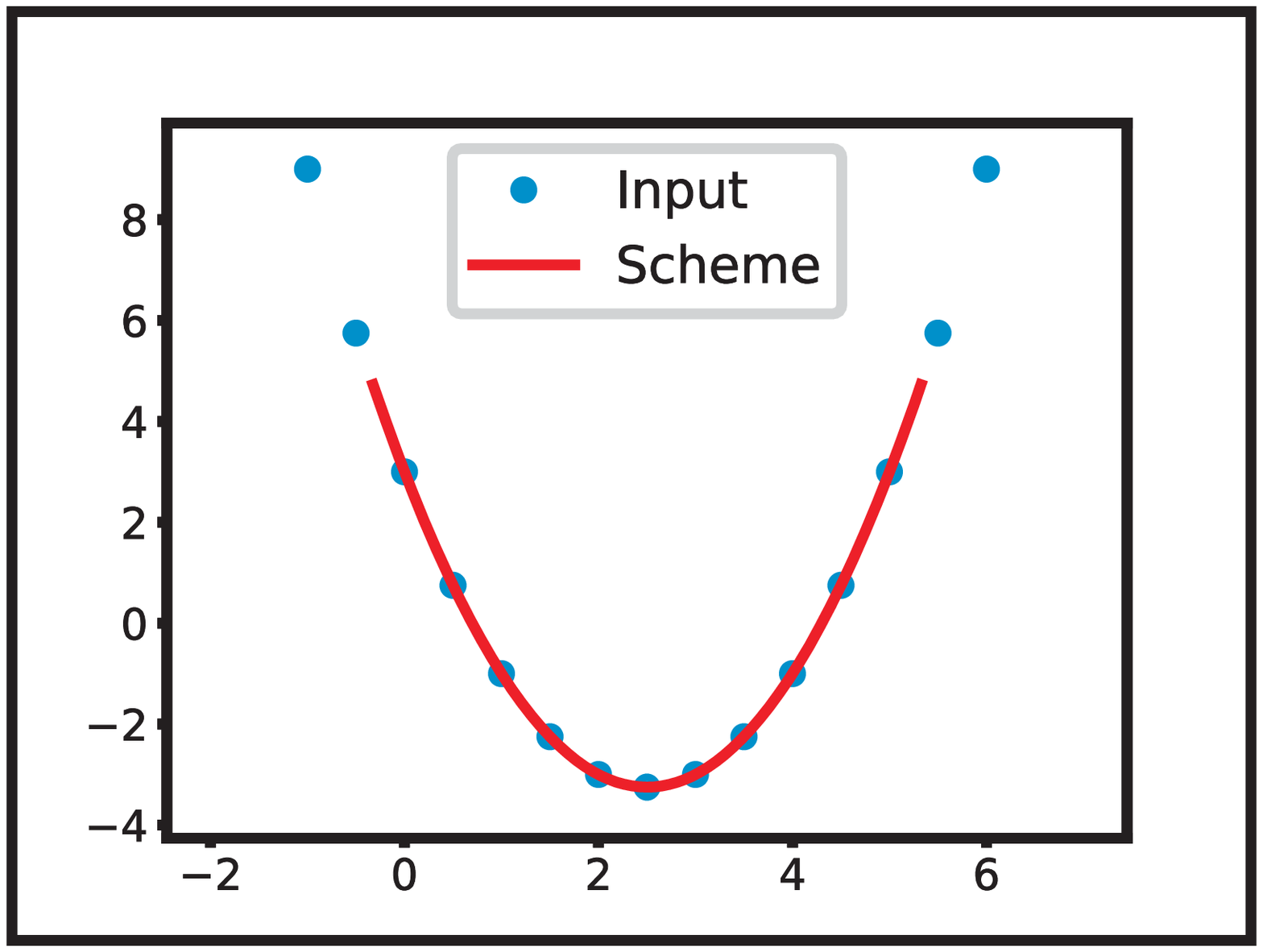, width=1.85 in}&
\epsfig{file=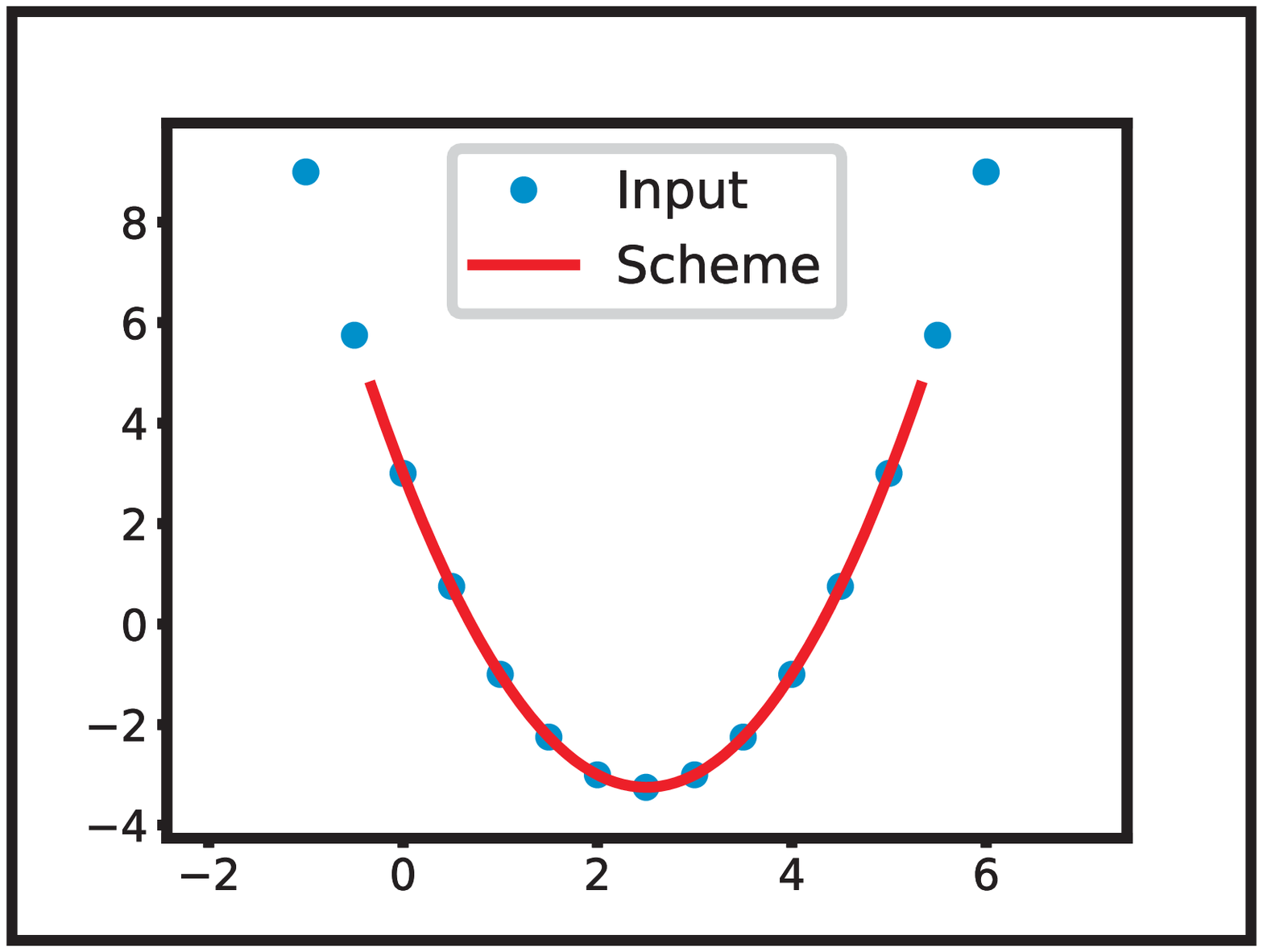, width=1.85 in}&\quad \\
(d) $D_{11,1}$ & (e) $D_{11,2}$ & (f) $D_{11,3}$
 \end{tabular}
\end{center}
\caption[Effects of the schemes $D_{h,d}$ on the polynomial data of degree 2.]{\label{quadpoly} \emph{Effects of the schemes $D_{h,d}$, where $h \in \{2n,2n+1:n=5\}$ and $1 \leqslant d \leqslant 3$, on the polynomial data of degree 2.}}
\end{figure}

\begin{figure}[htb!] 
 \begin{center}
\begin{tabular}{ccccccccccc}
\epsfig{file=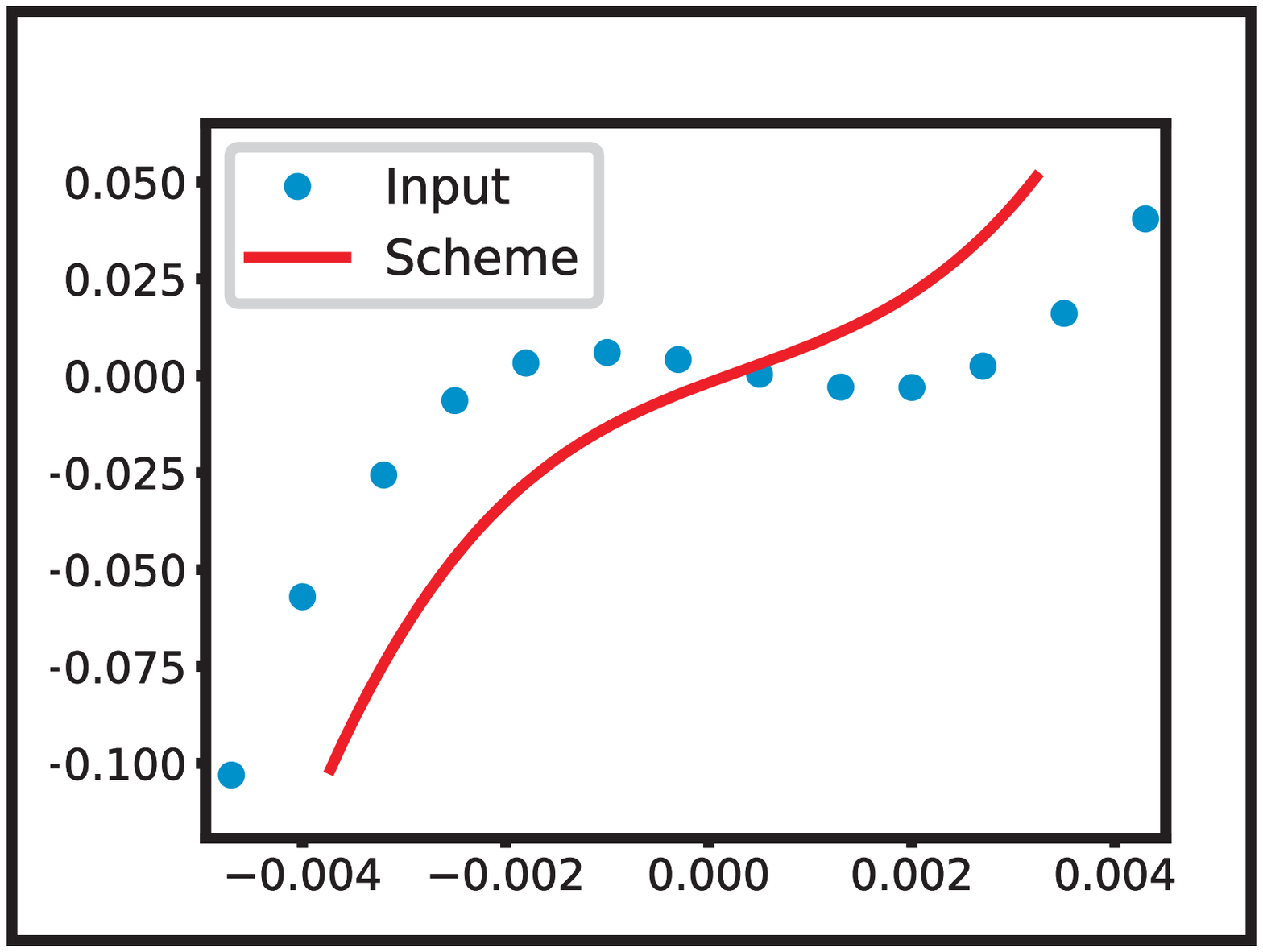, width=1.85 in}&
\epsfig{file=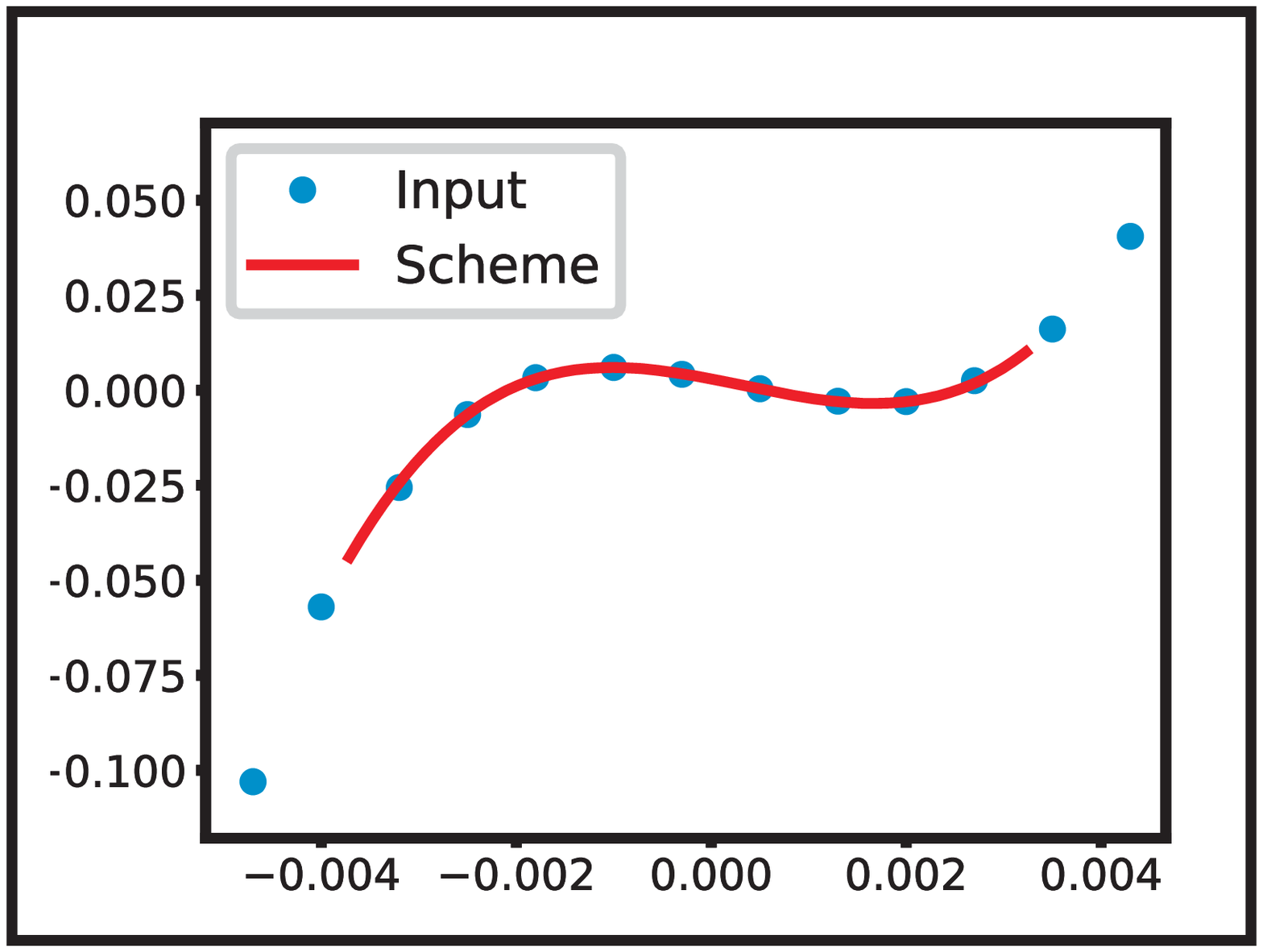, width=1.85 in}&
\epsfig{file=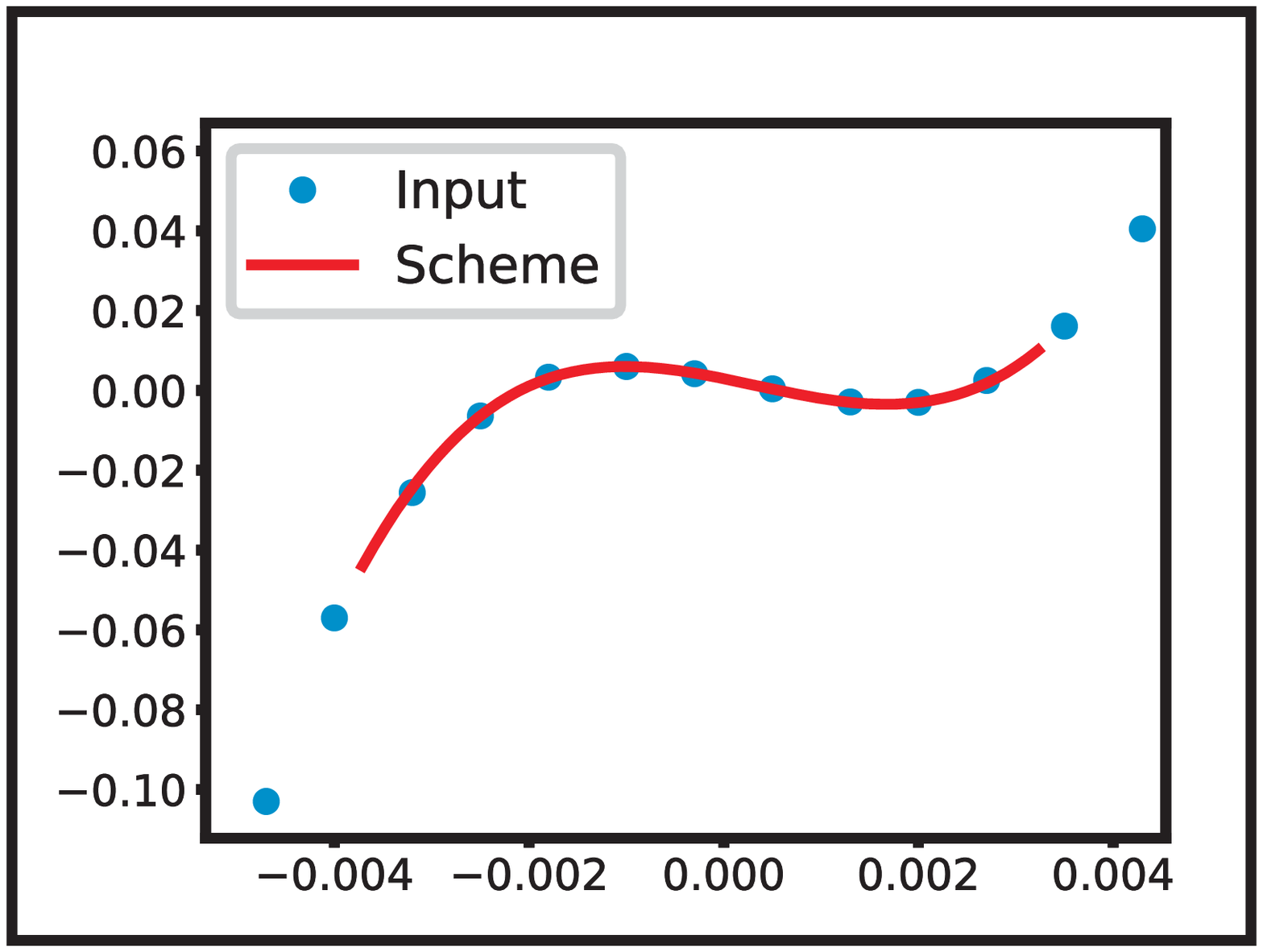, width=1.85 in}&\quad \\
(a) $D_{10,1}$ & (b) $D_{10,2}$ & (c) $D_{10,3}$
 \end{tabular}
\end{center}
 \begin{center}
\begin{tabular}{ccccccccccc}
\epsfig{file=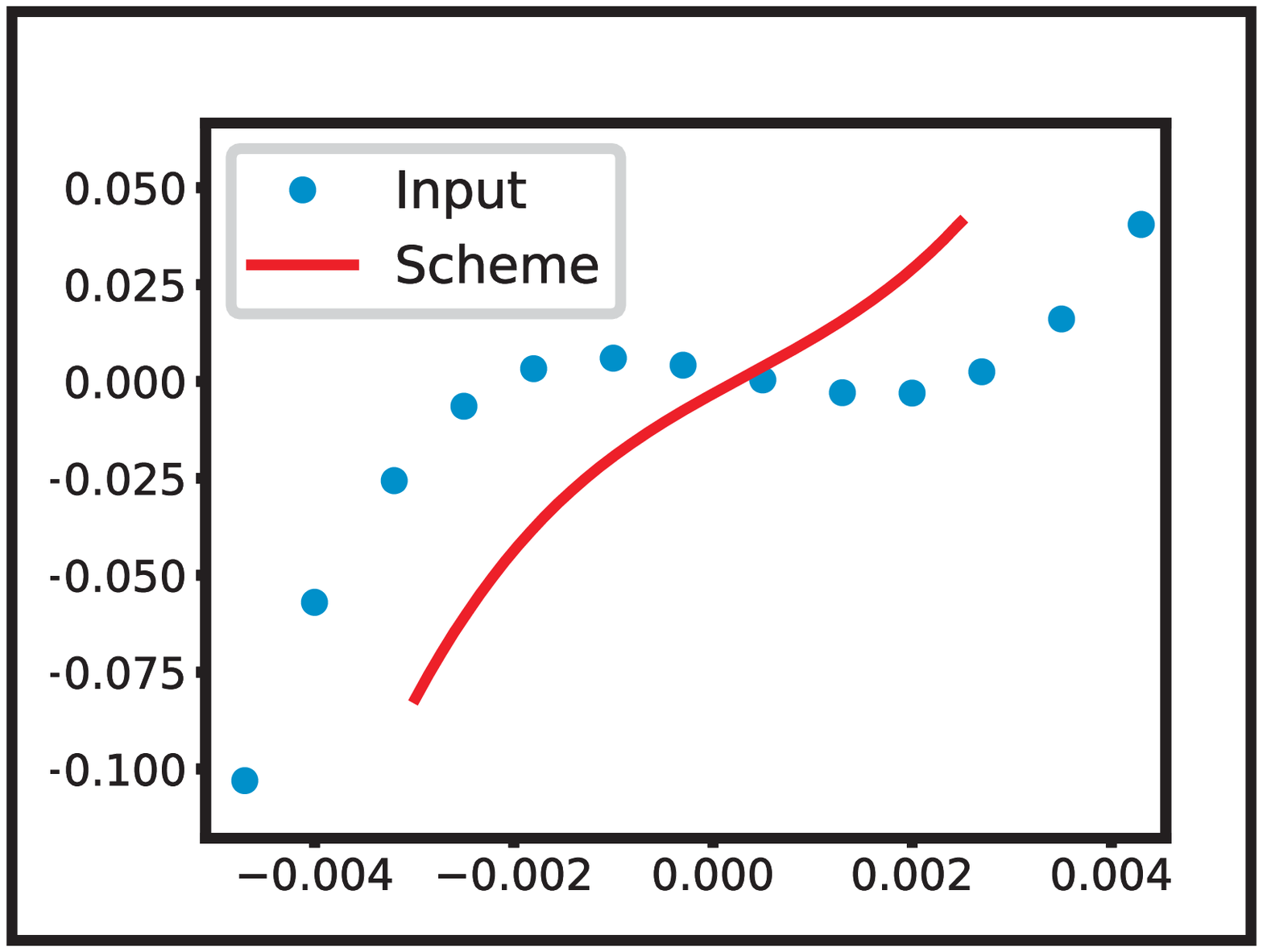, width=1.85 in}&
\epsfig{file=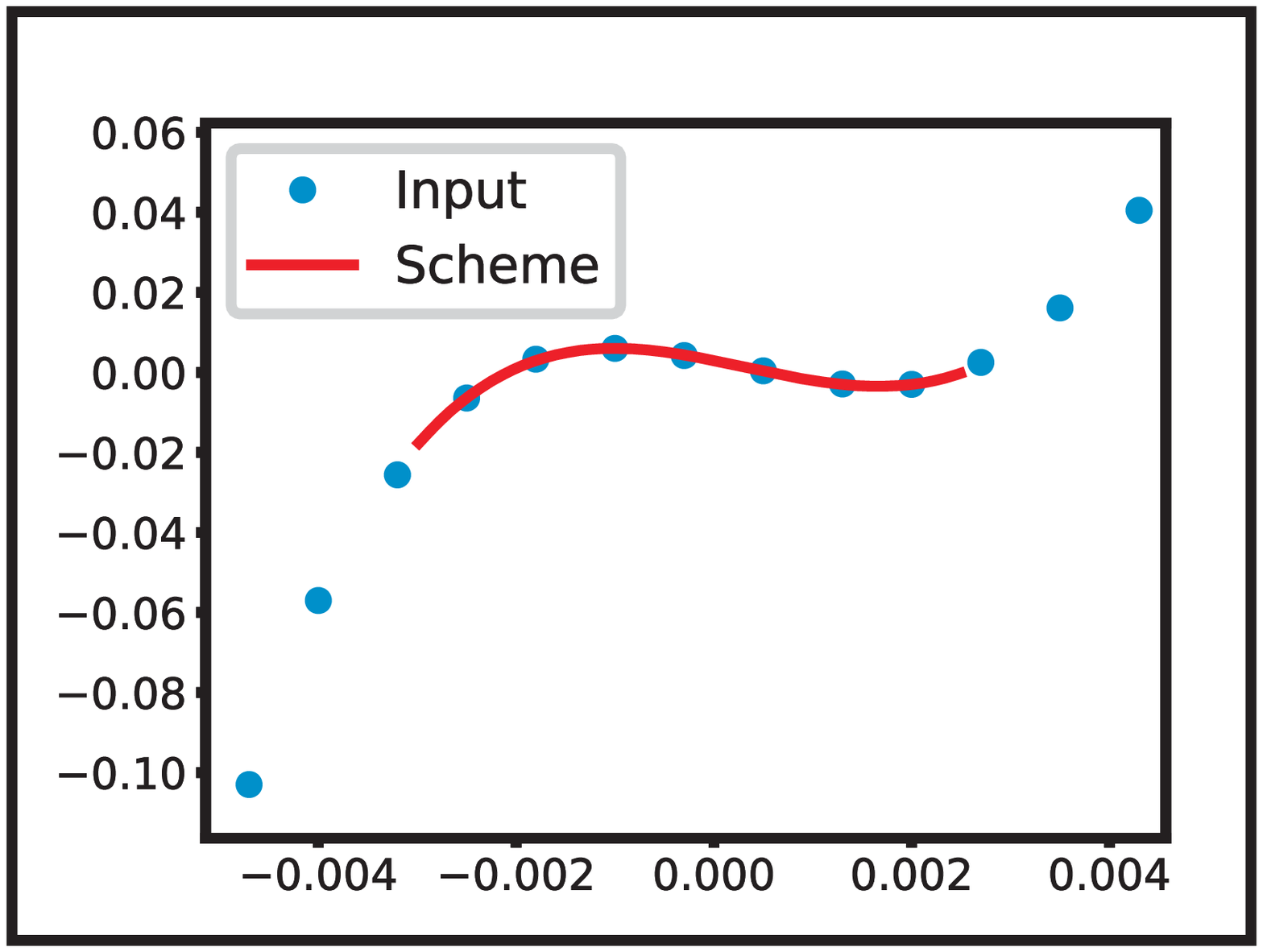, width=1.85 in}&
\epsfig{file=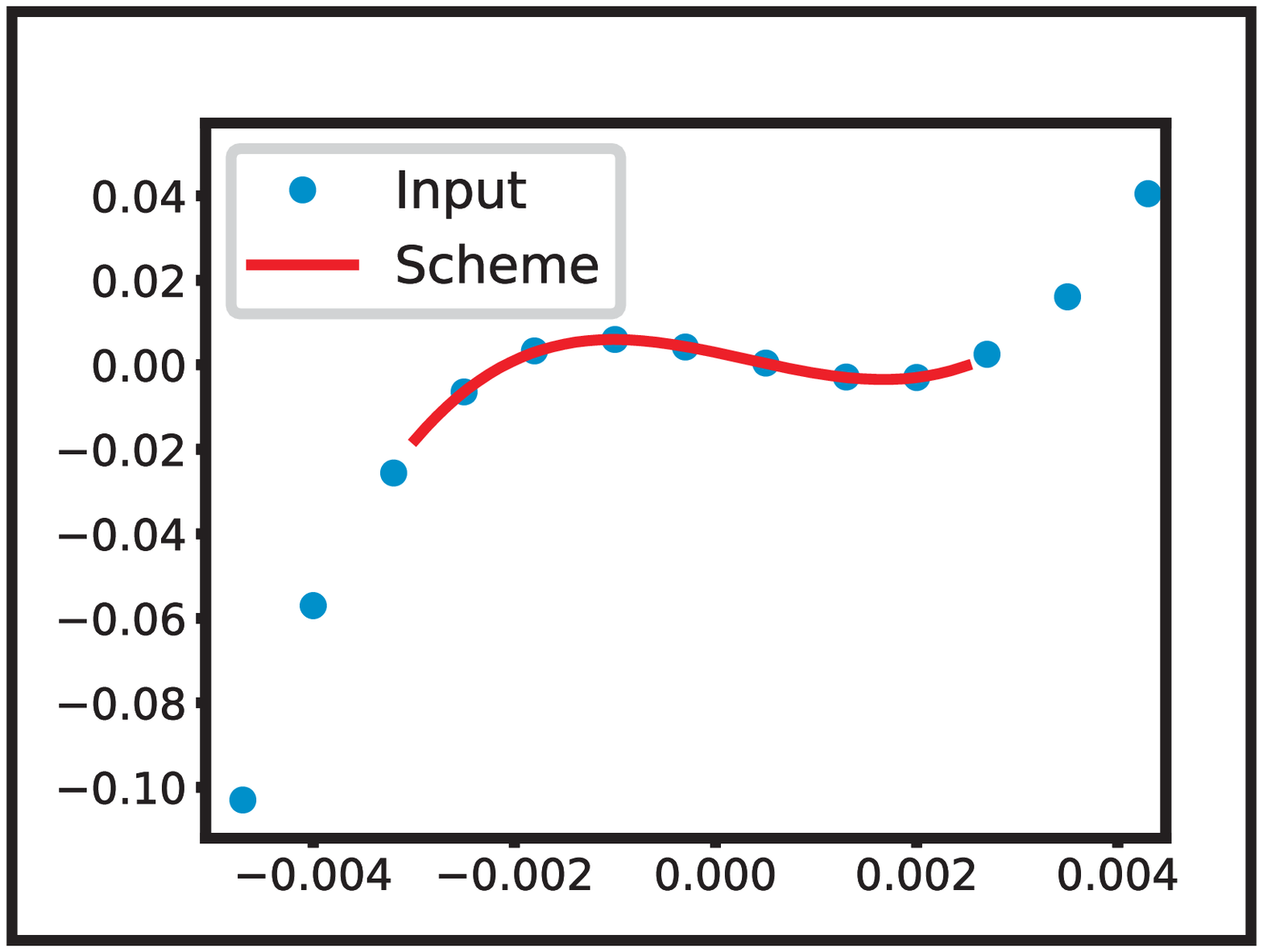, width=1.85 in}&\quad \\
(d) $D_{11,1}$ & (e) $D_{11,2}$ & (f) $D_{11,3}$
 \end{tabular}
\end{center}
\caption[Effects of the schemes $D_{h,d}$ on the polynomial data of degree 3.]{\label{cubepoly} \emph{Effects of the schemes $D_{h,d}$, where $h \in \{2n,2n+1:n=5\}$ and $1 \leqslant d \leqslant 3$, on the polynomial data of degree 3.}}
\end{figure}
\begin{figure}[htb!] 
 \begin{center}
\begin{tabular}{ccccccccccc}
\epsfig{file=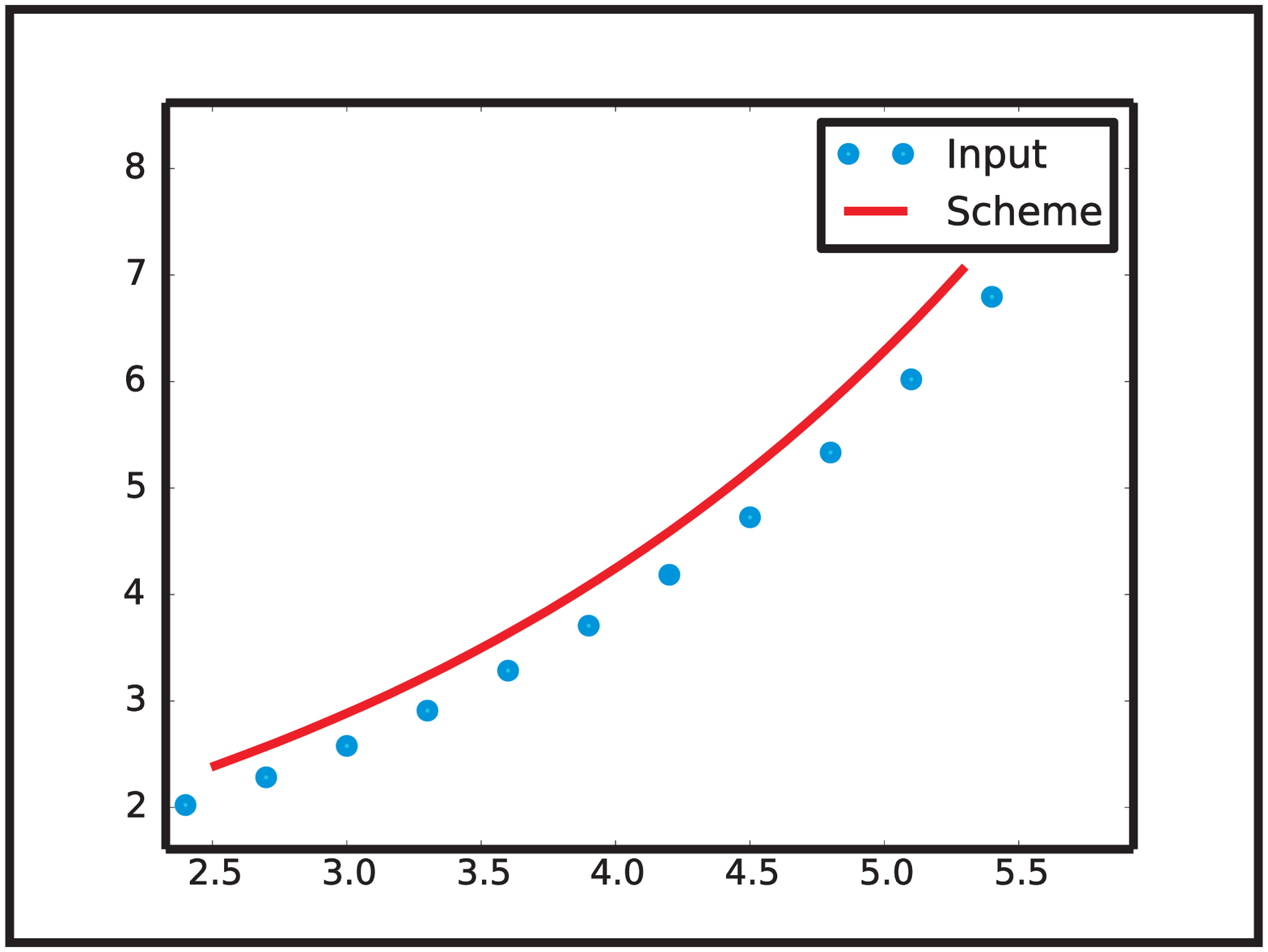, width=1.85 in}&
\epsfig{file=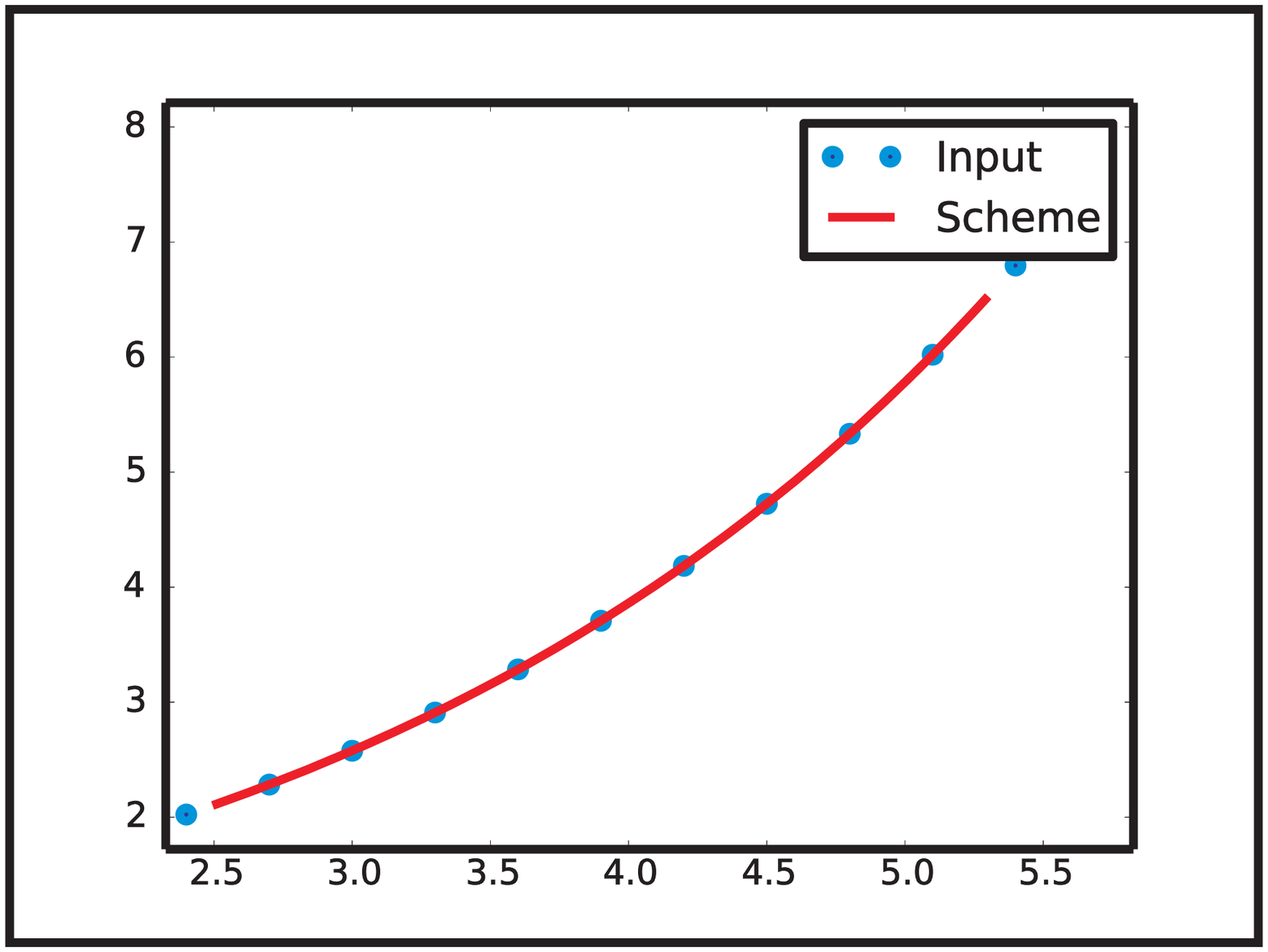, width=1.85 in}&
\epsfig{file=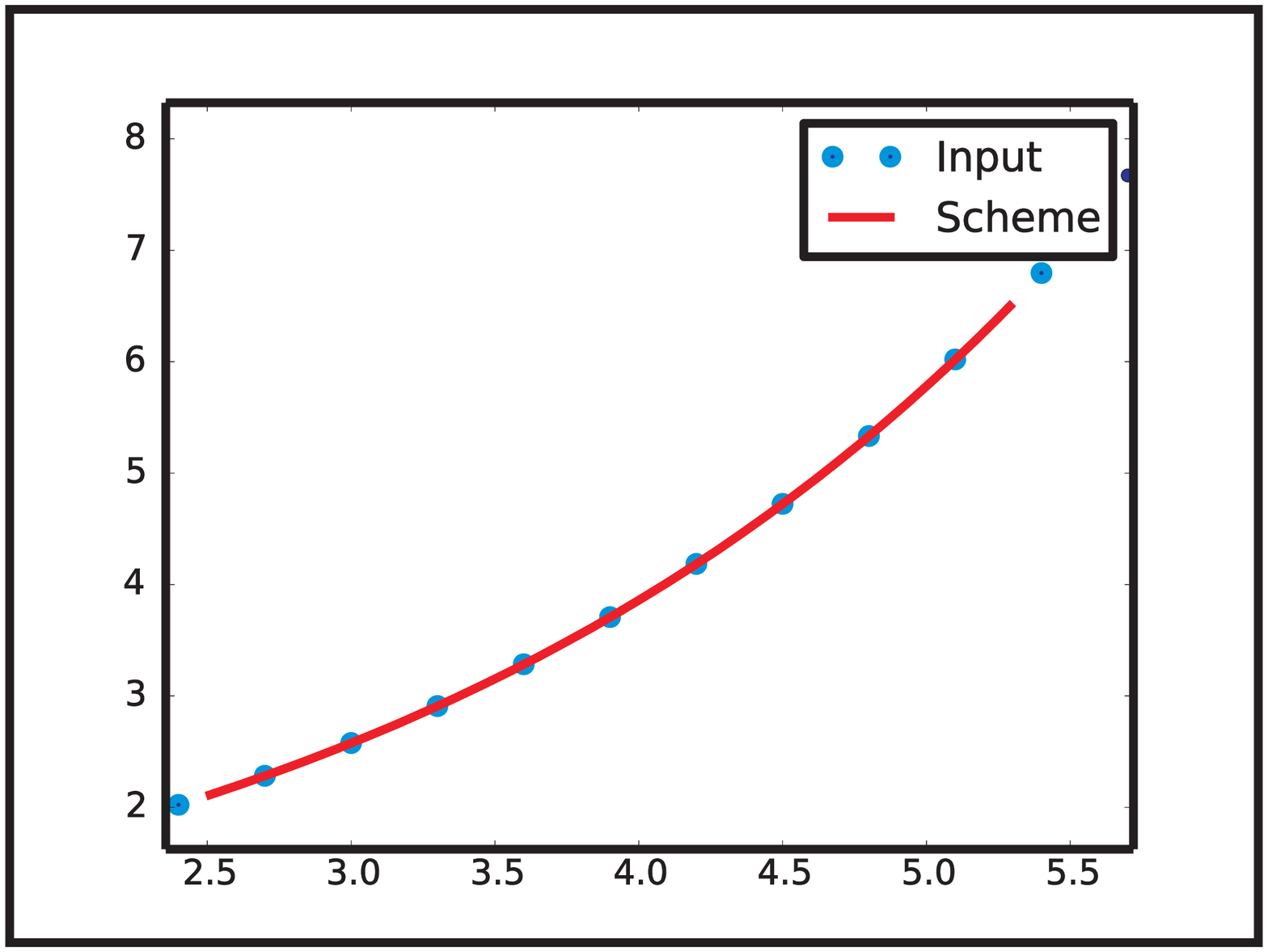, width=1.85 in}&\quad \\
(a) $D_{10,1}$ & (b) $D_{10,2}$ & (c) $D_{10,3}$
 \end{tabular}
\end{center}
 \begin{center}
\begin{tabular}{ccccccccccc}
\epsfig{file=11-point-exponential-linear.eps, width=1.85 in}&
\epsfig{file=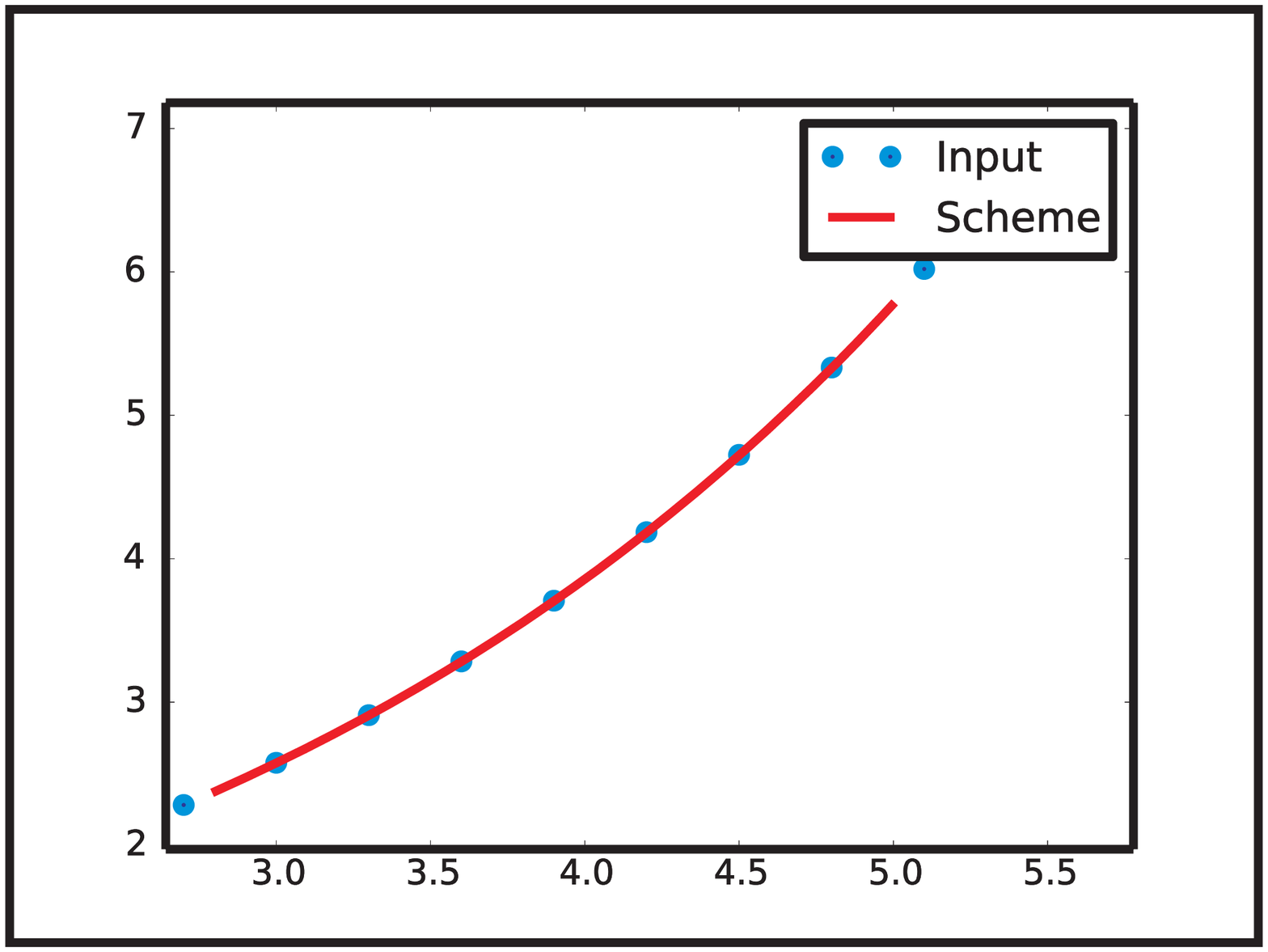, width=1.85 in}&
\epsfig{file=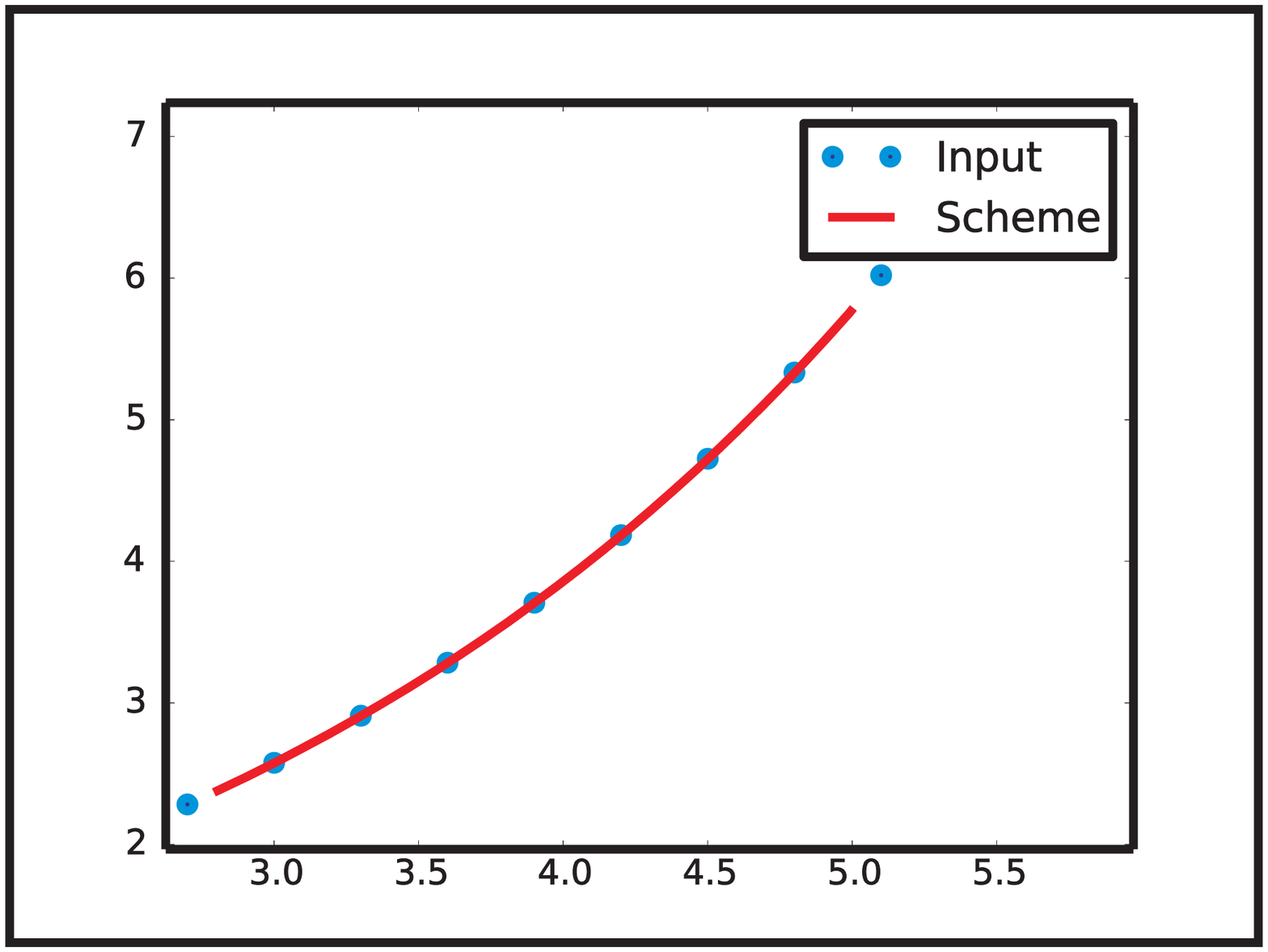, width=1.85 in}&\quad \\
(d) $D_{11,1}$ & (e) $D_{11,2}$ & (f) $D_{11,3}$
 \end{tabular}
\end{center}
\caption[Effects of the schemes $D_{h,d}$ on the exponential polynomial data.]{\label{exppoly} \emph{Effects of the schemes $D_{h,d}$, where $h \in \{2n,2n+1:n=5\}$ and $1 \leqslant d \leqslant 3$, on the exponential polynomial data.}}
\end{figure}
\end{Exp}
\begin{Exp}
\textbf{Response to the discontinuous data}\\
If we generate initial data from a discontinuous function
\begin{equation*}
 g_{4}(x)= \begin{cases}
-10  & \text{for $x \in ]-\infty,0]$}, \\
\,\ 10  & \text{for $x \in ]0,+\infty[$},
\end{cases}
\end{equation*}
then from Figure \ref{discontinuous-function} the fitted curves generated by the proposed schemes $D_{10,1}$, $D_{10,2}$, $D_{11,1}$ and $D_{11,2}$ do not suffer from Gibbs oscillations, whereas the curves generated by the proposed schemes $D_{10,3}$ and $D_{11,3}$ suffer from these oscillations. It means that the schemes based on higher degree polynomials are not suitable for fitting the data that comes from discontinuous functions.\emph{}
\begin{figure}[tbp] 
 \begin{center}
\begin{tabular}{ccccccccccc}
\epsfig{file=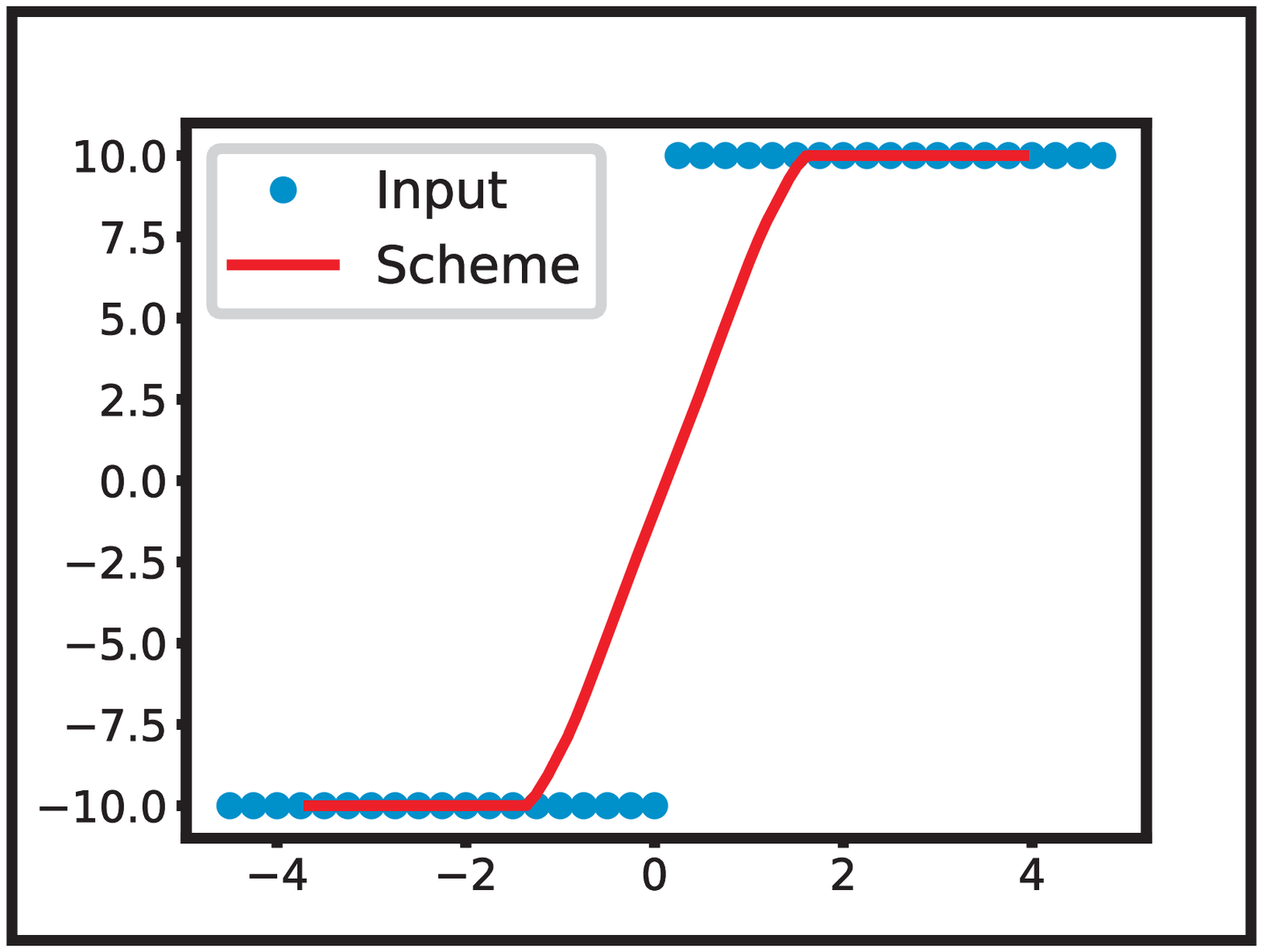, width=1.85 in}&
\epsfig{file=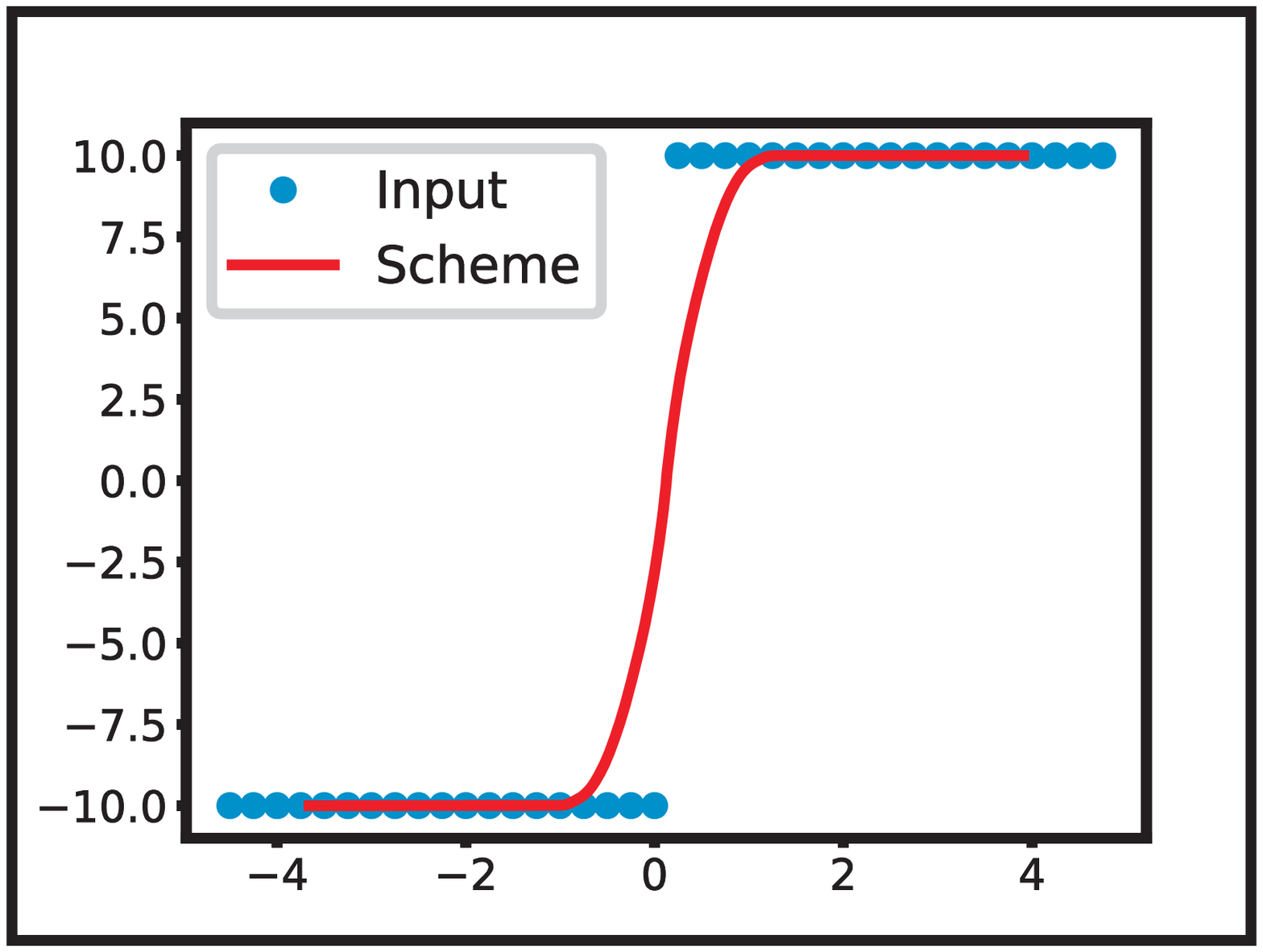, width=1.85 in}&
\epsfig{file=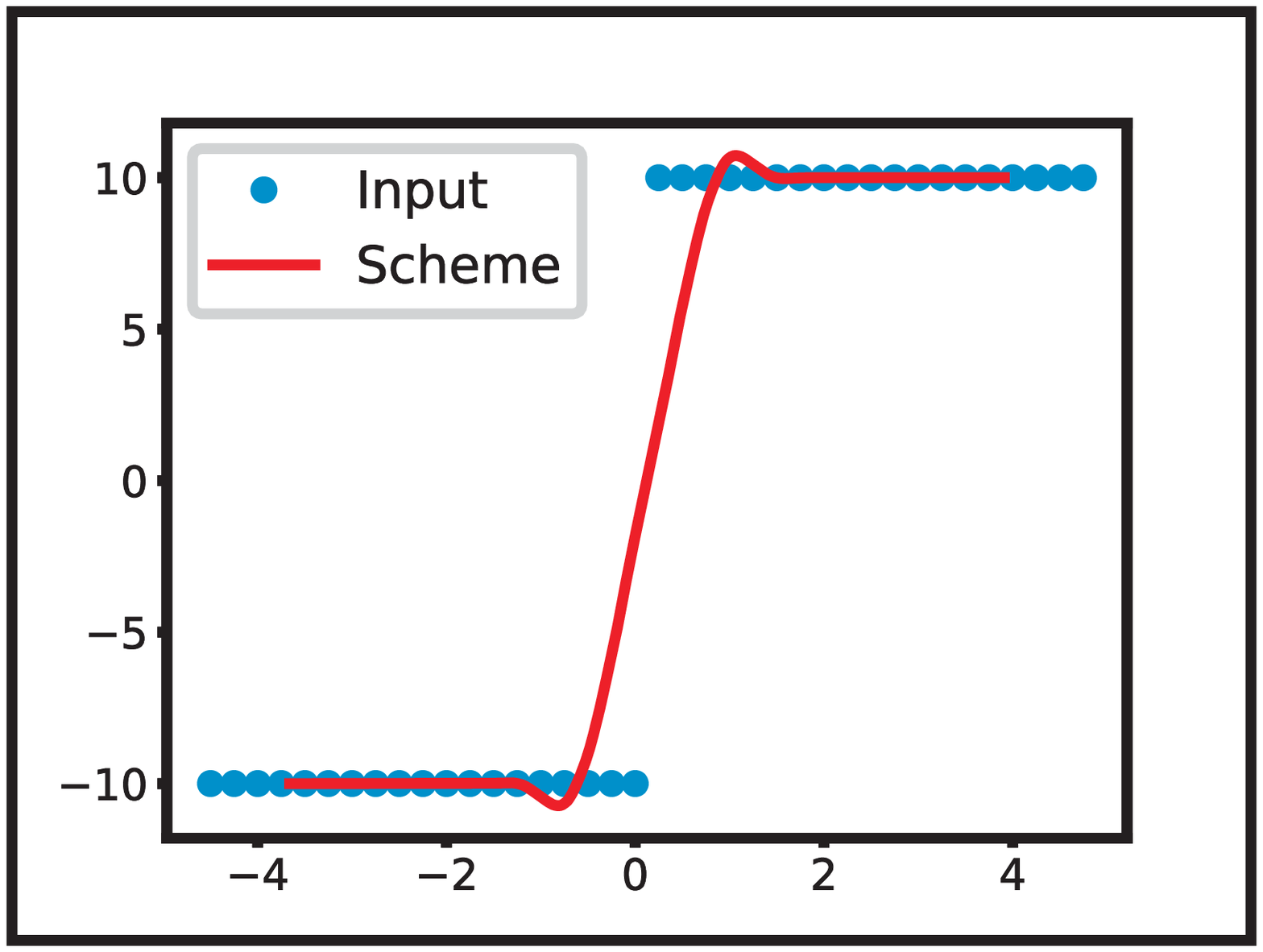, width=1.85 in}&\quad \\
(a) $D_{10,1}$ & (b) $D_{10,2}$ & (c) $D_{10,3}$
 \end{tabular}
\end{center}
 \begin{center}
\begin{tabular}{ccccccccccc}
\epsfig{file=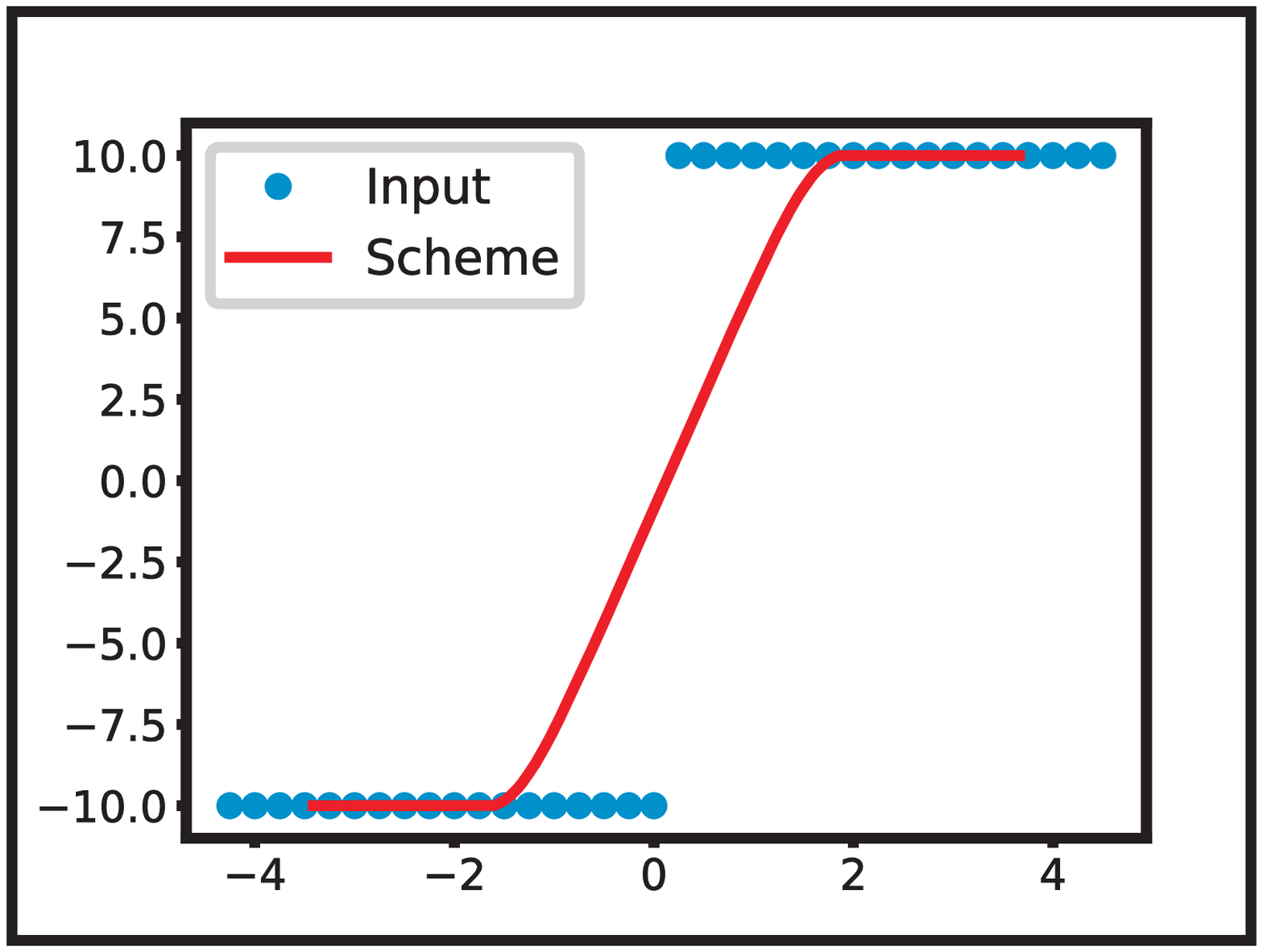, width=1.85 in}&
\epsfig{file=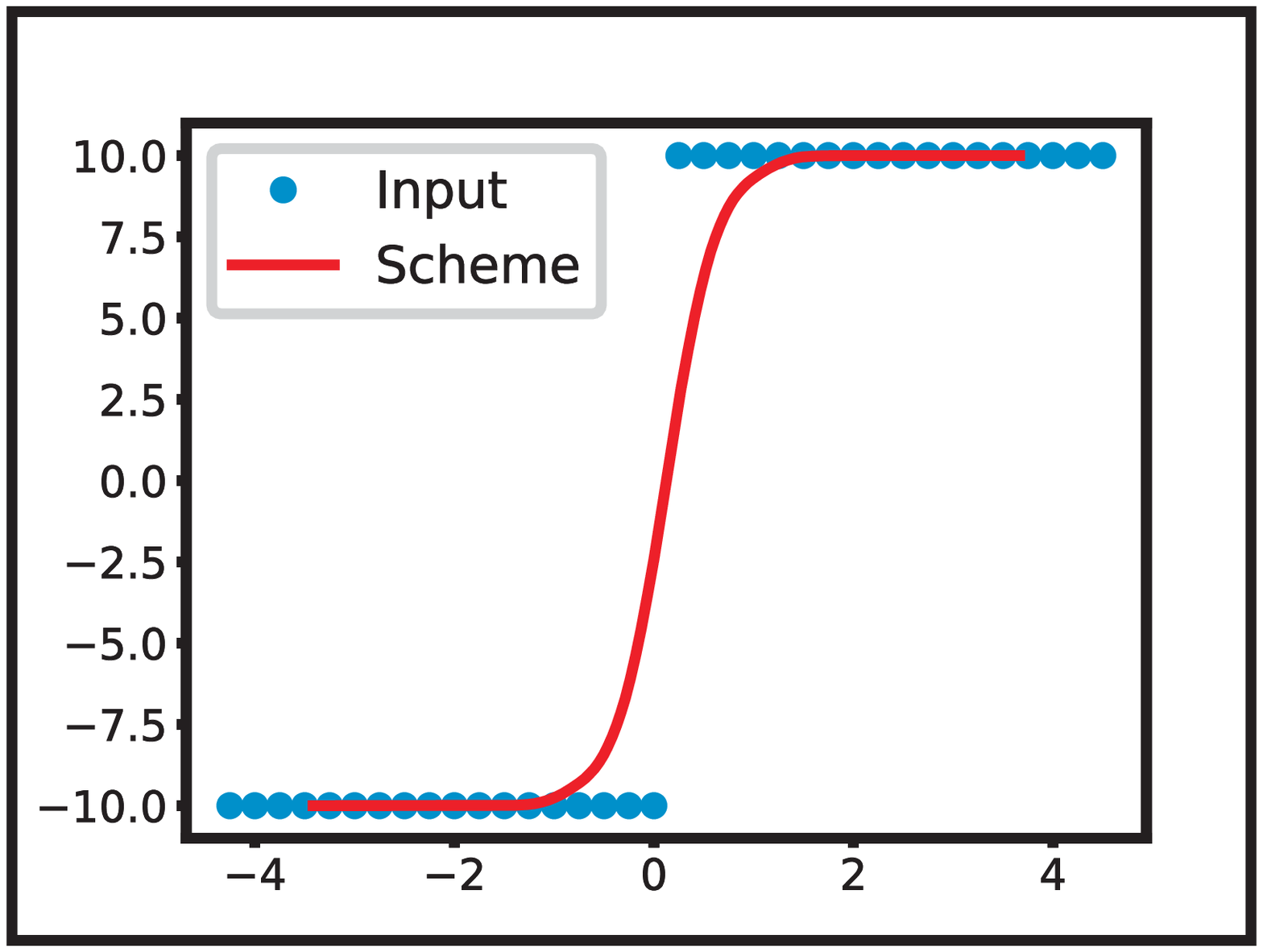, width=1.85 in}&
\epsfig{file=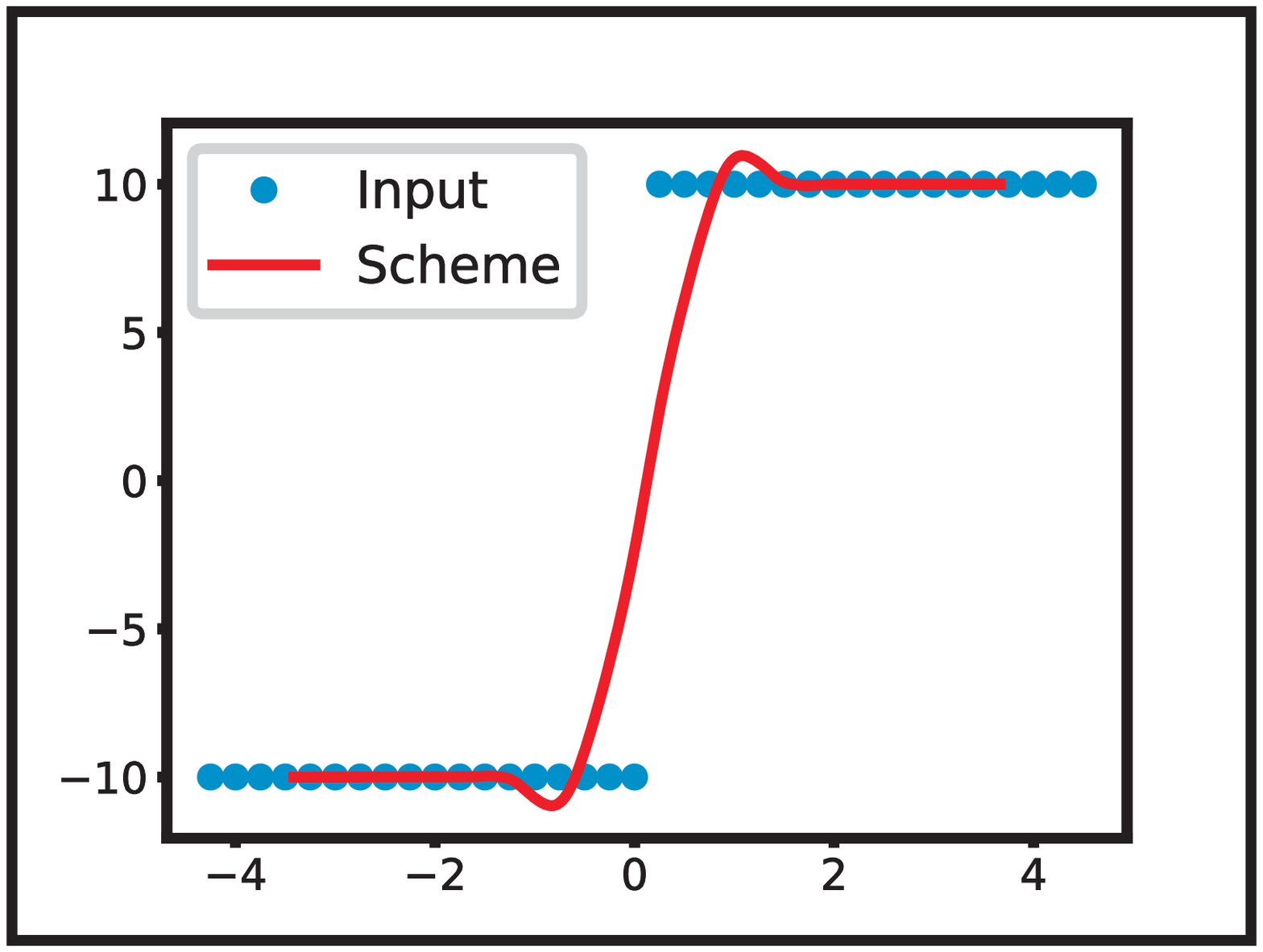, width=1.85 in}&\quad \\
(d) $D_{11,1}$ & (e) $D_{11,2}$ & (f) $D_{11,3}$
 \end{tabular}
\end{center}
\caption[Effects of the schemes $D_{h,d}$ on discontinuous data.]{\label{discontinuous-function} \emph{Effects of the schemes $D_{h,d}$, where $h \in \{2n,2n+1:n=6\}$ and $1 \leqslant d \leqslant 3$, on discontinuous data.}}
\end{figure}
\end{Exp}
\begin{Exp}
\textbf{Interpolating behavior to the non-noisy data}\\
The higher complexity schemes are more suitable to handle with noisy data \cite{Mustafa9}. It is observed that the schemes generated by the higher degree polynomials, i.e $D_{15,2}$, $D_{16,2}$, $D_{15,3}$ and $D_{16,3}$, show interpolatory behavior on non-noisy data points, while the schemes generated by the linear polynomial, i.e. $D_{15,1}$ and $D_{16,1}$, show approximating behavior (see Figure \ref{sharpe-corners}).
\begin{figure}[tbp] 
\begin{center}
\begin{tabular}{ccccccccccc}
\epsfig{file=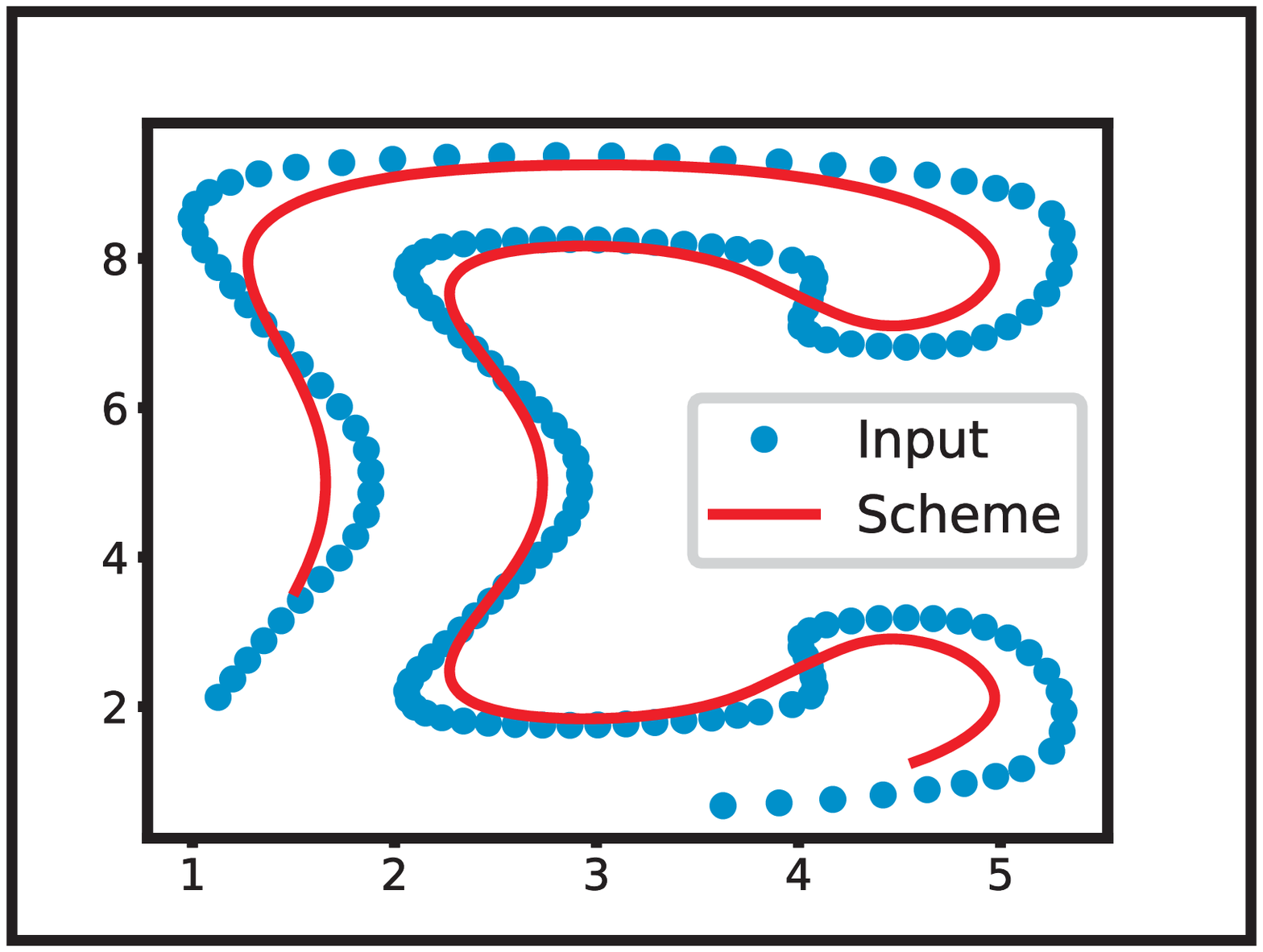, width=1.85 in}&
\epsfig{file=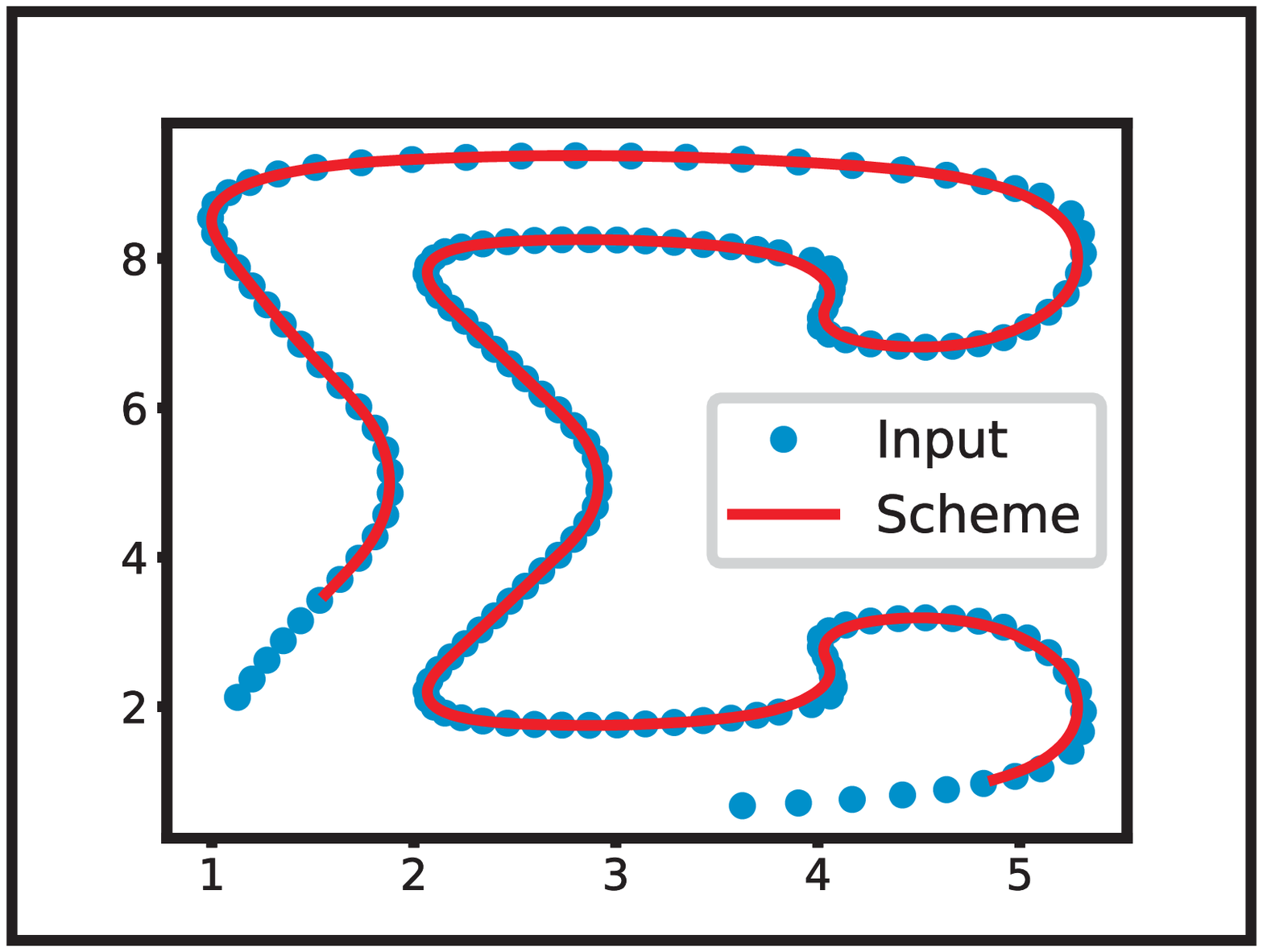, width=1.85 in}&
\epsfig{file=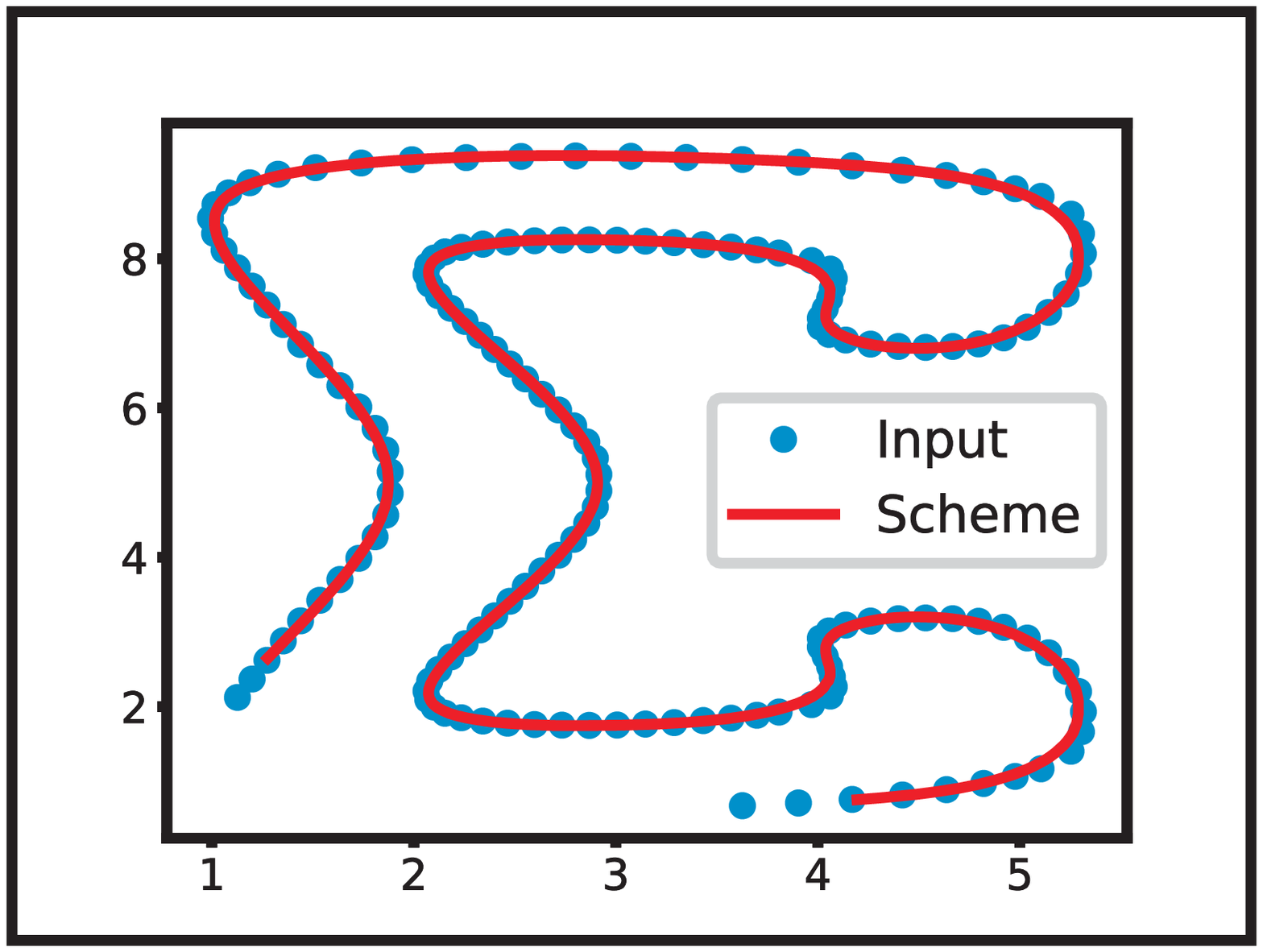, width=1.85 in}&\quad \\
(a) $D_{15,1}$ & (b) $D_{15,2}$ & (c) $D_{15,3}$
 \end{tabular}
\end{center}
 \begin{center}
\begin{tabular}{ccccccccccc}
\epsfig{file=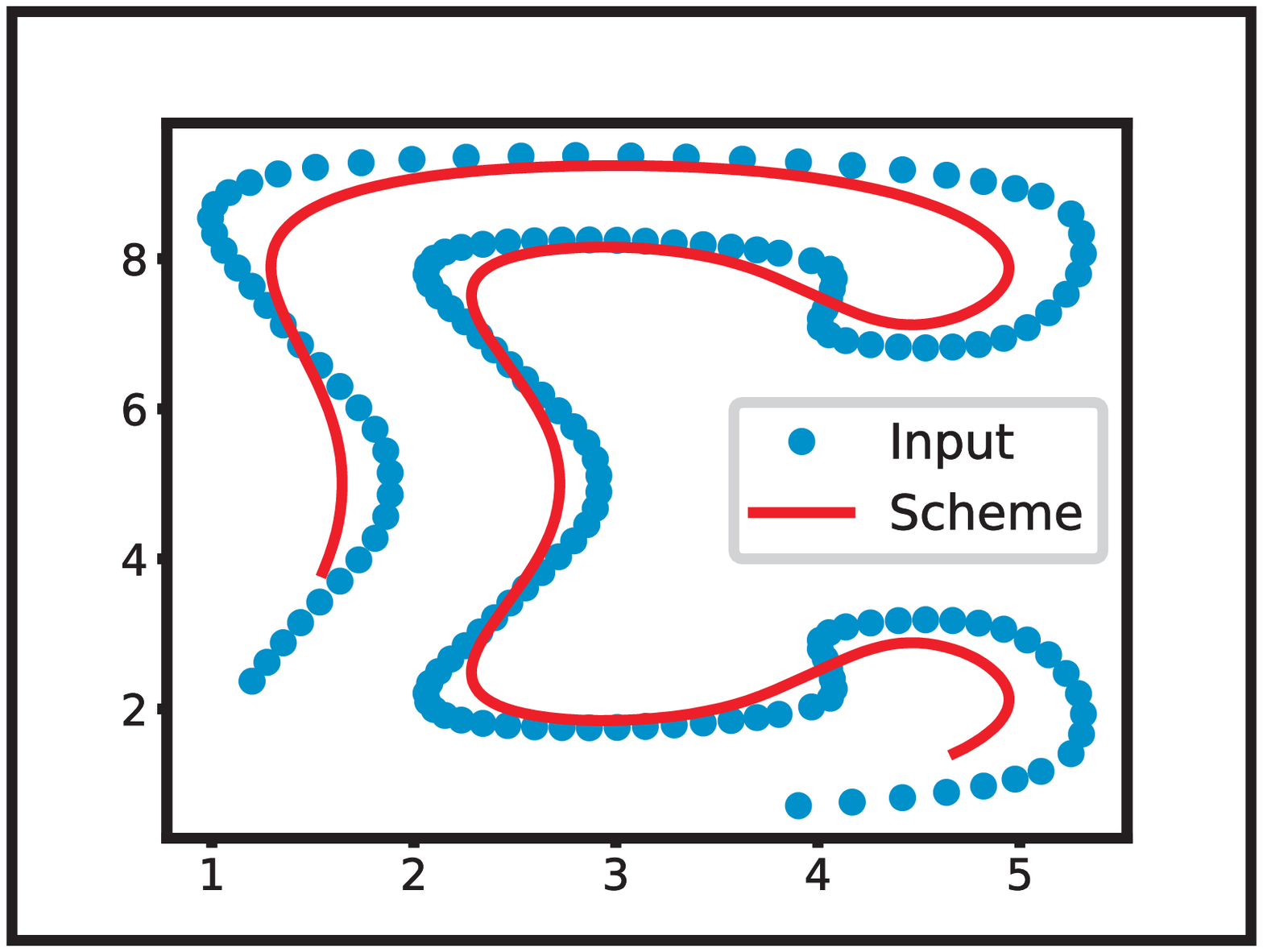, width=1.85 in}&
\epsfig{file=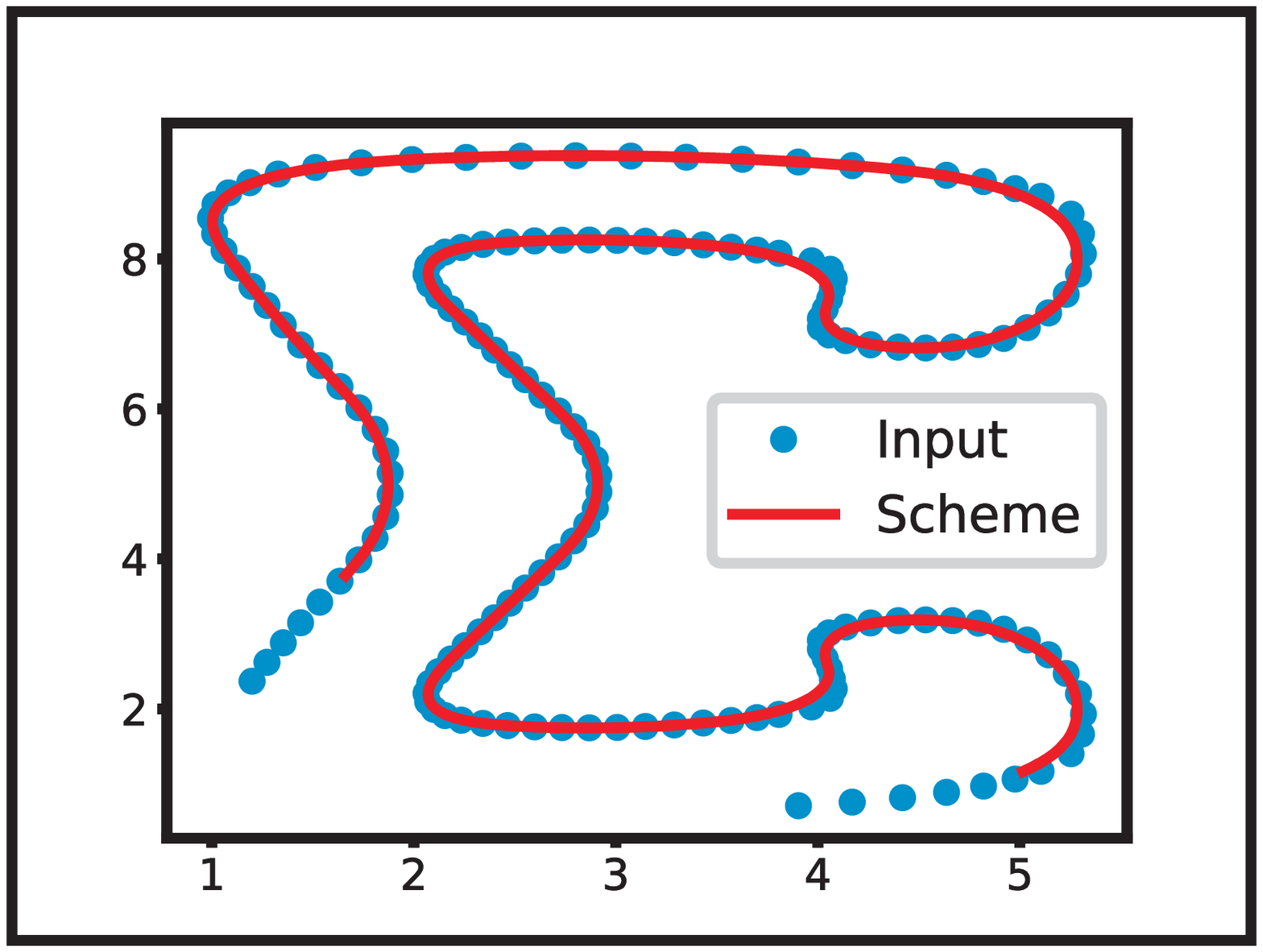, width=1.85 in}&
\epsfig{file=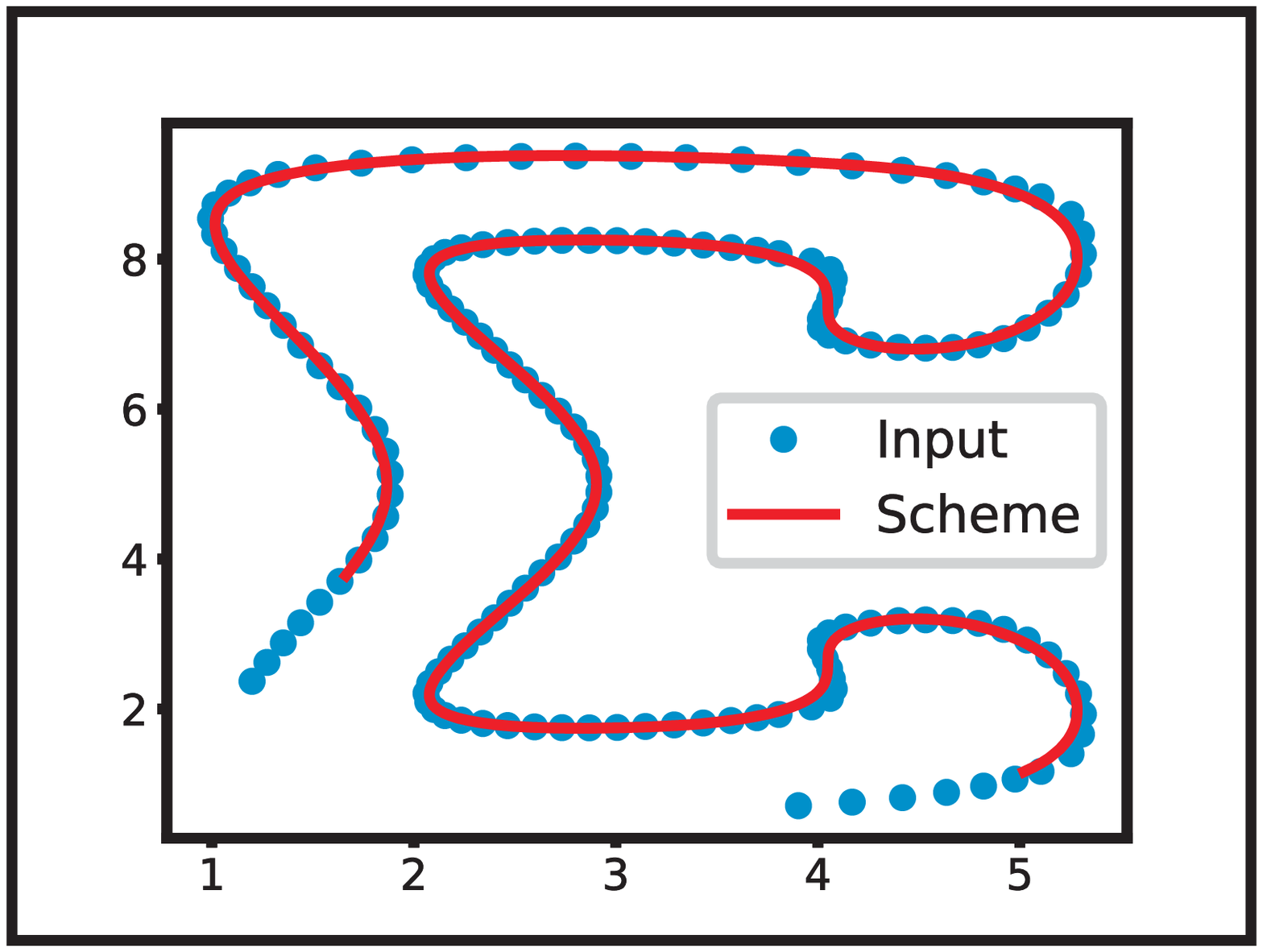, width=1.85 in}&\quad \\
(d) $D_{16,1}$ & (e) $D_{16,2}$ & (f) $D_{16,3}$
 \end{tabular}
\end{center}
\caption[Effects of the schemes $D_{h,d}$ on non-noisy data.]{\label{sharpe-corners} \emph{Curves generated by the proposed schemes.}}
\end{figure}
\end{Exp}

\begin{Exp}
\textbf{Response to the data that contains outliers}\\
In this example, we generate data from the oscillatory function $g_{5}(x)=\left(\frac{x}{40}-1\right)^{3}+cos\left(\frac{2x}{5}\right)$ along with twelve outliers. The limit curves generated by the proposed schemes are presented in Figure \ref{outlier1}. In this figure, blue bullets show the initial control points, while solid lines show fitted curves. This figure shows that the curves generated by the proposed schemes $D_{h,d}$ with $h \geqslant 5$ and $2 \leqslant d \leqslant 3$ are more consistent with the graph of $g_{5}(x)$ than the curves generated by $D_{h,1}$.

\begin{figure}[htb!] 
\begin{center}
\begin{tabular}{ccccccccccc}
\epsfig{file=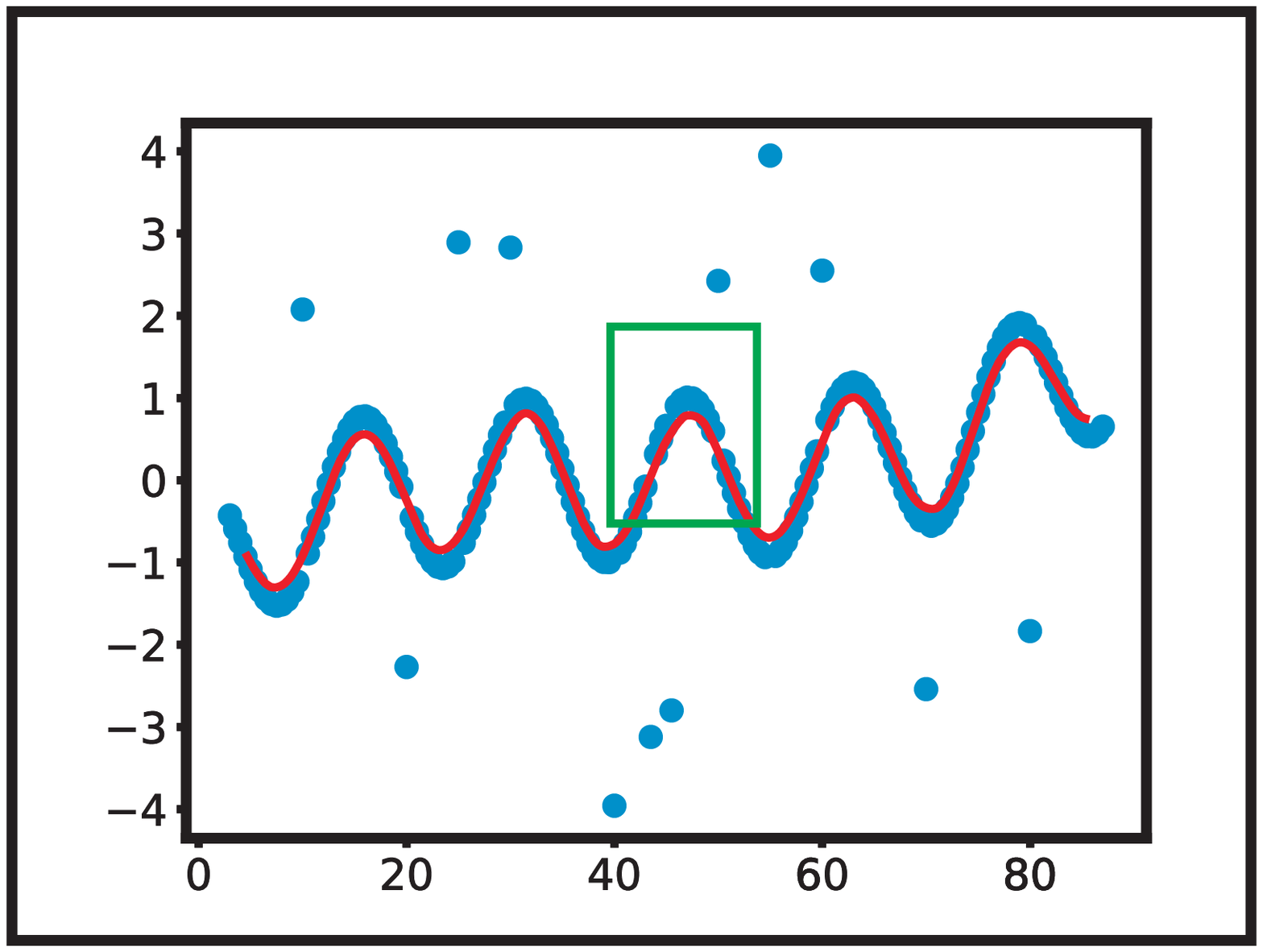, width=1.85 in}&
\epsfig{file=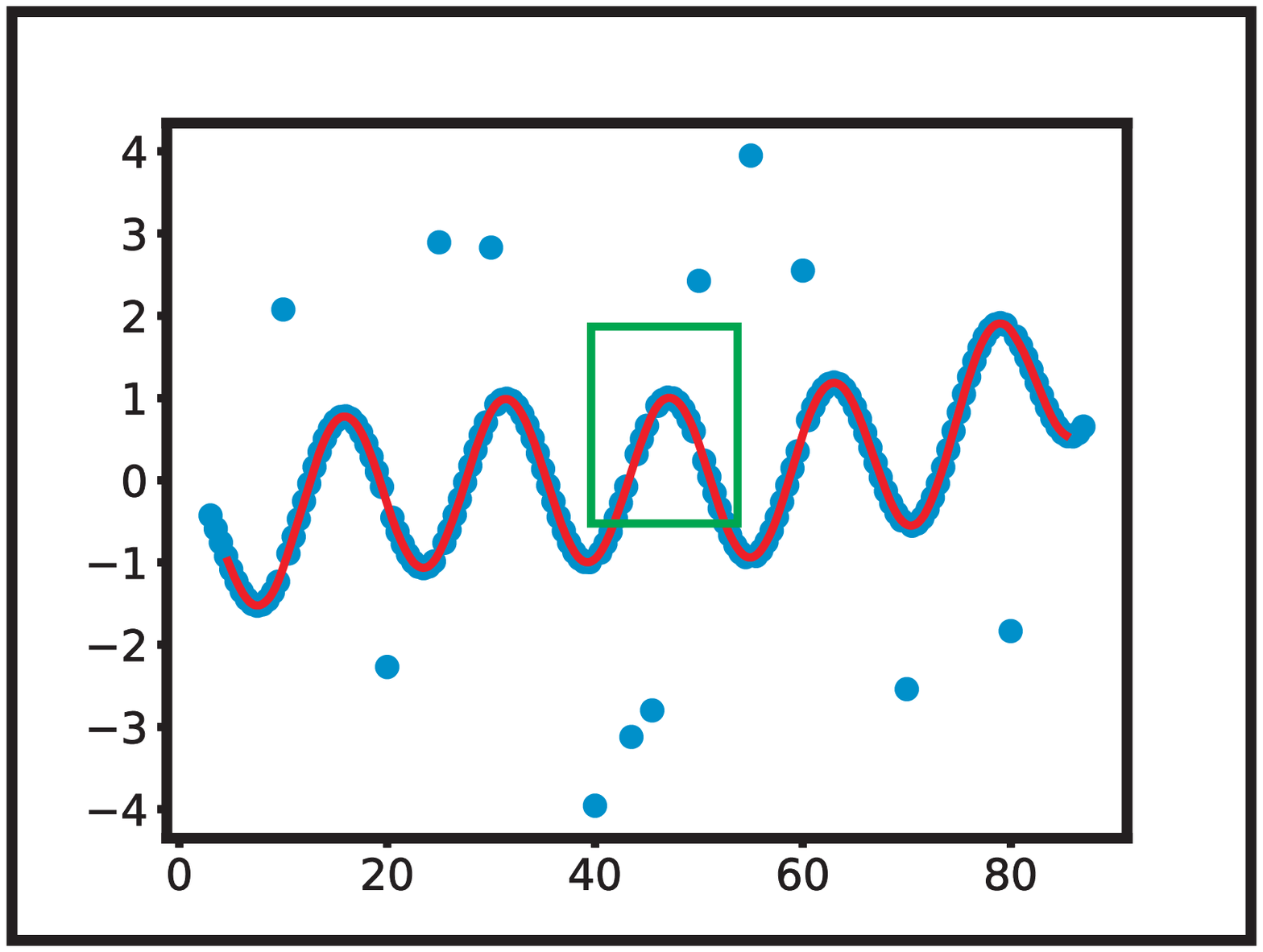, width=1.85 in}&
\epsfig{file=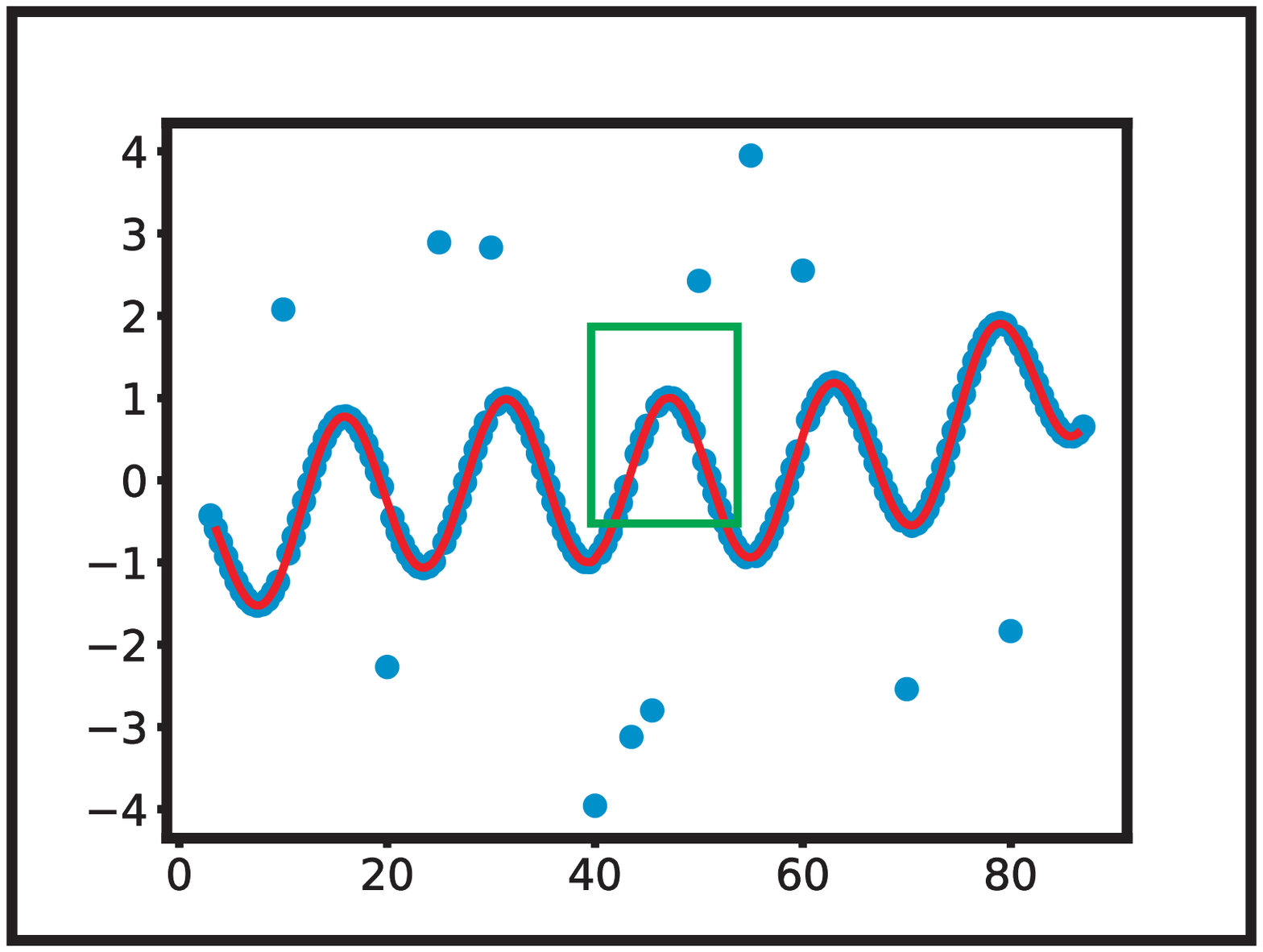, width=1.85 in}&\\
(a) $D_{10,1}$ & (b) $D_{10,2}$ & (c) $D_{10,3}$\\ \\
\epsfig{file=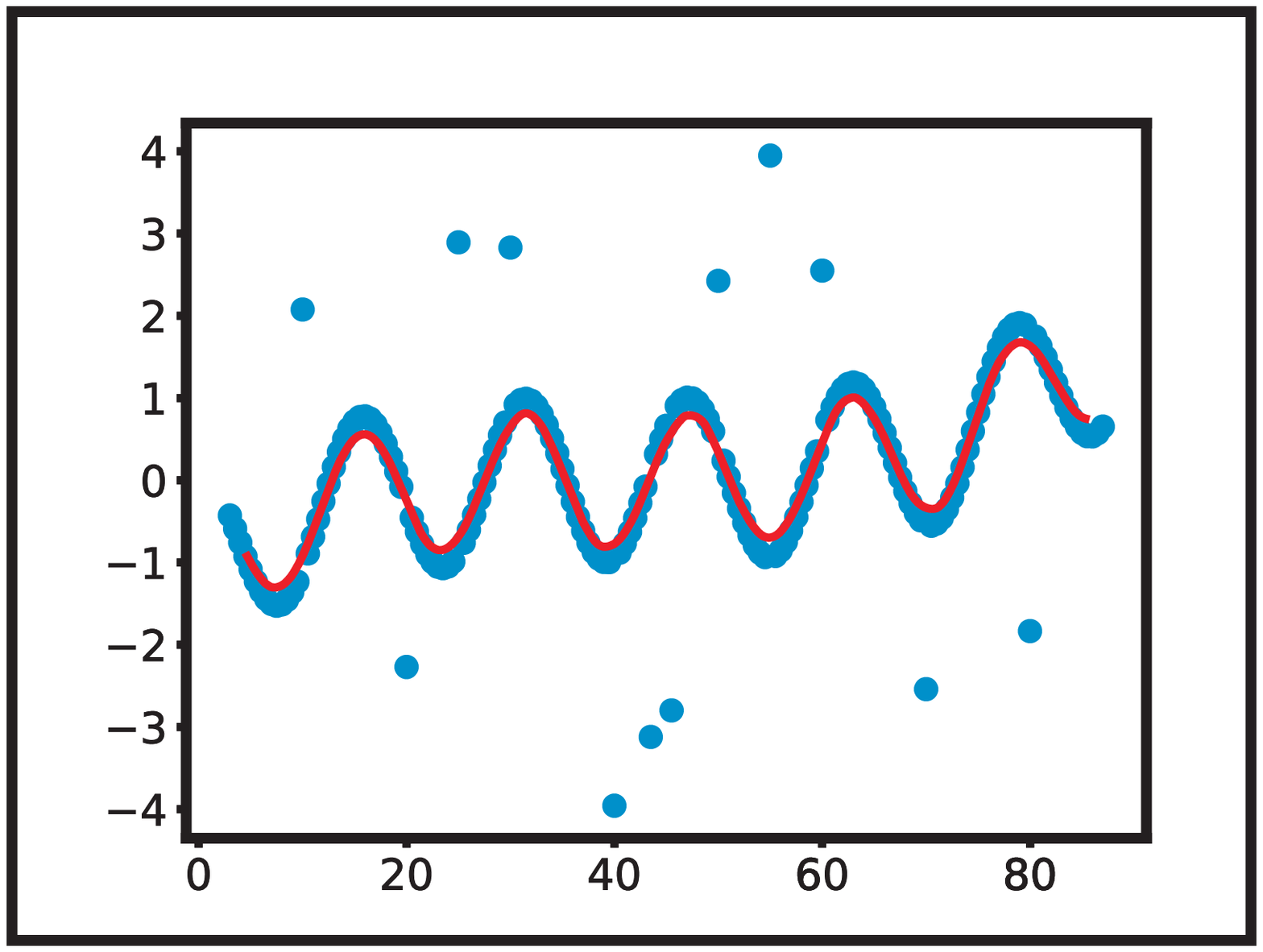, clip,trim=7.7cm 5.8cm 6.42cm 5.0cm, width=1.1 in}&
\epsfig{file=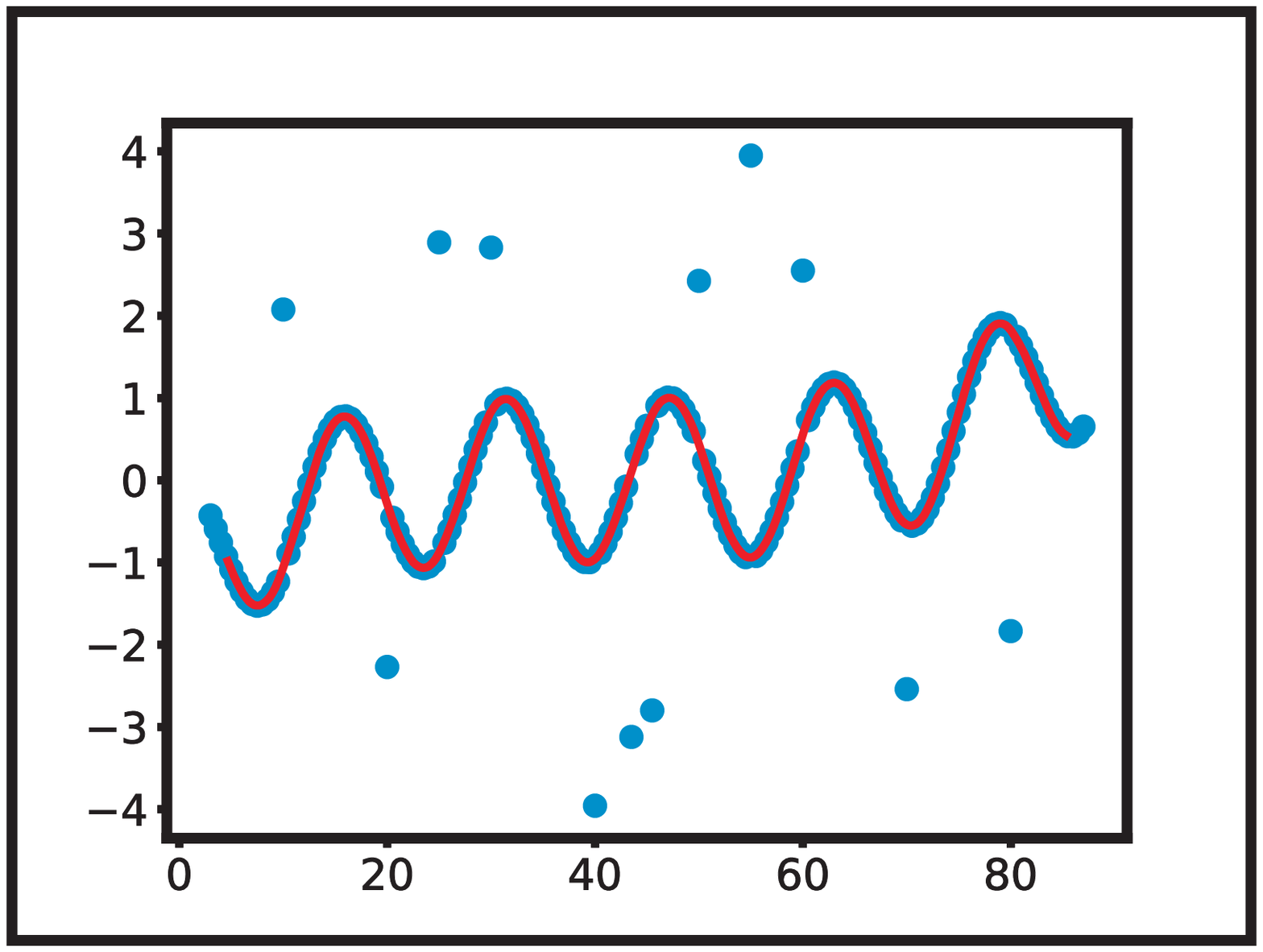, clip,trim=7.7cm 5.8cm 6.42cm 5.0cm, width=1.1 in}&
\epsfig{file=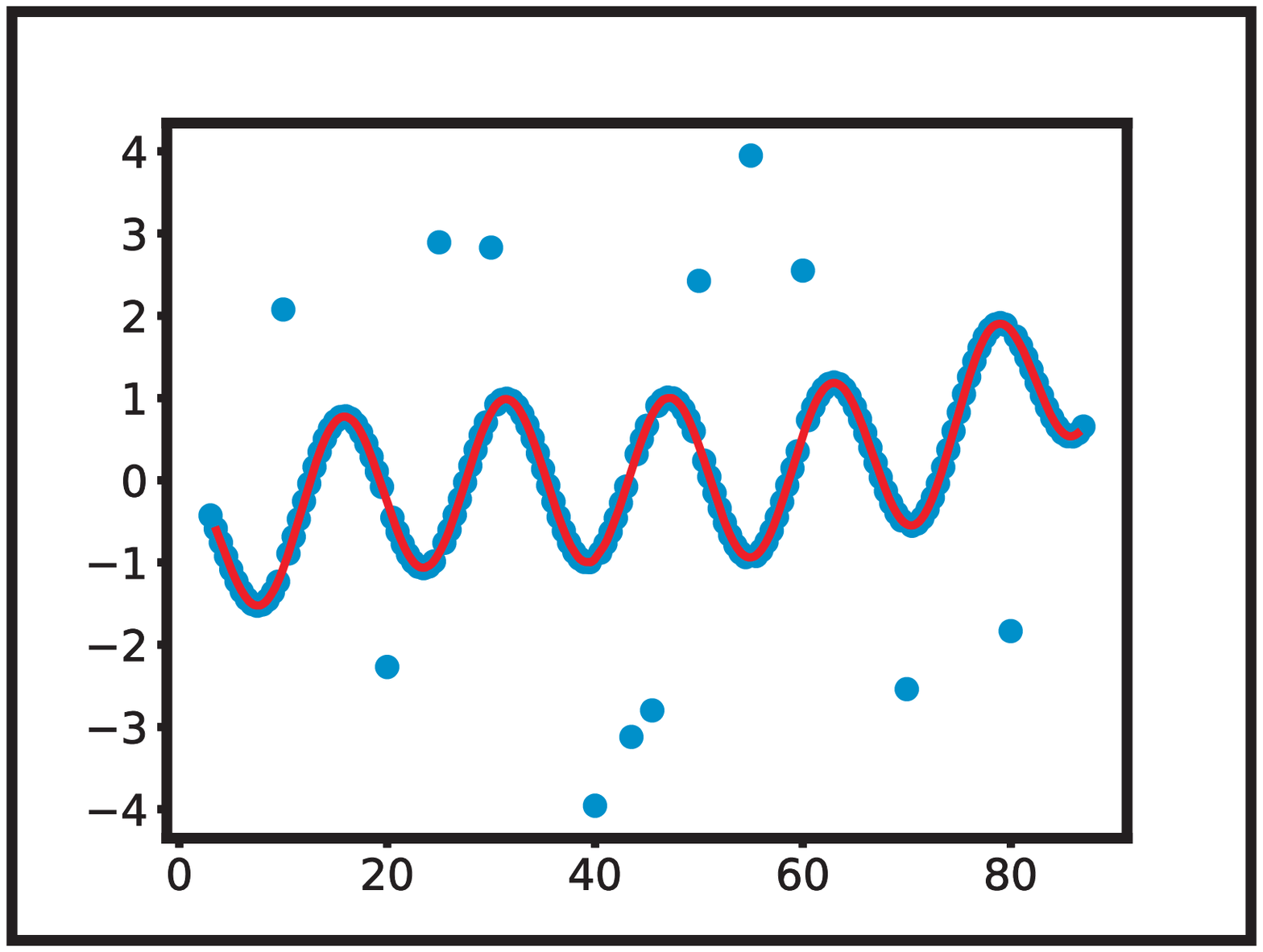, clip,trim=7.7cm 5.8cm 6.42cm 5.0cm, width=1.1 in}&\\
(d) & (e) & (f)
 \end{tabular}
\end{center}
\begin{center}
\begin{tabular}{ccccccccccc}
\epsfig{file=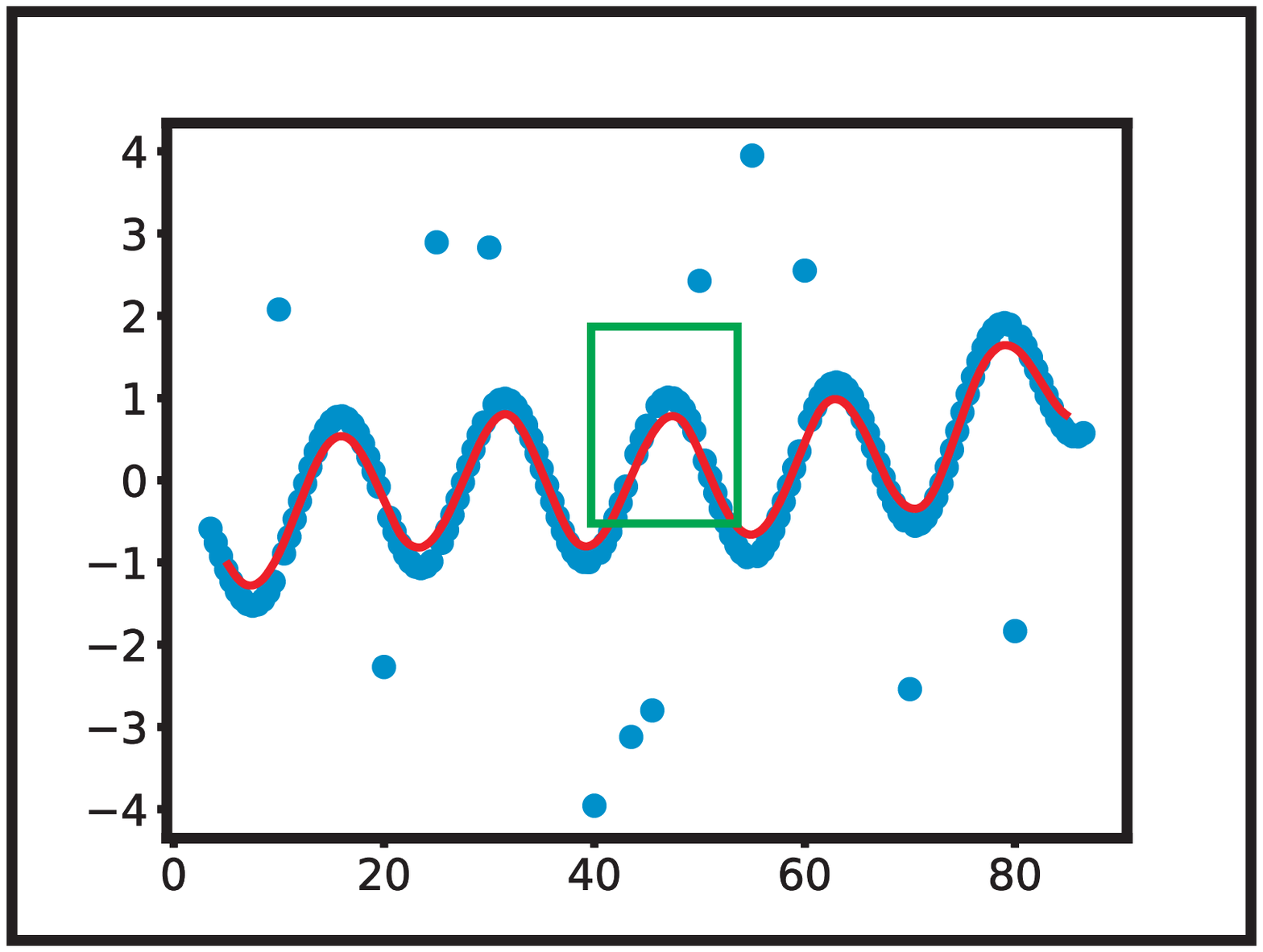, width=1.85 in}&
\epsfig{file=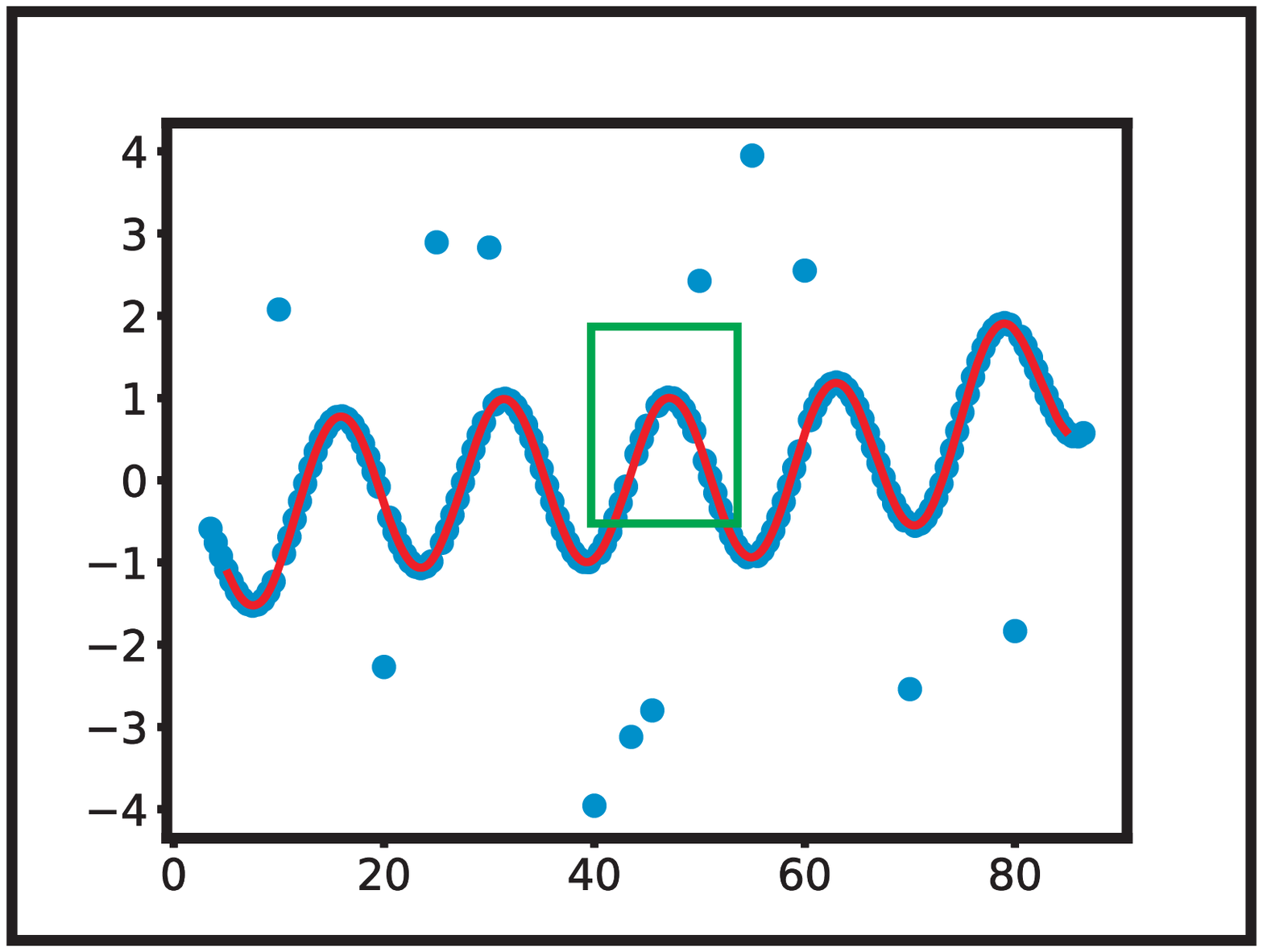, width=1.85 in}&
\epsfig{file=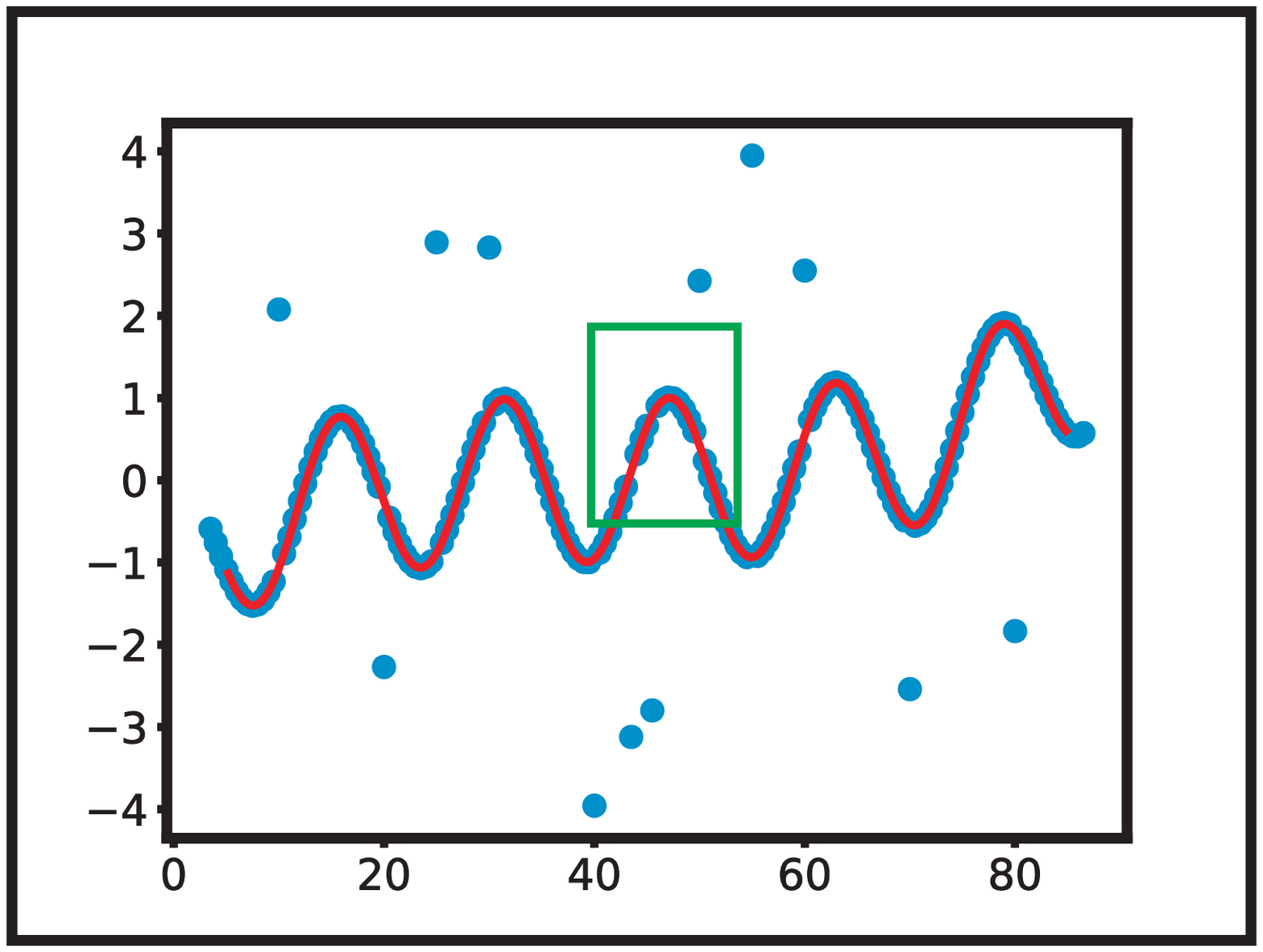, width=1.85 in}&\\
(g) $D_{11,1}$ & (h) $D_{11,2}$ & (i) $D_{11,3}$ \\ \\
\epsfig{file=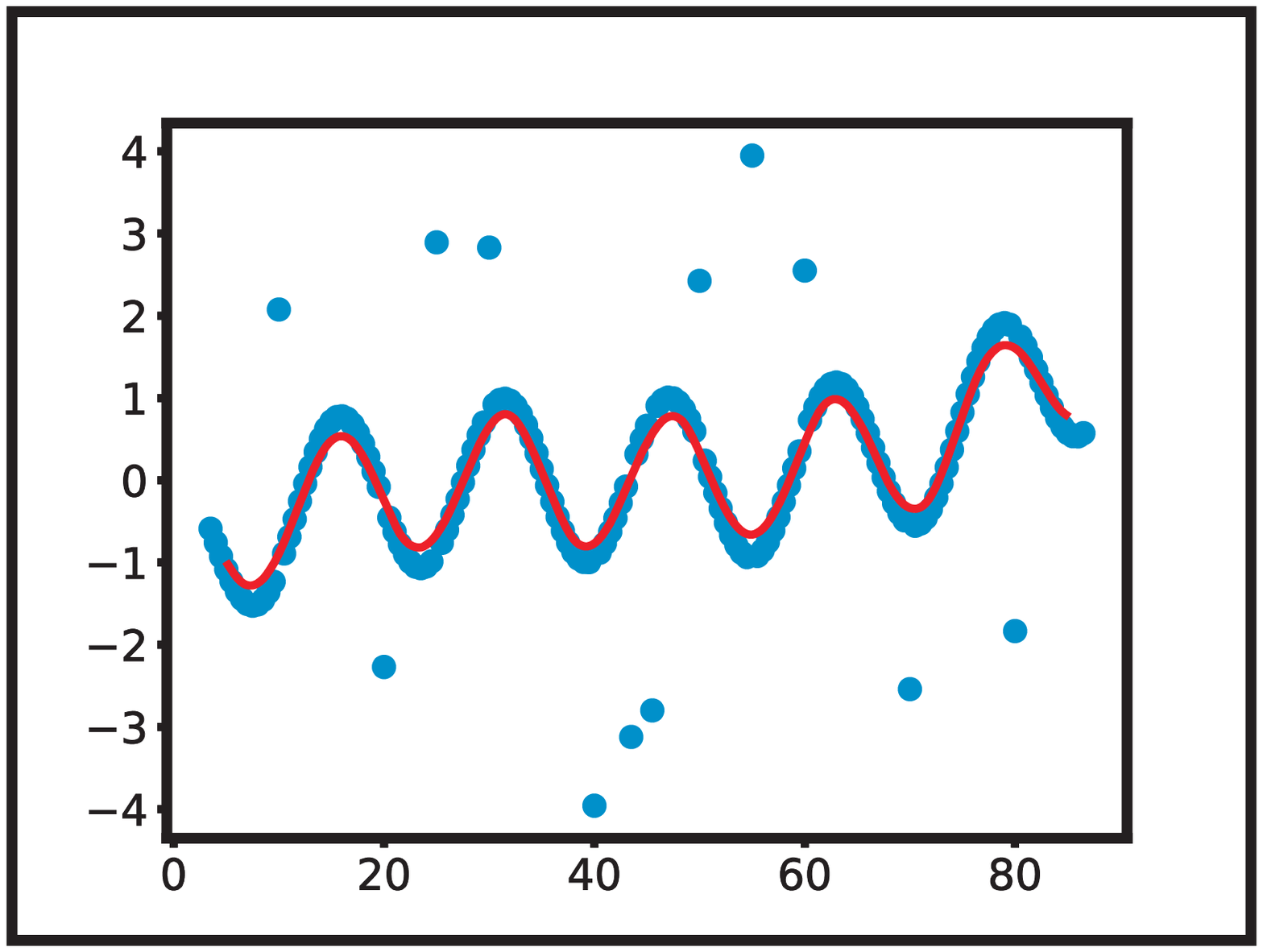, clip,trim=7.7cm 5.8cm 6.42cm 5.0cm, width=1.1 in}&
\epsfig{file=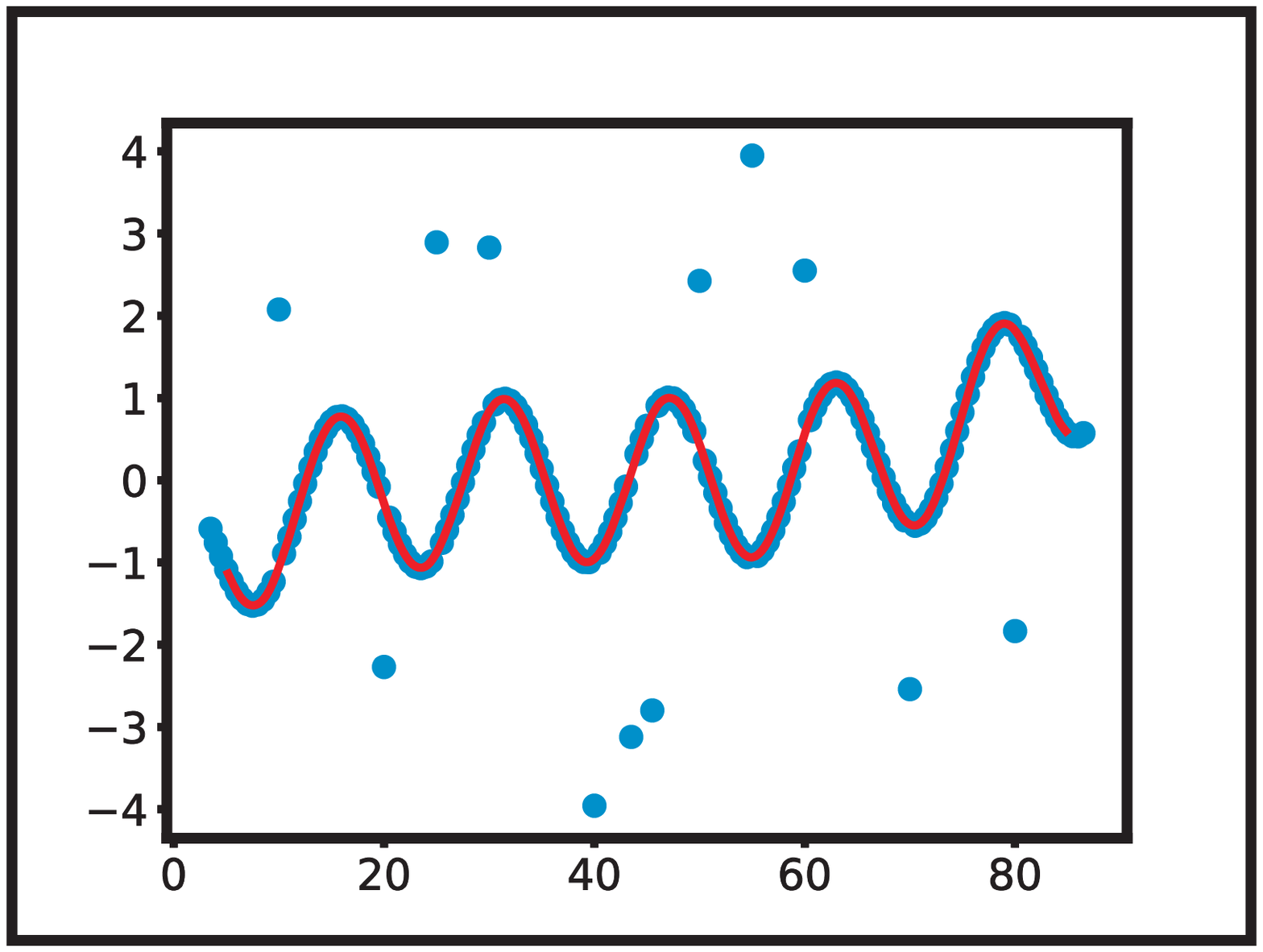, clip,trim=7.7cm 5.8cm 6.42cm 5.0cm, width=1.1 in}&
\epsfig{file=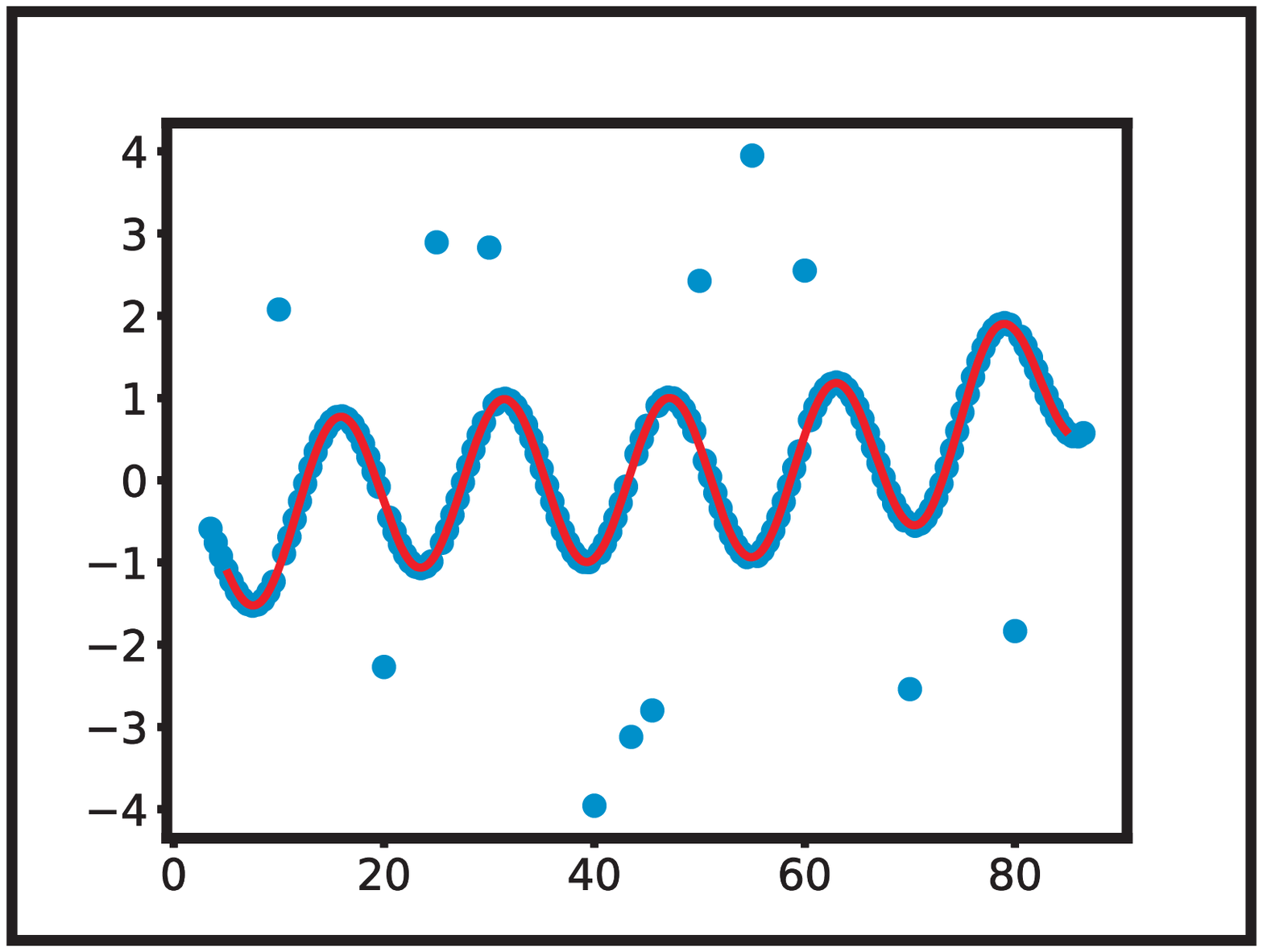, clip,trim=7.7cm 5.8cm 6.42cm 5.0cm, width=1.1 in}&\\
(j) & (k) & (l)
 \end{tabular}
\end{center}
\caption[Effects of the schemes $D_{h,d}$ on the data that contains outliers.]{\label{outlier1} \emph{Effects of the outliers on fitted curves. (d)-(f) and (j)-(l) are the mirror images of the central parts (inside the green rectangles) of (a)-(c) and (g)-(i) respectively.}}
\end{figure}
\end{Exp}

\begin{Exp}
\textbf{Response to the noisy data with outliers}\\
Let us take noisy data from the function $g_{6}(x)=e^{-\frac{x}{3}}sin(3x)$ with two outliers. We want to see the effects of this noisy data with outliers on the fitted curves. From Figure \ref{Noise-outlier1}, it is to be observed that the fitted curves generated by the proposed schemes $D_{19,1}$, $D_{19,2}$, $D_{19,3}$, $D_{20,1}$, $D_{20,2}$ and $D_{20,3}$ are not effected by outliers and noise.

\begin{figure}[htb!] 
\begin{center}
\begin{tabular}{ccccccccccc}
\epsfig{file=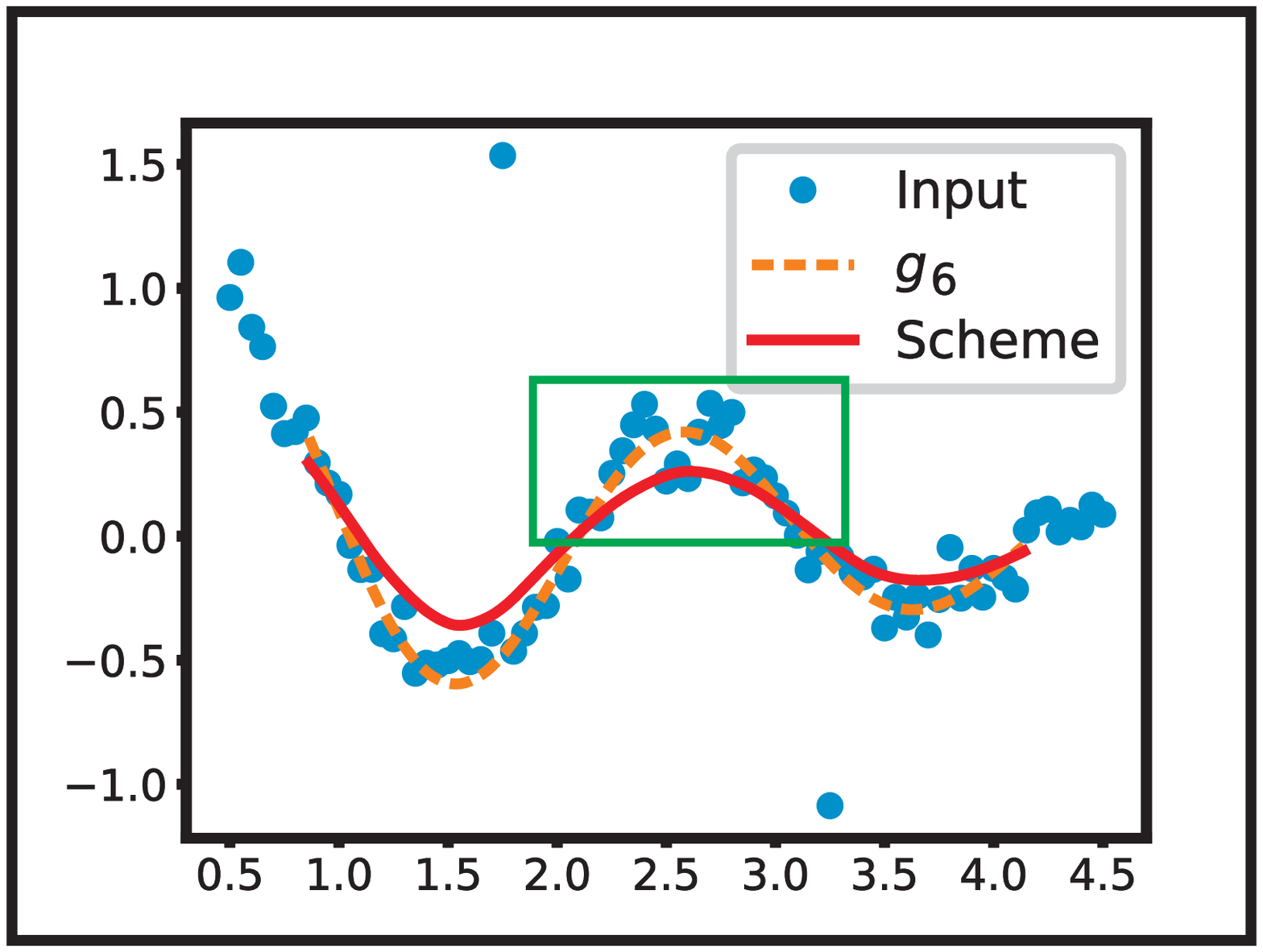, width=1.85 in}&
\epsfig{file=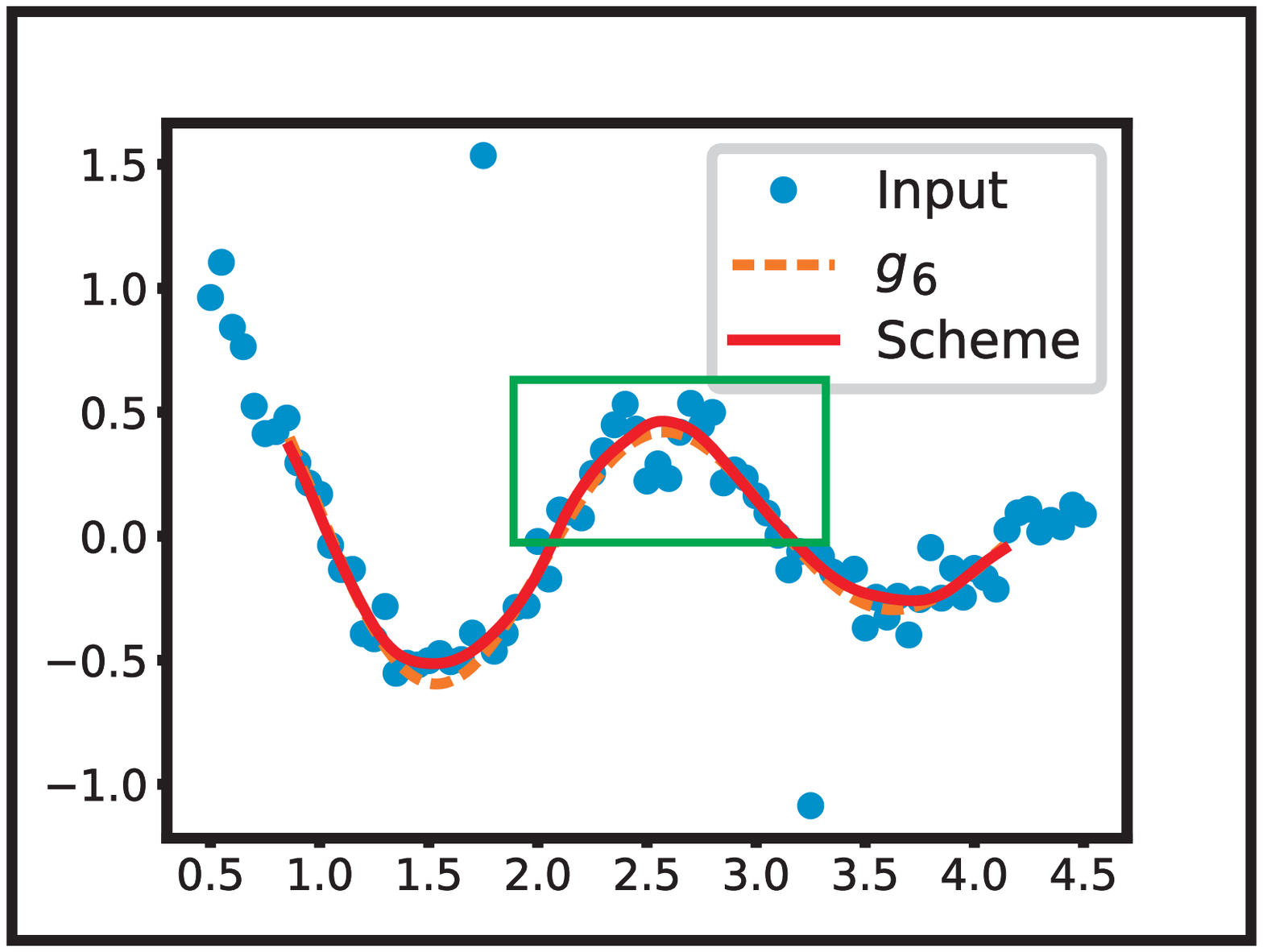, width=1.85 in}&
\epsfig{file=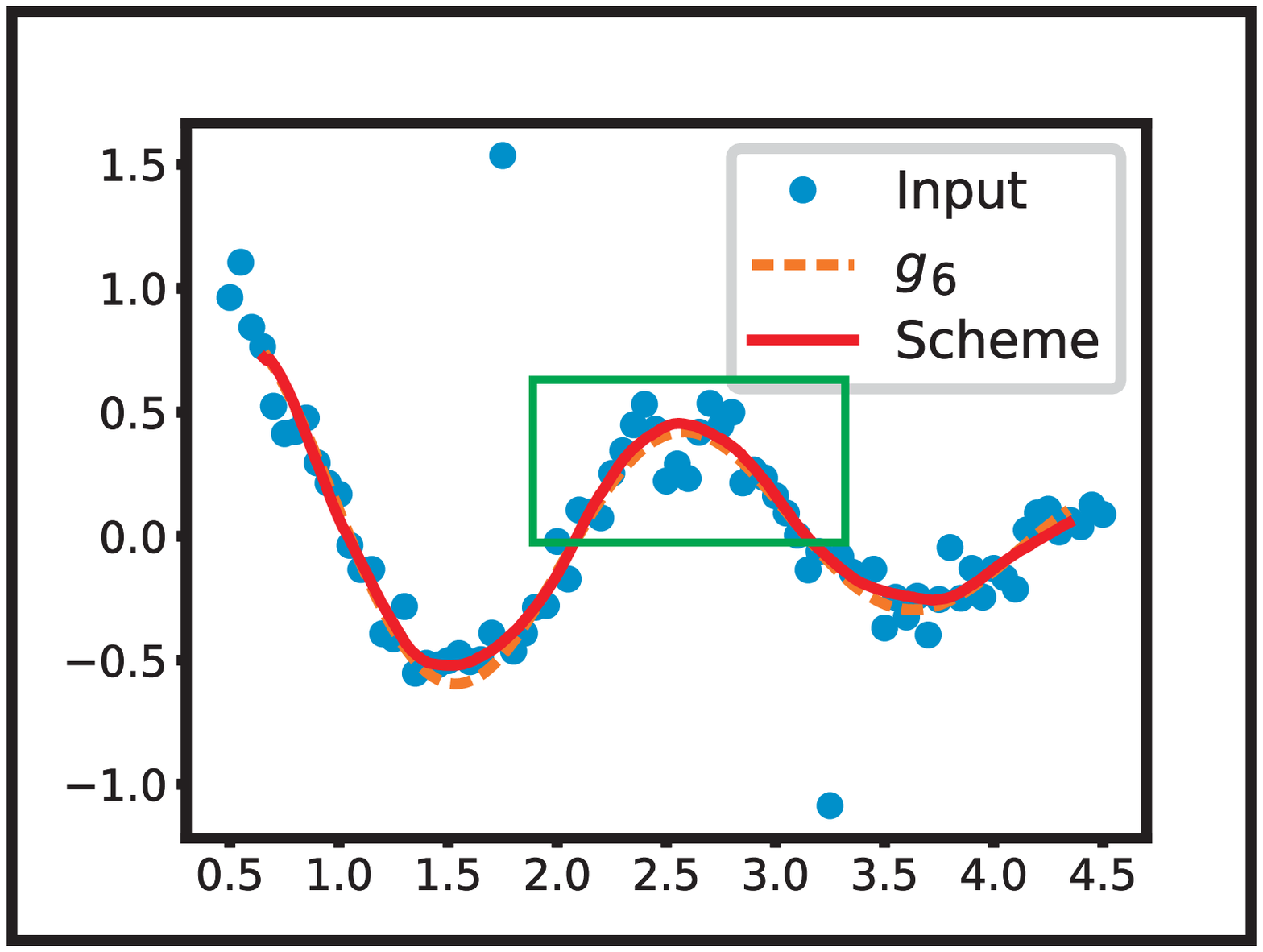, width=1.85 in}&\quad \\
(a) $D_{19,1}$ & (b) $D_{19,2}$ & (c) $D_{19,3}$\\  \\
\epsfig{file=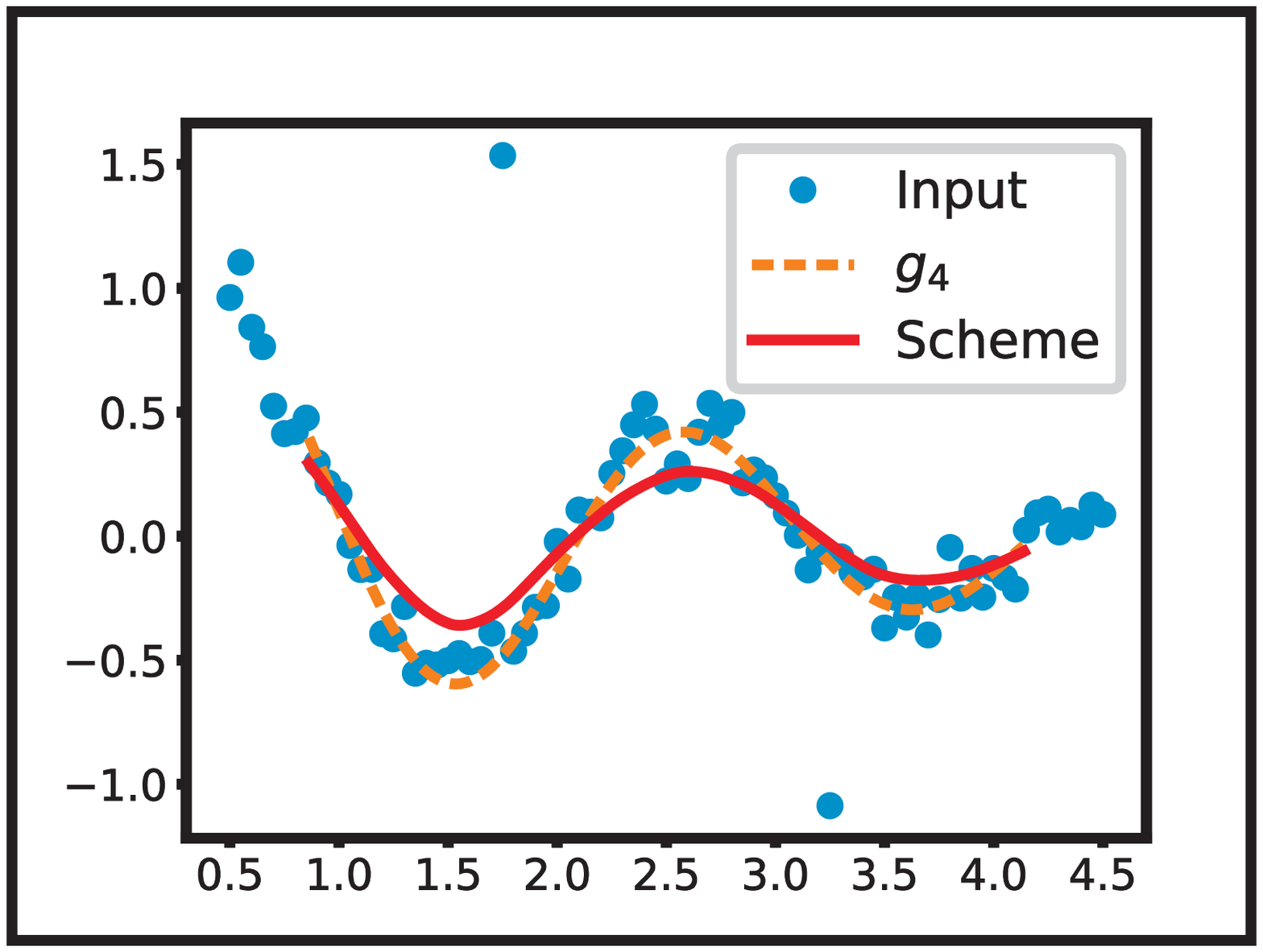, clip,trim=7cm 5.5cm 5.2cm 5.13cm, width=1.2 in}&
\epsfig{file=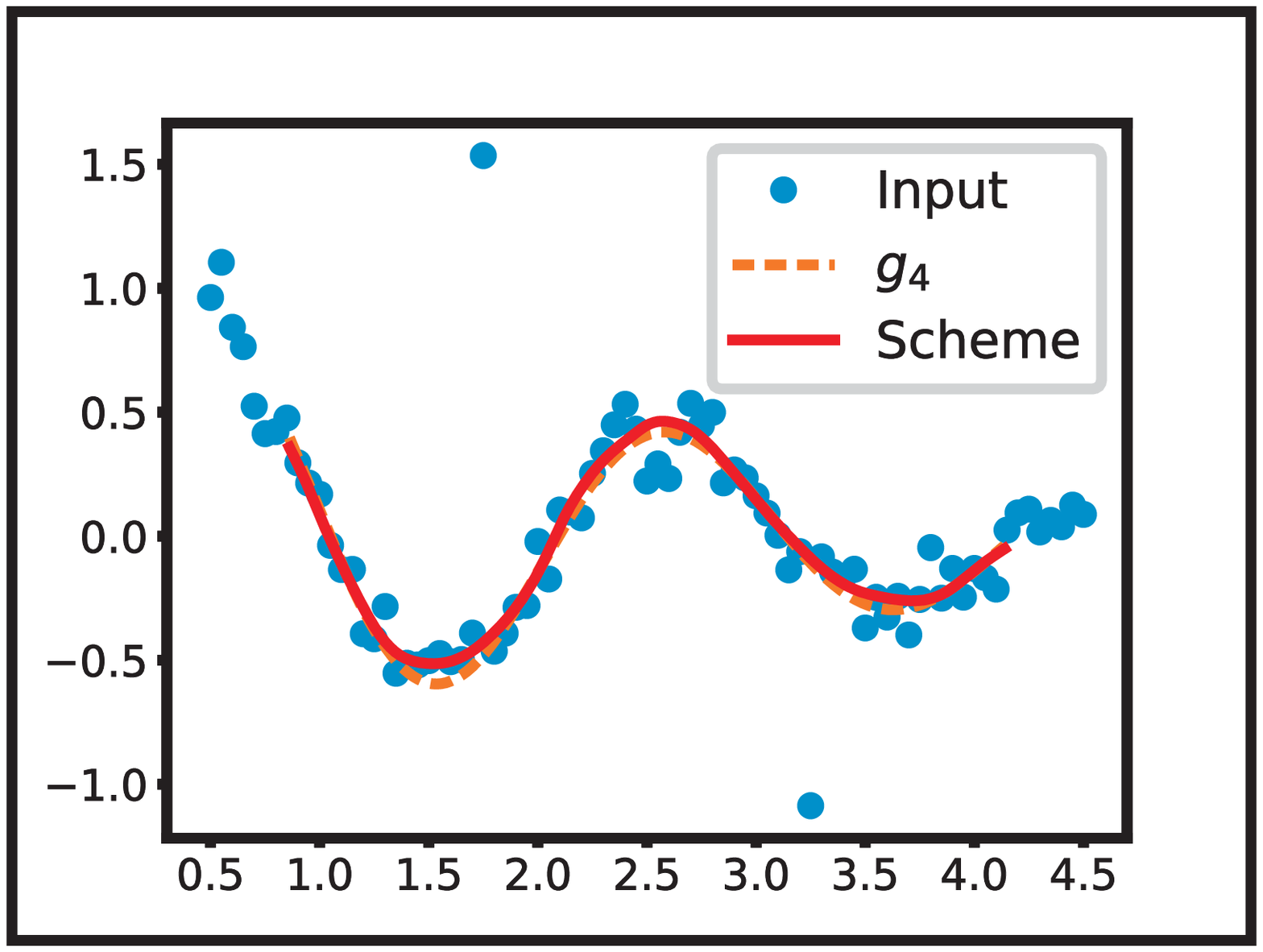, clip,trim=7cm 5.5cm 5.2cm 5.13cm, width=1.2 in}&
\epsfig{file=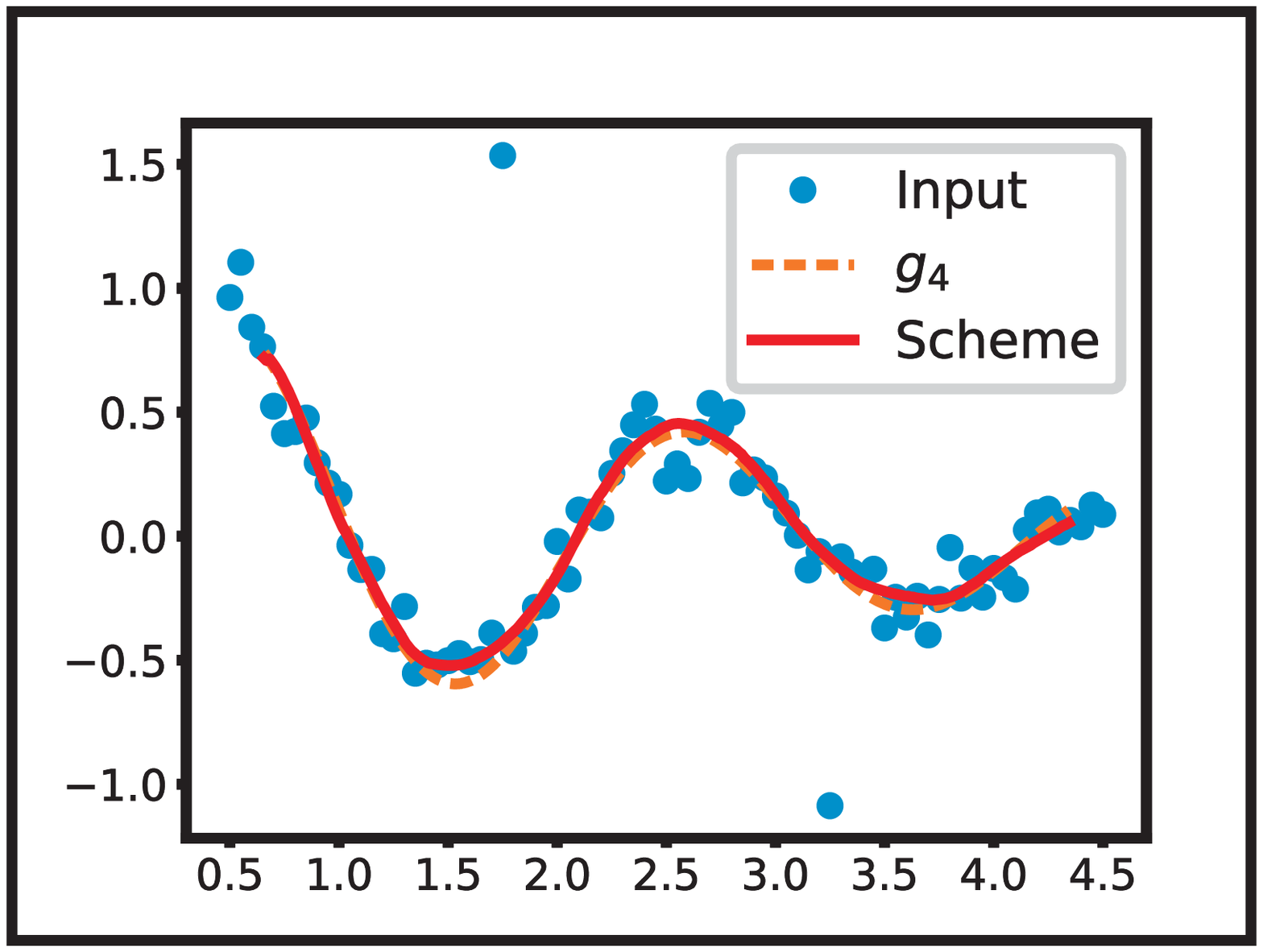, clip,trim=7cm 5.5cm 5.2cm 5.13cm, width=1.2 in}&\\
(d) & (e) & (f)
 \end{tabular}
\end{center}
\begin{center}
\begin{tabular}{ccccccccccc}
\epsfig{file=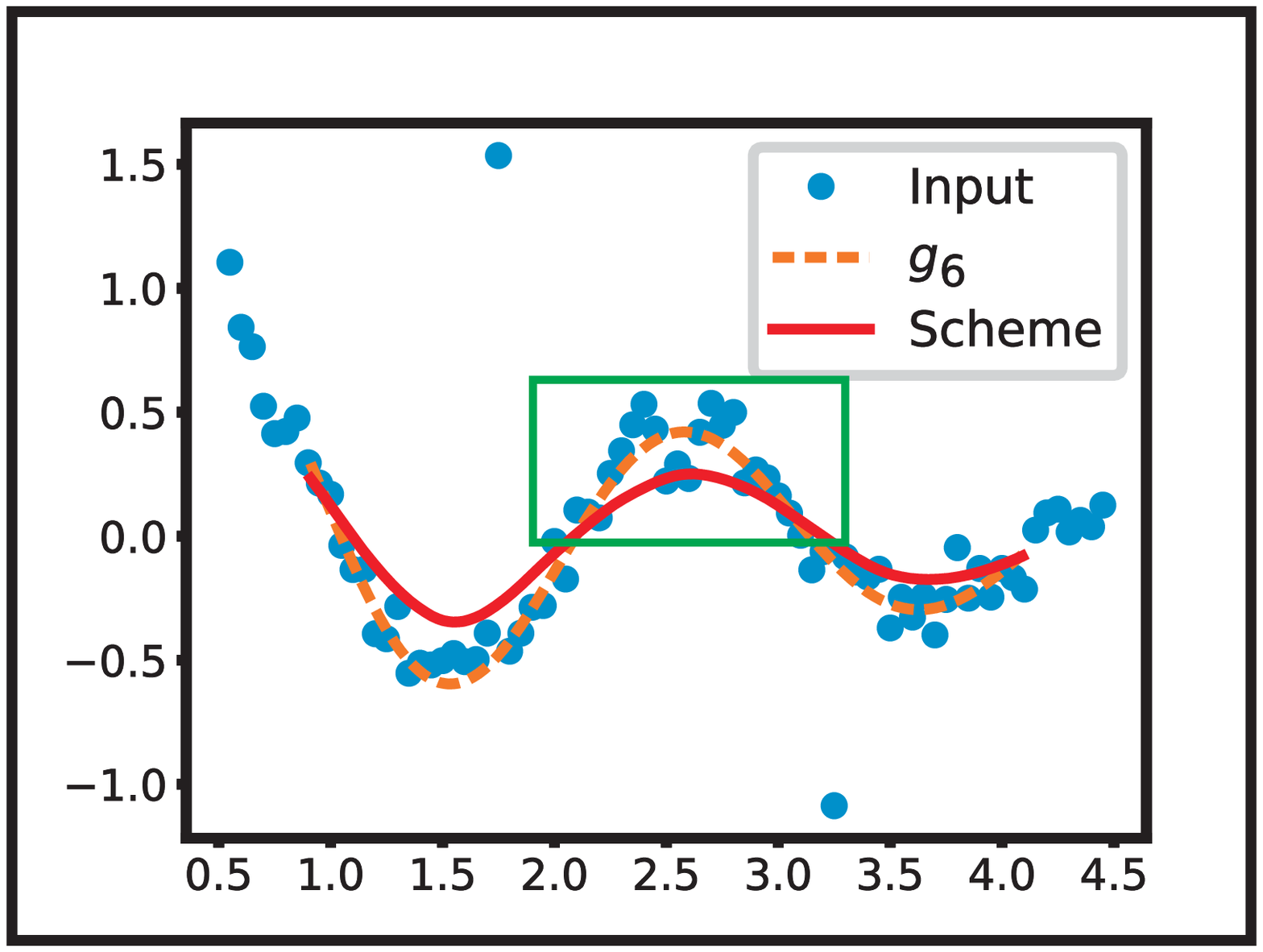, width=1.85 in}&
\epsfig{file=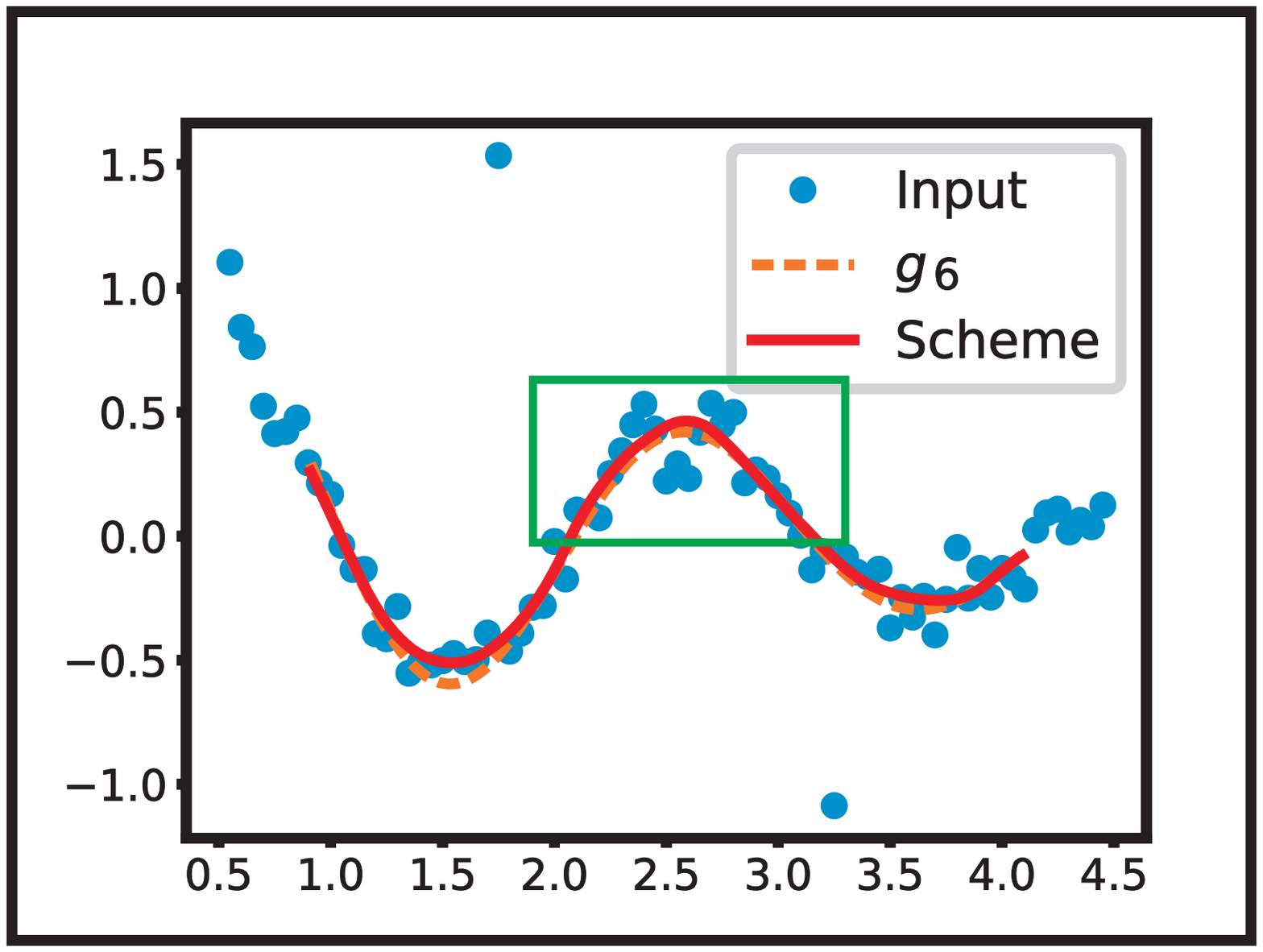, width=1.85 in}&
\epsfig{file=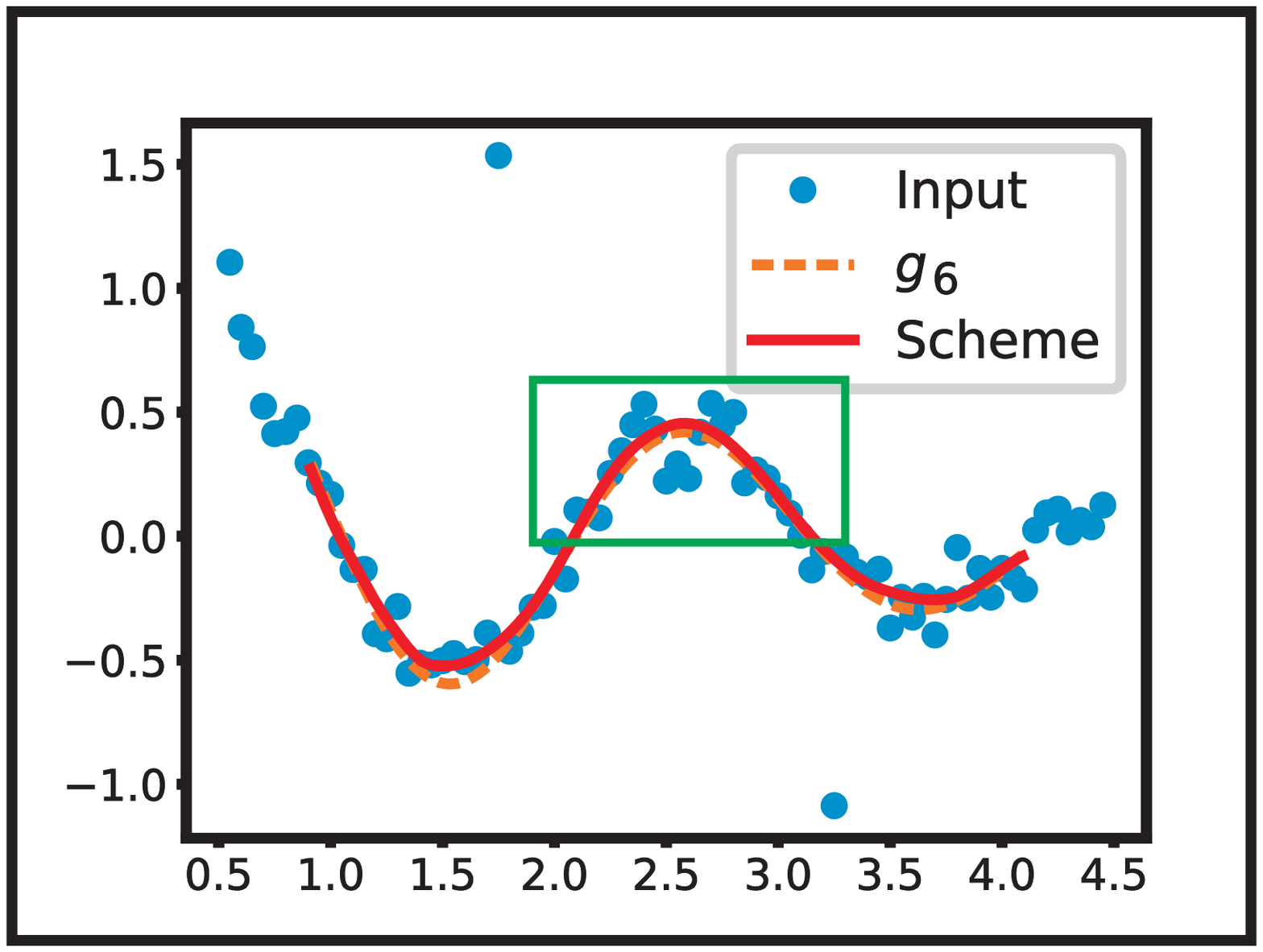, width=1.85 in}&\\
(g) $D_{20,1}$ & (h) $D_{20,2}$ & (i) $D_{20,3}$\\ \\
\epsfig{file=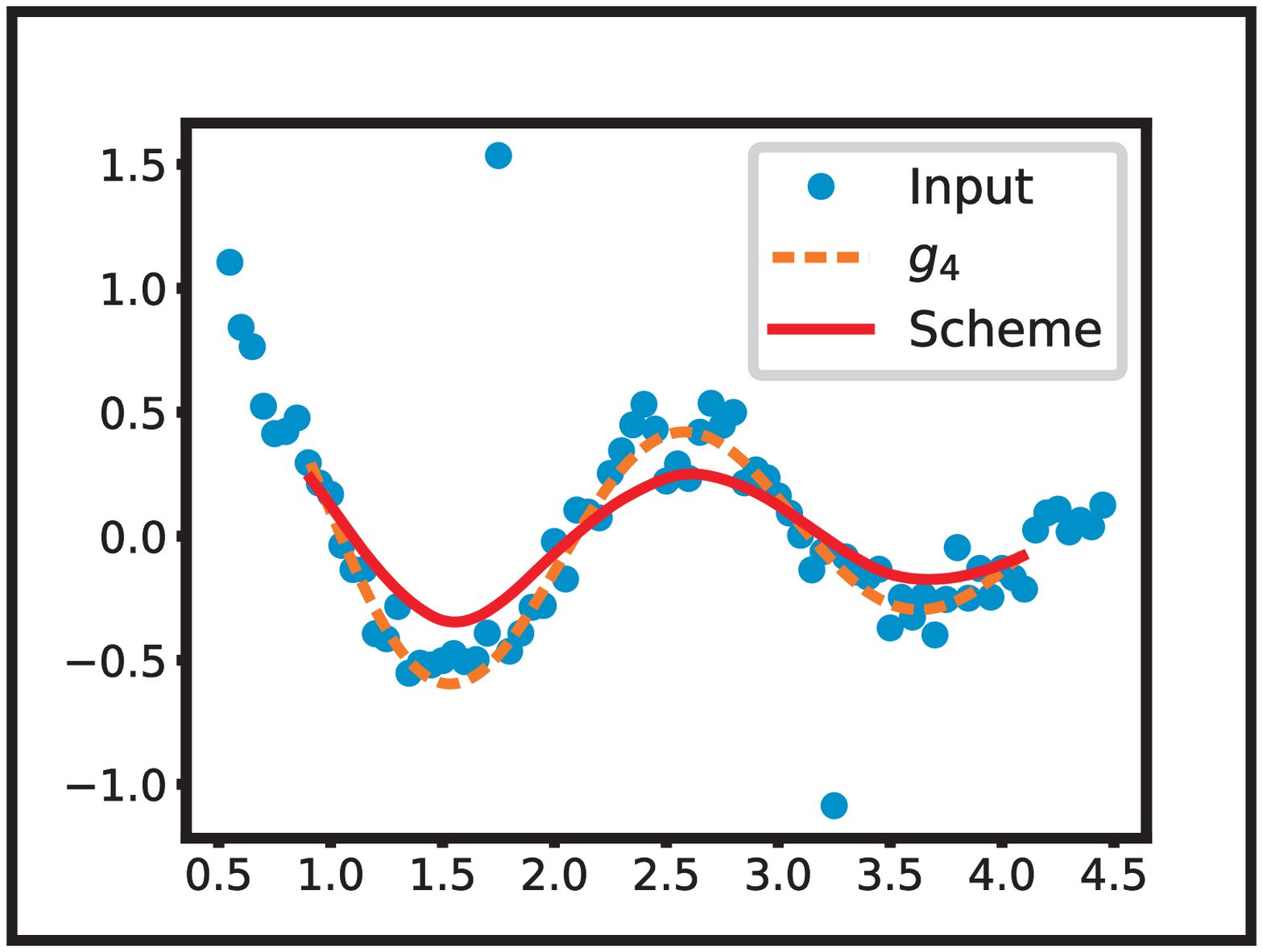, clip,trim=7cm 5.5cm 5.2cm 5.13cm, width=1.2 in}&
\epsfig{file=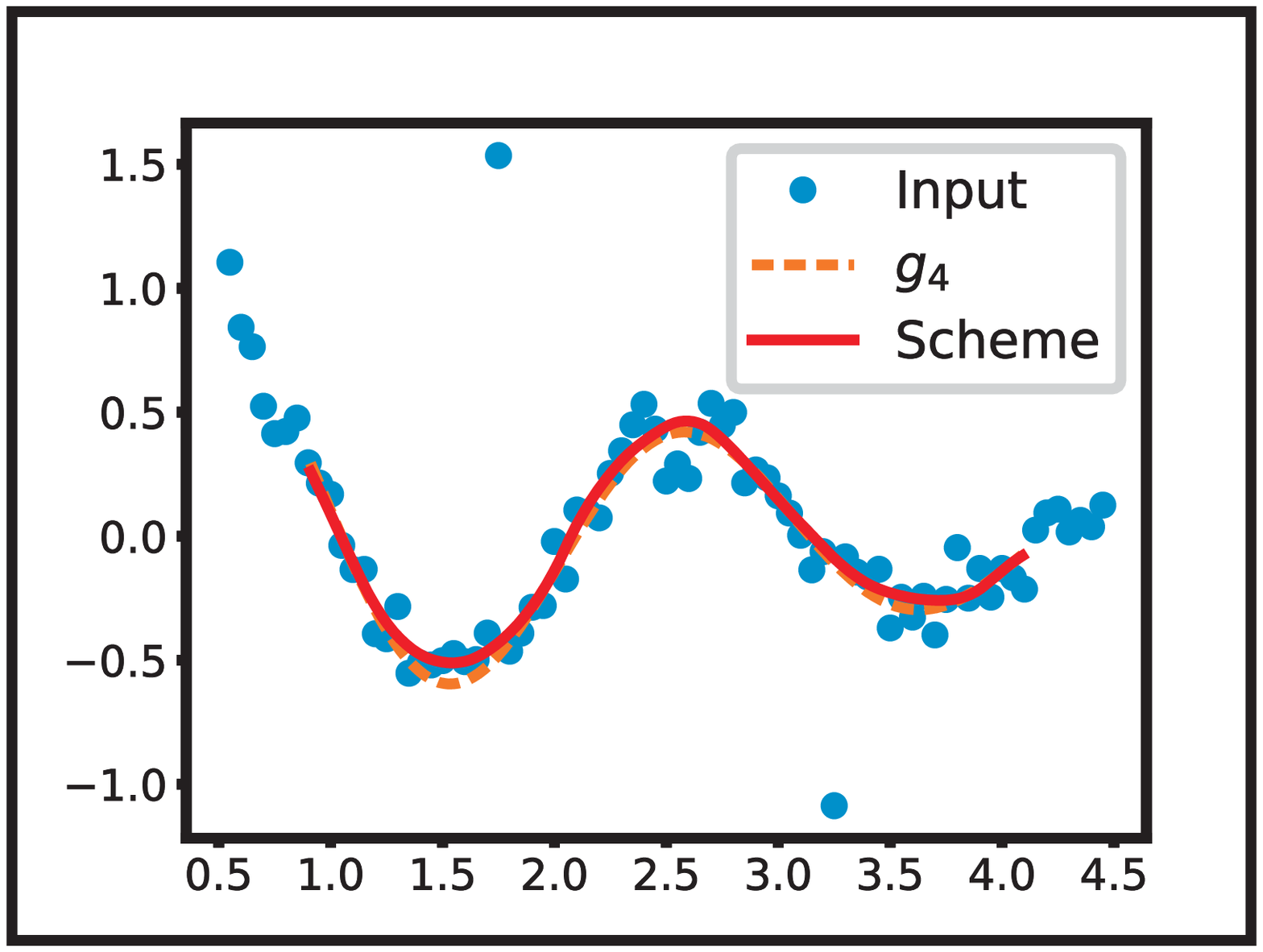, clip,trim=7cm 5.5cm 5.2cm 5.13cm, width=1.2 in}&
\epsfig{file=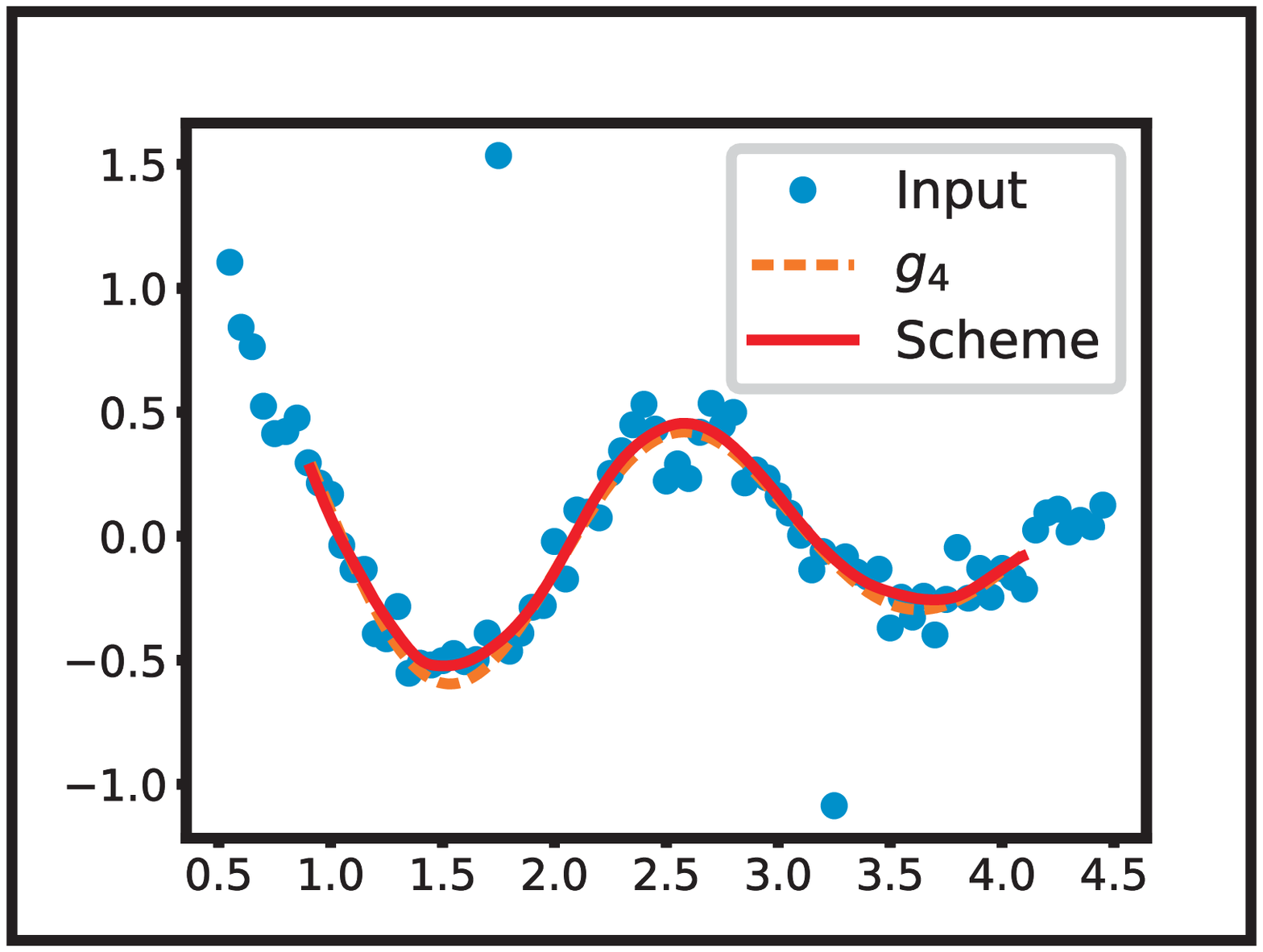, clip,trim=7cm 5.5cm 5.2cm 5.13cm, width=1.2 in}&\\
(j) & (k) & (l)
 \end{tabular}
\end{center}
\caption[Effects of the schemes $D_{h,d}$ on noisy data with outliers.]{\label{Noise-outlier1} \emph{Effects of the noisy data with outliers on fitted curves. (d)-(f) and (j)-(l) are the mirror images of the central parts (inside the green rectangles) of (a)-(c) and (g)-(i) respectively.}}
\end{figure}
\end{Exp}

\begin{rem}
Families of $2n$-point $N$-ary ($N \in \mathbb{N} \backslash \{1\}$) subdivision schemes can be obtained by substituting optimum $\beta^{(m+1)}_{\delta}$ of $\beta$ into degree three, two and one univariate polynomial functions and then evaluating these functions at $r=\frac{2M-1}{2N}$, $M=1,2,\ldots,N$. Similarly, by evaluating these functions at $r=\frac{M}{2N}$ with $M=-(N-1),\ldots,-5,-3,-1,1,3,5,\ldots,(N-1)$ and $r=\frac{M}{2N}$ with $M=-(N-1),\ldots,-5,-3,$ $-1,0,1,3,5,\ldots,(N-1)$, we get families of $(2n+1)$-point even N-ary subdivision schemes and $(2n+1)$-point odd N-ary subdivision schemes respectively.
\end{rem}

\section{Framework for the construction of bivariate schemes}\label{bivariate-construction-IRLS}
In this section, we generalize our representation of the 2-dimensional case to the 3-dimensional case to construct families of $(2n)^2$-point and $(2n+1)^2$-point non-tensor product binary schemes by using the IRLS technique as optimality criterion. Again the proposed framework contains two steps. These steps are presented in the following subsections.
\subsection{Step 1}
Consider the following bivariate polynomial function of degree $d$ to determine the best function to fit the data
\begin{eqnarray}\label{bipolyfunc}
f(x_{r},y_{s})&=&\sum\limits_{a=0}^{\frac{(d+1)(d+2)}{2}-1}B_{a}f_{a}(x_{r},y_{s}), \,\
\end{eqnarray}
where $f_{a}(x_{r},y_{s})=x_{r}^{a_{1}}y_{s}^{a_{2}}$ are the monomial functions with $a_{1}+a_{2} \leqslant d$ and $0 \leqslant a \leqslant \frac{(d+1)(d+2)}{2}-1$, i.e. there is one one correspondence between $a$$=$$\{0$, $1$, $2$, $\ldots$, $\frac{(d+1)(d+2)}{2}-2$, $\frac{(d+1)(d+2)}{2}-1\}$ and $(a_{1}$,$a_{2})$$=$$\{(0,0)$, $(1,0)$, $(0,1)$, $\ldots$, $(d,0)$, $(d-1,1)$, $\ldots$, $(1,d-1)$, $(0,d)\}$.\\
The polynomial function (\ref{bipolyfunc}) with respect to the observations $(x_{r}=r,y_{s}=s,f_{r,s})$ for $r=-n+1,\ldots,n$, $s=-n+1,\ldots,n$ where $n \in \mathbb{N}>1$ can be expressed as
\begin{eqnarray}\label{bipoly}
f_{r,s}&=&f(r,s)=\sum\limits_{a=0}^{\frac{(d+1)(d+2)}{2}-1}r^{a_{1}}s^{a_{2}}B_{a}, \,\ a_{1}+a_{2} \leqslant d.
\end{eqnarray}

The $\ell_{1}$ norm operator that is used for estimating $B=\{B_{a}:a=0,1,\ldots,\frac{(d+1)(d+2)}{2}-1 \}$ is defined as
\begin{eqnarray*}
&& \min_{B \in \mathbb{R}} \sum\limits_{r=-n+1}^{n}\sum\limits_{s=-n+1}^{n}\left|f_{r,s}-\sum\limits_{a=0}^{\frac{(d+1)(d+2)}{2}-1}r^{a_{1}}s^{a_{2}}B_{a}\right|, \,\ a_{1}+a_{2} \leqslant d.
\end{eqnarray*}
The optimum of $B$ is
\begin{eqnarray}\label{bi-optimization}
\nonumber B_{0},B_{1},\ldots,B_{\frac{(d+1)(d+2)}{2}-1}&=& \mbox{arg} \,\ \min_{B \in \mathbb{R}} \sum\limits_{r=-n+1}^{n}\sum\limits_{s=-n+1}^{n}\left|f_{r,s}-\sum\limits_{a=0}^{\frac{(d+1)(d+2)}{2}-1}r^{a_{1}}s^{a_{2}}B_{a}\right|\\
&&=\mbox{arg} \,\ \min_{B \in \mathbb{R}}F\left(B_{0},B_{1},\ldots,B_{\frac{(d+1)(d+2)}{2}-1}\right).
\end{eqnarray}
Use the IRLS method to approximate $F$ by $F_{\delta}$ ($\delta >0$) as
\begin{eqnarray}\label{bi-F-delta}
\nonumber F_{\delta}\left(B_{0},B_{1},\ldots,B_{\frac{(d+1)(d+2)}{2}-1}\right)&=&\sum\limits_{r=-n+1}^{n}\sum\limits_{s=-n+1}^{n}\left[\left(f_{r,s}-\sum\limits_{a=0}^{\frac{(d+1)(d+2)}{2}-1}r^{a_{1}}s^{a_{2}}B_{a}\right)^{2}+\delta\right]^{\frac{1}{2}}.\\
\end{eqnarray}
By regularizing (\ref{bi-optimization}) with (\ref{bi-F-delta}) and using the generalized version of (\cite{Bissantz}, Theorem 1), we get
\begin{eqnarray*}
B_{0,\delta},B_{1,\delta},\ldots,B_{\frac{(d+1)(d+2)}{2}-1,\delta}&=& \mbox{arg} \,\ \min_{B \in \mathbb{R}}F_{\delta}\left(B_{0},B_{1},\ldots,B_{\frac{(d+1)(d+2)}{2}-1}\right).
\end{eqnarray*}
Thus, we compute $B^{(m+1)}_{\delta}=\{B^{(m+1)}_{a,\delta}:a=0,1,\ldots,\frac{(d+1)(d+2)}{2}-1\}$ which is the approximation of $B$ by using the following iterative formula
\begin{eqnarray}\label{bi-I-Formulae}
\nonumber B^{(m+1)}_{0,\delta},B^{(m+1)}_{1,\delta},\ldots,B^{(m+1)}_{\frac{(d+1)(d+2)}{2}-1,\delta}&=& \mbox{arg} \,\ \min_{B \in \mathbb{R}} \sum\limits_{r=-n+1}^{n}\sum\limits_{s=-n+1}^{n}w^{(m)}_{r,s}\left(f_{r,s}-\right.\\&&\left.
\sum\limits_{a=0}^{\frac{(d+1)(d+2)}{2}-1}r^{a_{1}}s^{a_{2}}B_{a}\right)^{2},
\end{eqnarray}
where
\begin{eqnarray}\label{bi--weights}
w^{(m)}_{r,s}&=&\left[\left(f_{r,s}-\sum\limits_{a=0}^{\frac{(d+1)(d+2)}{2}-1}r^{a_{1}}s^{a_{2}}B_{a}\right)^{2}+\delta\right]^{-\frac{1}{2}}, \,\ a_{1}+a_{2} \leqslant d.
\end{eqnarray}
For starting values $B^{(0)}_{\delta}=\{B^{(0)}_{a,\delta}:a=0,1,\ldots,\frac{(d+1)(d+2)}{2}-1\}$, use the ordinary least squares method. \\
The sum of the squares of residuals for bivariate $d$ degree polynomial is
\begin{eqnarray}\label{bi-SSRM}
R=\sum\limits_{r=-n+1}^{n}\sum\limits_{s=-n+1}^{n}\left(f_{r,s}-\sum\limits_{a=0}^{\frac{(d+1)(d+2)}{2}-1}r^{a_{1}}s^{a_{2}}B^{(0)}_{a}\right)^{2}, \,\ a_{1}+a_{2} \leqslant d.
\end{eqnarray}
Differentiate (\ref{bi-I-Formulae}) and (\ref{bi-SSRM}) with respect to $B^{(0)}_{\delta}=B=\{B_{a}:a=0,1,\ldots,\frac{(d+1)(d+2)}{2}-1\}$ then solve these systems for unknowns values, we get $B^{(m+1)}_{\delta}=\{B^{(m+1)}_{a,\delta}:a=0,1,\ldots,\frac{(d+1)(d+2)}{2}-1\}$ and $B^{(0)}_{\delta}=\{B^{(0)}_{a,\delta}:a=0,1,\ldots,\frac{(d+1)(d+2)}{2}-1\}$ respectively.
Iterations will continue until
\begin{eqnarray*}
&&\mbox{max}\left(\left|B^{(m+1)}_{a,\delta}-B^{(m)}_{a,\delta}\right|:a=0,1,\ldots,\frac{(d+1)(d+2)}{2}-1\right)<\epsilon .
\end{eqnarray*}
Substitute $B^{(m+1)}_{\delta}=\{B^{(m+1)}_{a,\delta}:a=0,1,\ldots,\frac{(d+1)(d+2)}{2}-1\}$ in (\ref{bipoly}) to get the following best fitted $d$-degree polynomial to the $(2n)^{2}$-observations
\begin{eqnarray}\label{bi-beta-a-delta}
f(r,s)&=&\sum\limits_{a=0}^{\frac{(d+1)(d+2)}{2}-1}r^{a_{1}}s^{a_{2}}B^{(m+1)}_{a,\delta}, \,\ a_{1}+a_{2} \leqslant d.
\end{eqnarray}
This step is briefly defined in Algorithm \ref{bi-parameters-Estimating}.
\subsection{Step 2}
In this step, evaluate the bivariate $d$-degree polynomial (\ref{bi-beta-a-delta}) at $(r,s)=\left(\frac{1}{4},\frac{1}{4}\right)$, $\left(\frac{3}{4},\frac{1}{4}\right)$, $\left(\frac{1}{4},\frac{3}{4}\right)$ and $\left(\frac{3}{4},\frac{3}{4}\right)$.
Then, for the iterative point of view, set the following representations
\begin{eqnarray*}
&&f\left(\frac{1}{4},\frac{1}{4}\right)=f_{2i,2j}^{k+1}, \qquad \quad
f\left(\frac{3}{4},\frac{1}{4}\right)=f_{2i+1,2j}^{k+1},\\
&&f\left(\frac{1}{4},\frac{3}{4}\right)=f_{2i,2j+1}^{k+1}, \qquad
f\left(\frac{3}{4},\frac{3}{4}\right)=f_{2i+1,2j+1}^{k+1},\\
&&B^{(m+1)}_{a,\delta}=B^{(m+1)}_{a,\delta,i,j}, \qquad \quad
w_{r,s}^{(m)}=w_{i+r,j+s}^{(m)},
\end{eqnarray*}
and the 3D data values involved in $B^{(m+1)}_{a,\delta}$ are set as
\begin{eqnarray*}
f_{r,s} &=& f^{k}_{i+r,j+s}, \,\ i,j \in \mathbb{Z}
\end{eqnarray*}
where $f^{k+1}_{i,j}$ and $f^{k}_{i,j}$ are the control points at level $k+1$ and $k$ respectively. For different values of $d$ and $n$, the four refinement rules $f_{2i,2j}^{k+1}$ ,$f_{2i+1,2j}^{k+1}$, $f_{2i,2j+1}^{k+1}$ and $f_{2i+1,2j+1}^{k+1}$ make the $(2n)^{2}$-point scheme denoted by $D_{(2n)^{2},d}$ to fit a surface to the set of 3D data points.

\begin{algorithm}[htb!] 
   \caption[Estimation of the unknown parameters $B$ in (\ref{bipoly}).]{\label{bi-parameters-Estimating} \emph{Estimation of the unknown parameters $B$ in (\ref{bipoly}).}}
    \begin{algorithmic}[1]
    \State \textbf{input:}  $d$, $n$, $\epsilon$, $\sigma$, $m^{''}$
    \For{$m^{'} = 0$ to ${m^{''}}$}
            \If {$m^{'}=0$}
                \State calculate: $B^{(m^{'})}_{\delta}=\{B^{(m^{'})}_{a,\delta}:a=0,1,\ldots,\frac{(d+1)(d+2)}{2}-1\}$ by differentiating
                $\quad$ $\quad$ $\frac{}{}$ (\ref{bi-SSRM}) and  solving system of equations
            \Else
                \State calculate: $B^{(m^{'})}_{\delta}=\{B^{(m^{'})}_{a,\delta}:a=0,1,\ldots,\frac{(d+1)(d+2)}{2}-1\}$ by differentiating
                $\quad$ $\quad$ $\frac{}{}$ (\ref{bi-I-Formulae}) and solving system of equations
            \EndIf
            \State calculate: $w^{(m^{'})}_{r,s}$, $r=-n+1,\ldots,n$, $s=-n+1,\ldots,n$ by (\ref{bi--weights})
            \If {$m^{'} = 0$}
            \State \textbf{go to} 2
            \EndIf
            \State calculate: $M^{'}=\max \limits_{a}\left|B^{(m^{'})}_{a,\delta}-B^{(m^{'}-1)}_{a,\delta}\right|$, $a=0,1,\ldots,\frac{(d+1)(d+2)}{2}-1$.
            \If {$M^{'}<\epsilon$}
            \State \textbf{go to} 17
            \EndIf
        \EndFor
        \State set: $B^{(m+1)}_{a,\delta}=B^{(m^{'})}_{a,\delta}$
        \State \textbf{output:} $B^{(m+1)}_{a,\delta},a=0,1,\ldots,\frac{(d+1)(d+2)}{2}-1$
\end{algorithmic}
\end{algorithm}

\begin{algorithm}[htb!] 
   \caption[Subdivision of 3D data points.]{\label{bi-SS-A} \emph{Subdivision of 3D data points.}}
    \begin{algorithmic}[1]
    \State \textbf{input:}  $k^{''}$, $f^{0}=\{f^{0}_{i+r,j+s}:i,j \in \mathbb{Z};r=-n+1,\ldots,n;s=-n+1,\ldots,n\}$, $B^{(m+1)}_{\delta}=\{B^{(m+1)}_{a,\delta,i,j}:i,j \in \mathbb{Z};a=0,1,\ldots,\frac{(d+1)(d+2)}{2}-1\}$
    \State calculate: $f\left(\frac{1}{4},\frac{1}{4}\right)$, $f\left(\frac{3}{4},\frac{1}{4}\right)$, $f\left(\frac{1}{4},\frac{3}{4}\right)$ and $f\left(\frac{3}{4},\frac{3}{4}\right)$ from (\ref{bi-beta-a-delta})
    \State set: $f^{1}_{2i,2j}=f\left(\frac{1}{4},\frac{1}{4}\right)$, $f^{1}_{2i+1,2j}=f \left(\frac{3}{4},\frac{1}{4}\right)$, $f^{1}_{2i,2j+1}=f\left(\frac{1}{4},\frac{3}{4}\right)$ and $f^{1}_{2i+1,2j+1}=f \left(\frac{3}{4},\frac{3}{4}\right)$
    \State rewrite step 3: $f^{1}=S_{(2n)^{2},d}f^{0}$ \Comment{$S_{(2n)^{2},d}$ is the subdivision matrix/rule of the
    $\qquad$ $\qquad$ $\qquad$ $\qquad$ $\qquad$ $\qquad$ $\frac{}{}$ $\quad$ scheme $D_{(2n)^{2},d}$}
    \For{$k^{'}=1$ to $k^{''}$} \Comment{$k^{''} \in \mathbb{N}$ number of subdivision steps}
        \State $f^{k^{'}}=S_{(2n)^{2},d}f^{k^{'}-1}$
    \EndFor
    \State set: $f^{k^{''}}=f^{k+1}$
    \State \textbf{output:} $f^{k+1}=\{f^{k+1}_{i,j}:i,j \in \mathbb{Z}\}$ \Comment{$(k+1)$-th level subdivided data}
\end{algorithmic}
\end{algorithm}

\begin{prop}\label{bi-2n-betas-OLS-cor}
If we substitute $d=2$ in (\ref{bipoly}), then from Algorithm \ref{bi-parameters-Estimating} the starting values $B^{(0)}_{\delta}=\{B^{(0)}_{a,\delta}:a=0,1,\ldots,5\}$ are:
\begin{eqnarray}\label{bi-ols-B-5-2n}
 B^{(0)}_{5,\delta} &=& \sum\limits_{r=-n+1}^{n}\sum\limits_{s=-n+1}^{n} \frac{-15 (3 s- 1-3 s^{2} + n^{2})}{4 n^{2}(n^{2}-1)(4n^{2}-1)}f_{r,s},
\end{eqnarray}

\begin{eqnarray}\label{bi-ols-B-4-2n}
 B^{(0)}_{4,\delta} &=& \sum\limits_{r=-n+1}^{n}\sum\limits_{s=-n+1}^{n} \frac{9 (4 r s +1 -2 r - 2 s)}{n^{2}(4n^{2}-1)^{2}}f_{r,s},
\end{eqnarray}

\begin{eqnarray}\label{bi-ols-B-3-2n}
B^{(0)}_{3,\delta} &=& \sum\limits_{r=-n+1}^{n}\sum\limits_{s=-n+1}^{n} \frac{-15 (n^{2} + 3 r-3 r^{2} -1)}{4 n^{2}(n^{2}-1)(4n^{2}-1)}f_{r,s},
\end{eqnarray}

\begin{eqnarray}\label{bi-ols-B-2-2n}
B^{(0)}_{2,\delta} &=&\sum\limits_{r=-n+1}^{n}\sum\limits_{s=-n+1}^{n}  \frac{3 (2s-1)}{2 n^{2}(4n^{2}-1)}f_{r,s}-\frac{1}{2}B^{(0)}_{4,\delta}-
B^{(0)}_{5,\delta},
\end{eqnarray}
\begin{eqnarray}\label{bi-ols-B-1-2n}
B^{(0)}_{1,\delta} &=&\sum\limits_{r=-n+1}^{n}\sum\limits_{s=-n+1}^{n} \frac{3 (2r-1)}{2 n^{2}(4n^{2}-1)}f_{r,s}-B^{(0)}_{3,\delta}-\frac{1}{2}B^{(0)}_{4,\delta},
\end{eqnarray}
\begin{eqnarray}\label{bi-ols-B-0-2n}
\nonumber
B^{(0)}_{0,\delta}&=&\sum\limits_{r=-n+1}^{n}\sum\limits_{s=-n+1}^{n} \frac{1}{4 n^{2}} f_{r,s}-\frac{1}{2}B^{(0)}_{1,\delta}-\frac{1}{2}B^{(0)}_{2,\delta} -\frac{2n^{2}+1}{6} \\
&& \times B^{(0)}_{3,\delta}-\frac{1}{4}B^{(0)}_{4,\delta}-\frac{2n^{2}+1}{6} B^{(0)}_{5,\delta}.
\end{eqnarray}
\end{prop}

\begin{prop}\label{bi-betas-IRLS-cor}
If we substitute $d=2$ in (\ref{bipoly}), then from Algorithm \ref{bi-parameters-Estimating} the optimum $B^{(m+1)}_{\delta}=\{B^{(m+1)}_{a,\delta}:a=0,1,\ldots,5\}$ of unknown parameters $B$ are:

\begin{eqnarray}\label{bi-IRLS-betas-5-2n}
  B^{(m+1)}_{5,\delta} &=&-\frac{E^{(m)}_{4}}{E^{(m)}_{3}},
\end{eqnarray}

\begin{eqnarray}\label{bi-IRLS-betas-4-2n}
  B^{(m+1)}_{4,\delta} &=&-\frac{1}{E^{(m)}_{0}}\left(E^{(m)}_{1}+ E^{(m)}_{2} B^{(m+1)}_{5,\delta}\right),
\end{eqnarray}

\begin{eqnarray}\label{bi-IRLS-betas-3-2n}
  B^{(m+1)}_{3,\delta} &=&-\frac{1}{C^{(m)}_{5}}\left(C^{(m)}_{6}+C^{(m)}_{7} B^{(m+1)}_{4,\delta}+C^{(m)}_{8} B^{(m+1)}_{5,\delta}\right),
\end{eqnarray}

\begin{eqnarray}\label{bi-IRLS-betas-2-2n}
  B^{(m+1)}_{2,\delta} &=&-\frac{1}{C^{(m)}_{0}}\left(C^{(m)}_{1}+C^{(m)}_{2} B^{(m+1)}_{3,\delta}+C^{(m)}_{3} B^{(m+1)}_{4,\delta}+C^{(m)}_{4} B^{(m+1)}_{5,\delta}\right),
\end{eqnarray}

\begin{eqnarray}\label{bi-IRLS-betas-1-2n}
\nonumber  B^{(m+1)}_{1,\delta} &=& -\frac{1}{L^{(m)}_{0}} \left(K^{(m)}_{0}+L^{(m)}_{1} B^{(m+1)}_{2,\delta}+L^{(m)}_{2} B^{(m+1)}_{3,\delta}+L^{(m)} _{3} B^{(m+1)}_{4,\delta}+L^{(m)}_{4} \times \right. \\&&
\left. B^{(m+1)}_{5,\delta}\right),
\end{eqnarray}

\begin{eqnarray}\label{bi-IRLS-betas-0-2n}
\nonumber B^{(m+1)}_{0,\delta}&=&\frac{1}{\tau^{(m)}_{0}}\left(A^{(m)}_{0}-\tau^{(m)}_{1} B^{(m+1)}_{1,\delta}-\tau^{(m)}_{2} B^{(m+1)}_{2,\delta}-\tau^{(m)}_{3} B^{(m+1)}_{3,\delta}-\tau^{(m)}_{4} B^{(m+1)}_{4,\delta}- \right.\\&&
\left. \tau^{(m)}_{5} B^{(m+1)}_{5,\delta}\right),
\end{eqnarray}

where

\begin{eqnarray}\label{bi-weights}
\left\{\begin{array}{ccccccc}
&&E^{(m)}_{0}=C^{(m)}_{9}+C^{(m)}_{10}, \,\ E^{(m)}_{1}=C^{(m)}_{11}+C^{(m)}_{12}, \\ \\
&& E^{(m)}_{2}=C^{(m)}_{13}+C^{(m)}_{14},
\end{array}\right.
\end{eqnarray}

\begin{eqnarray}\label{bi-weights1}
\left\{\begin{array}{ccccccc}
&&E^{(m)}_{3}=C^{(m)}_{15}+C^{(m)}_{16}+C^{(m)}_{17}+C^{(m)}_{18}+C^{(m)}_{19}\\ \\
&& \qquad \quad \,\ +C^{(m)}_{20}+C^{(m)}_{21}+C^{(m)}_{22}+C^{(m)}_{23}+C^{(m)}_{24}, \\ \\
&&E^{(m)}_{4}=C^{(m)}_{25}+C^{(m)}_{26}+C^{(m)}_{27}+C^{(m)}_{28}+C^{(m)}_{29}\\ \\
&& \qquad \quad \,\ +C^{(m)}_{30}+C^{(m)}_{31}+C^{(m)}_{32}+C^{(m)}_{33}+C^{(m)}_{34},
\end{array}\right.
\end{eqnarray}

\begin{eqnarray}\label{bi-weights2}
\left\{\begin{array}{ccccccc}
&&C^{(m)}_{0}=L^{(m)}_{6} L^{(m)}_{0}-L^{(m)}_{5} L^{(m)}_{1}, \,\ C^{(m)}_{1}=K^{(m)}_{1} L^{(m)}_{0}-L^{(m)}_{5} K^{(m)}_{0}, \\ \\
&&C^{(m)}_{2}=L^{(m)}_{7} L^{(m)}_{0}-L^{(m)}_{5} L^{(m)}_{2}, \,\ C^{(m)}_{3}=L^{(m)}_{8} L^{(m)}_{0}-L^{(m)}_{5} L^{(m)}_{3}, \\ \\
&&C^{(m)}_{4}=L^{(m)}_{9} L^{(m)}_{0}-L^{(m)}_{5} L^{(m)}_{4},
\end{array}\right.
\end{eqnarray}

\begin{eqnarray}\label{bi-weights3}
\left\{\begin{array}{ccccccc}
&&C^{(m)}_{5}= -L^{(m)}_{12} L^{(m)}_{0} L^{(m)}_{6}-L^{(m)}_{10} L^{(m)}_{1} L^{(m)}_{7} + L^{(m)}_{11} L^{(m)}_{7} L^{(m)}_{0} \\ \\
&& \qquad \quad \,\ + L^{(m)}_{10} L^{(m)}_{2} L^{(m)}_{6}-L^{(m)}_{11} L^{(m)}_{5} L^{(m)}_{2} + L^{(m)}_{12} L^{(m)}_{5} L^{(m)}_{1}, \\ \\
&&C^{(m)}_{6}=K^{(m)}_{2} L^{(m)}_{5} L^{(m)}_{1}-K^{(m)}_{2} L^{(m)}_{0} L^{(m)}_{6}+L^{(m)}_{10} K^{(m)}_{0} L^{(m)}_{6} \\ \\
&& \qquad \quad \,\ -L^{(m)}_{10} L^{(m)}_{1} K^{(m)}_{1}+L^{(m)}_{11} K^{(m)}_{1} L^{(m)}_{0}-L^{(m)}_{11} L^{(m)}_{5} K^{(m)}_{0},
\end{array}\right.
\end{eqnarray}

\begin{eqnarray}\label{bi-weights4}
\left\{\begin{array}{ccccccc}
&&C^{(m)}_{7}= L^{(m)}_{13} L^{(m)}_{5} L^{(m)}_{1} +L^{(m)}_{11} L^{(m)}_{8} L^{(m)}_{0}-L^{(m)}_{10} L^{(m)}_{1} L^{(m)}_{8}  \\ \\
&& \qquad \quad \,\ -L^{(m)}_{13} L^{(m)}_{0} L^{(m)}_{6} +L^{(m)}_{10} L^{(m)}_{3} L^{(m)}_{6} -L^{(m)}_{11} L^{(m)}_{5} L^{(m)}_{3}, \\ \\
&&C^{(m)}_{8}=L^{(m)}_{11} L^{(m)}_{9} L^{(m)}_{0}-L^{(m)}_{10} L^{(m)}_{1} L^{(m)}_{9} -L^{(m)}_{11} L^{(m)}_{5} L^{(m)}_{4}  \\ \\
&& \qquad \quad \,\ -L^{(m)}_{14} L^{(m)}_{0} L^{(m)}_{6} +L^{(m)}_{10} L^{(m)}_{4} L^{(m)}_{6} +L^{(m)}_{14} L^{(m)}_{5} L^{(m)}_{1},
\end{array}\right.
\end{eqnarray}

\begin{eqnarray}\label{bi-weights5}
\nonumber C^{(m)}_{9}&=&L^{(m)}_{17} L^{(m)}_{13} L^{(m)}_{0} L^{(m)}_{6} - L^{(m)}_{0} L^{(m)}_{6} L^{(m)}_{18} L^{(m)}_{12} - L^{(m)}_{17} L^{(m)}_{11} L^{(m)}_{8} L^{(m)}_{0}   \\&&
\nonumber
+ L^{(m)}_{0} L^{(m)}_{18} L^{(m)}_{11} L^{(m)}_{7} + L^{(m)}_{0} L^{(m)}_{16} L^{(m)}_{8} L^{(m)}_{12} - L^{(m)}_{0} L^{(m)}_{16} L^{(m)}_{7} L^{(m)}_{13} \\&&
\nonumber
- L^{(m)}_{15} L^{(m)}_{2} L^{(m)}_{13} L^{(m)}_{6} + L^{(m)}_{6} L^{(m)}_{18} L^{(m)}_{10} L^{(m)}_{2} + L^{(m)}_{6} L^{(m)}_{15} L^{(m)}_{3} L^{(m)}_{12} \\&&
\nonumber
- L^{(m)}_{17} L^{(m)}_{10} L^{(m)}_{3} L^{(m)}_{6} - L^{(m)}_{18} L^{(m)}_{10} L^{(m)}_{1} L^{(m)}_{7} + L^{(m)}_{17} L^{(m)}_{11} L^{(m)}_{5} L^{(m)}_{3}, \\
\end{eqnarray}
\begin{eqnarray}\label{bi-weights6}
\nonumber C^{(m)}_{10}&=&
L^{(m)}_{16} L^{(m)}_{7} L^{(m)}_{10} L^{(m)}_{3} + L^{(m)}_{15} L^{(m)}_{2} L^{(m)}_{11} L^{(m)}_{8} + L^{(m)}_{18} L^{(m)}_{5} L^{(m)}_{1} L^{(m)}_{12}   \\ &&
\nonumber - L^{(m)}_{15} L^{(m)}_{1} L^{(m)}_{8} L^{(m)}_{12} - L^{(m)}_{16} L^{(m)}_{8} L^{(m)}_{10} L^{(m)}_{2} - L^{(m)}_{18} L^{(m)}_{11} L^{(m)}_{5} L^{(m)}_{2}   \\ &&
 \nonumber - L^{(m)}_{17} L^{(m)}_{13} L^{(m)}_{5} L^{(m)}_{1} - L^{(m)}_{16} L^{(m)}_{5} L^{(m)}_{3} L^{(m)}_{12} + L^{(m)}_{15} L^{(m)}_{1} L^{(m)}_{7} L^{(m)}_{13}   \\ &&
 \nonumber - L^{(m)}_{15} L^{(m)}_{3} L^{(m)}_{11} L^{(m)}_{7} + L^{(m)}_{17} L^{(m)}_{10} L^{(m)}_{1} L^{(m)}_{8} + L^{(m)}_{16} L^{(m)}_{5} L^{(m)}_{2} L^{(m)}_{13}, \\
\end{eqnarray}

\begin{eqnarray}\label{bi-weights7}
\nonumber C^{(m)}_{11}&=& K^{(m)}_{3} L^{(m)}_{11} L^{(m)}_{5} L^{(m)}_{2} - L^{(m)}_{15} L^{(m)}_{2} L^{(m)}_{11} K^{(m)}_{1} - L^{(m)}_{15} L^{(m)}_{1} L^{(m)}_{7} K^{(m)}_{2}    \\ &&
\nonumber - L^{(m)}_{16} L^{(m)}_{5} L^{(m)}_{2} K^{(m)}_{2} + K^{(m)}_{3} L^{(m)}_{10} L^{(m)}_{1} L^{(m)}_{7} - K^{(m)}_{3} L^{(m)}_{5} L^{(m)}_{1} L^{(m)}_{12}    \\ &&
\nonumber + L^{(m)}_{17} L^{(m)}_{11} K^{(m)}_{1} L^{(m)}_{0} - L^{(m)}_{0} K^{(m)}_{3} L^{(m)}_{11} L^{(m)}_{7} + L^{(m)}_{0} L^{(m)}_{16} L^{(m)}_{7} K^{(m)}_{2}    \\ &&
 \nonumber + L^{(m)}_{0} L^{(m)}_{6} K^{(m)}_{3} L^{(m)}_{12} - L^{(m)}_{0} L^{(m)}_{16} K^{(m)}_{1} L^{(m)}_{12} - L^{(m)}_{6} K^{(m)}_{3} L^{(m)}_{10} L^{(m)}_{2}, \\
\end{eqnarray}
\begin{eqnarray}\label{bi-weights8}
\nonumber C^{(m)}_{12}&=&
L^{(m)}_{15} L^{(m)}_{1} K^{(m)}_{1} L^{(m)}_{12} + L^{(m)}_{17} K^{(m)}_{2} L^{(m)}_{5} L^{(m)}_{1} - L^{(m)}_{17} K^{(m)}_{2} L^{(m)}_{0} L^{(m)}_{6}    \\ &&
\nonumber  + L^{(m)}_{17} L^{(m)}_{10} K^{(m)}_{0} L^{(m)}_{6} - L^{(m)}_{17} L^{(m)}_{10} L^{(m)}_{1} K^{(m)}_{1} + L^{(m)}_{16} K^{(m)}_{1} L^{(m)}_{10} L^{(m)}_{2}    \\ &&
 \nonumber + L^{(m)}_{16} L^{(m)}_{5} K^{(m)}_{0} L^{(m)}_{12} - L^{(m)}_{16} L^{(m)}_{7} L^{(m)}_{10} K^{(m)}_{0} - L^{(m)}_{6} L^{(m)}_{15} K^{(m)}_{0} L^{(m)}_{12}   \\ &&
  \nonumber - L^{(m)}_{17} L^{(m)}_{11} L^{(m)}_{5} K^{(m)}_{0} + L^{(m)}_{15} K^{(m)}_{0} L^{(m)}_{11} L^{(m)}_{7} + L^{(m)}_{15} L^{(m)}_{2} K^{(m)}_{2} L^{(m)}_{6}, \\
\end{eqnarray}

\begin{eqnarray}\label{bi-weights9}
\nonumber C^{(m)}_{13}&=& - L^{(m)}_{0} L^{(m)}_{19} L^{(m)}_{11} L^{(m)}_{7} + L^{(m)}_{16} L^{(m)}_{5} L^{(m)}_{4} L^{(m)}_{12} - L^{(m)}_{17} L^{(m)}_{10} L^{(m)}_{1} L^{(m)}_{9}  \\ &&
 \nonumber - L^{(m)}_{15} L^{(m)}_{1} L^{(m)}_{7} L^{(m)}_{14} - L^{(m)}_{17} L^{(m)}_{14} L^{(m)}_{0} L^{(m)}_{6} - L^{(m)}_{0} L^{(m)}_{16} L^{(m)}_{9} L^{(m)}_{12}  \\ &&
\nonumber - L^{(m)}_{15} L^{(m)}_{2} L^{(m)}_{11} L^{(m)}_{9} + L^{(m)}_{0} L^{(m)}_{16} L^{(m)}_{7} L^{(m)}_{14} + L^{(m)}_{17} L^{(m)}_{11} L^{(m)}_{9} L^{(m)}_{0}  \\ &&
\nonumber + L^{(m)}_{0} L^{(m)}_{6} L^{(m)}_{19} L^{(m)}_{12} + L^{(m)}_{17} L^{(m)}_{10} L^{(m)}_{4} L^{(m)}_{6} + L^{(m)}_{16} L^{(m)}_{9} L^{(m)}_{10} L^{(m)}_{2}, \\
\end{eqnarray}
\begin{eqnarray}\label{bi-weights10}
\nonumber C^{(m)}_{14}&=&
 - L^{(m)}_{17} L^{(m)}_{11} L^{(m)}_{5} L^{(m)}_{4} + L^{(m)}_{19} L^{(m)}_{11} L^{(m)}_{5} L^{(m)}_{2} + L^{(m)}_{15} L^{(m)}_{1} L^{(m)}_{9} L^{(m)}_{12}  \\ &&
\nonumber - L^{(m)}_{6} L^{(m)}_{19} L^{(m)}_{10} L^{(m)}_{2} - L^{(m)}_{6} L^{(m)}_{15} L^{(m)}_{4} L^{(m)}_{12} + L^{(m)}_{15} L^{(m)}_{4} L^{(m)}_{11} L^{(m)}_{7}  \\ &&
\nonumber + L^{(m)}_{15} L^{(m)}_{2} L^{(m)}_{14} L^{(m)}_{6} + L^{(m)}_{17} L^{(m)}_{14} L^{(m)}_{5} L^{(m)}_{1} - L^{(m)}_{16} L^{(m)}_{7} L^{(m)}_{10} L^{(m)}_{4}  \\ &&
\nonumber - L^{(m)}_{16} L^{(m)}_{5} L^{(m)}_{2} L^{(m)}_{14} - L^{(m)}_{19} L^{(m)}_{5} L^{(m)}_{1} L^{(m)}_{12} + L^{(m)}_{19} L^{(m)}_{10} L^{(m)}_{1} L^{(m)}_{7},\\
\end{eqnarray}

\begin{eqnarray}\label{bi-weights11}
 \nonumber C^{(m)}_{15}&=&L^{(m)}_{24} L^{(m)}_{16} L^{(m)}_{5} L^{(m)}_{2} L^{(m)}_{13} + L^{(m)}_{6} L^{(m)}_{22} L^{(m)}_{13} L^{(m)}_{15} L^{(m)}_{4} - L^{(m)}_{20} L^{(m)}_{1} L^{(m)}_{8} L^{(m)}_{17}   \\ &&
  \nonumber \times L^{(m)}_{14}+ L^{(m)}_{20} L^{(m)}_{1} L^{(m)}_{9} L^{(m)}_{17} L^{(m)}_{13} - L^{(m)}_{22} L^{(m)}_{8} L^{(m)}_{19} L^{(m)}_{10} L^{(m)}_{1} + L^{(m)}_{22} L^{(m)}_{13}   \\ &&
\nonumber \times L^{(m)}_{19} L^{(m)}_{5} L^{(m)}_{1}- L^{(m)}_{6} L^{(m)}_{23} L^{(m)}_{15} L^{(m)}_{4} L^{(m)}_{12} + L^{(m)}_{21} L^{(m)}_{9} L^{(m)}_{15} L^{(m)}_{2} L^{(m)}_{13} -    \\ &&
\nonumber L^{(m)}_{21} L^{(m)}_{5} L^{(m)}_{3} L^{(m)}_{17} L^{(m)}_{14} - L^{(m)}_{21} L^{(m)}_{8} L^{(m)}_{17} L^{(m)}_{10} L^{(m)}_{4} + L^{(m)}_{23} L^{(m)}_{17} L^{(m)}_{14} L^{(m)}_{5}    \\ &&
 \times L^{(m)}_{1} - L^{(m)}_{6} L^{(m)}_{20} L^{(m)}_{3} L^{(m)}_{19} L^{(m)}_{12},
\end{eqnarray}
\begin{eqnarray}\label{bi-weights12}
\nonumber C^{(m)}_{16}&=&L^{(m)}_{0} L^{(m)}_{22} L^{(m)}_{9} L^{(m)}_{16} L^{(m)}_{13} - L^{(m)}_{20} L^{(m)}_{2} L^{(m)}_{9} L^{(m)}_{16} L^{(m)}_{13} - L^{(m)}_{0} L^{(m)}_{11} L^{(m)}_{7} L^{(m)}_{23} \\&&
\nonumber \times L^{(m)}_{19} - L^{(m)}_{6} L^{(m)}_{24} L^{(m)}_{15} L^{(m)}_{2} L^{(m)}_{13} - L^{(m)}_{0} L^{(m)}_{11} L^{(m)}_{24} L^{(m)}_{17} L^{(m)}_{8} + L^{(m)}_{6} L^{(m)}_{24} \\&&
\nonumber \times  L^{(m)}_{12} L^{(m)}_{15} L^{(m)}_{3} - L^{(m)}_{0} L^{(m)}_{6} L^{(m)}_{24} L^{(m)}_{12} L^{(m)}_{18} - L^{(m)}_{0} L^{(m)}_{22} L^{(m)}_{8} L^{(m)}_{16} L^{(m)}_{14} +\\&&
\nonumber  L^{(m)}_{0} L^{(m)}_{21} L^{(m)}_{8} L^{(m)}_{17} L^{(m)}_{14} - L^{(m)}_{20} L^{(m)}_{1} L^{(m)}_{9} L^{(m)}_{18} L^{(m)}_{12} + L^{(m)}_{0} L^{(m)}_{6} L^{(m)}_{22} L^{(m)}_{14} \\&&
\times  L^{(m)}_{18} - L^{(m)}_{0} L^{(m)}_{23} L^{(m)}_{16} L^{(m)}_{9} L^{(m)}_{12},
\end{eqnarray}
\begin{eqnarray}\label{bi-weights13}
\nonumber C^{(m)}_{17}&=& L^{(m)}_{0} L^{(m)}_{21} L^{(m)}_{9} L^{(m)}_{18} L^{(m)}_{12} + L^{(m)}_{0} L^{(m)}_{24} L^{(m)}_{16} L^{(m)}_{8} L^{(m)}_{12} - L^{(m)}_{0} L^{(m)}_{21} L^{(m)}_{8} L^{(m)}_{19} \\&&
 \nonumber \times L^{(m)}_{12} - L^{(m)}_{21} L^{(m)}_{13} L^{(m)}_{19} L^{(m)}_{5} L^{(m)}_{2} + L^{(m)}_{22} L^{(m)}_{8} L^{(m)}_{16} L^{(m)}_{10} L^{(m)}_{4} - L^{(m)}_{22} L^{(m)}_{5} \\&&
 \nonumber \times L^{(m)}_{4} L^{(m)}_{16} L^{(m)}_{13} + L^{(m)}_{21} L^{(m)}_{8} L^{(m)}_{15} L^{(m)}_{4} L^{(m)}_{12} - L^{(m)}_{23} L^{(m)}_{16} L^{(m)}_{5} L^{(m)}_{2} L^{(m)}_{14} + \\&&
 \nonumber L^{(m)}_{0} L^{(m)}_{6} L^{(m)}_{24} L^{(m)}_{17} L^{(m)}_{13} - L^{(m)}_{22} L^{(m)}_{9} L^{(m)}_{15} L^{(m)}_{1} L^{(m)}_{13} - L^{(m)}_{23} L^{(m)}_{17} L^{(m)}_{9} L^{(m)}_{10}\\&&
\times L^{(m)}_{1} - L^{(m)}_{20} L^{(m)}_{4} L^{(m)}_{16} L^{(m)}_{8} L^{(m)}_{12},
\end{eqnarray}
\begin{eqnarray}\label{bi-weights14}
\nonumber C^{(m)}_{18}&=&
L^{(m)}_{6} L^{(m)}_{22} L^{(m)}_{10} L^{(m)}_{3} L^{(m)}_{19} + L^{(m)}_{6} L^{(m)}_{23} L^{(m)}_{17} L^{(m)}_{10} L^{(m)}_{4} + L^{(m)}_{22} L^{(m)}_{5} L^{(m)}_{3} L^{(m)}_{16} \\&&
\nonumber \times L^{(m)}_{14} + L^{(m)}_{6} L^{(m)}_{24} L^{(m)}_{18} L^{(m)}_{10} L^{(m)}_{2} - L^{(m)}_{22} L^{(m)}_{14} L^{(m)}_{18} L^{(m)}_{5} L^{(m)}_{1} + L^{(m)}_{6} L^{(m)}_{20} \\&&
\nonumber \times L^{(m)}_{2} L^{(m)}_{13} L^{(m)}_{19} - L^{(m)}_{6} L^{(m)}_{23} L^{(m)}_{19} L^{(m)}_{10} L^{(m)}_{2} + L^{(m)}_{23} L^{(m)}_{15} L^{(m)}_{1} L^{(m)}_{9} L^{(m)}_{12} + \\&&
\nonumber  L^{(m)}_{20} L^{(m)}_{1} L^{(m)}_{8} L^{(m)}_{19} L^{(m)}_{12} + L^{(m)}_{0} L^{(m)}_{11} L^{(m)}_{7} L^{(m)}_{24} L^{(m)}_{18} + L^{(m)}_{11} L^{(m)}_{7} L^{(m)}_{20} L^{(m)}_{3} \\&&
\times L^{(m)}_{19} - L^{(m)}_{11} L^{(m)}_{7} L^{(m)}_{20} L^{(m)}_{4} L^{(m)}_{18},
\end{eqnarray}
\begin{eqnarray}\label{bi-weights15}
\nonumber C^{(m)}_{19}&=&
L^{(m)}_{0} L^{(m)}_{11} L^{(m)}_{23} L^{(m)}_{17} L^{(m)}_{9} + L^{(m)}_{0} L^{(m)}_{11} L^{(m)}_{22} L^{(m)}_{8} L^{(m)}_{19} - L^{(m)}_{0} L^{(m)}_{11} L^{(m)}_{22} L^{(m)}_{9} \\&&
\nonumber \times L^{(m)}_{18} + L^{(m)}_{7} L^{(m)}_{24} L^{(m)}_{16} L^{(m)}_{10} L^{(m)}_{3} - L^{(m)}_{21} L^{(m)}_{8} L^{(m)}_{15} L^{(m)}_{2} L^{(m)}_{14} + L^{(m)}_{21} L^{(m)}_{5} \\&&
\nonumber \times L^{(m)}_{3} L^{(m)}_{19} L^{(m)}_{12} + L^{(m)}_{0} L^{(m)}_{7} L^{(m)}_{23} L^{(m)}_{16} L^{(m)}_{14} - L^{(m)}_{0} L^{(m)}_{7} L^{(m)}_{24} L^{(m)}_{16} L^{(m)}_{13} - \\&&
\nonumber  L^{(m)}_{0} L^{(m)}_{7} L^{(m)}_{21} L^{(m)}_{14} L^{(m)}_{18} + L^{(m)}_{6} L^{(m)}_{20} L^{(m)}_{4} L^{(m)}_{12} L^{(m)}_{18} - L^{(m)}_{0} L^{(m)}_{6} L^{(m)}_{22} L^{(m)}_{13} \\&&
\times L^{(m)}_{19} - L^{(m)}_{0} L^{(m)}_{6} L^{(m)}_{23} L^{(m)}_{17} L^{(m)}_{14},
\end{eqnarray}
\begin{eqnarray}\label{bi-weights16}
\nonumber C^{(m)}_{20}&=&
L^{(m)}_{0} L^{(m)}_{6} L^{(m)}_{23} L^{(m)}_{19} L^{(m)}_{12} + L^{(m)}_{0} L^{(m)}_{7} L^{(m)}_{21} L^{(m)}_{13} L^{(m)}_{19} - L^{(m)}_{7} L^{(m)}_{23} L^{(m)}_{16} L^{(m)}_{10} \\&&
\nonumber \times L^{(m)}_{4} - L^{(m)}_{7} L^{(m)}_{20} L^{(m)}_{3} L^{(m)}_{16} L^{(m)}_{14} - L^{(m)}_{0} L^{(m)}_{21} L^{(m)}_{9} L^{(m)}_{17} L^{(m)}_{13} - L^{(m)}_{7} L^{(m)}_{21} \\&&
\nonumber \times L^{(m)}_{10} L^{(m)}_{3} L^{(m)}_{19} - L^{(m)}_{7} L^{(m)}_{21} L^{(m)}_{13} L^{(m)}_{15} L^{(m)}_{4} - L^{(m)}_{7} L^{(m)}_{24} L^{(m)}_{10} L^{(m)}_{1} L^{(m)}_{18} + \\&&
\nonumber  L^{(m)}_{7} L^{(m)}_{24} L^{(m)}_{15} L^{(m)}_{1} L^{(m)}_{13} + L^{(m)}_{7} L^{(m)}_{23} L^{(m)}_{19} L^{(m)}_{10} L^{(m)}_{1} + L^{(m)}_{7} L^{(m)}_{20} L^{(m)}_{4} L^{(m)}_{16} \\&&
\times L^{(m)}_{13} + L^{(m)}_{7} L^{(m)}_{21} L^{(m)}_{10} L^{(m)}_{4} L^{(m)}_{18},
\end{eqnarray}
\begin{eqnarray}\label{bi-weights17}
\nonumber C^{(m)}_{21}&=&
L^{(m)}_{7} L^{(m)}_{20} L^{(m)}_{1} L^{(m)}_{14} L^{(m)}_{18} - L^{(m)}_{7} L^{(m)}_{23} L^{(m)}_{15} L^{(m)}_{1} L^{(m)}_{14} + L^{(m)}_{7} L^{(m)}_{21} L^{(m)}_{14} L^{(m)}_{15}  \\&&
\nonumber \times L^{(m)}_{3} - L^{(m)}_{7} L^{(m)}_{20} L^{(m)}_{1} L^{(m)}_{13} L^{(m)}_{19} + L^{(m)}_{11} L^{(m)}_{20} L^{(m)}_{2} L^{(m)}_{9} L^{(m)}_{18} + L^{(m)}_{11} L^{(m)}_{24} \\&&
\nonumber \times L^{(m)}_{15} L^{(m)}_{2} L^{(m)}_{8} - L^{(m)}_{11} L^{(m)}_{20} L^{(m)}_{2} L^{(m)}_{8} L^{(m)}_{19} - L^{(m)}_{11} L^{(m)}_{7} L^{(m)}_{24} L^{(m)}_{15} L^{(m)}_{3} + \\&&
\nonumber  L^{(m)}_{11} L^{(m)}_{7} L^{(m)}_{23} L^{(m)}_{15} L^{(m)}_{4} + L^{(m)}_{11} L^{(m)}_{23} L^{(m)}_{19} L^{(m)}_{5} L^{(m)}_{2} - L^{(m)}_{11} L^{(m)}_{22} L^{(m)}_{5} L^{(m)}_{3} \\&&
\times L^{(m)}_{19} - L^{(m)}_{11} L^{(m)}_{23} L^{(m)}_{17} L^{(m)}_{5} L^{(m)}_{4},
\end{eqnarray}
\begin{eqnarray}\label{bi-weights18}
\nonumber C^{(m)}_{22}&=&
L^{(m)}_{11} L^{(m)}_{22} L^{(m)}_{9} L^{(m)}_{15} L^{(m)}_{3} - L^{(m)}_{6} L^{(m)}_{24} L^{(m)}_{17} L^{(m)}_{10} L^{(m)}_{3} - L^{(m)}_{11} L^{(m)}_{20} L^{(m)}_{3} L^{(m)}_{17} \\&&
\nonumber \times L^{(m)}_{9} + L^{(m)}_{11} L^{(m)}_{20} L^{(m)}_{4} L^{(m)}_{17} L^{(m)}_{8} + L^{(m)}_{11} L^{(m)}_{24} L^{(m)}_{17} L^{(m)}_{5} L^{(m)}_{3} + L^{(m)}_{11} L^{(m)}_{22} \\&&
\nonumber \times L^{(m)}_{5} L^{(m)}_{4} L^{(m)}_{18} - L^{(m)}_{11} L^{(m)}_{23} L^{(m)}_{15} L^{(m)}_{2} L^{(m)}_{9} - L^{(m)}_{11} L^{(m)}_{22} L^{(m)}_{8} L^{(m)}_{15} L^{(m)}_{4} + \\&&
\nonumber  L^{(m)}_{22} L^{(m)}_{8} L^{(m)}_{15} L^{(m)}_{1} L^{(m)}_{14} - L^{(m)}_{21} L^{(m)}_{9} L^{(m)}_{15} L^{(m)}_{3} L^{(m)}_{12} - L^{(m)}_{11} L^{(m)}_{24} L^{(m)}_{18} L^{(m)}_{5} \\&&
\times L^{(m)}_{2} + L^{(m)}_{6} L^{(m)}_{23} L^{(m)}_{15} L^{(m)}_{2} L^{(m)}_{14},
\end{eqnarray}
\begin{eqnarray}\label{bi-weights19}
\nonumber C^{(m)}_{23}&=&
L^{(m)}_{6} L^{(m)}_{20} L^{(m)}_{3} L^{(m)}_{17} L^{(m)}_{14} + L^{(m)}_{20} L^{(m)}_{2} L^{(m)}_{8} L^{(m)}_{16} L^{(m)}_{14} + L^{(m)}_{21} L^{(m)}_{8} L^{(m)}_{19} L^{(m)}_{10} \\&&
\nonumber \times L^{(m)}_{2} + L^{(m)}_{23} L^{(m)}_{16} L^{(m)}_{9} L^{(m)}_{10} L^{(m)}_{2} + L^{(m)}_{23} L^{(m)}_{16} L^{(m)}_{5} L^{(m)}_{4} L^{(m)}_{12} - L^{(m)}_{6} L^{(m)}_{20} \\&&
\nonumber \times L^{(m)}_{2} L^{(m)}_{14} L^{(m)}_{18} - L^{(m)}_{24} L^{(m)}_{16} L^{(m)}_{8} L^{(m)}_{10} L^{(m)}_{2} + L^{(m)}_{21} L^{(m)}_{9} L^{(m)}_{17} L^{(m)}_{10} L^{(m)}_{3} + \\&&
\nonumber  L^{(m)}_{24} L^{(m)}_{5} L^{(m)}_{1} L^{(m)}_{12} L^{(m)}_{18} + L^{(m)}_{20} L^{(m)}_{3} L^{(m)}_{16} L^{(m)}_{9} L^{(m)}_{12} + L^{(m)}_{21} L^{(m)}_{14} L^{(m)}_{18} L^{(m)}_{5} \\&&
\times L^{(m)}_{2} - L^{(m)}_{24} L^{(m)}_{5} L^{(m)}_{1} L^{(m)}_{17} L^{(m)}_{13},
\end{eqnarray}
\begin{eqnarray}\label{bi-weights20}
\nonumber C^{(m)}_{24}&=&
L^{(m)}_{21} L^{(m)}_{5} L^{(m)}_{4} L^{(m)}_{17} L^{(m)}_{13} - L^{(m)}_{22} L^{(m)}_{9} L^{(m)}_{16} L^{(m)}_{10} L^{(m)}_{3} - L^{(m)}_{23} L^{(m)}_{19} L^{(m)}_{12} L^{(m)}_{5} \\&&
\nonumber \times L^{(m)}_{1} - L^{(m)}_{6} L^{(m)}_{20} L^{(m)}_{4} L^{(m)}_{17} L^{(m)}_{13} - L^{(m)}_{6} L^{(m)}_{22} L^{(m)}_{14} L^{(m)}_{15} L^{(m)}_{3} - L^{(m)}_{21} L^{(m)}_{9} \\&&
\nonumber \times L^{(m)}_{18} L^{(m)}_{10} L^{(m)}_{2} + L^{(m)}_{24} L^{(m)}_{10} L^{(m)}_{1} L^{(m)}_{17} L^{(m)}_{8} - L^{(m)}_{24} L^{(m)}_{16} L^{(m)}_{5} L^{(m)}_{3} L^{(m)}_{12} - \\&&
\nonumber L^{(m)}_{24} L^{(m)}_{15} L^{(m)}_{1} L^{(m)}_{8} L^{(m)}_{12} - L^{(m)}_{6} L^{(m)}_{22} L^{(m)}_{10} L^{(m)}_{4} L^{(m)}_{18} - L^{(m)}_{21} L^{(m)}_{5} L^{(m)}_{4} L^{(m)}_{12} \\&&
\times L^{(m)}_{18} + L^{(m)}_{22} L^{(m)}_{9} L^{(m)}_{18} L^{(m)}_{10} L^{(m)}_{1},
\end{eqnarray}
\begin{eqnarray}\label{bi-weights21}
\nonumber C^{(m)}_{25}&=& L^{(m)}_{11} K^{(m)}_{4} L^{(m)}_{15} L^{(m)}_{2} L^{(m)}_{8} - L^{(m)}_{11} L^{(m)}_{20} L^{(m)}_{2} L^{(m)}_{8} K^{(m)}_{3} + L^{(m)}_{0} L^{(m)}_{6} L^{(m)}_{23} K^{(m)}_{3} \\&&
\nonumber \times L^{(m)}_{12} - L^{(m)}_{0} L^{(m)}_{6} L^{(m)}_{23} L^{(m)}_{17} K^{(m)}_{2} - L^{(m)}_{0} L^{(m)}_{6} L^{(m)}_{22} L^{(m)}_{13} K^{(m)}_{3} - L^{(m)}_{22} L^{(m)}_{5} \\&&
\nonumber \times K^{(m)}_{0} L^{(m)}_{16} L^{(m)}_{13} + L^{(m)}_{21} L^{(m)}_{5} L^{(m)}_{3} K^{(m)}_{3} L^{(m)}_{12} - L^{(m)}_{0} L^{(m)}_{23} L^{(m)}_{16} K^{(m)}_{1} L^{(m)}_{12} + \\&&
\nonumber  L^{(m)}_{0} L^{(m)}_{6} L^{(m)}_{22} K^{(m)}_{2} L^{(m)}_{18} + L^{(m)}_{0} L^{(m)}_{6} K^{(m)}_{4} L^{(m)}_{17} L^{(m)}_{13} - L^{(m)}_{0} L^{(m)}_{21} K^{(m)}_{1} L^{(m)}_{17} \\&&
\times L^{(m)}_{13} + L^{(m)}_{0} L^{(m)}_{7} L^{(m)}_{23} L^{(m)}_{16} K^{(m)}_{2},
\end{eqnarray}
\begin{eqnarray}\label{bi-weights22}
\nonumber C^{(m)}_{26}&=&
L^{(m)}_{0} L^{(m)}_{7} L^{(m)}_{21} L^{(m)}_{13} K^{(m)}_{3} - L^{(m)}_{0} L^{(m)}_{7} L^{(m)}_{21} K^{(m)}_{2} L^{(m)}_{18} - L^{(m)}_{0} L^{(m)}_{7} K^{(m)}_{4} L^{(m)}_{16} \\&&
\nonumber \times L^{(m)}_{13} + L^{(m)}_{0} L^{(m)}_{11} L^{(m)}_{22} L^{(m)}_{8} K^{(m)}_{3} + L^{(m)}_{11} L^{(m)}_{7} L^{(m)}_{20} L^{(m)}_{3} K^{(m)}_{3} - L^{(m)}_{11} L^{(m)}_{7} \\&&
\nonumber \times K^{(m)}_{4} L^{(m)}_{15} L^{(m)}_{3} + L^{(m)}_{7} L^{(m)}_{21} L^{(m)}_{10} K^{(m)}_{0} L^{(m)}_{18} + L^{(m)}_{7} L^{(m)}_{20} L^{(m)}_{1} K^{(m)}_{2} L^{(m)}_{18} + \\&&
\nonumber L^{(m)}_{0} L^{(m)}_{22} K^{(m)}_{1} L^{(m)}_{16} L^{(m)}_{13} + L^{(m)}_{0} L^{(m)}_{21} K^{(m)}_{1} L^{(m)}_{12} L^{(m)}_{18} + L^{(m)}_{0} K^{(m)}_{4} L^{(m)}_{16} L^{(m)}_{8} \\&&
\times L^{(m)}_{12} - L^{(m)}_{20} L^{(m)}_{2} K^{(m)}_{1} L^{(m)}_{16} L^{(m)}_{13},
\end{eqnarray}
\begin{eqnarray}\label{bi-weights23}
\nonumber C^{(m)}_{27}&=&
L^{(m)}_{21} K^{(m)}_{1} L^{(m)}_{15} L^{(m)}_{2} L^{(m)}_{13} + L^{(m)}_{6} K^{(m)}_{4} L^{(m)}_{12} L^{(m)}_{15} L^{(m)}_{3} - L^{(m)}_{0} L^{(m)}_{22} L^{(m)}_{8} L^{(m)}_{16} \\&&
\nonumber \times K^{(m)}_{2} - L^{(m)}_{0} L^{(m)}_{6} K^{(m)}_{4} L^{(m)}_{12} L^{(m)}_{18} + L^{(m)}_{0} L^{(m)}_{21} L^{(m)}_{8} L^{(m)}_{17} K^{(m)}_{2} - L^{(m)}_{0} L^{(m)}_{21} \\&&
\nonumber \times L^{(m)}_{8} K^{(m)}_{3} L^{(m)}_{12} - L^{(m)}_{22} L^{(m)}_{8} K^{(m)}_{3} L^{(m)}_{10} L^{(m)}_{1} - L^{(m)}_{21} L^{(m)}_{13} K^{(m)}_{3} L^{(m)}_{5} L^{(m)}_{2} + \\&&
\nonumber  L^{(m)}_{6} L^{(m)}_{22} L^{(m)}_{10} L^{(m)}_{3} K^{(m)}_{3} + K^{(m)}_{4} L^{(m)}_{5} L^{(m)}_{1} L^{(m)}_{12} L^{(m)}_{18} - L^{(m)}_{11} L^{(m)}_{22} L^{(m)}_{8} L^{(m)}_{15} \\&&
\times K^{(m)}_{0} + L^{(m)}_{6} L^{(m)}_{22} L^{(m)}_{13} L^{(m)}_{15} K^{(m)}_{0},
\end{eqnarray}
\begin{eqnarray}\label{bi-weights24}
\nonumber C^{(m)}_{28}&=&
L^{(m)}_{7} K^{(m)}_{4} L^{(m)}_{16} L^{(m)}_{10} L^{(m)}_{3} + L^{(m)}_{7} K^{(m)}_{4} L^{(m)}_{15} L^{(m)}_{1} L^{(m)}_{13} - L^{(m)}_{7} L^{(m)}_{20} L^{(m)}_{1} L^{(m)}_{13}\\&&
\nonumber \times K^{(m)}_{3} - K^{(m)}_{4} L^{(m)}_{5} L^{(m)}_{1} L^{(m)}_{17} L^{(m)}_{13} - K^{(m)}_{4} L^{(m)}_{16} L^{(m)}_{8} L^{(m)}_{10} L^{(m)}_{2} - K^{(m)}_{4} L^{(m)}_{15} \\&&
\nonumber \times L^{(m)}_{1} L^{(m)}_{8} L^{(m)}_{12} + L^{(m)}_{21} L^{(m)}_{8} L^{(m)}_{15} K^{(m)}_{0} L^{(m)}_{12} + L^{(m)}_{20} L^{(m)}_{2} L^{(m)}_{8} L^{(m)}_{16} K^{(m)}_{2} - \\&&
\nonumber L^{(m)}_{6} L^{(m)}_{20} L^{(m)}_{3} K^{(m)}_{3} L^{(m)}_{12} + L^{(m)}_{6} L^{(m)}_{23} L^{(m)}_{17} L^{(m)}_{10} K^{(m)}_{0} + L^{(m)}_{0} L^{(m)}_{11} L^{(m)}_{7} K^{(m)}_{4} \\&&
\times L^{(m)}_{18} - L^{(m)}_{0} L^{(m)}_{11} K^{(m)}_{4} L^{(m)}_{17} L^{(m)}_{8},
\end{eqnarray}
\begin{eqnarray}\label{bi-weights25}
\nonumber C^{(m)}_{29}&=&
L^{(m)}_{0} L^{(m)}_{11} L^{(m)}_{23} L^{(m)}_{17} K^{(m)}_{1} - L^{(m)}_{0} L^{(m)}_{11} L^{(m)}_{22} K^{(m)}_{1} L^{(m)}_{18} - L^{(m)}_{20} K^{(m)}_{0} L^{(m)}_{16} L^{(m)}_{8} \\&&
\nonumber \times L^{(m)}_{12} - L^{(m)}_{23} L^{(m)}_{17} L^{(m)}_{10} L^{(m)}_{1} K^{(m)}_{1} + L^{(m)}_{22} L^{(m)}_{13} L^{(m)}_{5} L^{(m)}_{1} K^{(m)}_{3} - L^{(m)}_{23} L^{(m)}_{16} \\&&
\nonumber \times L^{(m)}_{5} L^{(m)}_{2} K^{(m)}_{2} + L^{(m)}_{23} L^{(m)}_{16} L^{(m)}_{5} K^{(m)}_{0} L^{(m)}_{12} - L^{(m)}_{21} K^{(m)}_{1} L^{(m)}_{12} L^{(m)}_{15} L^{(m)}_{3} + \\&&
\nonumber  L^{(m)}_{22} L^{(m)}_{10} L^{(m)}_{1} K^{(m)}_{1} L^{(m)}_{18} + L^{(m)}_{22} L^{(m)}_{8} L^{(m)}_{15} L^{(m)}_{1} K^{(m)}_{2} - L^{(m)}_{21} L^{(m)}_{8} L^{(m)}_{17} L^{(m)}_{10} \\&&
\times K^{(m)}_{0} - L^{(m)}_{22} K^{(m)}_{1} L^{(m)}_{15} L^{(m)}_{1} L^{(m)}_{13},
\end{eqnarray}
\begin{eqnarray}\label{bi-weights26}
\nonumber C^{(m)}_{30}&=&
L^{(m)}_{20} L^{(m)}_{1} L^{(m)}_{8} K^{(m)}_{3} L^{(m)}_{12} - L^{(m)}_{0} L^{(m)}_{11} L^{(m)}_{7} L^{(m)}_{23} K^{(m)}_{3} + L^{(m)}_{23} L^{(m)}_{15} L^{(m)}_{1} K^{(m)}_{1} \\&&
\nonumber \times L^{(m)}_{12} - L^{(m)}_{21} L^{(m)}_{8} L^{(m)}_{15} L^{(m)}_{2} K^{(m)}_{2} - L^{(m)}_{20} L^{(m)}_{1} K^{(m)}_{1} L^{(m)}_{12} L^{(m)}_{18} - L^{(m)}_{22} K^{(m)}_{1}\\&&
\nonumber \times L^{(m)}_{16} L^{(m)}_{10} L^{(m)}_{3} - K^{(m)}_{4} L^{(m)}_{16} L^{(m)}_{5} L^{(m)}_{3} L^{(m)}_{12} + L^{(m)}_{22} L^{(m)}_{8} L^{(m)}_{16} L^{(m)}_{10} K^{(m)}_{0} -\\&&
\nonumber L^{(m)}_{11} L^{(m)}_{23} L^{(m)}_{15} L^{(m)}_{2} K^{(m)}_{1} - L^{(m)}_{6} L^{(m)}_{22} K^{(m)}_{2} L^{(m)}_{15} L^{(m)}_{3} + L^{(m)}_{21} K^{(m)}_{1} L^{(m)}_{17} L^{(m)}_{10}\\&&
\times L^{(m)}_{3} - L^{(m)}_{6} L^{(m)}_{20} L^{(m)}_{2} K^{(m)}_{2} L^{(m)}_{18},
\end{eqnarray}
\begin{eqnarray}\label{bi-weights27}
\nonumber C^{(m)}_{31}&=&
L^{(m)}_{6} K^{(m)}_{4} L^{(m)}_{18} L^{(m)}_{10} L^{(m)}_{2} - L^{(m)}_{6} K^{(m)}_{4} L^{(m)}_{15} L^{(m)}_{2} L^{(m)}_{13} - L^{(m)}_{6} L^{(m)}_{22} L^{(m)}_{10} K^{(m)}_{0} \\&&
\nonumber \times L^{(m)}_{18} - L^{(m)}_{6} K^{(m)}_{4} L^{(m)}_{17} L^{(m)}_{10} L^{(m)}_{3} + L^{(m)}_{11} L^{(m)}_{22} K^{(m)}_{1} L^{(m)}_{15} L^{(m)}_{3} - L^{(m)}_{21} L^{(m)}_{5} \\&&
\nonumber \times L^{(m)}_{3} L^{(m)}_{17} K^{(m)}_{2} - L^{(m)}_{7} L^{(m)}_{21} L^{(m)}_{10} L^{(m)}_{3} K^{(m)}_{3} + L^{(m)}_{23} L^{(m)}_{17} K^{(m)}_{2} L^{(m)}_{5} L^{(m)}_{1} + \\&&
\nonumber L^{(m)}_{20} L^{(m)}_{3} L^{(m)}_{16} K^{(m)}_{1} L^{(m)}_{12} + L^{(m)}_{7} L^{(m)}_{21} K^{(m)}_{2} L^{(m)}_{15} L^{(m)}_{3} + K^{(m)}_{4} L^{(m)}_{10} L^{(m)}_{1} L^{(m)}_{17} \\&&
\times L^{(m)}_{8} + L^{(m)}_{6} L^{(m)}_{20} L^{(m)}_{2} L^{(m)}_{13} K^{(m)}_{3},
\end{eqnarray}
\begin{eqnarray}\label{bi-weights28}
\nonumber C^{(m)}_{32}&=&
L^{(m)}_{7} L^{(m)}_{23} K^{(m)}_{3} L^{(m)}_{10} L^{(m)}_{1} - L^{(m)}_{6} L^{(m)}_{20} K^{(m)}_{0} L^{(m)}_{17} L^{(m)}_{13} - L^{(m)}_{7} K^{(m)}_{4} L^{(m)}_{10} L^{(m)}_{1} \\&&
\nonumber \times L^{(m)}_{18} + L^{(m)}_{20} L^{(m)}_{1} K^{(m)}_{1} L^{(m)}_{17} L^{(m)}_{13} - L^{(m)}_{22} K^{(m)}_{2} L^{(m)}_{18} L^{(m)}_{5} L^{(m)}_{1} + L^{(m)}_{7} L^{(m)}_{20}\\&&
\nonumber \times K^{(m)}_{0} L^{(m)}_{16} L^{(m)}_{13} - L^{(m)}_{7} L^{(m)}_{23} L^{(m)}_{15} L^{(m)}_{1} K^{(m)}_{2} + L^{(m)}_{11} L^{(m)}_{22} L^{(m)}_{5} K^{(m)}_{0} L^{(m)}_{18} + \\&&
\nonumber L^{(m)}_{6} L^{(m)}_{23} L^{(m)}_{15} L^{(m)}_{2} K^{(m)}_{2} - L^{(m)}_{6} L^{(m)}_{23} L^{(m)}_{15} K^{(m)}_{0} L^{(m)}_{12} + L^{(m)}_{21} K^{(m)}_{2} L^{(m)}_{18} L^{(m)}_{5} \\&&
\times L^{(m)}_{2} - L^{(m)}_{11} L^{(m)}_{23} L^{(m)}_{17} L^{(m)}_{5} K^{(m)}_{0},
\end{eqnarray}
\begin{eqnarray}\label{bi-weights29}
\nonumber C^{(m)}_{33}&=&
L^{(m)}_{11} L^{(m)}_{20} K^{(m)}_{0} L^{(m)}_{17} L^{(m)}_{8} + L^{(m)}_{11} L^{(m)}_{23} K^{(m)}_{3} L^{(m)}_{5} L^{(m)}_{2} + L^{(m)}_{21} L^{(m)}_{5} K^{(m)}_{0} L^{(m)}_{17} \\&&
\nonumber \times L^{(m)}_{13} + L^{(m)}_{11} L^{(m)}_{20} L^{(m)}_{2} K^{(m)}_{1} L^{(m)}_{18} - L^{(m)}_{11} L^{(m)}_{20} L^{(m)}_{3} L^{(m)}_{17} K^{(m)}_{1} - L^{(m)}_{7} L^{(m)}_{21} \\&&
\nonumber \times L^{(m)}_{13} L^{(m)}_{15} K^{(m)}_{0} - L^{(m)}_{7} L^{(m)}_{20} L^{(m)}_{3} L^{(m)}_{16} K^{(m)}_{2} - L^{(m)}_{11} L^{(m)}_{7} L^{(m)}_{20} K^{(m)}_{0} L^{(m)}_{18} - \\&&
\nonumber  L^{(m)}_{21} K^{(m)}_{1} L^{(m)}_{18} L^{(m)}_{10} L^{(m)}_{2} + L^{(m)}_{22} L^{(m)}_{5} L^{(m)}_{3} L^{(m)}_{16} K^{(m)}_{2} + L^{(m)}_{23} L^{(m)}_{16} K^{(m)}_{1} L^{(m)}_{10} \\&&
\times L^{(m)}_{2} + L^{(m)}_{11} K^{(m)}_{4} L^{(m)}_{17} L^{(m)}_{5} L^{(m)}_{3},
\end{eqnarray}
\begin{eqnarray}\label{bi-weights30}
\nonumber C^{(m)}_{34}&=&
K^{(m)}_{4} L^{(m)}_{16} L^{(m)}_{5} L^{(m)}_{2} L^{(m)}_{13} - L^{(m)}_{20} L^{(m)}_{1} L^{(m)}_{8} L^{(m)}_{17} K^{(m)}_{2} - L^{(m)}_{11} L^{(m)}_{22} L^{(m)}_{5} L^{(m)}_{3} \\&&
\nonumber \times K^{(m)}_{3} - L^{(m)}_{23} K^{(m)}_{3} L^{(m)}_{5} L^{(m)}_{1} L^{(m)}_{12} - L^{(m)}_{11} K^{(m)}_{4} L^{(m)}_{18} L^{(m)}_{5} L^{(m)}_{2} - L^{(m)}_{7} L^{(m)}_{23} \\&&
\nonumber \times L^{(m)}_{16} L^{(m)}_{10} K^{(m)}_{0} + L^{(m)}_{21} L^{(m)}_{8} K^{(m)}_{3} L^{(m)}_{10} L^{(m)}_{2} + L^{(m)}_{11} L^{(m)}_{7} L^{(m)}_{23} L^{(m)}_{15} K^{(m)}_{0} + \\&&
\nonumber  L^{(m)}_{6} L^{(m)}_{20} L^{(m)}_{3} L^{(m)}_{17} K^{(m)}_{2} - L^{(m)}_{6} L^{(m)}_{23} K^{(m)}_{3} L^{(m)}_{10} L^{(m)}_{2} - L^{(m)}_{21} L^{(m)}_{5} K^{(m)}_{0} L^{(m)}_{12} \\&&
\times L^{(m)}_{18} + L^{(m)}_{6} L^{(m)}_{20} K^{(m)}_{0} L^{(m)}_{12} L^{(m)}_{18},
\end{eqnarray}


\begin{eqnarray}\label{bi-weights31}
\left\{\begin{array}{ccccccc}
&&L^{(m)}_{0}=\tau^{(m)}_{0} \tau^{(m)}_{3}-\tau^{(m)}_{1} \tau^{(m)}_{1},\,\ L^{(m)}_{1}= \tau^{(m)}_{0} \tau^{(m)}_{4}-\tau^{(m)}_{1} \tau^{(m)}_{2},\\ \\ &&L^{(m)}_{2}=\tau^{(m)}_{0} \tau^{(m)}_{6}-\tau^{(m)}_{1} \tau^{(m)}_{3}, \,\ L^{(m)}_{3}=\tau^{(m)}_{0} \tau^{(m)}_{7}-\tau^{(m)}_{1} \tau^{(m)}_{4}, \\ \\
&&L^{(m)}_{4}=\tau^{(m)}_{0} \tau^{(m)}_{8}-\tau^{(m)}_{1} \tau^{(m)}_{5}, \,\ L^{(m)}_{5}=\tau^{(m)}_{1} \tau^{(m)}_{4}-\tau^{(m)}_{2} \tau^{(m)}_{3},
\end{array}\right.
\end{eqnarray}

\begin{eqnarray}\label{bi-weights32}
\left\{\begin{array}{ccccccc}
&&L^{(m)}_{6}=\tau^{(m)}_{1} \tau^{(m)}_{5}-\tau^{(m)}_{2} \tau^{(m)}_{4},\,\ L^{(m)}_{7}= \tau^{(m)}_{1} \tau^{(m)}_{7}-\tau^{(m)}_{2} \tau^{(m)}_{6},\\ \\ &&L^{(m)}_{8}=\tau^{(m)}_{1} \tau^{(m)}_{8}-\tau^{(m)}_{2} \tau^{(m)}_{7}, \,\ L^{(m)}_{9}=\tau^{(m)}_{1} \tau^{(m)}_{9}-\tau^{(m)}_{2} \tau^{(m)}_{8}, \\ \\
&&L^{(m)}_{10}=\tau^{(m)}_{2} \tau^{(m)}_{6}-\tau^{(m)}_{3} \tau^{(m)}_{4}, \,\ L^{(m)}_{11}=\tau^{(m)}_{2} \tau^{(m)}_{7}-\tau^{(m)}_{3} \tau^{(m)}_{5},
\end{array}\right.
\end{eqnarray}

\begin{eqnarray}\label{bi-weights33}
\left\{\begin{array}{ccccccc}
&&L^{(m)}_{12}=\tau^{(m)}_{2} \tau^{(m)}_{10}-\tau^{(m)}_{3} \tau^{(m)}_{7},\,\ L^{(m)}_{13}= \tau^{(m)}_{2} \tau^{(m)}_{11}-\tau^{(m)}_{3} \tau^{(m)}_{8},\\ \\ &&L^{(m)}_{14}=\tau^{(m)}_{2} \tau^{(m)}_{12}-\tau^{(m)}_{3} \tau^{(m)}_{9}, \,\ L^{(m)}_{15}=\tau^{(m)}_{3} \tau^{(m)}_{7}-\tau^{(m)}_{4} \tau^{(m)}_{6}, \\ \\
&&L^{(m)}_{16}=\tau^{(m)}_{3} \tau^{(m)}_{8}-\tau^{(m)}_{4} \tau^{(m)}_{7}, \,\ L^{(m)}_{17}=\tau^{(m)}_{3} \tau^{(m)}_{11}-\tau^{(m)}_{4} \tau^{(m)}_{10},
\end{array}\right.
\end{eqnarray}

\begin{eqnarray}\label{bi-weights34}
\left\{\begin{array}{ccccccc}
&&L^{(m)}_{18}=\tau^{(m)}_{3} \tau^{(m)}_{12}-\tau^{(m)}_{4} \tau^{(m)}_{11},\,\ L^{(m)}_{19}= \tau^{(m)}_{3} \tau^{(m)}_{13}-\tau^{(m)}_{4} \tau^{(m)}_{12},\\ \\ &&L^{(m)}_{20}=\tau^{(m)}_{4} \tau^{(m)}_{8}-\tau^{(m)}_{5} \tau^{(m)}_{7}, \,\ L^{(m)}_{21}=\tau^{(m)}_{4} \tau^{(m)}_{9}-\tau^{(m)}_{5} \tau^{(m)}_{8}, \\ \\
&&L^{(m)}_{22}=\tau^{(m)}_{4} \tau^{(m)}_{12}-\tau^{(m)}_{5} \tau^{(m)}_{11}, \,\ L^{(m)}_{23}=\tau^{(m)}_{4} \tau^{(m)}_{13}-\tau^{(m)}_{5} \tau^{(m)}_{12},\\ \\
&& L^{(m)}_{24}=\tau^{(m)}_{4} \tau^{(m)}_{14}-\tau^{(m)}_{5} \tau^{(m)}_{13},
\end{array}\right.
\end{eqnarray}

\begin{eqnarray}\label{bi-weights35}
\left\{\begin{array}{ccccccc}
&&  K^{(m)}_{b} = \tau^{(m)}_{b+1}A^{(m)}_{b}-\tau^{(m)}_{b}A^{(m)}_{b+1}, \,\ 0 \leqslant b \leqslant 4, \\ \\
&&  A^{(m)}_{b} =  \sum\limits_{r=-n+1}^{n}\sum\limits_{s=-n+1}^{n}r^{b_{1}}s^{b_{2}}w^{(m)}_{r,s}f_{r,s}, \,\ 0 \leqslant b \leqslant 5, \,\ b_{1}+b_{2} \leqslant 2, \\ \\
&&  \tau^{(m)}_{b} =  \sum\limits_{r=-n+1}^{n}\sum\limits_{s=-n+1}^{n}r^{b_{1}}s^{b_{2}}w^{(m)}_{r,s}, \,\ 0 \leqslant b \leqslant 14, \,\ b_{1}+b_{2} \leqslant 4.
\end{array}\right.
\end{eqnarray}
\end{prop}

\begin{prop}\label{bi-2n+1-B-OLS-cor}
If we put $d=2$ in (\ref{bipoly}) and replace $r=-n+1,\ldots,n$ $\&$ $s=-n+1,\ldots,n$ by $r=-n,\ldots,n$ $\&$ $s=-n,\ldots,n$ respectively, then from Algorithm \ref{bi-parameters-Estimating} the starting values $B^{(0)}_{\delta}=\{B^{(0)}_{a,\delta}:a=0,1,\ldots,5\}$ are:

\begin{eqnarray}\label{bi-ols-betas-5-2n+1}
B^{(0)}_{5,\delta} &=&\sum\limits_{r=-n}^{n} \sum\limits_{s=-n}^{n} \frac{-15 (n^{2} + n - 3 s^{2})}{n(n+1)(2n-1)(2n+1)^{2}(2n+3)} f_{r,s},
\end{eqnarray}

\begin{eqnarray}\label{bi-ols-betas-4-2n+1}
B^{(0)}_{4,\delta} &=&\sum\limits_{r=-n}^{n} \sum\limits_{s=-n}^{n} \frac{9 r s}{n^{2}(n+1)^{2}(2n+1)^{2}} f_{r,s},
\end{eqnarray}

\begin{eqnarray}\label{bi-ols-betas-3-2n+1}
B^{(0)}_{3,\delta} &=&\sum\limits_{r=-n}^{n} \sum\limits_{s=-n}^{n} \frac{-15 (n^{2} + n - 3 r^{2 })}{n(n+1)(2n-1)(2n+1)^{2}(2n+3)} f_{r,s},
\end{eqnarray}

\begin{eqnarray}\label{bi-ols-betas-2-2n+1}
B^{(0)}_{2,\delta} &=&\sum\limits_{r=-n}^{n} \sum\limits_{s=-n}^{n} \frac{3 s}{n(n+1)(2n+1)^{2}} f_{r,s},
\end{eqnarray}

\begin{eqnarray}\label{bi-ols-betas-1-2n+1}
B^{(0)}_{1,\delta} &=& \sum\limits_{r=-n}^{n} \sum\limits_{s=-n}^{n} \frac{3 r}{n(n+1)(2n+1)^{2}} f_{r,s},
\end{eqnarray}

\begin{eqnarray}\label{bi-ols-betas-0-2n+1}
B^{(0)}_{0,\delta}&=&\sum\limits_{r=-n}^{n} \sum\limits_{s=-n}^{n} \frac{1}{(2n+1)^{2}} f_{r,s} - \frac{n(n+1)}{3} B^{(0)}_{3,\delta} - \frac{n(n+1)}{3} B^{(0)}_{5,\delta}.
\end{eqnarray}
\end{prop}

\begin{rem}
Since in Step 2 of Section \ref{bivariate-construction-IRLS}, we have set $B^{(m+1)}_{a,\delta}=B^{(m+1)}_{a,\delta,i,j}$ and $w^{(m)}_{r,s}=w^{(m)}_{i+r,j+s}$. Furthermore, set the symbols used in Proposition \ref{bi-betas-IRLS-cor} as $E^{(m)}_{b}=E^{(m)}_{b,i,j}$ for $b=0,1,\ldots,4$; $C^{(m)}_{b}=C^{(m)}_{b,i,j}$ for $b=0,1,\ldots,34$; $L^{(m)}_{b}=L^{(m)}_{b,i,j}$ for $b=0,1,\ldots,24$; $K^{(m)}_{b}=K^{(m)}_{b,i,j}$ for $b=0,1,\ldots,4$; $A^{(m)}_{b}=A^{(m)}_{b,i,j}$ for $b=0,1,\ldots,5$ and $\tau^{(m)}_{b}=\tau^{(m)}_{b,i,j}$ for $b=0,1,\ldots,14$.
\end{rem}

\begin{thm}
If first we take values of $B^{(0)}_{\delta}$ and $B^{(m+1)}_{\delta}$ from Propositions \ref{bi-2n-betas-OLS-cor}-\ref{bi-betas-IRLS-cor} respectively and then apply Algorithm \ref{bi-SS-A} for $d=2$, we get the family of $(2n)^{2}$-point subdivision schemes $D_{(2n)^{2},2}$ based on fitting bivariate polynomial of degree 2.
\end{thm}

\begin{thm}
 We get the family of $(2n)^{2}$-point subdivision schemes $D_{(2n)^{2},1}$ based on fitting bivariate polynomial of degree 1, if first we put $B^{(0)}_{5,\delta}=B^{(0)}_{4,\delta}=B^{(0)}_{3,\delta}=0$ and $B^{(m+1)}_{5,\delta}=B^{(m+1)}_{4,\delta}=B^{(m+1)}_{3,\delta}=0$ in Propositions \ref{bi-2n-betas-OLS-cor}-\ref{bi-betas-IRLS-cor} respectively and take values of $B^{(0)}_{2,\delta}$, $B^{(0)}_{1,\delta}$, $B^{(0)}_{0,\delta}$, $B^{(m+1)}_{2,\delta}$, $B^{(m+1)}_{1,\delta}$ and $B^{(m+1)}_{0,\delta}$ from Propositions \ref{bi-2n-betas-OLS-cor}-\ref{bi-betas-IRLS-cor} and then use Algorithm \ref{bi-SS-A} for $d=1$.
\end{thm}

\begin{thm}
If first we replace $r=-n+1,\ldots,n$ $\&$ $s=-n+1,\ldots,n$ by $r=-n,\ldots,n$ $\&$ $s=-n,\ldots,n$ respectively in Algorithm \ref{bi-SS-A} and (\ref{bi-IRLS-betas-5-2n})-(\ref{bi-weights35}) of Proposition \ref{bi-betas-IRLS-cor}, then use Propositions \ref{bi-betas-IRLS-cor}-\ref{bi-2n+1-B-OLS-cor} and then apply Algorithm \ref{bi-SS-A} for $d=2$ and $(r,s)=\left(-\frac{1}{4},-\frac{1}{4}\right)$, $\left(\frac{1}{4},-\frac{1}{4}\right)$, $\left(-\frac{1}{4},\frac{1}{4}\right)$, $\left(\frac{1}{4},\frac{1}{4}\right)$, we get the four refinement rules for the family of $(2n+1)^{2}$-point subdivision schemes $D_{(2n+1)^{2},2}$ based on fitting bivariate polynomial of degree 2.
\end{thm}

\begin{thm}
Let us replace $r=-n+1,\ldots,n$ $\&$ $s=-n+1,\ldots,n$ by $r=-n,\ldots,n$ $\&$ $s=-n,\ldots,n$ respectively in Algorithm \ref{bi-SS-A} and (\ref{bi-IRLS-betas-5-2n})-(\ref{bi-weights35}) of Proposition \ref{bi-betas-IRLS-cor}. Now put $B^{(m+1)}_{5,\delta}=B^{(m+1)}_{4,\delta}=B^{(m+1)}_{3,\delta}=0$ and $B^{(0)}_{5,\delta}=B^{(0)}_{4,\delta}=B^{(0)}_{3,\delta}=0$ in Propositions \ref{bi-betas-IRLS-cor}-\ref{bi-2n+1-B-OLS-cor} respectively and take the values of $B^{(m+1)}_{2,\delta}$, $B^{(m+1)}_{1,\delta}$, $B^{(m+1)}_{0,\delta}$, $B^{(0)}_{2,\delta}$, $B^{(0)}_{1,\delta}$ and $B^{(0)}_{0,\delta}$ from Propositions \ref{bi-betas-IRLS-cor}-\ref{bi-2n+1-B-OLS-cor}. Then we apply Algorithm \ref{bi-SS-A} for $d=1$ and $(r,s)=\left(-\frac{1}{4},-\frac{1}{4}\right)$, $\left(\frac{1}{4},-\frac{1}{4}\right)$, $\left(-\frac{1}{4},\frac{1}{4}\right)$, $\left(\frac{1}{4},\frac{1}{4}\right)$, we get the $(2n+1)^{2}$-point schemes $D_{(2n+1)^{2},1}$ based on fitting linear bivariate polynomial.
\end{thm}
\begin{figure}[htb!] 
 \begin{center}
(a) \quad (b)\\
Figure will be supplied on demand. Exceeded size limits.
\end{center}
\caption[Original surface with its data points and noisy surface with its data points.]{\label{bi-original-surf} \emph{(a) is original surface with original data points while (b) is noisy surface with outliers and its data points.}}
\end{figure}

\begin{figure}[htb!] 
\begin{center}
(a) 1st subdivision level \quad (b) 2nd subdivision level \quad (c) 3rd subdivision level \\
\end{center}
\begin{center}
(d) 1st subdivision level \quad (e) 2nd subdivision level \quad (f) 3rd subdivision level\\
Figure will be supplied on demand. Exceeded size limits.
\end{center}
\caption[Effects of the schemes $D_{h^{2},1}$ on 3D noisy data with outliers.]{\label{bi-linear-surf} \emph{Effects of the noisy data with outliers on fitted surfaces.}}
\end{figure}

\begin{figure}[htb!] 
\begin{center}
(a) 1st subdivision level \quad (b) 2nd subdivision level \quad (c) 3rd subdivision level\\
\end{center}
\begin{center}
(d) 1st subdivision level \quad (e) 2nd subdivision level \quad (f) 3rd subdivision level\\
Figure will be supplied on demand. Exceeded size limits.
\end{center}
\caption[Effects of the schemes $D_{h^{2},2}$ on 3D noisy data with outliers.]{\label{bi-quad-surf} \emph{Effects of the noisy data with outliers on fitted surfaces.}}
\end{figure}

\subsection{Numerical example}
 Here we present numerical example to show the performance of the proposed non-tensor product bivariate schemes. In this example, we examine a parametric torus surface expressed by the following parametric equations
\begin{eqnarray*}
x(\nu_{1},\nu_{2})&=& \left(c_{2}+c_{1} cos (\nu_{2})\right)cos (\nu_{1}),\\
y(\nu_{1},\nu_{2})&=& \left(c_{2}+c_{1} cos (\nu_{2})\right)sin(\nu_{1}),\\
z(\nu_{1},\nu_{2})&=& c_{1} sin(\nu_{2}).
\end{eqnarray*}
The surface and data shown in Figure \ref{bi-original-surf}(a) have been generated by taking $c_{1}=2$, $c_{2}=5$ and $0 \leq \nu_{1}, \nu_{2} \leq 2 \pi$, whereas the surface and data shown in Figure \ref{bi-original-surf}(b) have been generated by adding noise and some outliers in the data that has been shown in \ref{bi-original-surf}(a). Figures \ref{bi-linear-surf} and \ref{bi-quad-surf} show the performance of the proposed bivariate schemes that are based on linear and quadratic bivariate polynomials respectively. It is clear from these figures that the proposed bivariate schemes are best for fitting 3D data that contains noise and outliers. Moreover, from Figures \ref{bi-linear-surf}-\ref{bi-quad-surf}, it is also observed that the bivariate schemes based on bi-quadratic polynomials smooth data/mesh more efficiently than the schemes based on bi-linear polynomials.

\section{Conclusion}

In this article, we have presented six families of univariate binary subdivision schemes which are based on univariate polynomials up to degree three and $\ell_{1}$-regression. These schemes have been constructed for noisy data with outliers and have the ability to remove under and over fitting of the data. Since in this article we have extended the technique of Mustafa et al. \cite{Mustafa9} for the construction of the families of schemes, their schemes have become special cases of the proposed schemes. The schemes of \cite{Mustafa11} are also the special cases of our schemes. We have given several numerical examples to see the performance of the proposed schemes. Our schemes reproduce non-linear polynomial data and have the ability to show the interpolatory and approximating behavior, according to the nature of initial data points. We have also given the suggestion to increase the arity of these schemes. Visual performance of the families of $\ell_{1}$- schemes based on bivariate polynomials up to degree two have also been presented.

\subsection*{Acknowledgement}
This work is supported by NRPU Project. No. 3183, Pakistan.

\end{document}